%% file: main.tex
\documentclass[a4paper,11pt,leqno]{amsart}

\input{Preamble/Preamble}

\usepackage{subfiles}

\title{Path and cycle decompositions of dense graphs}

\begin{document}
\renewcommand{\onlyinsubfile}[1]{}
\renewcommand{\notinsubfile}[1]{#1}

\begin{abstract}

\subfile{Abstract}
\end{abstract}

\maketitle

\section{Introduction}

\subfile{Introduction}

\subfile{Organisation}

\section{Proof overview of the main theorems}\label{sec:sketch}
\subfile{Sketch_of_Proof}

\section{Notation, definitions, and probabilistic tools}\label{sec:notation}
\subfile{Notation_and_Probabilistic_Tools}

\section{Preliminary results}\label{sec:preliminaries}
In this section, we introduce some preliminary results which will be useful in the proof of our main theorems. In \cref{sec:regularity,sec:tyingpaths}, we collect some useful properties of $\varepsilon$-regular pairs and prove some lemmas for tying paths together. These results will be used repeatedly in the rest of the paper. In \cref{sec:regularising,sec:robust}, we introduce some tools for regularising superregular pairs and state the robust decomposition lemma of \cite{kuhn2013hamilton}, which will be needed in \cref{sec:Gamma}.

\subfile{Regularity}

\subfile{Tying_Paths}

\subfile{Regularising}

\subfile{Robust_Decomposition_Lemma}

\section{Proof of the main theorems}\label{sec:proofs}

In \cref{sec:Szemeredi,sec:insideV0,sec:exceptionaledges,sec:bad,sec:Gamma}, we prove the main lemmas that will be needed for the proof of our theorems. These intermediate results are organised according to the structure of the proof overview.
\cref{thm:Delta/2,thm:n/2,thm:epspathdecomp} are proved in \cref{sec:thms}.

\subfile{Szemeredi_and_Cleaning}

\subfile{Inside_V0}

\subfile{Exceptional_Edges}

\subfile{Bad_Graph}

\subfile{Decomposing_Gamma}

\subfile{Proof_of_Main_Theorems}

\section{Concluding remarks}\label{sec:conclusion}
\subfile{Concluding_Remarks}

\NEW{\section*{Acknowledgements}
We are grateful to the referee for detailed comments on an earlier version of this paper.}

\bibliographystyle{abbrv}
\bibliography{Bibliography/papers}

\APPENDIX{\appendix
\section*{Appendix: Proof of Lemma~\ref{lm:goodgraph}}

\subfile{Appendix}}

\end{document}

%% file: Preamble/Preamble.tex
\usepackage[shortlabels]{enumitem}
\usepackage{xr-hyper, zref}
\usepackage{hyperref}
\usepackage{etoolbox}
\usepackage{bm}
\usepackage{amsthm,amssymb,amsmath}
\usepackage{graphicx,subcaption,placeins}
\usepackage{psfrag,enumitem,mathrsfs}
\usepackage{mathtools}
\usepackage[foot]{amsaddr}
\usepackage[british]{babel}
\usepackage[babel]{microtype}
\usepackage[margin=1in]{geometry}
\usepackage[capitalise]{cleveref}
\usepackage[textsize=footnotesize]{todonotes}
\usepackage{environ}
\usepackage{float}
\usepackage{bigfoot}
\usepackage[noadjust]{cite}

\date{}
\address[A. Gir\~{a}o]{Institut f\"{u}r Informatik, Universit\"{a}t Heidelberg, 69120 Heidelberg, Deutschland}
\address[B. Granet, D. K\"{u}hn, and D. Osthus]{School of Mathematics, University of Birmingham, Edgbaston, Birmingham, B15 2TT, United Kingdom}

\author[Ant\'{o}nio~Gir\~{a}o]{Ant\'{o}nio Gir\~{a}o}
\email{tzgirao@gmail.com}

\author[Bertille~Granet]{Bertille Granet}
\email{bxg855@bham.ac.uk}

\author[Daniela~K\"{u}hn]{Daniela K\"{u}hn}
\email{d.kuhn@bham.ac.uk}

\author[Deryk~Osthus]{Deryk Osthus}
\email{d.osthus@bham.ac.uk}

\thanks{This project has received partial funding from the European Research 
	Council (ERC) under the European Union's Horizon 2020 research and innovation programme (grant agreement no. 786198, D.~K\"{u}hn and D.~Osthus).
	The research leading to these results was also partially supported by the EPSRC, grant nos. EP/N019504/1 (A.~Gir\~{a}o and D.~K\"{u}hn) and EP/S00100X/1 (D.~Osthus), as well as the Royal Society and the Wolfson Foundation (D.~K\"{u}hn).
}

\subjclass[2010]{\NEW{05C70, 05C38, 05D40}} 

\newcommand{\COMMENT}[1]{}
\newcommand{\NEW}[1]{#1}
\newcommand{\OLD}[1]{}
\newcommand{\APPENDIX}[1]{}
\newcommand{\NOAPPENDIX}[1]{#1}
\renewcommand{\APPENDIX}[1]{#1}\renewcommand{\NOAPPENDIX}[1]{}% comment out for the appendix-free version

\newcommand{\onlyinsubfile}[1]{#1}
\newcommand{\notinsubfile}[1]{}

\hypersetup{colorlinks=true,
	citecolor=blue,
	filecolor=blue,
	linkcolor=blue,
}

\SetEnumitemKey{longlabel}{wide=0.5cm, leftmargin=*, align=right}
\setlist[enumerate]{itemsep=3pt, topsep=5pt,leftmargin=1.2cm}
\setlist[itemize]{itemsep=3pt, topsep=2pt}

\newlist{steps}{enumerate}{1}
\setlist[steps,1]{
	label=\textbf{Step \arabic*:},
	ref=\arabic*, % if to be used with \cref, don't provide label string or parentheses
	wide,
	parsep=0pt,
	itemsep=10pt,
	topsep=10pt}
\Crefname{stepsi}{Step}{Steps}
\crefname{stepsi}{Step}{Steps}

\newlist{case}{enumerate}{1}
\setlist[case,1]{
	label=\textbf{Case \arabic*:},
	ref=\arabic*, % if to be used with \cref, don't provide label string or parentheses
	wide,
	parsep=0pt,
	itemsep=10pt,
	topsep=10pt}
\Crefname{casei}{Case}{Cases}
\crefname{casei}{Case}{Cases}

\AtBeginEnvironment{thm}{\setlist[enumerate,1]{label=\upshape(\roman*),ref=(\roman*)}}
\AtBeginEnvironment{lm}{\setlist[enumerate,1]{label=\upshape(\roman*)}}
\AtBeginEnvironment{prop}{\setlist[enumerate,1]{label=\upshape(\roman*)}}
\AtBeginEnvironment{cor}{\setlist[enumerate,1]{label=\upshape(\roman*)}}
\AtBeginEnvironment{claim}{\setlist[enumerate,1]{label=\upshape(\roman*)}}
\AtBeginEnvironment{thm}{\setlist[enumerate,2]{label=\upshape(\alph*)}}
\AtBeginEnvironment{lm}{\setlist[enumerate,2]{label=\upshape(\alph*)}}
\AtBeginEnvironment{prop}{\setlist[enumerate,2]{label=\upshape(\alph*)}}
\AtBeginEnvironment{cor}{\setlist[enumerate,2]{label=\upshape(\alph*)}}
\AtBeginEnvironment{claim}{\setlist[enumerate,2]{label=\upshape(\alph*)}}

\Crefname{enumi}{}{}
\Crefname{thm}{Theorem}{Theorems}
\Crefname{lm}{Lemma}{Lemmas}
\Crefname{cor}{Corollary}{Corollaries}
\Crefname{prop}{Proposition}{Propositions}
\Crefname{claim}{Claim}{Claims}
\Crefname{equation}{}{}
\Crefname{conjecture}{Conjecture}{Conjectures}
\Crefname{figure}{Figure}{Figures}

\newtheorem{definition}{Definition}[section]
\newtheorem{claim}{Claim}
\newtheorem{prop}[definition]{Proposition}
\newtheorem{thm}[definition]{Theorem}
\newtheorem{cor}[definition]{Corollary}
\newtheorem{lm}[definition]{Lemma}
\newtheorem{fact}[definition]{Fact}
\newtheorem{conjecture}[definition]{Conjecture}

\NewEnviron{property}[1]{\begin{equation}
	\tag{\(#1\)}
	\parbox{0.9\linewidth}{{\itshape\BODY}}
	\end{equation}}

\AtBeginEnvironment{lm}{\setcounter{claim}{0}}
\AtBeginEnvironment{thm}{\setcounter{claim}{0}}
\AtBeginEnvironment{prop}{\setcounter{claim}{0}}
\AtBeginEnvironment{cor}{\setcounter{claim}{0}}

\newenvironment{proofclaim}{\begin{proof}[Proof of Claim]}{\end{proof}}

\numberwithin{equation}{section}

\newcommand{\eps}{\varepsilon}
\renewcommand{\epsilon}{\varepsilon}

\DeclareMathOperator{\Bin}{Bin}
\DeclareMathOperator{\HGeom}{Hyp}
\DeclareMathOperator{\odd}{odd}

\DeclareMathOperator{\Fict}{Fict}

\newcommand{\cC}{\mathcal{C}}
\newcommand{\cD}{\mathcal{D}}
\newcommand{\cE}{\mathcal{E}}
\newcommand{\cF}{\mathcal{F}}

\newcommand{\cH}{\mathcal{H}}
\newcommand{\cI}{\mathcal{I}}

\newcommand{\cN}{\mathcal{N}}
\newcommand{\cO}{\mathcal{O}}
\newcommand{\cP}{\mathcal{P}}
\newcommand{\cQ}{\mathcal{Q}}

\newcommand{\cS}{\mathcal{S}}

\newcommand{\tG}{\widetilde{G}}

\newcommand{\tH}{\widetilde{H}}

\newcommand{\tR}{\widetilde{R}}

\newcommand{\tGamma}{\widetilde{\Gamma}}

\newcommand{\hG}{\widehat{G}}
\newcommand{\hR}{\widehat{R}}

\newcommand{\sC}{\mathscr{C}}

\newcommand{\sP}{\mathscr{P}}

%% file: Abstract.tex
	We make progress on three long standing conjectures from the 1960s about path and cycle decompositions of graphs. Gallai conjectured that any connected graph on~$n$ vertices can be decomposed into at most~$\left\lceil \frac{n}{2}\right\rceil$ paths, while a conjecture of Haj\'{o}s states that any Eulerian graph on~$n$ vertices can be decomposed into at most~$\left\lfloor \frac{n-1}{2}\right\rfloor$ cycles. The Erd\H{o}s-Gallai conjecture states that any graph on~$n$ vertices can be decomposed into~$O(n)$ cycles and edges.
	
	We show that if~$G$ is a sufficiently large graph on~$n$ vertices with linear minimum degree, then the following hold.
	\begin{enumerate}[label=\upshape(\roman*)]
		\item $G$ can be decomposed into at most $\frac{n}{2}+o(n)$ paths.\label{abs:pathdecomp}
		\item If~$G$ is Eulerian, then it can be decomposed into at most $\frac{n}{2}+o(n)$ cycles.\label{abs:cycledecomp}
		\item $G$ can be decomposed into at most $\frac{3 n}{2}+o(n)$ cycles and edges. \label{abs:cycleedgedecomp}
	\end{enumerate}
	If in addition~$G$ satisfies a weak expansion property, we asymptotically determine the required number of paths/cycles for each such~$G$.
	\begin{enumerate}[resume,label=\upshape(\roman*)]
		\item $G$ can be decomposed into $\max \left\{\frac{\odd(G)}{2},\frac{\Delta(G)}{2}\right\}+o(n)$ paths, where~$\odd(G)$ is the number of odd-degree vertices of~$G$.\label{abs:quasirandompathdecomp}
	 	\item If~$G$ is Eulerian, then it can be decomposed into~$\frac{\Delta(G)}{2}+o(n)$ cycles.\label{abs:quasirandomcycledecomp}
	\end{enumerate}
	All bounds in \cref{abs:cycledecomp,abs:cycleedgedecomp,abs:pathdecomp,abs:quasirandomcycledecomp,abs:quasirandompathdecomp}
	are asymptotically best possible.

%% file: Introduction.tex
	\onlyinsubfile{
		\setcounter{section}{1}
		\addtocounter{section}{-1}
	\section{Introduction}}

\subsection{Background}

Graph decomposition is a central field of graph theory, which encompasses some of the oldest and most famous problems in combinatorics. For example, the decomposition of complete graphs into Hamilton cycles or Hamilton paths was attributed to Walecki and dates back to 1883 \cite{lucas1883recreationsII} (see \cite{alspach2008wonderful} for a description in English of Walecki's construction).  
Extensive research has also been done on decompositions of graphs into (not necessarily Hamiltonian) paths and/or cycles. 
One of the most famous results in this area is due to Lov\'{a}sz.

\begin{thm}[\cite{lovasz1968covering}]\label{thm:Lovasz}
	Let~$G$ be a graph on~$n$ vertices. Then~$G$ can be decomposed into at most $\left\lfloor\frac{n}{2}\right\rfloor$ paths and cycles.
\end{thm}

We observe that this result is sharp. Indeed, a vertex of odd degree in a graph~$G$ must be the endpoint of at least one path in a path and cycle decomposition of~$G$. Thus,~$n$-vertex graphs with at most one vertex of even degree cannot be decomposed into fewer than $\left\lfloor\frac{n}{2}\right\rfloor$ paths and cycles.

The result of Lov\'{a}sz was inspired by the following conjecture of Gallai (see \cite{lovasz1968covering}).

\begin{conjecture}[Gallai]\label{conj:Gallai}
	Any connected graph on~$n$ vertices can be decomposed into at most $\left\lceil\frac{n}{2}\right\rceil$ paths.
\end{conjecture}

Complete graphs show that the conjecture of Gallai would be best possible. 
Lov\'{a}sz \cite{lovasz1968covering} observed that \cref{thm:Lovasz} implies that any graph can be decomposed into at most~$n-1$ paths. This was later improved by Donald \cite{donald1980upper} who showed that~$\left\lfloor\frac{3n}{4}\right\rfloor$ paths are sufficient. It was subsequently shown by Dean and Kouider \cite{dean2000gallai} and independently by Yan \cite{yan1998path} that~$\left\lfloor\frac{2n}{3}\right \rfloor$ paths suffice.  The covering version of Gallai's conjecture (where the paths are not necessarily edge-disjoint) was solved by Fan \cite{fan2002subgraph}.

Although the conjecture of Gallai remains open, it has been verified for several classes of graphs. 
We direct the readers to \cite{fan2005path,harding2014gallai,botler2017path,jimenez2017path,bonamy2019gallai,pyber1996covering,lovasz1968covering} for some of these results.%

The analogous problem for cycle decompositions was posed by Haj\'{o}s (see \cite{lovasz1968covering}). We note that the original problem suggested by Haj\'os asked for a decomposition of Eulerian~$n$-vertex graphs into at most $\left\lfloor\frac{n}{2}\right\rfloor$ cycles, but Dean \cite{dean1986smallest} observed that this is equivalent to the following.

\begin{conjecture}[Haj\'{o}s]\label{conj:Hajos}
	Any Eulerian graph on~$n$ vertices can be decomposed into at most $\left\lfloor\frac{n-1}{2}\right\rfloor$ cycles.
\end{conjecture}

Eulerian graphs with maximum degree~$n-1$ demonstrate that the conjecture of Haj\'{o}s would be best possible. 
\cref{conj:Hajos} has only been verified for specific classes of graphs. See \cite{fan2002hajos} for some of these results.%
\COMMENT{Heinrich, Natale, Streicher \cite{heinrich2017haj} verified the conjecture for graphs of order at most~$12$. Granville and Moisiadis \cite{granville1987hajos}, Favaron and Kouider \cite{favaron1988path}, and Botler, Sambinelli, Coelho and Lee \cite{botler2017gallai} independently proved Haj\'{o}s conjecture for graphs of maximum degree at most~$4$. Favaron and Kouider \cite{favaron1988path} also showed Haj\'{o}s conjecture is true for minimally~$2$-connected or minimally~$2$-edge-connected Eulerian graphs, while Botler, Sambinelli, Coelho and Lee \cite{botler2017gallai} also showed Haj\'{o}s conjecture for graphs of treewidth at most~$3$. 
	Fuchs, Gellert, Heinrich \cite{fuchs2017cycle} verified the conjecture for graphs of pathwidth at most~$6$. 
	Fan and Xu \cite{fan2002hajos} showed that Haj\'{o}s conjecture holds for projective graphs and~$K_6^-$-minor free graphs ($K_6$ with an edge deleted).
	Tao \cite{jiang1984hajos} and Seyffarth \cite{seyffarth1992hajos} verified Haj\'{o}s' conjecture for planar graphs.}
Again, the analogous covering problem was resolved by Fan \cite{fan2003covers}.

Jackson \cite{jackson1980decompositions} conjectured the analogue of \cref{conj:Hajos} for Eulerian oriented graphs. However, Dean \cite{dean1986smallest} observed that this conjecture is false and conjectured instead that any Eulerian \NEW{oriented graph on $n$ vertices}\OLD{digraph} can be decomposed into $\left\lfloor\frac{2n}{3}\right\rfloor$ dicycles \NEW{and any Eulerian digraph on $n>1$ vertices can be decomposed into $\frac{8n}{3}-3$ dicycles.}

Very little progress has been made on \cref{conj:Hajos} for general graphs. In particular, the related problem of decomposing Eulerian graphs into~$O(n)$ cycles is still open and is equivalent to a problem posed in  \cite{erdos1966representation} which is known as the Erd\H{o}s-Gallai conjecture (see \cite{erdos1983some}). 

\begin{conjecture}[Erd\H{o}s-Gallai]\label{conj:Erdos-Gallai}
	Any graph on~$n$ vertices can be decomposed into~$O(n)$ cycles and edges.
\end{conjecture}

Observe that given any~$n$-vertex graph~$G$, by repeatedly removing cycles until no longer possible, we obtain a forest~$F$ such that~$G\setminus F$ is Eulerian. Since this forest contains at most~$n-1$ edges, the problem of decomposing graphs into~$O(n)$ cycles and edges reduces to decomposing Eulerian graphs into~$O(n)$ cycles. Conversely, given a decomposition of an Eulerian graph~$G$ into~$O(n)$ cycles and edges, one can easily obtain a decomposition of~$G$ into~$O(n)$ cycles. Thus, \cref{conj:Erdos-Gallai} is equivalent to the problem of decomposing Eulerian graphs into~$O(n)$ cycles. 
Also observe that \cref{conj:Hajos} would imply that any graph can be decomposed into at most~$\frac{3(n-1)}{2}$ cycles and edges. 
Thus, the Erd\H{o}s-Gallai conjecture holds for all classes of graphs for which \cref{conj:Hajos} has been verified. 
Additionally, the Erd\H{o}s-Gallai conjecture was verified for graphs of linear minimum degree by Conlon, Fox, and Sudakov \cite{conlon2014cycle}. More precisely, they showed the following.

\begin{thm}[{\cite{conlon2014cycle}}]\label{thm:linmindegErdosGallai}
	For any~$\alpha >0$, if~$G$ is a graph on~$n$ vertices with minimum degree~$\delta(G)\geq\alpha n$, then~$G$ can be decomposed into~$O(\alpha^{-12}n)$ cycles and edges.
\end{thm} 

\Cref{conj:Erdos-Gallai} remains open for general graphs, while the covering version was proved by Pyber \cite{pyber1985erdHos}.

It is not hard to show that $\frac{n\log(n)}{2}+O(n)$ cycles and edges are sufficient to decompose any graph (see \cite{erdos1983some}). 
An example of Erd\H{o}s \cite{erdos1983some} shows that at least~$(\frac{3}{2}-o(1))n$ cycles and edges are necessary for some graphs. It was recently shown in \cite{conlon2014cycle} that~$O(n\log\log n)$ cycles and edges are sufficient, and this is currently the best known result for general graphs. More precisely, they proved the following. 

\begin{thm}[{\cite{conlon2014cycle}}]\label{thm:loglognErdosGallai}
	Let~$G$ be a graph on~$n$ vertices with average degree~$d$. Then~$G$ can be decomposed into~$O(n\log\log d)$ cycles and edges. 
\end{thm} 

Progress on \cref{conj:Gallai,conj:Erdos-Gallai} was also made for random graphs. As proved in \cite{conlon2014cycle}, for any edge probability~$p\coloneqq p(n)$, the binomial random graph~$G(n,p)$ satisfies \cref{conj:Erdos-Gallai} asymptotically almost surely. More details about decompositions of random graphs into cycles and edges can be found in \cite{korandi2015decomposing}, where Kor\'{a}ndi, Krivelevich and Sudakov provided an asymptotically tight result for a large range of edge probabilities~$p(n)$.

For constant edge probability~$0<p<1$, Glock, K\"uhn, and Osthus \cite{glock2016optimal} strengthened the bounds of \cite{conlon2014cycle,korandi2015decomposing} to obtain precise results for decompositions of~$G(n,p)$ into paths or cycles and into matchings. In fact, they obtained their results for quasirandom graphs. More precisely, they used the following notion of quasirandomness. An~$n$-vertex graph~$G$ is \emph{lower-$(\varepsilon,p)$-regular} if for any disjoint~$S,T\subseteq V(G)$ with $|S|, |T|\geq \varepsilon n$, we have $e_G(S,T)\geq (p-\varepsilon)|S||T|$.
Given a graph~$G$, we denote by~$\odd(G)$ the number of odd-degree vertices of~$G$.

\begin{thm}[{\cite{glock2016optimal}}]\label{thm:quasirandom}
	For any~$0<p<1$, there exist $\varepsilon, \eta, n_0>0$ such that for any~$n\geq n_0$ the following hold. Let~$G$ be a lower-$(\varepsilon, p)$-regular graph on~$n$ vertices with $\Delta(G)-\delta(G)\leq \eta n$. Then,
	\begin{enumerate}
		\item $G$ can be decomposed into $\max\left\{\frac{\odd(G)}{2}, \left\lceil\frac{\Delta(G)+1}{2}\right\rceil\right\}$ paths, and
		\item if~$G$ is Eulerian, then it can be decomposed into~$\frac{\Delta(G)}{2}$ cycles.
	\end{enumerate}
\end{thm}

These bounds are best possible for each~$G$, but do not hold in general (some examples can be found in \cref{sec:conclusion}).

\NEW{Bienia and Meyniel \cite{bienia1986partitions} conjectured the analogue of \cref{conj:Erdos-Gallai} for Eulerian digraphs.
\begin{conjecture}[Bienia and Meyniel]\label{conj:O(n)digraph}
	There exists $\alpha\in \mathbb{R}$ such that any Eulerian digraph on~$n$ vertices can be decomposed into at most $\alpha n$ dicycles. 
\end{conjecture}
As mentioned in \cite{bienia1986partitions,dean1986smallest}, unions of complete symmetric digraphs $K_4^*$ which are all sharing a common vertex show that, if \cref{conj:O(n)digraph} is true, then $\alpha\geq \frac{4}{3}$.
\cref{conj:O(n)digraph} is also discussed in \cite{bollobas1996proof}.
%	\COMMENT{For each Eulerian digraph $D$, denote by $\cn(D)$ the minimum $k$ such that $D$ can be decomposed into $k$ dicycles. Assume \cref{conj:O(n)digraph} holds. Suppose for a contradiction that there exist a digraph $D$ on $n$ vertices such that $\alpha(n-1)<\cn(D)\leq \alpha n$.
%	Pick an integer $k>\frac{\alpha}{\lceil\alpha(n-1)\rceil-\alpha(n-1)-1}$ and observe that $\frac{\alpha}{k}<\lceil\alpha(n-1)\rceil-\alpha(n-1)-1$. Let $D'$ be obtained by taking $k$ copies of $D$, all sharing a vertex. Then, $\cn(D')=k\cn(D)$. Thus,
%	\[\cn(D)=\frac{\cn(D')}{k}\leq \left\lfloor\frac{\alpha(k(n-1)+1)}{k}\right\rfloor=\left\lfloor \alpha(n-1)+\frac{\alpha}{k}\right\rfloor<\lfloor\alpha(n-1)+\lceil\alpha(n-1)\rceil-\alpha(n-1)\rfloor=\lceil\alpha(n-1)\rceil,\]
%	a contradiction.}.
It is still open but some progress was recently made by Knierim, Larcher, Martinsson, and Noever \cite{knierim2019long}.
\begin{thm}[{\cite{knierim2019long}}]
	Let $D$ be an Eulerian digraph on $n$ vertices and with maximum degree $\Delta$. Then $D$ can be decomposed into $O(n\log\Delta)$ dicycles.
\end{thm}}

\subsection{New results}

First, we prove approximate versions of \cref{conj:Gallai,conj:Hajos} for sufficiently large graphs of linear minimum degree (see \namecrefs{thm:n/2}~\labelcref{thm:n/2}\cref{thm:n/2pathdecomp} and~\labelcref{thm:n/2}\labelcref{thm:n/2cycledecomp}). \cref{thm:n/2}\cref{thm:n/2cycledecomp} easily implies \cref{thm:n/2}\cref{thm:ErdosGallai}, which improves \cref{thm:linmindegErdosGallai} and gives (asymptotically) the best possible constant.

\begin{thm}\label{thm:n/2}
	For any~$\alpha,\delta >0$, there exists~$n_0$ such that if~$G$ is a graph on~$n\geq n_0$ vertices with $\delta(G)\geq \alpha n$, then the following hold.
	\begin{enumerate}
		\item $G$ can be decomposed into at most~$\frac{n}{2}+\delta n$ paths. \label{thm:n/2pathdecomp}
		\item If~$G$ is Eulerian, then it can be decomposed into at most~$\frac{n}{2} +\delta n$ cycles.\label{thm:n/2cycledecomp}
		\item $G$ can be decomposed into at most~$\frac{3n}{2}+\delta n$ cycles and edges. \label{thm:ErdosGallai}
	\end{enumerate}
\end{thm}

Secondly, we prove approximate versions of the bounds in \cref{thm:quasirandom} for sufficiently large graphs with linear minimum degree which satisfy a weak version of quasirandomness. More precisely, we say an~$n$-vertex graph~$G$ is \emph{weakly-$(\varepsilon, p)$-quasirandom} if for any partition~$A\cup B$ of~$V(G)$ with $|A|, |B|\geq \varepsilon n$ we have $e_G(A, B)\geq p|A||B|$. \NEW{This notion of weak quasirandomness implies that the reduced graph obtained after applying the regularity lemma to a dense graph is connected. This is the only property required to obtain the bounds in the following theorem.}

\begin{thm}\label{thm:Delta/2}
	For any~$\alpha,\delta, p>0$, there exists~$n_0$ such that if~$G$ is a weakly-$(\frac{\alpha}{2},p)$-quasirandom graph on~$n\geq n_0$ vertices with $\delta(G)\geq \alpha n$, then the following hold. 
	\begin{enumerate}
		\item $G$ can be decomposed into at most~$\max\left\{ \frac{\odd(G)}{2},\frac{\Delta(G)}{2}\right\} +\delta n$ paths.\label{thm:Delta/2pathdecomp}
		\item If~$G$ is Eulerian, then it can be decomposed into at most~$\frac{\Delta(G)}{2} +\delta n$ cycles.\label{thm:Delta/2cycledecomp}
	\end{enumerate}
\end{thm}

In particular, the following holds.

\begin{cor}\label{cor:largemindegree}
	For any $\delta,\varepsilon >0$, there exists~$n_0$ such that if~$G$ is a graph on~$n\geq n_0$ vertices with $\delta(G)\geq \frac{n}{2}+\varepsilon n$ then the following hold.
	\begin{enumerate}
		\item $G$ can be decomposed into at most $\max\left\{ \frac{\odd(G)}{2},\frac{\Delta(G)}{2}\right\} +\delta n$ paths.
		\item If~$G$ is Eulerian, then it can be decomposed into at most~$\frac{\Delta(G)}{2} +\delta n$ cycles.\label{cor:largemindegreecycledecomp}
	\end{enumerate}
\end{cor}

Note that, if in addition~$G$ is regular, then the error terms of $\varepsilon n$ and~$\delta n$ can be removed in \cref{cor:largemindegree}\cref{cor:largemindegreecycledecomp}, see \cite{csaba2016proof}.

The next result shows that one can drop the linear minimum degree condition in \cref{thm:Delta/2}\cref{thm:Delta/2pathdecomp} if the quasirandomness covers a larger range of partition class sizes.

\begin{thm}\label{thm:epspathdecomp}
	For any~$p,\delta>0$, there exist~$\varepsilon,n_0>0$ such that the following holds. If~$G$ is a weakly-$(\varepsilon,p)$-quasirandom graph on~$n\geq n_0$ vertices, then~$G$ can be decomposed into at most $\max\left\{ \frac{\odd(G)}{2},\frac{\Delta(G)}{2}\right\} +\delta n$ paths.
\end{thm}

For \cref{thm:n/2}, the linear minimum degree condition is likely to be an artefact of our proof.
On the other hand, in \cref{sec:conclusion}, we will give some examples to show that neither the linear minimum degree condition (or even the stronger assumption of linear connectivity), nor the weakly-$\left(\frac{\alpha}{2}, p\right)$-quasirandom property is sufficient on its own to obtain the bounds in \cref{thm:Delta/2}.
However, \cref{thm:epspathdecomp} shows that, in the case of path decompositions, the linear minimum degree condition can be dropped if we assume~$G$ to be weakly-$(\varepsilon,p)$-quasirandom for a sufficiently small constant~$\varepsilon>0$. 
Surprisingly, it turns out that the Erd\H{o}s-Gallai conjecture is equivalent to the following analogue of \cref{thm:epspathdecomp} for cycle decompositions of Eulerian graphs (see \cref{prop:equivErdosGallai}).
\COMMENT{We will also see that, using \cref{thm:loglognErdosGallai}, we can actually show the analogue of \cref{thm:epspathdecomp} for decomposing weakly-$(\frac{\varepsilon}{\log\log(n)},p)$-quasirandom  Eulerian graphs into cycles, where~$\varepsilon>0$ is a sufficiently small constant.}

\begin{conjecture}\label{conj:epscycledecomp}
	For any~$\delta, p>0$, there exist~$\varepsilon,n_0>0$ such that the following holds. If~$G$ is an Eulerian weakly-$(\varepsilon,p)$-quasirandom graph on~$n\geq n_0$ vertices, then~$G$ can be decomposed into at most $\frac{\Delta(G)}{2} +\delta n$ cycles.
\end{conjecture}

We can prove \cref{conj:epscycledecomp} if weak-$(\varepsilon, p)$-quasirandomness is replaced by weak-$(\frac{\varepsilon}{\log\log n},p)$-quasirandomness (see \cref{prop:loglogn}).

We note that \cref{thm:Delta/2,thm:epspathdecomp} differ from \cref{thm:quasirandom} in the following way. Firstly, we have no restriction on the difference between the maximum and minimum degree. Secondly, weak-$(\varepsilon, p)$-quasirandomness is a significantly weaker property than lower-$(\varepsilon,p)$-regularity. More\-over, the {$\varepsilon$-parameter} in \cref{thm:quasirandom} is much smaller than the {$p$-parameter}. We do not require this in \cref{thm:Delta/2}, and while this is necessary in \cref{thm:epspathdecomp}, there we do not require the minimum degree to be linear. On the other hand, \cref{thm:Delta/2,thm:epspathdecomp} have an additional~$o(n)$ term in the number of paths/cycles compared to \cref{thm:quasirandom}.

Finally, we observe that the following is immediately implied by \cref{cor:largemindegree}.

\begin{cor}\label{cor:boundedDelta}
	For any~$\varepsilon>0$, there exists~$n_0$ such that the conjecture of Haj\'{o}s is true for all Eulerian graphs~$G$ on~$n\geq n_0$ vertices with $\delta(G)\geq \frac{n}{2}+\varepsilon n$ and $\Delta(G)\leq n-\varepsilon n$.
\end{cor}

We remark that by \cref{thm:Delta/2}, \cref{cor:boundedDelta} holds more generally for sufficiently large weakly-quasirandom graphs with maximum degree bounded away from~$n$.

A key tool in our proofs will be the main technical result of \cite{kuhn2013hamilton}, which generates a Hamilton decomposition of a graph satisfying certain robust expansion properties (see \cref{sec:robust} for the statement). This was developed originally in \cite{kuhn2013hamilton} to give a proof of Kelly's conjecture (which states that every large regular tournament has a Hamilton decomposition), and applied e.g.\ in \cite{csaba2016proof} to prove the~$1$-factorisation conjecture (see also \cite{kuhn2014hamilton} for some early applications).

\onlyinsubfile{\bibliographystyle{abbrv}
	\bibliography{Bibliography/papers}}

%% file: Organisation.tex
\onlyinsubfile{
\setcounter{section}{1}
}

\subsection{Organisation of the paper}

The paper is organised as follows. We start by providing a proof overview of our main theorems in \cref{sec:sketch}. Notation and probabilistic tools are introduced in \cref{sec:notation}, and preliminary results are collected in \cref{sec:preliminaries}. Theorems \labelcref{thm:n/2}\cref{thm:n/2pathdecomp}, \labelcref{thm:n/2}\cref{thm:n/2cycledecomp}, \labelcref{thm:Delta/2}, and \labelcref{thm:epspathdecomp} are proved in \cref{sec:proofs}. Finally, we derive \cref{thm:n/2}\cref{thm:ErdosGallai} and make some concluding remarks in \cref{sec:conclusion}.

%% file: Sketch_of_Proof.tex
\onlyinsubfile{
\setcounter{section}{2}
\addtocounter{section}{-1}

\section{Proof overview of the main theorems}}

The proofs of Theorems~\labelcref{thm:n/2}\cref{thm:n/2pathdecomp}, \labelcref{thm:n/2}\cref{thm:n/2cycledecomp}, \labelcref{thm:Delta/2,thm:epspathdecomp} follow a similar strategy, and so, for simplicity, we only sketch the proof of \cref{thm:n/2}\cref{thm:n/2cycledecomp}.

Fix additional constants~$\varepsilon, \zeta, \beta$, and~$n_0$ such that $0<\frac{1}{n_0} \ll \varepsilon\ll \zeta \ll \beta \ll \alpha, \delta \leq 1$. Let~$G$ be a graph on~$n\geq n_0$ vertices with $\delta(G)\geq \alpha n$. We decompose~$G$ by repeatedly constructing cycles. For simplicity, whenever edges are used to form a cycle, they are implicitly deleted from the graph (so all the cycles constructed below are edge-disjoint, as desired). We obtain the bulk of our cycles in \cref{step:sketchexceptional}, all other cycles will contribute to the error term. In \cref{step:sketchexceptional}, we need to be very efficient (i.e.\ the average length of the cycles needs to be large), while there is room to spare in the other steps.

\begin{steps}
	\item \textbf{Applying Szemer\'{e}di's regularity lemma and setting aside \NEW{some} random subgraphs~$\Gamma$ and~$\Gamma'$.}\label{step:sketchcleaning}
	We start by applying Szemer\'{e}di's regularity lemma and a cleaning procedure similar to the one used to prove the degree form of the regularity lemma. 
	We will thus obtain a subgraph~$H\subseteq G$ of small maximum degree and a partition of~$V(G)$ into clusters~$V_1, \dots, V_k$ and an exceptional set~$V_0$. 
	Moreover, in each non-empty pair of clusters of~$G\setminus H$, almost all vertices have degree close to the density of the pair, while the few other vertices are isolated. Moreover, in each pair, the vertices of positive degree span an~$\varepsilon$-regular bipartite graph.
	
	We also set aside two sparse edge-disjoint random spanning subgraphs $\Gamma, \Gamma'\subseteq G\setminus H$ such that, in~$\Gamma$, each non-empty pair of clusters has density close to~$\beta$, while in~$\Gamma'$ each such pair has density close to~$\zeta$.
	By \cref{thm:Lovasz} and by splitting clusters if necessary, we may assume that the reduced graph~$R'$ of~$\Gamma$ can be decomposed into at most~$\frac{|R'|}{2}=\frac{k}{2}$ cycles of even length (this will be needed in \cref{step:sketchGamma}).
	Let $G^*\coloneqq G\setminus (H \cup \Gamma\cup \Gamma')$. Denote by~$G^*_{ij}$ the~$\varepsilon$-regular (almost spanning) subgraph of the pair~$G^*[V_i,V_j]$, and define~$\Gamma_{ij}$ similarly.~$\Gamma$ and~$\Gamma'$ will be used to tie together given sets of paths of~$G^*$ into cycles.
	
	\item \textbf{Covering the edges of~$G[V_0]$.}\label{step:sketchinsideV0}
	Apply \cref{thm:Lovasz} to~$G[V_0]$. The paths obtained are extended to paths with endpoints in~$V(G)\setminus V_0$ and then closed into cycles using edges of~$\Gamma$. Since~$V_0$ is small, this results in only a few cycles and we can use edges of~$\Gamma$ sparingly so that its properties are not destroyed. 
	
	\item \textbf{Covering most of~$G^*$ with at most roughly~$\frac{n}{2}$ cycles.}\label{step:sketchexceptional}
	The idea is to decompose the edges of~$G^*$ into paths and then link some of these paths together using the edges in~$\Gamma\cup \Gamma'$ to form cycles. 
	The bipartite graph~$G^*[V_0, V(G)\setminus V_0]$ is decomposed into paths of length~$2$ with midpoints in~$V_0$, called \emph{exceptional paths}, while~$\varepsilon$-regular pairs~$G_{ij}^*$ are approximately decomposed into long but not spanning paths, so that a few vertices are set aside for tying up paths. 
	We then use edges of~$\Gamma\cup \Gamma'$ to link these paths into cycles.
	More precisely, \NEW{we proceed as follows. Suppose first that the reduced graph $R$ of $G$ is connected.}	
	We construct an auxiliary reduced graph~$\hR$ such that the multiplicity of the edges between~$V_i$ and~$V_j$ in~$\hR$ is proportional to the density of corresponding pair~$G_{ij}^*$ of~$G^*$. We optimally decompose~$\hR$ into matchings.
	Given a matching~$M$ of~$\hR$, we form sets~$\cP$ of paths consisting of exactly one path of~$G^*_{ij}$ for each~$V_iV_j\in M$, and of exceptional paths which cover vertices of~$V_0$ with highest degree.
	Since~$M$ is a matching of clusters and our non-exceptional paths do not span entire clusters, we can ensure that each set~$\cP$ of paths obtained in this way consists of vertex-disjoint paths and does not span entire clusters. Thus, after this step, we still have some uncovered vertices, called \emph{reservoir vertices}, which can be used to link the paths in each set~$\cP$ into a cycle using edges of~$\Gamma\cup \Gamma'$.
	
	Since the edge multiplicity between two clusters in~$\hR$ is proportional to the density of the corresponding pair of~$G^*$ and at each stage we cover exceptional vertices of highest degree, we obtain \NEW{an upper bound of roughly~$\frac{\Delta(G^*)}{2}$}\OLD{at most roughly~$\frac{n}{2}$} cycles in total. \NEW{In general, $R$ may be disconnected and, by construction, $\Gamma\cup \Gamma'$ contains no edges between the different components of $R$. Thus, we cannot tie together paths from different components and we need to apply the above argument to each component of $R$ separately. But, if a component of $R$ contains $n'$ vertices of $G^*$ (say), then the subgraph of $G^*$ induced by this component has maximum degree at most $n'$ and we obtain at most roughly $\frac{n'}{2}$ cycles from that component.
	Thus, we get an upper bound of roughly~$\frac{n}{2}$ cycles in total.}
	
	By alternating which vertices are used as reservoir vertices, we ensure that the leftover graph~$H'$ has small maximum degree. Moreover, we use edges of~$\Gamma$ sparingly so that the properties of~$\Gamma$ are maintained. Since the density~$\zeta$ of~$\Gamma'$ is small, we can add the remaining edges of~$\Gamma'$ to~$H'$ without significantly increasing the maximum degree of~$H'$.
	
	We remark that in \cref{step:sketchinsideV0} it was possible to tie together paths using only~$\Gamma$ because we had some room to spare (in the sense that the number of cycles produced might be fairly large compared to the number of edges covered). But in \cref{step:sketchexceptional}, we need to use edges of both~$\Gamma$ and~$\Gamma'$ in order to be efficient and obtain the desired number of cycles. (The reason that using~$\Gamma\cup \Gamma'$ is more efficient is that the reduced graph of~$\Gamma\cup \Gamma'$ equals that of~$G^*$. We cannot guarantee this property for~$\Gamma$ alone since for \cref{step:sketchbad} the non-empty pairs~$\Gamma_{ij}$ of~$\Gamma$ need to be fairly dense.)
	
	\item \textbf{Covering the leftovers.}\label{step:sketchbad}
	By construction,~$H\cup H'$ has small maximum degree and so can be decomposed into few small matchings. We tie the edges of each matching into a cycle using edges of~$\Gamma$. Once again, we make sure that the relevant properties of~$\Gamma$ are preserved.
	
	\item \textbf{Fully decomposing~$\Gamma$.}\label{step:sketchGamma}
	\begin{figure}[htb]
		\centering
		\begin{subfigure}{0.4\textwidth}
			\centering
			\includegraphics[width=0.8\textwidth]{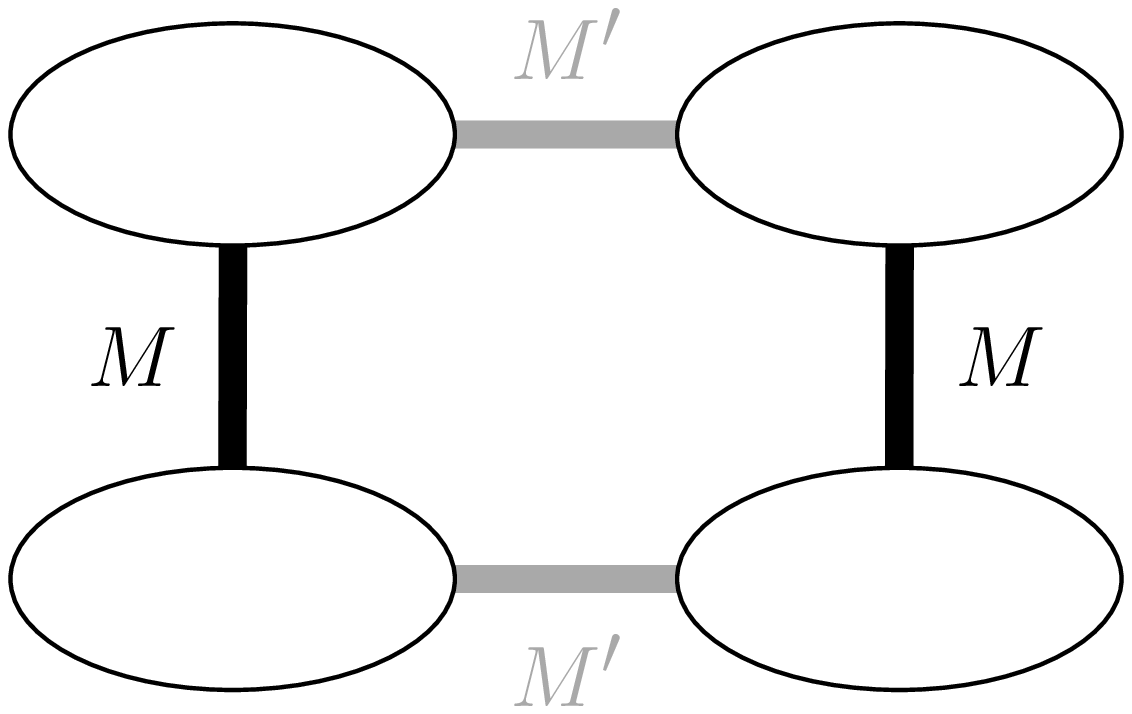}
			\caption{Pair of matchings~$(M, M')$ in the reduced graph. \\   \\}
		\end{subfigure}
		\hspace{0.1\textwidth}
		\begin{subfigure}{0.4\textwidth}
			\centering
			\includegraphics[width=0.8\textwidth]{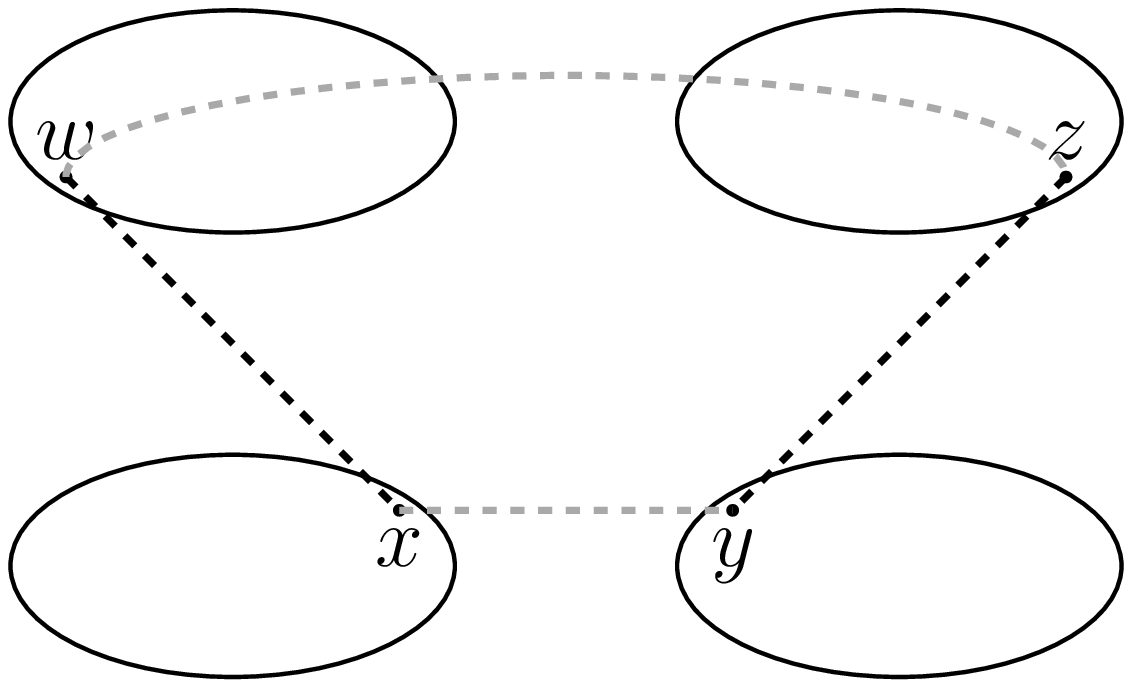}
			\caption{We set aside an edge from each pair in~$M'$ (dashed grey) and replace them by a fictive edge in each pair of~$M$ (dashed black).}
		\end{subfigure}
		
		\vspace{0.05\textwidth}
		
		\begin{subfigure}{0.4\textwidth}
			\centering
			\includegraphics[width=0.8\textwidth]{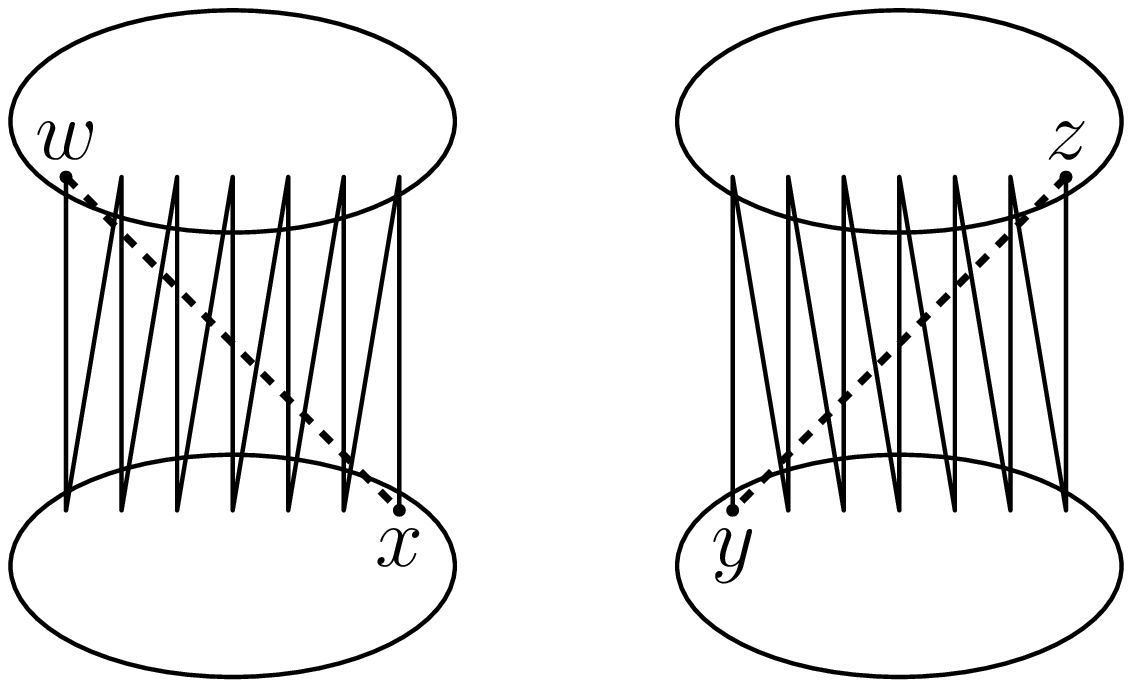}
			\caption{We find a Hamilton cycle of each pair of~$M$ containing a single fictive edge (dashed black).\\}
		\end{subfigure}
		\hspace{0.1\textwidth}
		\begin{subfigure}{0.4\textwidth}
			\centering
			\includegraphics[width=0.8\textwidth]{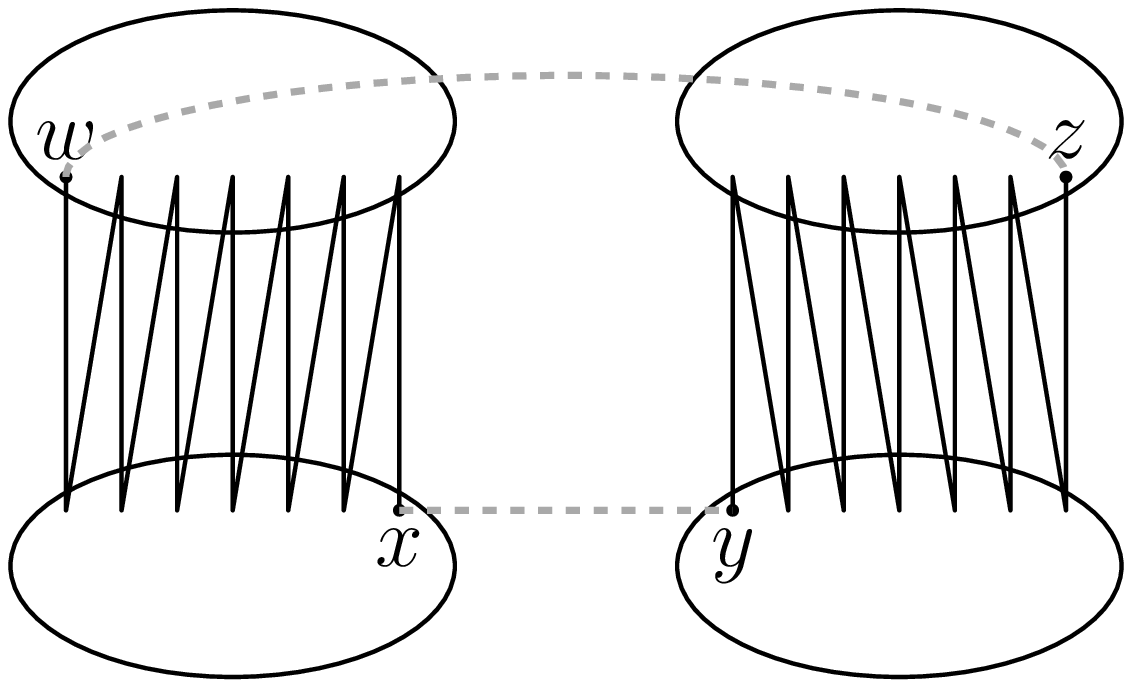}
			\caption{We remove the fictive edges from the cycles of pairs of~$M$ and insert back the edges set aside from pairs of~$M'$ (dashed grey).}
		\end{subfigure}
		\caption{Construction of a cycle of~$\Gamma$.}
		\label{fig:sketch}
	\end{figure}
	It only remains to decompose (the remainder of)~$\Gamma$. The idea is to initially decompose the reduced graph of~$\Gamma$ into~$\frac{k}{2}$ cycles of even length (as discussed in \cref{step:sketchcleaning}). For each such cycle~$C$, the subgraph~$\Gamma_C$ of~$\Gamma$ corresponding to the blow-up of~$C$ is first approximately decomposed into Hamilton cycles of~$\Gamma_C$ that ``wind around"~$C$. The leftover is then decomposed using the main technical result of \cite{kuhn2013hamilton} as follows.	
	
	The cycle~$C$ is initially decomposed into a pair~$(M, M')$ of matchings. For each~$V_iV_j\in M\cup M'$, we first set aside a small set~$\cE_{ij}$ of edges of~$\Gamma_{ij}$ and then decompose the remaining edges into a set~$\cH_{ij}$ of Hamilton paths. We make sure the set of endpoints of the paths in~$\bigcup_{V_iV_j\in M}\cH_{ij}$ equals the set of endpoints of the edges in~$\bigcup_{V_iV_j\in M'}\cE_{ij}$, and similarly for~$M$ and~$M'$ exchanged. Thus we can tie together a path of~$\cH_{ij}$ for each~$V_iV_j\in M$ using exactly one edge of~$\cE_{i'j'}$ for each~$V_{i'}V_{j'}\in M'$. We proceed similarly to tie paths of~$\bigcup_{V_iV_j\in M'}\cH_{ij}$ into cycles. 
	We thus obtain a Hamilton decomposition of~$\Gamma_C$.
	
	In order to prescribe the endpoints of the Hamilton paths, we add some suitable edges to~$\Gamma_C$, called \emph{fictive edges}, and then actually find a Hamilton decomposition of each pair~$\Gamma_{ij}\setminus \cE_{ij}$ such that each cycle in the decomposition contains exactly one fictive edge (see \cref{fig:sketch}).
	Such decompositions are guaranteed by the ``robust decomposition lemma" of \cite{kuhn2013hamilton}. 
	Since by construction all pairs of~$\Gamma$ have density close to~$\beta$, we obtain, in total, about~$\frac{\beta n}{2}\ll \delta n$ cycles.
\end{steps}

\FloatBarrier
\onlyinsubfile{\bibliographystyle{abbrv}
	\bibliography{Bibliography/papers}}

%% file: Notation_and_Probabilistic_Tools.tex
\onlyinsubfile{
\setcounter{section}{3}
\addtocounter{section}{-1}

\section{Notation, definitions, and probabilistic tools}}

\subsection{Notation}

Let~$G$ be a graph. If~$X\subseteq V(G)$ is a set of vertices of~$G$ we write~$G[X]$ for the subgraph of~$G$ induced by~$X$ and~$G-X$ for~$G[V(G)\setminus X]$. Given a set~$F\subseteq E(G)$ of edges of~$G$, we write~$G\setminus F$ for the graph obtained from~$G$ by deleting all edges in~$F$. Similarly, given a subgraph~$H\subseteq G$, we write~$G\setminus H$ for~$G\setminus E(H)$. If~$F$ is a set of non-edges of~$G$, we write~$G\cup F$ for the graph obtained from~$G$ by adding all edges in~$F$. If~$G$ and~$H$ are edge-disjoint graphs we write~$G\cup H$ for the graph with vertex set~$V(G)\cup V(H)$ and edge set~$E(G)\cup E(H)$.

Assume~$G$ is a graph. For any~$x\in V(G)$, we denote by~$N_G(x)$ the set of neighbours of~$x$ and by~$d_G(x)$ the degree of~$x$ in~$G$. Given~$x,y\in V(G)$, we define $d_G(x, y)\coloneqq |N_G(x)\cap N_G(y)|$. The subscripts may be omitted if this is unambiguous. We say~$G$ is \emph{Eulerian} if all its vertices have even degree. (Note that~$G$ is not necessarily connected.) 

Given a graph~$G$ and~$A,B\subseteq V(G)$, we write~$e_G(A,B)$ for the number of edges of~$G$ which have an endpoint in~$A$ and an endpoint~$B$. If~$A,B$ are disjoint then we write~$G[A,B]$ for the bipartite subgraph of~$G$ with vertex classes~$A$ and~$B$ and all edges of~$G$ with an endpoint in~$A$ and an endpoint in~$B$.

Let~$\overrightarrow{G}$ be a digraph. Given vertices~$x,y\in V(\overrightarrow{G})$, we write~$xy$ for the edge directed from~$x$ to~$y$. The vertex~$x$ is called the \emph{initial vertex} of~$xy$ and~$y$ the \emph{final vertex} of~$xy$.
Given a vertex~$x\in V(\overrightarrow{G})$, the \emph{outneighbourhood} of~$x$, denoted~$N_{\overrightarrow{G}}^+(x)$, is the set of vertices~$y$ such that~$xy\in E(\overrightarrow{G})$. Similarly, the \emph{inneighbourhood}~$N_{\overrightarrow{G}}^-(x)$ of a vertex~$x\in V(\overrightarrow{G})$ is the set of vertices~$y$ such that~$yx\in V(\overrightarrow{G})$. We say~$\overrightarrow{G}$ is~$r$-regular if for any vertex~$x\in V(G)$, we have $|N_{\overrightarrow{G}}^+(x)|=|N_{\overrightarrow{G}}^-(x)|=r$.
For any $A,B\subseteq V(\overrightarrow{G})$, we write~$e_{\overrightarrow{G}}(A, B)$ for the number of edges of~$\overrightarrow{G}$ whose initial vertex belongs to~$A$ and whose final vertex belongs to~$B$.
For any disjoint $A,B\subseteq V(\overrightarrow{G})$, we write~$\overrightarrow{G}[A,B]$ for the bipartite subdigraph of~$\overrightarrow{G}$ with vertex classes~$A$ and~$B$ and whose edges are all the edges of~$\overrightarrow{G}$ whose initial vertex belongs to~$A$ and whose final vertex belongs to~$B$.

The \emph{length} of a path is number of edges it contains. An \emph{$(x,y)$-path} is a path whose endpoints are~$x$ and~$y$. Given a path~$P$ and~$x,y\in V(P)$, we write~$xPy$ for the~$(x,y)$-path induced by~$P$. We use the terms \emph{set of vertex-disjoint paths} and \emph{linear forest} interchangeably. In particular, by slightly abusing notation, given a set~$\cP$ of vertex-disjoint paths, we write~$V(\cP)$ for the set of vertices of the paths in~$\cP$ and define~$E(\cP)$ similarly.

We write~$\mathbb{N}$ for the set of natural numbers (including~$0$) and~$\mathbb{N}^*$ for the set of positive natural numbers. For any~$k\in \mathbb{N}^*$, we write~$[k]\coloneqq \{1, 2,\dots, k\}$, $[k]_{\rm odd}\coloneqq \{i\in [k]\mid k\text{ is odd}\}$, and, similarly, $[k]_{\rm even}\coloneqq \{i\in [k]\mid k\text{ is even}\}$.

Let~$a,b,c\in \mathbb{R}$. We write~$a=b\pm c$ if~$b-c\leq a\leq b+c$. For simplicity, we use hierarchies instead of explicitly calculating the values of constants for which statements hold. Namely, if we write~$0<a\ll b\ll c\leq 1$ in a statement, we mean that there exist non-decreasing functions $f\colon (0,1]\longrightarrow (0,1]$ and $g\colon (0,1]\longrightarrow (0,1]$ such that the statement holds for all~$0<a,b,c\leq 1$ satisfying~$b\leq f(c)$ and~$a\leq g(b)$. Hierarchies with more constants are defined in a similar way.
We assume large numbers to be integers and omit floors and ceilings, provided this does not affect the argument. 

Let~$G$ be a graph. A \emph{decomposition of~$G$} is a set~$\cD$ of edge-disjoint subgraphs of~$G$ such each edge of~$G$ belongs to exactly one subgraph in~$\cD$. A \emph{path decomposition} (respectively, \emph{cycle decomposition}) is a decomposition~$\cD$ of~$G$ such that each subgraph in~$\cD$ is a path (respectively, a cycle). We say~$G$ can be \emph{decomposed into~$d$ paths} (respectively, \emph{decomposed into~$d$ cycles}) if~$G$ has a path (respectively, cycle) decomposition~$\cD$ of size~$d$. Similarly, we say~$G$ can be \emph{decomposed into~$d$ paths and cycles} (respectively, \emph{decomposed into~$d$ cycles and edges}) if~$G$ has a decomposition~$\cD$ of size~$d$ such that each subgraph in~$\cD$ is either a path or a cycle (a cycle or an edge, respectively).

\subsection{Regularity}

Let~$G$ be a bipartite graph on vertex classes $A, B$. The \textit{density} of~$G$ is $d_G(A,B) \coloneqq  \frac{e_G(A,B)}{|A| |B|}$.
We may write $d(A,B)$ instead of $d_G(A,B)$ if this is unambiguous. 
For any $\varepsilon>0$, we say~$G$ is \textit{$\varepsilon$-regular} if, for any $A'\subseteq A$ and~$B'\subseteq B$ with~$|A'|\geq \varepsilon |A|$ and~$|B'| \geq \varepsilon |B|$, we have~$|d(A',B') - d(A,B)|< \varepsilon$.

Let~$d\in [0,1]$. 
We say~$G$ is \textit{$(\varepsilon, d)$-regular} if~$G$ is $\varepsilon$-regular and has density~$d$. We write~$G$ is \textit{$(\varepsilon, \geq d)$-regular} is~$G$ \NEW{if} $\varepsilon$-regular of density at least~$d$. 
We say~$G$ is \textit{$[\varepsilon, d]$-superregular} if~$G$ is $\varepsilon$-regular, for all $a\in A$, $d(a)=(d\pm \varepsilon)|B|$, and, for all $b\in B$, $d(b)=(d\pm \varepsilon)|A|$. 
We say~$G$ is \textit{$[\varepsilon, \geq d]$-superregular} if there exists $d'\geq d$ such that~$G$ is $[\varepsilon, d']$-superregular.

We also define a sparse version of~$\varepsilon$-(super)regularity to allow for~$d<\varepsilon$. Let~$G$ be a bipartite graph on vertex classes~$A, B$ of size~$m$. We say~$G$ is \textit{$\{\varepsilon, d\}$-regular} if for any~$A'\subseteq A$ and~$B'\subseteq B$ with~$|A'|,|B'| \geq \varepsilon m$, we have~$d(A',B')= (1\pm \varepsilon) d$. For any~$0<c<1$, we say~$G$ is \textit{$(\varepsilon, d, c)$-regular} if the following hold:
\newcounter{saveenum}
\begin{enumerate} [longlabel, label=(Reg~\arabic*)]
	\item $G$ is~$\{\varepsilon, d\}$-regular; \label{def:reg1}
	\item For any distinct~$a,a'\in A$ we have~$|N(a)\cap N(a')|\leq c^2 m$, and similarly~$|N(b)\cap N(b')|\leq c^2 m$ for any distinct~$b,b'\in B$;\label{def:reg2}
	\item $\Delta(G)\leq cm$. \label{def:reg3}
	\setcounter{saveenum}{\value{enumi}}
\end{enumerate}
For any~$0<d^*<1$, we say that~$G$ is \textit{$(\varepsilon, d, d^*, c)$-superregular} if it is~$(\varepsilon, d, c)$-regular and the following holds:
\begin{enumerate} [longlabel,label=(Reg~\arabic*)]
	\setcounter{enumi}{\value{saveenum}}
	\item $\delta(G)\geq d^* m$. \label{def:reg4}
\end{enumerate}

Given a bipartite digraph~$\overrightarrow{G}$ with vertex classes~$A, B$, recall that~$\overrightarrow{G}[A, B]$ denotes the bipartite subgraph of~$\overrightarrow{G}$ whose edges are all the edges directed from~$A$ to~$B$ in~$\overrightarrow{G}$. We often view~$\overrightarrow{G}[A,B]$ as an undirected bipartite graph. In particular, we say~$\overrightarrow{G}[A,B]$ is~$\varepsilon$-regular if this holds when~$\overrightarrow{G}[A,B]$ is viewed as an undirected graph. We define $(\varepsilon, d)$-regularity, $(\varepsilon, \geq d)$-regularity, $[\varepsilon, d]$-superregularity, and $[\varepsilon, \geq d]$-superregularity for directed bipartite graphs similarly.

Let~$G$ be a graph and~$V_0, V_1, \dots, V_k$ be a partition of~$V(G)$ into~$k$ \emph{clusters}~$V_1, \dots, V_k$ and an \emph{exceptional set}~$V_0$. The vertices in~$V_0$ are called the \emph{exceptional vertices} of~$G$ and an edge of~$G$ is called \emph{exceptional} if it has an endpoint in~$V_0$. The \emph{reduced graph of~$G$} (with respect to the partition~$V_0, V_1, \dots, V_k$) is the graph~$R$ with~$V(R)\coloneqq \{V_1, \dots, V_k\}$ and~$E(R)\coloneqq \{V_iV_j\mid e(G[V_i,V_j])>0\}$. For clarity, we sometimes abuse notation and denote by~$1, \dots, k$ the vertices of~$R$.
If~$C$ is a connected component of~$R$, we let~$V_G(C)\coloneqq\bigcup C$, i.e.~$V_G(C)$ is the set of vertices~$x\in V(G)$ such that~$x\in V_i$ for some~$V_i\in V(C)$.
The \emph{reduced digraph}~$\overrightarrow{R}$ of a digraph~$\overrightarrow{G}$ is defined similarly.

Let~$G$ be an~$n$-vertex graph. Let~$V_0, V_1, \dots, V_k$ be a partition of~$V(G)$ and~$R$ be the corresponding reduced graph. For any distinct~$i,j\in[k]$, the \emph{support cluster of~$V_i$ with respect to~$V_j$} is the set~$V_{ij}\coloneqq \{x\in V_i \mid N_G(x)\cap V_j\neq \emptyset\}$. We also say~$V_{ij}$ and~$V_{ji}$ are the \emph{support clusters of the pair~$G[V_i, V_j]$}. Let~$ij\in E(R)$ and $x\in V_i$. We say~\emph{$x$ belongs to the superregular pair~$G[V_i, V_j]$} if~$x$ belongs to the support cluster~$V_{ij}$.
Let $V_0', V_1', \dots, V_{k'}'$ be a partition of~$V(G)$ such that, for all~$i\in [k']$, there exists~$j\in [k]$ such that $V_i'\subseteq V_j$. We say the support clusters of the partition $V_0', V_1', \dots, V_{k'}'$ are \emph{induced by} the partition $V_0, V_1, \dots, V_k$ if, for all $i',j'\in [k']$, the support cluster~$V_{i'j'}'$ of~$V_{i'}'$ with respect to~$V_{j'}'$ satisfies $V_{i'j'}'=V_{ij}\cap V_{i'}'$, where~$i,j\in [k]$ are such that $V_{i'}'\subseteq V_i$, $V_{j'}'\subseteq V_j$, and $V_{ij}\coloneqq \emptyset$ if~$i=j$.
Let~$G'$ be a graph on~$V(G)$ with reduced graph~$R'$ (with respect to the partition $V_0, V_1, \dots, V_k$). We say~$G$ and~$G'$ \emph{have the same support clusters} if for any~$ij\in E(R)\cap E(R')$, the support clusters of the pairs~$G[V_i, V_j]$ and~$G'[V_i, V_j]$ are the same.

We say~$V_0, V_1, \dots, V_k$ is an \emph{$(\varepsilon, \geq d, k, m, R)$-superregular partition of~$G$} if the following hold.
\begin{enumerate}[longlabel, label=\upshape(\(\text{SRP}\arabic*\))]
	\item $|V_1|=\dots=|V_k|=m$. \label{def:supregpartitionkm}
	\item $|V_0|\leq \varepsilon n$. \label{def:supregpartitionV0}
	\item $G[V_i]$ is empty for all~$i\in [k]$. \label{def:supregpartitioninside}
	\item $R$ is the reduced graph of~$G$.\label{def:supregpartitionreduced}
	\newcounter{supregpartitionsupreg}
	\setcounter{supregpartitionsupreg}{\value{enumi}}
	\item For any~$ij\in E(R)$, let~$V_{ij}, V_{ji}$ be the support clusters of~$G[V_i, V_j]$. Then,~$G[V_{ij}, V_{ji}]$ is~$[\varepsilon, \geq d]$-superregular and~$|V_{ij}|,|V_{ji}|\geq (1-\varepsilon)m$.\label{def:supregpartitionsupreg}	
\end{enumerate}
We say~$V_0, V_1, \dots, V_k$ is an \emph{$(\varepsilon, \geq d, k, m, m', R)$-superregular equalised partition of~$G$} if \cref{def:supregpartitionkm,def:supregpartitionV0,def:supregpartitionreduced,def:supregpartitioninside,def:supregpartitionreduced,def:supregpartitionsupreg} are satisfied and, moreover, the following holds.
\begin{enumerate}[resume,label=\upshape(\(\text{SRP}\arabic*\)), leftmargin=\widthof{SRP1'0000}, align=right,]
	\item $m'\geq (1-\varepsilon)m$ and, for any~$ij\in E(R)$,~$|V_{ij}|=|V_{ji}|=m'$.\label{def:supregpartitionequi}
\end{enumerate}
We say~$V_0, V_1, \dots, V_k$ is an \emph{$(\varepsilon, d, k, m, R)$-superregular partition of~$G$} if \cref{def:supregpartitionkm,def:supregpartitioninside,def:supregpartitionV0,def:supregpartitionreduced,def:supregpartitionsupreg} hold, except that $[\varepsilon, \geq d]$-superregularity is replaced by $[\varepsilon, d]$-superregularity in \cref{def:supregpartitionsupreg}.
We define an \emph{$(\varepsilon, d, k, m, m', R)$-superregular equalised partition of~$G$} analogously. 

We say a graph~$G$ \textit{admits a superregular (equalised) partition} if there exist $V_0, V_1, \dots, V_k$, $\varepsilon, d,$~$k, m, R$ (and~$m'$) such that~$V_0, V_1, \dots, V_k$ is an~$(\varepsilon, \geq d, k, m, R)$-superregular (equalised) partition.

\subsection{Probabilistic estimates}\label{sec:prob}

Let~$X$ be a random variable. We write~$X\sim \Bin(n,p)$ if~$X$ follows  a binomial distribution with parameters~$n, p$. Let~$N,n,m\in \mathbb{N}$ be such that~$\max\{n, m\} \leq N$. Let~$\Gamma$ be a set of size~$N$ and~$\Gamma'\subseteq \Gamma$ be of size~$m$. Recall that~$X$ has a \textit{hypergeometric distribution with parameters~$N, n, m$} if~$X=|\Gamma_n\cap \Gamma'|$, where~$\Gamma_n$ is a random subset of~$\Gamma$ with~$|\Gamma_n|=n$ (i.e.~$\Gamma_n$ is obtained by drawing~$n$ elements of~$\Gamma$ without replacement). We will denote this by~$X\sim \HGeom(N,n,m)$. \COMMENT{Note that if~$X\sim \HGeom(N,n,m)$ then~$\mathbb{E}[X]=\frac{nm}{N}$.}

\par We will use the following Chernoff-type bound.

\begin{lm}[{see e.g.\ \cite[Theorem 2.1 and Theorem 2.10]{janson2011random}}]\label{lm:Chernoff}
	Assume~$X\sim \Bin(n,p)$ or~$X\sim \HGeom(N,n,m)$. Then the following hold for any~$0< \varepsilon < 1$:
	\begin{enumerate}
		\item $\mathbb{P}\left[X\leq (1-\varepsilon)\mathbb{E}[X]\right] \leq \exp \left(-\frac{\varepsilon^2}{3}\mathbb{E}[X]\right)$;
		\item $\mathbb{P}\left[X\geq (1+\varepsilon)\mathbb{E}[X]\right] \leq \exp \left(-\frac{\varepsilon^2}{3}\mathbb{E}[X]\right)$.
	\end{enumerate}
\end{lm}

\onlyinsubfile{\bibliographystyle{plain}
	\bibliography{Bibliography/papers}}

%% file: Regularity.tex
\onlyinsubfile{
\setcounter{section}{4}
\setcounter{subsection}{1}
\addtocounter{subsection}{-1}
}

\subsection{Regularity}\label{sec:regularity}

The first \lcnamecref{lm:subgraph} \NEW{follows} easily from the definition of~$\varepsilon$-regularity.

\begin{lm} \label{lm:subgraph}
	Let~$0<\frac{1}{m}\ll \varepsilon\leq d <1$ and assume~$ \varepsilon\leq \eta \leq 1$. Let~$G$ be a~$(\varepsilon, d)$-regular bipartite graph on vertex classes~$A, B$ of size m. If~$A'\subseteq A$,~$B'\subseteq B$ have size at least~$\eta m$, then~\NEW{$G[A', B']$}\OLD{$G[A'\cup B']$} is~$\frac{\varepsilon}{\eta}$-regular of density~$\geq d-\varepsilon$.
\end{lm}

\COMMENT{\begin{proof}
	This follows by definition of~$\varepsilon$-regularity. Indeed, since~$\eta\geq \varepsilon$, we have~$d(A',B')\geq d-\varepsilon$. Moreover, for~$A''\subseteq A', B''\subseteq B$ with~$|A''|\geq \frac{\varepsilon}{\eta} |A'|\geq \varepsilon m$ and~$|B''|\geq \frac{\varepsilon}{\eta} |B'|\geq \varepsilon m$, we have~$d_{G[A', B']}(A'', B'')=d_G(A'', B'')< d+ \varepsilon\leq d+\frac{\varepsilon}{\eta}$. Similarly,~$d_{G[A',B']}(A'', B'')> d-\frac{\varepsilon}{\eta}$.
\end{proof}}

The following \lcnamecref{lm:verticesedgesremoval} states that~$\varepsilon$-regularity is preserved if only few vertices and edges are removed. This result will be used repeatedly in the rest of the paper.

\begin{lm}[{\cite[Proposition 4.3]{kuhn2013hamilton}}]\label{lm:verticesedgesremoval}
	Let~$0< \frac{1}{m} \ll \varepsilon\leq d' \leq d \leq 1$ and~$G$ be a bipartite graph on vertex classes of size~$m$. Suppose~$G'$ is obtained from~$G$ by removing at most~$d'm$ vertices from each vertex class and at most~$d'm$ edges incident to each vertex from~$G$. 
	\begin{enumerate}
		\item If~$G$ is~$(\varepsilon, d)$-regular, then~$G'$ is~$(2\sqrt{d'}, \geq d-2\sqrt{d'})$-regular. \label{lm:verticesedgesremovalreg}
		\item If~$G$ is~$[\varepsilon, d]$-superregular, then~$G'$ is~$[2\sqrt{d'}, d]$-superregular.
	\end{enumerate}
\end{lm}

An analogous result holds for the sparse version of regularity.

\begin{lm}[{\cite[Proposition 4.8]{kuhn2013hamilton}}]\label{lm:edgesremovalsparse}
	Let~$0< \frac{1}{m} \ll d'\ll \varepsilon, d, d^*, c \leq 1$. Let~$G$ be an~$(\varepsilon, d, d^*, c)$-superregular bipartite graph on vertex classes of size~$m$. Suppose~$G'$ is obtained from~$G$ by removing at most~$d'm$ edges incident to each vertex from~$G$. Then~$G'$ is~$\left(2\varepsilon, d, d^*-d', c\right)$-superregular.	
\end{lm}
\COMMENT{This is actually a special case of \cite[Proposition 4.8]{kuhn2013hamilton}.}
\COMMENT{\begin{proof}
	Clearly \ref{def:reg2} and \ref{def:reg3} are satisfied. Also~$\delta(G')\geq \delta(G)-d'm\geq (d^*-d')m$ so \ref{def:reg4} holds. To show \ref{def:reg1}, let~$A'\subseteq A$ and~$B'\subseteq B$ be such that~$|A'|, |B'|\geq 2\varepsilon m\geq \varepsilon m$. Then we have~$d_{G'}(A', B')\leq d_G(A', B')\leq (1+\varepsilon)d\leq (1+2\varepsilon)d$. Moreover,
	\begin{equation*}
	\begin{aligned}
	d_{G'}(A', B')&\geq d_G(A', B')-\frac{|A'| d' m}{|A'||B'|}\\
	&\geq (1-\varepsilon)d-\frac{d'}{\varepsilon}\\
	&\geq (1-2\varepsilon)d.
	\end{aligned}
	\end{equation*}
\end{proof}}

\NEW{The following \lcnamecref{prop:supreg4} states an $\varepsilon$-regular bipartite graph of linear minimum degree has small diameter.
\begin{prop}\label{prop:supreg4}
	Let $0<\frac{1}{m_A}, \frac{1}{m_B}\ll \varepsilon \leq d\leq 1$. Let $G$ be an $\varepsilon$-regular bipartite graph on vertex classes $A$ and $B$ of size $m_A$ and $m_B$.
	Suppose that each $x\in A$ satisfies $d_G(x)\geq dm_B$ and each $y\in B$ satisfies $d_G(y)\geq dm_A$.
	Then, for any $x\in A$ and $y\in B$, $G$ contains an $(x,y)$-path of length at most $3$. In particular, for any distinct $x,y\in V(G)$, $G$ contains an $(x, y)$-path of length at most $4$.
\end{prop}}

\COMMENT{\begin{proof}
		Suppose first that $x\in A$ and $y\in B$. Then, $N_G(x)\subseteq B$ and $N_G(y)\subseteq A$. Moreover, $|N_G(x)|\geq dm_B\geq \varepsilon m_B$ and, similarly, $|N_G(y)|\geq \varepsilon m_A$. Therefore, there exist $x'\in N_G(x)$ and $y'\in N_G(y)$ such that $x'y'\in E(G)$. Thus, $xx'y'y$ is an $(x,y)$-path of length at most $3$.
		If both $x,y\in B$, let $x'\in N_G(x)\subseteq A$ and find an $(x', y)$-path $P$ of length $3$. Then, $xx'Py$ is an $(x, y)$-path $P$ of length $4$.
	\end{proof}}

The following \lcnamecref{lm:regHamcycle} states that balanced~$\varepsilon$-regular bipartite graphs of large minimum degree are Hamiltonian.

\begin{lm}[{see for instance \cite[Lemma 3.3]{glock2016optimal}}]\label{lm:regHamcycle}
	Let~$0< \frac{1}{m} \ll \varepsilon \ll \alpha \leq 1$. If~$G$ is an~$\varepsilon$-regular bipartite graph on vertex classes of size~$m$ such that~$\delta(G)\geq \alpha m$, then~$G$ contains a Hamilton cycle. 
\end{lm}

\begin{cor}\label{lm:regperfectmatching}
	Let~$0< \frac{1}{m} \ll \varepsilon \ll \alpha \leq 1$. If~$G$ is an~$\varepsilon$-regular bipartite graph on vertex classes of size~$m$ such that~$\delta(G)\geq \alpha m$, then~$G$ contains a perfect matching.
\end{cor}

The next \lcnamecref{lm:sparsesubgraph} states that any superregular pair contains a sparse superregular pair as a subgraph.

\begin{lm}[{\cite[Lemma 4.10]{kuhn2013hamilton}}]\label{lm:sparsesubgraph}
	Let~$0< \frac{1}{m} \ll \varepsilon, d' \leq d \leq 1$ and suppose~$\varepsilon \ll d$. Let~$G$ be an~$[\varepsilon, d]$-superregular bipartite graph on vertex classes of size~$m$. Then~$G$ contains an $(\varepsilon^{\frac{1}{12}}, d', \frac{d'}{2}, \frac{3d'}{2d})$-superregular spanning subgraph.
\end{lm}

By considering a random partition of the edges, one can show that the edges of an~$\varepsilon$-regular pair can be partitioned without destroying the~$\varepsilon$-regularity (see e.g.\ the proof of \cite[Lemma 4.10]{kuhn2013hamilton}). 

\begin{lm}[Partitioning the edges of a regular pair]\label{lm:edgepartition}
	Assume~$0<\frac{1}{m}\ll \varepsilon\ll d_1, \dots, d_\ell\leq  d\leq 1$ with~$\sum_{i=1}^\ell d_i \leq d$. Let~$G$ be a bipartite graph on vertex classes~$A,B$ of size~$m$. Then~$G$ can be decomposed into edge-disjoint spanning subgraphs~$G_0,G_1, \dots, G_\ell\subseteq G$ such that~$G_0$ is empty if $\sum_{i=1}^\ell d_i =d$, and the following hold for each~$i\in[\ell]$.
	\begin{enumerate}
		\item If~$G$ is~$(\varepsilon, d)$-regular, then~$G_i$ is~$(\varepsilon^{\frac{1}{12}}, d_i\pm \varepsilon^{\frac{1}{12}})$-regular.
		\item If~$G$ is~$[\varepsilon, d]$-superregular, then~$G_i$ is~$[\varepsilon^{\frac{1}{12}}, d_i]$-superregular.
	\end{enumerate}
\end{lm}

\COMMENT{\begin{proof}
	Randomly partition the edges of~$G$ to form subgraphs~$G_1, \dots, G_\ell$, where each edge is assigned to~$G_i$ with probability~$p_i\coloneqq \frac{d_i}{d}$. Let~$i\in[\ell]$, by the same arguments as in the proof of \cite[Lemma 4.10]{kuhn2013hamilton}, the following hold with high probability:
	\begin{itemize}
		\item if~$G$ is~$(\varepsilon, d)$-regular, then~$G_i$ is~$(\varepsilon^{\frac{1}{12}}, d_i\pm \varepsilon^{\frac{1}{12}})$-regular;
		\item if~$G$ is~$[\varepsilon, d]$-superregular, then~$G_i$ is~$[\varepsilon^{\frac{1}{12}}, d_i]$-superregular.
	\end{itemize}
	So a union bound over~$i\in[\ell]$ gives that~$G_1, \dots, G_\ell$ all satisfy the \lcnamecref{lm:edgepartition} with positive probability.
\end{proof}}

\begin{cor}\label{cor:supregdigraph}
	Suppose $0<\frac{1}{m}\ll \varepsilon \ll d\leq 1$. Let~$G$ be an~$[\varepsilon,d]$-superregular bipartite graph on vertex classes~$A$ and~$B$ of size~$m$. Then, there exists an orientation~$\overrightarrow{G}$ of the edges of~$G$ such that both~$\overrightarrow{G}[A,B]$ and~$\overrightarrow{G}[B,A]$ are $[\varepsilon^{\frac{1}{12}},\frac{d}{2}]$-superregular.
\end{cor}

Using \cref{lm:Chernoff} and a result of \cite{alon1994algorithmic} which characterises~$\varepsilon$-regularity in terms of co-degree%
\COMMENT{Let~$0< \frac{1}{m_A}, \frac{1}{m_B} < \varepsilon \leq d \leq 1$ and~$G$ be an~$(\varepsilon, d)$-regular bipartite graph on vertex classes~$A$ and~$B$ of size~$m_A$ and~$m_B$, respectively. Then, fewer than~$\varepsilon m_A$ vertices in~$A$ have degree at least~$(d+\varepsilon)m_B$ and fewer than~$\varepsilon m_A$ vertices in~$A$ have degree at most~$(d-\varepsilon)m_B$.\label{comment:badvertices}
%This lemma is in the appendix
}%
\COMMENT{\label{lm:badpairs}Let~$0< \frac{1}{m_A}, \frac{1}{m_B} < \varepsilon \ll d \leq 1$ and~$G$ be an~$(\varepsilon, d)$-regular bipartite graph on vertex classes~$A$ and~$B$ of size~$m$. Then, fewer than~$2\varepsilon m^2$ pairs of distinct vertices~$a, a'\in A$ satisfy~$|N(a) \cap N(a')|\geq (d +\varepsilon)^2m_B$ and fewer than~$2\varepsilon m^2$ pairs of distinct vertices~$a, a'\in A$ satisfy~$|N(a) \cap N(a')|\leq (d-\varepsilon^2)m_B$.
\begin{proof}
	We show that fewer than~$2\varepsilon m_A^2$ pairs of distinct vertices~$a, a'\in A$ satisfy~$|N(a) \cap N(a')|\geq (d +\varepsilon)^2m_B$. 
	Let~$X\coloneqq \{(a,a')\in A^2\mid d(a,a')\geq (d+\varepsilon)^2m_B\}$. We will show that~$|X| < 2\varepsilon m_A^2$. Assume~$(a,a')\in X$. Then one of the following hold:
	\begin{enumerate}[label=(\roman*)]
		\item $d(a)\geq (d+\varepsilon)m_B$;\label{lm:badpaira}
		\item $\varepsilon m_B\leq d(a) < (d+\varepsilon)m_B$ and~$d(a, a')> (d+\varepsilon)d(a)$.\label{lm:badpairb}
	\end{enumerate}
	By \cref{comment:badvertices}, there are fewer than~$\varepsilon m_A \cdot m_A$ pairs satisfying \cref{lm:badpaira}. Let~$a\in A$ be such that~$\varepsilon m_B\leq d(a) < (d+\varepsilon)m_B$. Let~$A'\coloneqq \{a'\in A\mid d(a, a')> (d+\varepsilon)d(a)\}$. Then we have~$d(A', N(a))>d+\varepsilon$. Since~$|N(a)|\geq \varepsilon m_B$, we must have~$|A'|<\varepsilon m_A$. Thus there are less than~$m_A\cdot \varepsilon m_A$ pairs satisfying \cref{lm:badpairb}. Thus we have~$|X|<2\varepsilon m_A^2$, as desired. The second part of the statements holds by similar arguments.
\end{proof}}%
,
one can show that the vertex classes of superregular pairs can partitioned into superregular subpairs.

\begin{lm}[Partitioning the vertices of a regular pair]\label{lm:vertexpartition}
	Assume~$0<\frac{1}{m}\ll \varepsilon \ll d \leq 1$ and~$\frac{1}{m}\ll \frac{1}{r}, \eta$. Let~$G$ be a bipartite graph on vertex classes~$A$ and~$B$ of size~$m$. Let~$A'$ and~$B'$ be the support clusters of~$A$ and~$B$, respectively. Assume that $m'\coloneqq |A'|=|B'|\geq (1-\varepsilon)m$ and~$G[A', B']$ is~$[\varepsilon, d]$-superregular.
	Let $m_1, \dots, m_r\in \mathbb{N}^*$ be such that $\sum_{i\in [r]}m_i=m$ and, for each $i\in [r]$, $m_i=\frac{m}{r}\pm 1$.
	Assume~$A$ and~$B$ are randomly partitioned into~$r$ subsets~$A_1, \dots, A_r$ and~$B_1, \dots, B_r$ such that for all $i\in [r]$, $|A_i|=|B_i|=m_i$.
	Then, with high probability, all of the following hold.
	\begin{enumerate}
		\item For any~$i\in [r]$, we have $|A'\cap A_i|=(1\pm \varepsilon)\frac{m'}{r}$ and similarly $|B'\cap B_i|=(1\pm \varepsilon)\frac{m'}{r}$.
		\item $G[X,Y]$ is $[\varepsilon^{\frac{1}{7}}, d]$-superregular for any $X\in \{A'\cap A_i, A'\setminus A_i \mid i\in [r]\}$ and $Y\in \{B'\cap B_i, B'\setminus B_i \mid i\in [r]\}$.
		\item For any $X\in \{A'\cap A_i \mid i\in [r]\}$ and $Y\in \{B'\cap B_i \mid i\in [r]\}$, $\Delta(G[X,Y])\leq \frac{\Delta(G)+\eta m'}{r}$ and $\delta(G[X,Y])\geq \frac{\delta(G)-\eta m'}{r}$.\label{lm:vertexpartitiondegree}	\end{enumerate}
\end{lm}

Finally, the following simple fact will be needed in \cref{sec:Gamma}.

\begin{prop}\label{lm:regcycle}
	Suppose $0\leq \frac{1}{m}\ll \varepsilon\ll d\leq 1$ and~$k\in \mathbb{N}^*$. Let~$G$ be a graph and~$V_1, \dots, V_k$ be a partition of~$V(G)$ into~$k$ clusters of size~$m$. Let~$R$ be the corresponding reduced graph of~$G$ and assume that for each~$ij\in E(R)$, the pair~$G[V_i, V_j]$ is $[\varepsilon, \geq d]$-superregular. If~$R$ is a cycle of length~$k$, then~$G$ contains~$\varepsilon m$ vertex-disjoint cycles of length~$k$ which intersect each of the clusters~$V_1, \dots, V_k$.
\end{prop}

\COMMENT{\begin{proof}
	Assume without loss of generality that $E(R)=\{V_iV_{i+1}\mid i\in [k]\}$, where~$V_{k+1}\coloneqq V_1$. Assume inductively that for some $0\leq t<\varepsilon m$,~$C_1, \dots, C_t$ are~$t$ vertex-disjoint cycles of~$G$ satisfying the desired properties. Select a vertex $x_1\in V_1\setminus V(C_1\cup \dots \cup C_t)$. Since~$G[V_1, V_2]$ is $[\varepsilon, \geq d]$-superregular and~$t\leq \varepsilon m$, we can find a neighbour $x_2\in V_2\setminus V(C_1\cup \dots \cup C_t)$ of~$x_1$. Similarly, we can find $x_i\in  V_i\setminus V(C_1\cup \dots \cup C_t)$ for each~$i=3, \dots, k-2$ such that~$x_1x_2 \dots x_{k-2}$ is a path in $G-V(C_1\cup \dots \cup C_t)$. Let $S\coloneqq N_{G[V_{k-2},V_{k-1}]}(x_{k-2})\setminus  V(C_1\cup \dots\cup C_t)$ and $T\coloneqq N_{G[V_1,V_k]}(x_1)\setminus  V(C_1\cup \dots \cup C_t)$. Then $|S|, |T|\geq \varepsilon m$ so there exist~$x_{k-1}\in S$ and~$x_k\in T$ which are adjacent in~$G[V_{k-1}, V_{k}]$. Let $C_{t+1}\coloneqq x_1\dots x_k$.
\end{proof} }

\onlyinsubfile{\bibliographystyle{plain}
	\bibliography{Bibliography/papers}}

%% file: Tying_Paths.tex
\onlyinsubfile{
\setcounter{section}{4}
\setcounter{subsection}{2}
\addtocounter{subsection}{-1}
\setcounter{definition}{10}
}

\subsection{Tying paths together}\label{sec:tyingpaths}
	
	Throughout the proof of our main theorems, we will form linear forests and aim to tie together some of the paths in each forest to form cycles. This section gathers several tools to achieve this.
	\cref{lm:tyingmanypaths} will be used to efficiently reduce the number of components of linear forests (i.e.\ to merge paths), from a linear number of components to bounded number, while
	\cref{lm:tyingfewpaths} will be used to further reduce the number of components, from a large constant to a smaller one.
	\cref{lm:closingcycles} will be used to turn linear forests with few components into small sets of vertex-disjoint cycles. Finally, we will use \cref{lm:tyingfewpathscycle,lm:closingcycle,lm:matchingtying} to turn small linear forests into a cycle each.
	
	Let~$\Gamma$ be a graph and~$P_1, \dots, P_\ell$ be vertex-disjoint paths with endpoints in~$V(\Gamma)$. By \emph{tying the paths~$P_1, \dots, P_\ell$ together into a path~$P$ (a cycle~$C$) using the edges of~$\Gamma$}, we mean forming a path~$P$ (a cycle~$C$) such that for each~$i\in [\ell]$, the path~$P_i$ is a subpath of~$P$ (of~$C$), the other edges of~$P$ (of~$C$) are edges of~$\Gamma$ and the endpoints of~$P$ are in~$\bigcup_{i\in [\ell]} V(P_i)$. 
	A subpath~$P'$ of~$P$ (of~$C$) is called a \emph{link path} if~$E(P')\cap E(P_i)=\emptyset$ for each~$i\in [\ell]$ and the endpoints of~$P'$ are in~$\bigcup_{i\in [\ell]} V(P_i)$. 
	In particular, we say~$P'$ \emph{links}~$P_i$ and~$P_j$ if the endpoints of~$P'$ are an endpoint of~$P_i$ and an endpoint of~$P_j$.
	Moreover, if~$A, B$ are distinct clusters, we say~$P'$ is an \emph{$(A,B)$-link path} if~$E(P')\subseteq E(\Gamma[A,B])$ and both endpoints of~$P'$ belong to~$A$.

	The idea behind the next \lcnamecref{lm:tyingmanypaths} is to iteratively tie two paths which have an endpoint in a common cluster using a single superregular pair of~$\Gamma$.
	
	\begin{lm}\label{lm:tyingmanypaths}
		Suppose $0<\frac{1}{n}\ll \frac{1}{k}\ll \varepsilon\leq\zeta\ll \beta \leq 1$.
		Let~$\Gamma$ be a graph on vertex set~$V$ of size~$n$ such that the following hold.
		\begin{enumerate}
			\item $V_0, V_1, \dots, V_k$ is an~$(\varepsilon, \beta, k, m, R)$-superregular partition of~$\Gamma$.\label{lm:tyingmanypathssupreg}
			\item Any~$x\in V\setminus V_0$ belongs to at least~$\beta k$ superregular pairs of~$\Gamma$. \label{lm:tyingmanypathsbadinfewpairs}
		\end{enumerate}
		Let~$\cP_1, \dots, \cP_\ell$ be sets of paths\COMMENT{Note that the internal vertices of~$\cP_1, \dots, \cP_\ell$ are not necessarily in~$V$} satisfying the following.
		\begin{enumerate}[resume]
			%\item $\cP_1\cup \dots\cup \cP_\ell$ is edge-disjoint from~$\Gamma$.\label{lm:tyingmanypathsedgedisjoint}
			\item For each~$i\in[\ell]$,~$\cP_i$ is a set of vertex-disjoint paths with endpoints in~$V\setminus V_0$.\label{lm:tyingmanypathsvertexdisjoint}
			\item For each~$i\in[\ell]$ and~$j\in[k]$,~$|V(\cP_i)\cap V_j|\leq \zeta m$. In particular, $|\cP_i|\leq \zeta n$. \label{lm:tyingmanypathsroom}
			\newcounter{tyinglemma0}
			\setcounter{tyinglemma0}{\value{enumi}}
			\item For any~$x\in V$, there are at most~$\varepsilon n$ paths in~$\cP_1\cup \dots\cup \cP_\ell$ which have~$x$ as an endpoint. \label{lm:tyingmanypathsendpoint}
			\newcounter{tyinglemma1}
			\setcounter{tyinglemma1}{\value{enumi}}
		\end{enumerate}	
		Then, there exist disjoint $E_1, \dots, E_\ell \subseteq E(\Gamma)$ such that the following hold.
		\begin{enumerate}[label=\upshape(\alph*)]
			\item For any~$i\in [\ell]$, by using each edge in~$E_i$ exactly once, we can tie together some of the paths in~$\cP_i$ to form a set~$\cQ_i$ of vertex-disjoint paths such that, for any~$j\in [k]$, at most~$2\beta^{-2}$ paths in~$\cQ_i$ have an endpoint in~$V_j$.\label{lm:tyingmanypathsdecomp}
			
			\item For any distinct~$i,j\in [k]$ and~$x\in V_i$,~$E_1\cup\dots\cup E_\ell$ contains at most~$3\varepsilon^{\frac{1}{4}}m$ edges of~$\Gamma[V_i, V_j]$ which are incident to~$x$.\label{lm:tyingmanypathsneighbours}

			\item For any~$i\in [\ell]$ and~$j\in[k]$, $|V(\cP_i\cup E_i)\cap V_j|\leq \sqrt{\zeta} m$.\label{lm:tyingmanypathsroomQ}
		\end{enumerate}
	\end{lm}

	To prove \cref{lm:tyingmanypaths}, we will use edges of~$\Gamma$ to tie together some of the paths in~$\cP_i$, for each~$i\in[\ell]$. We will only tie together paths which have an endpoint in a common cluster and use a single superregular pair of~$\Gamma$ to do so.
	
	\begin{proof}
		Let $E_1, \dots, E_\ell\subseteq E(\Gamma)$ be (possibly empty) disjoint sets of edges of~$\Gamma$ and assume inductively that for each~$i\in [\ell]$, by using each edge in~$E_i$ exactly once, we can tie together some of the paths in~$\cP_i$ to form a set~$\cQ_i$ of vertex-disjoint paths such that the following is satisfied.
		\begin{enumerate}[label=(\arabic*)]
			\item If $P\in \cQ_1\cup\dots\cup \cQ_\ell$ and~$P'$ is a link path of~$P$, then~$P'$ is an~$(A,B)$-link path of length at most~$4$, for some clusters~$A,B$.\label{lm:tyingmanypathslinkpath}
			\item For any clusters~$A, B$ and any~$x\in A$, there are at most~$\varepsilon^{\frac{1}{2}}m$~$(A, B)$-link paths in~$\cQ_1\cup\dots\cup \cQ_\ell$ which have~$x$ as an endpoint.\label{lm:tyingmanypathslinkedge}
			\item For any clusters~$A, B$ and any~$x\in A\cup B$,~$\cQ_1\cup\dots\cup \cQ_\ell$ contains at most~$\varepsilon^{\frac{1}{4}}m$~$(A,B)$-link paths which have~$x$ as an internal vertex.\label{lm:tyingmanypathslinkvertex}
			\item For any cluster~$A$ and~$i\in [\ell]$, there are at most~$\frac{\sqrt{\zeta} m}{4}$ tuples~$(B,P)$ such that~$B\neq A$ is a cluster and~$P$ is a~$(B,A)$-link path in~$\cQ_i$.\label{lm:tyingmanypathsroomcluster}
		\end{enumerate}
		(In \cref{lm:tyingmanypathslinkedge,lm:tyingmanypathslinkvertex,lm:tyingmanypathsroomcluster} and below, by a link path in~$\cQ_i$, we mean a link path of some path in~$\cQ_i$.)
		
		If for any~$i\in [\ell]$ and~$j\in[k]$, the set~$\cQ_i$ contains at most~$2\beta^{-2}$ paths with an endpoint in~$V_j$, then \cref{lm:tyingmanypathsdecomp} holds. Moreover, \cref{lm:tyingmanypathslinkedge} and \cref{lm:tyingmanypathslinkvertex} imply \cref{lm:tyingmanypathsneighbours}, while \cref{lm:tyingmanypathsroomQ} follows from \cref{lm:tyingmanypathslinkpath,lm:tyingmanypathsroomcluster,lm:tyingmanypathsroom}%
		\COMMENT{for any~$i\in [\ell]$ and~$j\in[k]$,~$\cP_i$ originally contained at most~$\zeta m$ vertices of~$V_j$, then, by \cref{lm:tyingmanypathslinkpath}, we add at most one vertex when we tie two paths with endpoints in~$V_j$ (this occurs at most~$\frac{\zeta m}{2}$ times by \cref{lm:tyingmanypathsroom}) and at most two vertices when we use~$V_j$ to tie paths with endpoints in~$V(\Gamma)\setminus V_j$ so by \cref{lm:tyingmanypathsroomcluster} $|V(\cP_i)\cap V_j|\leq \zeta m+\frac{\zeta m}{2}+2\frac{\sqrt{\zeta}m}{4} \leq\sqrt{\zeta} m$},
		and we are done.
		
		We may therefore assume that there exist~$i\in [\ell]$ and~$j\in [k]$ such that~$\cQ_i$ contains more than~$2\beta^{-2}$ paths with an endpoint in~$V_j$. Then, we claim that there exist distinct~$P, P'\in \cQ_i$, each with an endpoint in~$V_j$, such that the following hold. There exists~$j'\in [k]$ such that
		\begin{enumerate}[label=(\Roman*)]
			\item $x\in V_j$ is an endpoint of~$P$ and~$x'\in V_j$ is an endpoint of~$P'$;\label{lm:tyingmanypathsx}
			\item $x, x'\in V_{jj'}$, where~$V_{jj'}$ is the support cluster of~$V_j$ with respect to~$V_{j'}$; \label{lm:tyingmanypathsxintV}
			\item $\cQ_1\cup\dots\cup \cQ_\ell$ contains fewer than~$\varepsilon^{\frac{1}{2}}m$~$(V_j, V_{j'})$-link paths which have~$x$ as an endpoint, and similarly for~$x'$;\label{lm:tyingmanypathsxlinkedge}
			\item there are fewer that~$\frac{\sqrt{\zeta} m}{4}$ tuples~$(A, Q)$ such that~$A\neq V_{j'}$ is a cluster and~$Q$ is an~$(A,V_{j'})$-link path in~$\cQ_i$.\label{lm:tyingmanypathsxroomcluster}
		\end{enumerate}
		Indeed, for any~$x\in V_j$, there are at least~$\beta k$ indices~$j'\in [k]$ such that~$x\in V_{jj'}$ (by \cref{lm:tyingmanypathsbadinfewpairs}). By \cref{lm:tyingmanypathsendpoint}, there are at most~$2\varepsilon^{\frac{1}{2}}k$ such indices~$j'\in[k]$ such that~$\cQ_1\cup\dots\cup \cQ_\ell$ contains~$\varepsilon^{\frac{1}{2}}m$~$(V_j, V_{j'})$-link paths with~$x$ as an endpoint\COMMENT{At most $\frac{\varepsilon n}{\varepsilon^{\frac{1}{2}}m}\leq 2\varepsilon^{\frac{1}{2}}k$.}.
		Moreover, by \cref{lm:tyingmanypathsroom}, there are at most%
		\COMMENT{$\frac{\zeta n}{\frac{\sqrt{\zeta} m}{4}}\leq 5\sqrt{\zeta} k$}%
		~$5\sqrt{\zeta} k$ indices~$j'\in[k]$ such that there exist~$\frac{\sqrt{\zeta} m}{4}$ tuples~$(A, Q)$ where~$A\neq V_{j'}$ is a cluster and~$Q$ is an~$(A,V_{j'})$-link path in~$\cQ_i$.
		Thus, for any~$x\in V_j$, there are at least~$\frac{\beta^2 k}{2}$ indices~$j'$ such that~$x$ satisfies \cref{lm:tyingmanypathsxintV,lm:tyingmanypathsxlinkedge,lm:tyingmanypathsxroomcluster}. Therefore, since by assumption~$\cQ_i$ contains more than~$2\beta^{-2}$ paths with an endpoint in~$V_j$, we can find~$P,P', x, x'$, and~$j'$ satisfying \cref{lm:tyingmanypathsx,lm:tyingmanypathsxintV,lm:tyingmanypathsxlinkedge,lm:tyingmanypathsxroomcluster}. 
		
		We can now find an~$(x,x')$-path in~$\Gamma[V_{jj'}, V_{j'j}]$ to tie~$P$ and~$P'$ together as follows.		
	 	Let $\Gamma'\coloneqq\Gamma\setminus (E_1\cup\dots\cup E_\ell)$. By \cref{lm:tyingmanypathslinkedge,lm:tyingmanypathslinkvertex}, \cref{lm:verticesedgesremoval} implies that~$\Gamma'[V_{jj'}, V_{j'j}]$ is still~$[\varepsilon^{\frac{1}{9}}, \beta]$-superregular%
		\COMMENT{By \cref{lm:tyingmanypathslinkedge,lm:tyingmanypathslinkvertex}, we have deleted at most~$3\varepsilon^{\frac{1}{4}}m$ edges incident to each vertex in~$\Gamma[V_{jj'}, V_{j'j}]$.}. 
		Let~$V_j'$ be obtained from~$V_{jj'}$ by deleting the following vertices:
		\begin{itemize}
			\item vertices in~$V(\cQ_i)\setminus\{x,x'\}$ (by \cref{lm:tyingmanypathsroom,lm:tyingmanypathslinkpath,lm:tyingmanypathsroomcluster}, there are at most~$\sqrt{\zeta} m$ such vertices\COMMENT{at most~$\zeta m$ original vertices, plus one vertex when we use~$(V_j,V_{j''})$ (for some~$j''$) to tie two paths with endpoints in~$V_j$, plus two vertices when we use~$(V_j,V_{j''})$ to tie two paths with endpoints in~$V_{j''}$. So at most $\zeta m+\frac{\zeta m}{2}+\frac{\sqrt{\zeta} m}{2}$.});
			\item vertices in $V_{jj'}\setminus \{x, x'\}$ which are an internal vertex of~$\varepsilon^{\frac{1}{4}}m$~$(V_j, V_{j'})$-link paths of~$\cQ_1\cup\dots\cup \cQ_\ell$ (by \cref{lm:tyingmanypathslinkpath,lm:tyingmanypathslinkedge}, there are at most~$\varepsilon^{\frac{1}{4}}m$ such vertices\COMMENT{by \cref{lm:tyingmanypathslinkedge}, there are at most~$\frac{\varepsilon^{\frac{1}{2}}m^2}{2}$~$(V_j, V_{j'})$-link paths in total. By \cref{lm:tyingmanypathslinkpath}, these link paths have at most~$2$ internal vertices in~$V_{j'}$ and at most~$1$ internal vertex in~$V_j$.}).
		\end{itemize}
		Note that the number of deleted vertices is $|V_{jj'}\setminus V_j'|\leq \sqrt{\zeta} m+\varepsilon^{\frac{1}{4}}m\leq 2\sqrt{\zeta}|V_{jj'}|$. Define~$V_{j'}'$ similarly. Then, by \cref{lm:verticesedgesremoval}, $\Gamma'[V_j', V_{j'}']$ is~$[\zeta^{\frac{1}{5}}, \beta]$-superregular. Thus, \NEW{by \cref{prop:supreg4}, $\Gamma'[V_j', V_{j'}']$} contains an~$(x,x')$-path~$P''$ of length at most~$4$. Add the edges of~$P''$ to~$E_i$ and  replace in~$\cQ_i$ the paths~$P$ and~$P'$ by the concatenation of~$P, P''$, and~$P'$. By construction, \cref{lm:tyingmanypathslinkpath,lm:tyingmanypathslinkedge,lm:tyingmanypathslinkvertex,lm:tyingmanypathsroomcluster} are still satisfied, as desired for the induction step.
	\end{proof}
	
	After applying \cref{lm:tyingmanypaths}, we obtain linear forests with few components. One can then be less economical and use several superregular pairs of~$\Gamma$ to tie paths together. This is achieved in the next \lcnamecref{lm:tyingfewpaths}.

	\begin{lm}\label{lm:tyingfewpaths}
		Suppose $0<\frac{1}{n}\ll \frac{1}{k}\ll \varepsilon\leq \zeta\ll \beta \leq 1$. Let~$\Gamma$ be a graph on vertex set~$V$ of size~$n$ and~$\cP_1, \dots, \cP_\ell$ be sets of paths. Assume~$\Gamma$ and~$\cP_1, \dots, \cP_\ell$ satisfy properties \cref{lm:tyingmanypathssupreg,lm:tyingmanypathsbadinfewpairs,%
		%lm:tyingmanypathsedgedisjoint,
		lm:tyingmanypathsvertexdisjoint,lm:tyingmanypathsroom,lm:tyingmanypathsendpoint} of \cref{lm:tyingmanypaths}, as well as the following.
		\begin{enumerate}
			\setcounter{enumi}{\value{tyinglemma1}}
			\item $\ell\leq n$.\label{lm:tyingfewpathsnumbersets}\COMMENT{Could take any linear upper bound. Later we will apply the lemma with~$\ell \leq \frac{n}{2}+o(n)$.}
			\newcounter{tyinglemma2}
			\setcounter{tyinglemma2}{\value{enumi}}
			\item For any~$i\in [\ell]$ and~$j\in [k]$,~$\cP_i$ contains at most~$2\beta^{-2}$ paths with an endpoint in~$V_j$ (and thus at most~$2\beta^{-2}k$ paths in total).\label{lm:tyingfewpathscluster}\COMMENT{In particular, this implies~$|\cP_i|\leq 2\beta^{-2}k$.}
			\newcounter{tyinglemma3}
			\setcounter{tyinglemma3}{\value{enumi}}
		\end{enumerate}
		Then, there exist disjoint $E_1, \dots, E_\ell \subseteq E(\Gamma)$ such that the following hold.
		\begin{enumerate}[label=\upshape(\alph*)]
						
			\item For any~$i\in [\ell]$, by using each edge in~$E_i$ exactly once, we can tie together some of the paths in~$\cP_i$ to form a set~$\cQ_i$ of vertex-disjoint paths such that for any connected component~$C$ of~$R$,~$\cQ_i$ contains at most one path with an endpoint in~$V_\Gamma(C)$.\label{lm:tyingfewpathsdecomp}
			
			\item For any distinct~$i,j\in [k]$ and~$x\in V_i$,~$E_1\cup\dots\cup E_\ell$ contains at most~$3\varepsilon^{\frac{1}{4}}m$ edges of~$\Gamma[V_i, V_j]$ which are incident to~$x$.\label{lm:tyingfewpathsneighbours}

			\item For any~$i\in [\ell]$ and~$j\in[k]$, $|V(\cP_i\cup E_i)\cap V_j|\leq \sqrt{\zeta} m$.\label{lm:tyingfewpathsroomQ}
		\end{enumerate}
	\end{lm}

	This is proved similarly to \cref{lm:tyingmanypaths} but since we now have fewer paths to link, we can use several superregular pairs of~$\Gamma$ to tie together paths whose endpoints are not necessarily in a same cluster. Thus, the main difference to the proof of \cref{lm:tyingmanypaths} is that, in order to link two paths, we no longer need to find a suitable superregular pair of~$\Gamma$ but a suitable walk in the reduced graph of~$\Gamma$.
	Moreover, since we have few paths to tie together, we no longer need to ensure that no superregular pair is used too many times (condition \cref{lm:tyingmanypathsroomcluster} in the proof of \cref{lm:tyingmanypaths}).		
	Finally, note that since link paths may now intersect several superregular pairs of~$\Gamma$, it no longer makes sense to talk about~$(A,B)$-link paths, so we only use the generic term \emph{link path} (defined at the beginning of \cref{sec:tyingpaths}). 
	
	\begin{proof}
		Let $E_1, \dots, E_\ell\subseteq E(\Gamma)$ be (possibly empty) disjoint sets of edges of~$\Gamma$ and assume inductively that for each~$i\in [\ell]$, by using each edge in~$E_i$ exactly once, we can tie together some of the paths in~$\cP_i$ to form a set~$\cQ_i$ of vertex-disjoint paths such that the following is satisfied.
		\begin{enumerate}[label=(\arabic*)]
			\item If $P\in \cQ_1\cup\dots\cup \cQ_\ell$ and~$P'$ is a link path of~$P$, then~$P'$ contains at most~$3$ vertices from each cluster and at most~$4$ edges from each superregular pair of~$\Gamma$.\label{lm:tyingfewpathslinkpath}
			\item For any clusters~$A$ and~$B$, and any~$x\in A$, the set~$\cQ_1\cup \dots\cup\cQ_\ell$ contains at most~$\varepsilon^{\frac{1}{2}}m$ link paths which have~$x$ as an endpoint and whose edge incident to~$x$ belongs to~$\Gamma[A, B]$.\label{lm:tyingfewpathslinkedge}
			\item For any~$x\in V(\Gamma)$, there are at most~$\varepsilon^{\frac{1}{4}}m$ link paths in~$\cQ_1\cup\dots\cup \cQ_\ell$ which contain~$x$ as an internal vertex.\label{lm:tyingfewpathslinkvertex}
		\end{enumerate}
		If for any~$i\in [\ell]$ and any connected component~$C$ of~$R$, the set~$\cQ_i$ contains at most one path with an endpoint in~$V_\Gamma(C)$, then \cref{lm:tyingfewpathsdecomp} holds. Moreover, \cref{lm:tyingfewpathslinkedge} and \cref{lm:tyingfewpathslinkvertex} imply \cref{lm:tyingfewpathsneighbours}, while \cref{lm:tyingfewpathsroomQ} follows from \cref{lm:tyingmanypathsroom,lm:tyingfewpathslinkpath,lm:tyingfewpathscluster}%
		\COMMENT{Let~$i\in [\ell]$ and~$j\in [k]$.~$\cP_i$ contains at most~$\zeta m$ vertices of~$V_j$ by \cref{lm:tyingmanypathsroom}. Moreover, by \cref{lm:tyingfewpathslinkpath}, each link path in~$\cQ_i$ contains at most~$3$ vertices of~$V_j$. By \cref{lm:tyingfewpathscluster},~$\cQ_i$ contains at most~$2\beta^{-2}k$ link paths. Therefore, $|V(\cP_i) \cup E_i\cap V_j|\leq \zeta m+ 6\beta^{-2}k\leq 2\zeta m$.},
		and we are done. 
		
		We may therefore assume that there exist~$i\in[\ell]$, a component~$C$ of~$R$, distinct paths~$P, P'\in \cQ_i$ and distinct vertices~$x,x'\in V_\Gamma(C)$ such that~$x$ and~$x'$ are endpoints of~$P$ and~$P'$, respectively.
		We find an~$(x, x')$-path in~$\Gamma$ to link~$P$ and~$P'$ as follows. Let $\Gamma'\coloneqq\Gamma\setminus (E_1\cup\dots\cup E_\ell)$. 
		By \cref{lm:tyingfewpathslinkedge,lm:tyingfewpathslinkvertex}, \cref{lm:verticesedgesremoval} implies that for any~$jj'\in E(R)$,~$\Gamma'[V_{jj'}, V_{j'j}]$ is still $[\varepsilon^{\frac{1}{9}}, \beta]$-superregular, where~$V_{jj'}$ and~$V_{j'j}$ are the support clusters of~$\Gamma[V_j, V_{j'}]$.
		
		Let~$i',i''\in [k]$ be such that~$x\in V_{i'}$ and~$x'\in V_{i''}$. 
		Choose~$j'\in [k]$ such that~$x\in V_{i'j'}$ and~$\cQ_1\cup\dots\cup \cQ_\ell$ contains fewer than~$\varepsilon^{\frac{1}{2}}m$ link paths which have~$x$ as an endpoint and whose edge incident to~$x$ belongs to~$\Gamma[V_{i'}, V_{j'}]$.
		The existence of such an index~$j'$ is guaranteed by \cref{lm:tyingmanypathsendpoint,lm:tyingmanypathsbadinfewpairs}. Indeed, by \cref{lm:tyingmanypathsendpoint}, there are at most\COMMENT{$\leq \frac{\varepsilon n}{\varepsilon^{\frac{1}{2}}m}\leq 2\varepsilon^{\frac{1}{2}}k$} $2\varepsilon^{\frac{1}{2}}k<\beta k$ indices~$j'$ such that~$\cQ_1\cup\dots\cup \cQ_\ell$ contains~$\varepsilon^{\frac{1}{2}}m$ link paths which have~$x$ as an endpoint and whose edge incident to~$x$ belongs to~$\Gamma[V_{i'}, V_{j'}]$.
		The existence of the desired index~$j'$ now follows from \cref{lm:tyingmanypathsbadinfewpairs}. 
		Similarly, pick~$j''\in [k]$ such that~$x'\in V_{i''j''}$ and~$\cQ_1\cup\dots\cup \cQ_\ell$ contains fewer than~$\varepsilon^{\frac{1}{2}}m$ link paths which have~$x'$ as an endpoint and whose edge incident to~$x'$ belongs to~$\Gamma[V_{i''}, V_{j''}]$.
		
		Let~$V_{i_1} \dots V_{i_r}$ be a~$(V_{j'}, V_{j''})$-path in~$R$, where $i_1\coloneqq j'$ and~$i_r\coloneqq j''$. Let~$i_0\coloneqq i'$ and~$i_{r+1}\coloneqq i''$.
		Then, $V_{i_0}V_{i_1}\dots V_{i_r}V_{i_{r+1}}$ is a~$(V_{i'}, V_{i''})$-walk in~$R$ where the clusters~$V_{i'}$ and~$V_{i''}$ appear at most twice and all other clusters occur at most once.
		For~$0\leq s\leq r+1$, let~$V_{i_s}'$ be obtained from~$V_{i_s}$ by deleting the following vertices:
		\begin{itemize}
			\item vertices in $V_{i_{s}}\setminus (V_{i_{s-1}i_s}\cap V_{i_si_{s+1}})$ (by \cref{lm:tyingmanypathssupreg}, there are at most~$2\varepsilon m$ such vertices);
			\item vertices in~$V(\cQ_i)\setminus\{x, x'\}$ (by \cref{lm:tyingmanypathsroom,lm:tyingfewpathscluster,lm:tyingfewpathslinkpath}, there are at most~$\frac{3\zeta m}{2}$ such vertices%
			\COMMENT{at most~$\zeta m$ original vertices in~$\cP_i$, plus at most~$3|\cP_i|\leq 3\cdot 2\beta^{-2}k$ additional vertices for tying up});
			\item vertices in $V_{i_s}\setminus \{x, x'\}$ which are an internal vertex of~$\varepsilon^{\frac{1}{4}}m$ link paths in~$\cQ_1\cup\dots\cup \cQ_\ell$ (by \cref{lm:tyingfewpathsnumbersets,lm:tyingfewpathslinkpath}, there are at most~$\varepsilon m$ such vertices%
			\COMMENT{at most $\frac{3\cdot \ell \cdot 2\beta^{-2}k}{\varepsilon^{\frac{1}{4}}m}\leq \varepsilon m$}).
		\end{itemize}
		Then, $|V_{i_s}'|\geq m-2\zeta m$. So for any~$s\in [r+1]$, by \cref{lm:verticesedgesremoval}, $\Gamma'[V_{i_{s-1}}', V_{i_s}']$ is still~$[\zeta^{\frac{1}{3}}, \beta]$-superregular%
		\COMMENT{We remove at most $2\zeta m\leq 3\zeta|V_{i_{s-1}i_s}|$ vertices from~$V_{i_{s-1}i_s}$ and similarly, at most~$3\zeta|V_{i_si_{s-1}}|$ vertices from~$V_{i_si_{s-1}}$.}. 
		We can therefore find an~$(x,x')$-path~$P''$ in~$\Gamma'$ containing exactly one edge of $\Gamma'[V_{i_{s-1}}', V_{i_s}']$ for each~$s\in [r]$ and at most~$3$ edges of $\Gamma'[V_{i_{r}}', V_{i_{r+1}}']$. 
		We \NEW{add the edges of~$P''$ to~$E_i$ and }replace in~$\cQ_i$ the paths~$P,P'$ by the concatenation of~$P, P''$, and~$P'$. By construction, \cref{lm:tyingfewpathslinkpath,lm:tyingfewpathslinkvertex,lm:tyingfewpathslinkedge} are still satisfied%
		\COMMENT{We get~$3$ vertices and~$4$ edges in \cref{lm:tyingfewpathslinkpath} if~$i_{r-1}=i_{r+1}$; \cref{lm:tyingfewpathslinkedge} holds by our choice of~$i_1$ and~$i_r$; \cref{lm:tyingfewpathslinkvertex} is satisfied by construction of the sets~$V_{i_s}'$.},
		as desired.
	\end{proof}

	The methods used to prove the previous lemma can be used to close a path~$P$ into a cycle provided the endpoints of~$P$ lie in a same connected component of~$\Gamma$. More generally, one can show the following. 
	
	\begin{lm}\label{lm:tyingfewpathscycle}
		Suppose $0<\frac{1}{n}\ll \frac{1}{k}\ll \varepsilon\leq \zeta \ll \beta \leq 1$. Let~$\Gamma$ be a graph on vertex set~$V$ of size~$n$ and~$\cP_1, \dots, \cP_\ell$ be sets of paths. Assume~$\Gamma$ and~$\cP_1, \dots, \cP_\ell$ satisfy properties \cref{lm:tyingmanypathssupreg,lm:tyingmanypathsbadinfewpairs,%
		%lm:tyingmanypathsedgedisjoint,
		lm:tyingmanypathsvertexdisjoint,lm:tyingmanypathsroom,lm:tyingmanypathsendpoint,lm:tyingfewpathsnumbersets,lm:tyingfewpathscluster} of \cref{lm:tyingfewpaths,lm:tyingmanypaths}. Suppose moreover that the following holds.
		\begin{enumerate}
			\setcounter{enumi}{\value{tyinglemma3}}
			\item For each~$i\in [\ell]$, there exists an ordering~$P_{i,1}, \dots, P_{i, \ell_i}$ of the paths in~$\cP_i$, and, for each~$j\in [\ell_i]$, an ordering~$x_{i,j}, x'_{i,j}$ of the endpoints of~$P_{i,j}$ such that the following holds. For each~$i\in [\ell]$ and~$j\in [\ell_i]$, there exists a component~$C$ of~$R$ such that $x'_{i,j}, x_{i,j+1}\in V_\Gamma(C)$, where $x_{i,\ell_i+1}\coloneqq x_{i,1}$.\label{lm:tyingfewpathscycleordering}
			\newcounter{tyinglemma4}
			\setcounter{tyinglemma4}{\value{enumi}}
		\end{enumerate}
		Then, there exist disjoint $E_1, \dots, E_\ell \subseteq E(\Gamma)$ such that the following hold.
		\begin{enumerate}[label=\upshape(\alph*)]
			\item For any~$i\in [\ell]$,~$\cP_i\cup E_i$ forms a cycle.\label{lm:tyingfewpathsdecompcycles}
			
			\item For any distinct~$i,j\in [k]$ and~$x\in V_i$,~$E_1\cup \dots \cup E_\ell$ contains at most~$3\varepsilon^{\frac{1}{4}}m$ edges of~$\Gamma[V_i, V_j]$ which are incident to~$x$.\label{lm:tyingfewpathsneighbourscycles}
		\end{enumerate}
	\end{lm}

	\begin{proof}
		The idea is to link the paths~$P_{i,j}$ and~$P_{i,j+1}$ together for each~$i\in [\ell]$ and~$j\in [\ell_i]$, where $P_{i, \ell_i+1}\coloneqq P_{i,1}$. This can be done by using the arguments of \cref{lm:tyingfewpaths} to find an~$(x_{i,j}', x_{i,j+1})$-path in~$\Gamma$ for each~$i\in [\ell]$ and~$j\in [\ell_i]$.
	\end{proof}

	In general, our sets of paths will not satisfy property \cref{lm:tyingfewpathscycleordering} of \cref{lm:tyingfewpathscycle}. In that case, we need to add suitable edges to our sets of paths before applying \cref{lm:tyingfewpathscycle}. This is achieved in the next lemma. 
	
	\begin{lm}\label{lm:closingcycle}
		Suppose $0<\frac{1}{n}\ll \frac{1}{k}\ll \varepsilon\ll\zeta \ll \beta \leq 1$. Let~$\Gamma$ be a graph on vertex set~$V$ of size~$n$ and~$\cP_1, \dots, \cP_\ell$ be sets of paths. Assume~$\Gamma$ and~$\cP_1, \dots, \cP_\ell$ satisfy properties \cref{lm:tyingmanypathssupreg,lm:tyingmanypathsbadinfewpairs,%
		%lm:tyingmanypathsedgedisjoint,
		lm:tyingmanypathsvertexdisjoint,lm:tyingmanypathsroom,lm:tyingmanypathsendpoint,lm:tyingfewpathsnumbersets} of \cref{lm:tyingmanypaths,lm:tyingfewpaths}, as well as the following.
		\begin{enumerate}[label=\upshape(\roman*$'$)]
			\setcounter{enumi}{\value{tyinglemma2}}
			\item For any~$i\in[\ell]$ and any connected component~$C$ of~$R$,~$\cP_i$ contains at most one path with an endpoint in~$V_\Gamma(C)$.\label{lm:closingcyclecomponent}
		\end{enumerate}
		Let~$\Gamma'$ be a graph on~$V$ such that~$\Gamma$ and~$\Gamma'$ are edge-disjoint and the following hold.
		\begin{enumerate}
			\setcounter{enumi}{\value{tyinglemma4}}
			\item $V_0,V_1,\dots, V_k$ is an~$(\varepsilon, \zeta, k, m, R')$-superregular partition of~$\Gamma'$. \label{lm:closingcycleGamma'}
			\item $R\cup R'$ is connected.\label{lm:closingcycleconnected}
			\newcounter{tyinglemma5}
			\setcounter{tyinglemma5}{\value{enumi}}
		\end{enumerate}
		Then, there exist disjoint $E_1, \dots, E_\ell\subseteq E(\Gamma)$ and $E_1',\dots, E_\ell' \subseteq E(\Gamma')$ such that the following hold.
		\begin{enumerate}[label=\upshape(\alph*)]
			\item For any~$i\in [\ell]$,~$\cP_i\cup E_i\cup E_i'$ forms a cycle.\label{lm:closingcycledecomp}
			
			\item For any distinct~$i,j\in [k]$ and~$x\in V_i$,~$E_1\cup\dots\cup E_\ell$ contains at most~$3\varepsilon^{\frac{1}{4}}m$ edges of~$\Gamma[V_i, V_j]$ which are incident to~$x$ and~$E_1'\cup\dots\cup E_\ell'$ contains at most~$\varepsilon m$ edges of~$\Gamma'[V_i, V_j]$ which are incident to~$x$. \label{lm:closingcycleneighbours}
		\end{enumerate}
	\end{lm}
	
	\begin{proof}
		We add some edges of~$\Gamma'$ to each~$\cP_i$ in order to satisfy property \cref{lm:tyingfewpathscycleordering} of \cref{lm:tyingfewpathscycle} as follows.
		For any~$i\in [\ell]$, denote $\cP_i\coloneqq \{P_{i,1}, \dots, P_{i,\ell_i}\}$ and, for each~$j\in [\ell_i]$, denote by~$x_{i,j}, x_{i,j}'$ the endpoints of~$P_{i,j}$.
		For each~$i\in[\ell]$, indices ranging in~$[\ell_i]$ are taken modulo~$\ell_i$, in particular~$\ell_i+1\coloneqq 1$.
		 
		Assume inductively that for some~$0\leq i\leq \ell$, $E_1', \dots, E_i'\subseteq E(\Gamma')$ are disjoint and satisfy the following.
		\begin{enumerate}[label=(\arabic*)]
			\item For each~$j\in[i]$, the edges in~$E_j'$ are vertex-disjoint from each other and from paths in~$\cP_j$. In particular,~$\cP_j'\coloneqq \cP_j\cup E_j'$ is a set of vertex-disjoint paths. \label{vertexdisjoint}
			\item For any distinct~$j,j'\in [k]$ and~$x\in V_j$,~$E_1'\cup \dots\cup E_i'$ contains at most~$\varepsilon m$ edges of~$\Gamma'[V_j, V_{j'}]$ which are incident to~$x$.\label{fewedgesofGamma'}
			\item For any~$x\in V$, there are at most~$\varepsilon n$ paths in $\cP_1'\cup \dots\cup \cP_i'\cup \cP_{i+1}\cup \dots \cup \cP_\ell$ which have~$x$ as an endpoint.\label{endpoint}
			\item For each~$j\in[i]$, there exists a partition $E_{j,1}' \cup \dots\cup E_{j,\ell_j}'$ of~$E_j'$ such that the following holds. 			
			For each~$j'\in [\ell_j]$, there exist an ordering~$y_1y_1', \dots, y_ty_t'$ of the edges in~$E_{j,j'}'$, and, distinct connected components~$C_0, \dots, C_t$ of~$R$ such that~$x_{j,j'}'\in V_\Gamma(C_0)$,~$x_{j,j'+1}\in V_\Gamma(C_t)$ and, for each~$s\in [t]$,~$y_s\in V_\Gamma(C_{s-1})$ and~$y_s'\in V_\Gamma(C_s)$. \label{componentpath}	
		\end{enumerate}
		Assume~$i=\ell$. Then, by \cref{fewedgesofGamma'}, the second part of property \cref{lm:closingcycleneighbours} holds. Also note that~$\cP_1', \dots, \cP_\ell'$ satisfy conditions \cref{lm:tyingmanypathssupreg,lm:tyingmanypathsbadinfewpairs,%
		%lm:tyingmanypathsedgedisjoint,
		lm:tyingmanypathsvertexdisjoint,lm:tyingmanypathsroom,lm:tyingmanypathsendpoint,lm:tyingfewpathsnumbersets,lm:tyingfewpathscluster,lm:tyingfewpathscycleordering} of \cref{lm:tyingfewpathscycle}, with~$2\zeta$ playing the role of~$\zeta$. Indeed, \cref{lm:tyingmanypathssupreg,%
		%lm:tyingmanypathsedgedisjoint,
		lm:tyingmanypathsbadinfewpairs,lm:tyingfewpathsnumbersets,lm:tyingmanypathsvertexdisjoint} are clearly satisfied. 
		Moreover, by \cref{lm:tyingmanypathsbadinfewpairs},~$R$ has at most~$\beta^{-1}$ connected components and thus \cref{lm:closingcyclecomponent,componentpath} imply%
		\COMMENT{For each~$j\in [\ell_i]$,~$|E_{i,j}'|\leq \beta^{-1}$ so $|\cP_i'|=|\cP_i|+|E_i|\leq (1+\beta^{-1})|\cP_i|\leq 2\beta^{-2}$ since, by \cref{lm:closingcyclecomponent},~$|\cP_i|\leq \beta^{-1}$.}
		$|\cP_i'|\leq (1+\beta^{-1})|\cP_i|\leq 2\beta^{-2}$.
		Therefore, \cref{lm:tyingfewpathscluster} holds. By \cref{componentpath}, for each~$i\in [\ell]$ and~$j\in [k]$, we have $|V(\cP_i')\cap V_j|\leq |V(\cP_i)\cap V_j|+2|\cP_i|\leq 2\zeta m$, so \cref{lm:tyingmanypathsroom} holds with~$2\zeta$ playing the role of~$\zeta$. Finally, \cref{lm:tyingmanypathsendpoint} holds by \cref{endpoint} and \cref{lm:tyingfewpathscycleordering} follows from \cref{componentpath}. 
		Thus, we can apply \cref{lm:tyingfewpathscycle} and we are done. We may therefore assume that~$i<\ell$.
		
		We construct~$E_{i+1}'$ as follows.
		Consider the auxiliary reduced graph~$\tR$ with the connected components of~$R$ as vertices and an edge between~$C$ and~$C'$ if~$R'$ contains an edge between~$C$ and~$C'$. 
		Note that by \cref{lm:closingcycleconnected},~$\tR$ is connected. 
		For each~$j\in [\ell_{i+1}]$, let~$C_j, C_j'$ be the connected components of~$R$ such that~$x_{i+1,j}\in V_\Gamma(C_j)$ and~$x_{i+1,j}'\in V_\Gamma(C_j')$. 
		For each~$j\in[\ell_{i+1}]$, fix a~$(C_j',C_{j+1})$-path~$Q_j$ in~$\tR$.
		
		Let~$\Gamma''$ be obtained from~$\Gamma'$ by deleting the following edges:
		\begin{itemize}
			\item edges in~$E_1'\cup \dots \cup E_i'$ (by \cref{fewedgesofGamma'}, we delete at most~$\varepsilon m^2$ such edges from each superregular pair of~$\Gamma'$);
			\item edges incident to some vertex~$x$ such that $\cP_1'\cup \dots\cup \cP_i'\cup \cP_{i+1}\cup \dots \cup \cP_\ell$ contains~$\varepsilon n$ paths which have~$x$ as an endpoint (by \cref{lm:tyingfewpathsnumbersets} and since for all~$i\in [\ell]$, $|\cP_i|\leq |\cP_i'|\leq 2\beta^{-2}$, we delete at most~$\varepsilon m^2$ such edges per superregular pair of~$\Gamma$%
			\COMMENT{There are at most $\frac{2\beta^{-2}\ell}{\varepsilon n}\leq \frac{2}{\varepsilon\beta^2}~$ such vertices so at most $\frac{2}{\varepsilon\beta} m\leq \varepsilon m^2$ edges per pair.});
			\item edges incident to some vertex in~$V(\cP_{i+1})$ (by \cref{lm:tyingmanypathsroom,lm:closingcycleGamma'}, we delete at most~$3\zeta^2 m^2$ such edges from each superregular pair of~$\Gamma'$\COMMENT{at most $2\zeta m \cdot (\zeta+\varepsilon)m\leq 3\zeta^2 m^2$});
			\item for each~$i',i''\in [k]$, edges of~$\Gamma'[V_{i'}, V_{i''}]$ which are incident to some vertex~$x$ such that $E_1'\cup \dots\cup E_i'$ contains~$\varepsilon m$ edges of~$\Gamma'[V_{i'}, V_{i''}]$ incident to~$x$ (by the fact that~$R$ has at most~$\beta^{-1}$ connected components (by \cref{lm:tyingmanypathsbadinfewpairs}) as well as \cref{componentpath,lm:closingcyclecomponent,lm:closingcycleGamma',lm:tyingfewpathsnumbersets}, we delete at most~$2\varepsilon\zeta m^2$ such edges from each superregular pair of~$\Gamma'$\COMMENT{by \cref{componentpath,lm:closingcyclecomponent,lm:tyingmanypathsbadinfewpairs,lm:tyingfewpathsnumbersets},~$E_1'\cup \dots\cup E_i'$ contains at most~$\beta^{-1}n$ edges from each pair of~$\Gamma'$. So, for each~$i',i''\in [k]$, there are at most $\frac{2\beta^{-1}n}{\varepsilon m}\leq \varepsilon m$ vertices which are incident to~$\varepsilon m$ edges in $E(\Gamma'[V_{i'}, V_{i''}])\cap (E_1'\cup\dots\cup E_i')$. Therefore, by \cref{lm:closingcycleGamma'}, we delete at most $\varepsilon m \cdot (\zeta +\varepsilon)m\leq 2\varepsilon\zeta m^2$ edges from each pair of~$\Gamma'$.}).
		\end{itemize} 
		Then, note that by \cref{lm:closingcycleGamma'}, for any~$i'i''\in E(R')$, $e(\Gamma''[V_{i'},V_{i''}])\geq (\zeta -\varepsilon)(1-\varepsilon)^2m^2-\varepsilon m^2-\varepsilon m^2-3 \zeta^2 m^2-2\varepsilon\zeta m^2 \geq \varepsilon m^2$.
		Thus, there exists a set~$E_{i+1}'\subseteq E(\Gamma'')$ of vertex-disjoint edges of~$\Gamma''$ such that \cref{componentpath,vertexdisjoint,fewedgesofGamma',endpoint} are still satisfied for~$i=i+1$, where, for~$j=i+1$ and each~$j'\in [\ell_{i+1}]$, the components in~$Q_{j'}$ play the roles of~$C_0, \dots, C_t$ in \cref{componentpath}%
		\COMMENT{By construction of~$\tR$ and since we need to take at most~$|\cP_{i+1}|\leq \beta^{-1}$ vertex-disjoint edges from each non-empty pair of~$\Gamma''$.}.
	\end{proof}

We will not always be able to add suitable edges to our sets of paths in order to apply \cref{lm:tyingfewpathscycle}. This problem can be circumvented by splitting paths and forming new sets of paths. This is achieved in the next \lcnamecref{lm:closingcycles}. Note that the cost of this approach is that we may obtain more cycles than in \cref{lm:closingcycle}, as well as a few leftover edges. Thus, \cref{lm:closingcycles} will only be used when we have some room to spare. 

\begin{lm}\label{lm:closingcycles}
	Suppose $0<\frac{1}{m}\ll \frac{1}{k}\ll \varepsilon\leq \zeta \ll \beta \leq 1$. Let~$\Gamma$ be a graph on vertex set~$V$ of size~$n$. Let~$\cP_1, \dots, \cP_\ell$ be sets of paths on~$V$\COMMENT{Note that this time we also require internal vertices to be in~$V$}. Assume~$\Gamma$ and~$\cP_1, \dots, \cP_\ell$ are all pairwise edge-disjoint and satisfy properties \cref{lm:tyingmanypathssupreg,lm:tyingmanypathsbadinfewpairs,%
	%lm:tyingmanypathsedgedisjoint,
	lm:tyingmanypathsvertexdisjoint,lm:tyingmanypathsroom} of \cref{lm:tyingmanypaths}, property \cref{lm:closingcyclecomponent} of \cref{lm:closingcycle}, as well as the following.
	\begin{enumerate}[label=\upshape(\roman*$'$)]
		\setcounter{enumi}{\value{tyinglemma0}}
		\item For any~$x\in V\setminus V_0$, $E(\cP_1\cup\dots\cup\cP_\ell)$ contains at most~$\varepsilon n$ edges incident with~$x$. \label{lm:closingcyclesmaxdegree}
		\setcounter{enumi}{\value{tyinglemma1}}
		\item ~$\ell\leq \zeta n$.\label{lm:closingcyclesnumbersets}
	\end{enumerate}
	\begin{enumerate}
		\setcounter{enumi}{\value{tyinglemma5}}
		\item For any~$x\in V_0$, if~$xy$ and~$xy'$ are distinct edges in~$E(\cP_1\cup\dots\cup\cP_\ell)$ then~$y\in V_\Gamma(C)$ and~$y'\in V_\Gamma(C')$ for some distinct components~$C$ and~$C'$ of~$R$.\label{lm:closingcyclesV0}
	\end{enumerate}
	Then, there exists~$E\subseteq E(\Gamma)$ such that the following hold.
	\begin{enumerate}[label=\upshape(\alph*)]
		
		\item $(\cP_1\cup\dots\cup\cP_\ell)\cup E$ can be decomposed into a set~$\cC$ of at most~$\beta n$ edge-disjoint cycles and a set~$E'$ of at most~$\beta^{-2}$ edges.\label{lm:closingcyclesdecomp}
		
		\item For any distinct~$i,j\in [k]$, and~$x\in V_i$,~$E$ contains at most~$\varepsilon^{\frac{1}{73}}m$ edges of~$\Gamma[V_i, V_j]$ which are incident to~$x$.\label{lm:closingcyclesneighbours}
	\end{enumerate}
\end{lm}

To prove \cref{lm:closingcycles}, we need the following results.

\begin{thm}[Vizing's theorem (see e.g.\ {\cite[Theorem 17.5]{bondy2008graph}})]\label{thm:Vizing}
	Let~$G$ be a multigraph with multiplicity~$\mu(G)$. Then the edge-chromatic number~$\chi'(G)$ of~$G$ satisfies~$\chi'(G)\leq \Delta(G)+\mu(G)$. In particular, if~$G$ is simple, then~$\chi'(G)\leq \Delta(G)+1$.
\end{thm}

\begin{lm}\label{lm:evenmatchings}
	Assume~$G$ is a multigraph with maximum degree~$\Delta$, multiplicity~$\mu$, and~$|E(G)|$ even. Then,~$G$ can be decomposed into at most~$\frac{3(\Delta+\mu)}{2}$ matchings of even size and at most~$\frac{\Delta+\mu}{2}$ paths and cycles of length~$2$.
\end{lm}

\begin{proof}
	Let~$M_1, \dots, M_r$ be an optimal matching decomposition of~$G$ \COMMENT{i.e.\ with~$r$ least possible}. By \cref{thm:Vizing},~$r\leq \Delta +\mu$. Let~$S$ be the set of indices~$i\in[r]$ such that~$|M_i|=1$ and~$T$ be the set of indices~$i\in [r]$ such that~$|M_i|$ is odd and at least~$3$. If~$|S|$ is odd, remove some~$i\in S$ and add it to~$T$ so that~$|S|$ is now even. Note that since~$|E(G)|$ is even,~$|T|$ must also be even. 
	
	For any distinct~$i,j\in S$, by minimality of~$r$,~$M_i\cup M_j$ is either a path of length~$2$ or a pair of parallel edges. Therefore,~$\bigcup_{i\in S} M_i$ can be decomposed into at most $\frac{\Delta + \mu}{2}$ paths and cycles of length~$2$. 
	For each distinct~$i, j\in T$, since~$|M_i|$ and~$|M_j|$ are odd and at most one of~$|M_i|$ and~$|M_j|$ is equal to~$1$, we can find vertex-disjoint~$e_i\in M_i$ and~$e_j\in M_j$, and thus decompose~$M_i\cup M_j$ into at most~$3$ matchings of even size: $M_i\setminus e_i, M_j\setminus e_j$, and~$\{e_i, e_j\}$. 
	Thus $\bigcup_{i\in [r]\setminus S} M_i$ can be decomposed into at most $\frac{3(\Delta +\mu)}{2}$ matchings of even size. 
\end{proof}	

\begin{proof}[Proof of \cref{lm:closingcycles}]
	Start with $\cC\coloneqq \emptyset$. 
	Note that by \cref{lm:tyingmanypathsbadinfewpairs,lm:tyingmanypathssupreg},~$R$ has~$c\leq \beta^{-1}$ connected components~$C_1, \dots, C_c$. For each~$i\in [c]$, colour each~$x\in V_\Gamma(C_i)$ with colour~$i$.
	We say a path in some~$\cP_j$ is \emph{monochromatic} if its endpoints are coloured with the same colour, and \emph{bichromatic} otherwise.
	We say a monochromatic path is \emph{coloured~$i$} if its endpoints are coloured~$i$, and we say a bichromatic path is \emph{coloured with~$\{i,i'\}$} if one of its endpoint is coloured~$i$ and the other is coloured~$i'$.
	Observe that exceptional vertices are left uncoloured but, by \cref{lm:tyingmanypathsvertexdisjoint}, all paths in $\cP_1\cup \dots \cup \cP_\ell$ have coloured endpoints.
	A path of length~$2$ with internal vertex in~$V_0$ is called an \emph{exceptional path}.
	
	By \cref{lm:closingcyclecomponent}, for each~$i\in [\ell]$,~$|\cP_i|\leq c\leq \beta^{-1}$. 
	Moreover, \cref{lm:tyingmanypathsvertexdisjoint,lm:closingcyclesV0} imply that no path in $\cP_1\cup \dots \cup \cP_\ell$ contains an edge inside~$V_0$. Thus, by repeatedly taking maximal monochromatic subpaths, each path in~$\cP_1\cup\dots\cup\cP_\ell$ admits a decomposition $\cD_{\rm mono}\cup \cD_{\rm bi}$,  where~$\cD_{\rm mono}$ is a set of at most~$c\leq\beta^{-1}$ monochromatic paths of distinct colours and~$\cD_{\rm bi}$ is a set of bichromatic edges and exceptional paths such that, if~$P,P'\in \cD_{\rm bi}$ are distinct, then they are coloured with distinct pairs of colours. This induces a decomposition of~$\cP_1\cup \dots \cup \cP_\ell$ into
	\begin{itemize}
		\item $\ell'\leq \beta^{-2}\ell$ monochromatic subpaths~$P_1, \dots, P_{\ell'}$; and
		\item for each~$1\leq i<i'\leq c$, a set~$\cQ_{ii'}$ of at most~$\beta^{-1}\ell$ bichromatic edges and exceptional paths coloured with~$\{i,i'\}$.
	\end{itemize}
	\NEW{Each monochromatic path $P_i$ will be tied into a cycle. In order to cover the bichromatic paths, we partition each $\cQ_{ii'}$ as follows.}	
	Observe that, by \cref{lm:closingcyclesV0}, the following holds.
	\begin{property}{\dagger}\label{prop:vertexdisjointpaths}
		For any $1\leq i<i'\leq c$, the exceptional paths in~$\cQ_{ii'}$ have distinct internal vertices.
	\end{property}

	By removing at most one edge or exceptional path from each~$\cQ_{ii'}$, we may assume~$|\cQ_{ii'}|$ is even for any~$1\leq i<i'\leq c$. Let~$E'$ be the set of deleted edges. Then, $|E'|\leq\beta^{-2}$, as desired for \cref{lm:closingcyclesdecomp}.\COMMENT{$|E'|\leq2\binom{c}{2}$ since exceptional paths have length~$2$.}	

	Let $1\leq i<i'\leq c$. Define a multiset~$\cQ_{ii'}^*$ of bicoloured edges coloured with~$\{i,i'\}$ by replacing each exceptional path in~$\cQ_{ii'}$ by a \emph{fictive edge} between its endpoints. Since $|V_0|\leq \varepsilon n$, each edge in~$\cQ_{ii'}^*$ has multiplicity at most $\varepsilon n+1$.
	Then, by \cref{lm:closingcyclesmaxdegree}, we can apply \cref{lm:evenmatchings} with~$\cQ_{ii'}^*$ playing the role of~$G$, $\Delta \leq \varepsilon n$, and $\mu \leq \varepsilon n +1$ to obtain $\ell_{ii'}^*\leq 4\varepsilon n$ matchings of even size, $\ell_{ii'}\leq 2\varepsilon n$ monochromatic paths of length~$2$, and $\ell_{ii'}'\leq 2\varepsilon n$ cycles of length~$2$.	
	Denote the matchings by~$M_{ii's}^*$, with $s\in [\ell_{ii'}^*]$.
	Replace, in the paths and cycles of length~$2$, the fictive edges by their corresponding exceptional paths. By \cref{prop:vertexdisjointpaths}, we thus obtain~$\ell_{ii'}$ edge-disjoint monochromatic paths, which we denote by~$P_{ii's}$, with $s\in [\ell_{ii'}]$, and~$\ell_{ii'}'$ edge-disjoint cycles which we add to~$\cC$. Note that $|\cC|=\sum_{1\leq i<i'\leq c}\ell_{ii'}'\leq \sqrt{\varepsilon}n$.
	\NEW{Each monochromatic path $P_{ii's}$ will be tied into a cycle.}
	
	For each~$1\leq i < i' \leq c$ and~$j\in [\ell_{ii'}^*]$, if there exists~$j'\in [k]$ such that $|V(M_{ii'j}^*)\cap V_{j'}|>\zeta m$, then randomly partition~$M_{ii'j}^*$ into~$2\zeta^{-1}$ submatchings whose sizes are even and approximately equal.
	By \cref{lm:Chernoff}, we may assume that each of the submatchings obtained contains at most~$\zeta m$ edges with an endpoint in~$V_{j'}$, for each~$j'\in [k]$%
	\COMMENT{Let~$M=M_{ii'j}^*$ and consider a random subset~$M'$ of~$M$ as above. Note that $\frac{|M|}{4\zeta^{-1}}\leq |M'|\leq \frac{3|M|}{4\zeta^{-1}}$. Let~$j'\in[k]$. We have $\mathbb{E}[|V(M')\cap V_{j'}|]\leq \frac{3\zeta |V(M)\cap V_{j'}|}{4}\leq \frac{3\zeta m}{4}$ and $\mathbb{E}[|V(M')\cap V_{j'}|]\geq \frac{\zeta |V(M)\cap V_{j'}|}{4}\leq \frac{\zeta^2 m}{4}$. 
	By \cref{lm:Chernoff}
		\begin{equation*}
		\mathbb{P}\left[|V(M')\cap V_{j'}|>\zeta m\right]
		\leq \mathbb{P}\left[|V(M')\cap V_{j'}|>\frac{4}{3}\mathbb{E}[|V(M')\cap V_{j'}|]\right]
		\leq \exp\left(-\frac{\zeta^2 m}{108}\right).
		\end{equation*}
	A union bound gives that~$|V(M')\cap V_{j'}|\leq\zeta m$ for each random subset~$M'$ of~$M$ and each~$j'\in [k]$ with positive probability.}.
	Denote by~$\cQ_{ii's}$, with $s\in [\ell_{ii'}'']$, the $\ell_{ii'}''\leq \frac{8 \varepsilon n}{\zeta}$ sets of paths obtained from these submatchings by replacing the fictive edges by their corresponding exceptional paths. By construction and \cref{prop:vertexdisjointpaths}, the following hold.	
	\begin{enumerate}[label=(\arabic*)]
		\item For any~$j'\in[k]$, $|V(\cQ_{ii'j})\cap V_{j'}|\leq \zeta m$. \label{lm:closingcyclesroom}
		\item $|\cQ_{ii'j}|$ is even. \label{lm:closingcycleseven}		
		\item All paths in~$\cQ_{ii'j}$ are pairwise vertex-disjoint, bichromatic, and coloured with~$\{i,i'\}$.	\label{lm:closingcyclescrosscomponents}
		\item The paths in~$\cQ_{ii's}$, for all $s\in [\ell_{ii'}'']$, and the paths~$P_{ii's}$, for all $s\in [\ell_{ii'}]$, are all pairwise edge-disjoint.
	\end{enumerate}
	\NEW{Each set $\cQ_{ii's}$ will be tied into a cycle.}
	
	\NEW{We now aim to apply \cref{lm:tyingmanypaths,lm:tyingfewpathscycle}. Note that all the edges that we still need to cover with cycles belong to one the monochromatic paths $P_i$, or one of the monochromatic paths $P_{ii's}$, or one of the sets $\cQ_{ii's}$. As mentioned above, the goal is tie each of these into a cycle, i.e.\ the sets $\cP_1, \dots, \cP_\ell$ in \cref{lm:tyingmanypaths,lm:tyingfewpathscycle} will consist of the sets of the form $\{P_i\}$, $\{P_{ii's}\}$, or $\cQ_{ii's}$. Formally, proceed as follows.}
	Let $\ell''\coloneqq \ell' +\sum _{1\leq i<i'\leq c} (\ell_{ii'}+\ell_{ii'}'')$. Then, by \cref{lm:closingcyclesnumbersets}, $\ell''\leq \ell\beta^{-2} + \beta^{-2}\left(2\varepsilon n+ \frac{8\varepsilon n}{\zeta}\right)\leq \beta^2 n$.
	Denote by~$\cP_1', \dots, \cP'_{\ell''}$ the sets in
	\[\left\{\{P_i\} \mid i\in [\ell']\right\}\cup \left\{ \{P_{ii'j}\} \mid 1\leq i<i'\leq c,\, j\in [\ell_{ii'}']\right\}\cup  \left\{ \cQ_{ii'j} \mid 1\leq i<i'\leq c, \, j\in [\ell_{ii'}'']\right\}.\]

	By construction, we can successively apply \cref{lm:tyingmanypaths,lm:tyingfewpathscycle} to tie up the paths in each~$\cP_i'$ into a cycle as follows. 
	First, let~$E_1, \dots, E_{\ell''}$ be the sets of edges of~$\Gamma$ obtained after applying \cref{lm:tyingmanypaths} with~$\cP_1', \dots, \cP_{\ell''}'$, and~$\ell''$ playing the roles of~$\cP_1, \dots, \cP_\ell$, and~$\ell$, respectively%
	\COMMENT{Conditions \cref{lm:tyingmanypathssupreg,lm:tyingmanypathsbadinfewpairs%
		%,lm:tyingmanypathsedgedisjoint
		,lm:tyingmanypathsvertexdisjoint} clearly hold, \cref{lm:tyingmanypathsroom} is satisfied by \cref{lm:closingcyclesroom} and since our original sets of paths~$\cP_1, \dots, \cP_\ell$ satisfied \cref{lm:tyingmanypathsroom}, while \cref{lm:tyingmanypathsendpoint} is entailed by \cref{lm:closingcyclesmaxdegree}.}. 
	Let~$\cQ_1, \dots, \cQ_{\ell''}$ be the sets of paths as in part \cref{lm:tyingmanypathsdecomp} of \cref{lm:tyingmanypaths}. Note that, for any~$i\in[\ell'']$ and~$j\in [k]$, by part \cref{lm:tyingmanypathsroomQ} of \cref{lm:tyingmanypaths}, $|V(\cQ_i)\cap V_j|\leq \sqrt{\zeta} m$.
	Moreover, condition \cref{lm:tyingfewpathscycleordering} of \cref{lm:tyingfewpathscycle} holds for the sets~$\cQ_1, \dots, \cQ_{\ell''}$ since, by construction, each~$\cP_i'$ either contains a single monochromatic path or, an even number of bichromatic paths coloured with the same pair of colours.	
	Let $\Gamma'\coloneqq\Gamma \setminus (E_1\cup \dots \cup E_{\ell''})$ and note that, by \cref{lm:verticesedgesremoval} and part \cref{lm:tyingmanypathsneighbours} of \cref{lm:tyingmanypaths},~$V_0, V_1, \dots, V_k$ is an $(\varepsilon^{\frac{1}{9}}, \beta, k, m, R)$-superregular partition of~$\Gamma'$.
	Thus, we can now apply \cref{lm:tyingfewpathscycle} with $\cQ_1, \dots, \cQ_{\ell''},\ell'', \Gamma',\varepsilon^{\frac{1}{9}}$, and~$\sqrt{\zeta}$ playing the roles of $\cP_1,\ldots, \cP_\ell, \ell, \Gamma, \varepsilon$, and~$\zeta$, respectively. Add all cycles obtained to~$\cC$ and note that $|\cC|\leq \sqrt{\varepsilon}n+\beta^2 n \leq \beta n$. Denote by~$E_1', \dots, E_{\ell''}'$ the sets of edges of~$\Gamma'$ obtained. Define $E\coloneqq E_1\cup\dots\cup E_{\ell''}\cup E_1'\cup\dots \cup E_{\ell''}'$ and observe that property \cref{lm:closingcyclesneighbours} of \cref{lm:closingcycles} holds by part \cref{lm:tyingmanypathsneighbours,lm:tyingfewpathsneighbourscycles} of \cref{lm:tyingmanypaths,lm:tyingfewpathscycle}\COMMENT{at most $3\varepsilon^{\frac{1}{4}}m+3(\varepsilon^{\frac{1}{9}})^{\frac{1}{4}}m\leq \varepsilon^{\frac{1}{37}}m$ edges incident to each vertex}. This completes the proof.		
\end{proof}

Finally, in \cref{step:sketchGamma} of the proof of \cref{thm:n/2}\cref{thm:n/2cycledecomp} (see the proof overview), we will need to cover a few excess edges. This will be achieved using \cref{lm:matchingtying}. The idea is similar to the approach described in \cref{fig:sketch}. An even cycle in the reduced graph can \NEW{be} decomposed into two matchings~$M$ and~$M'$. We can then use a few of the edges of the pairs in~$M'$ to tie together a path from each pair in~$M$, and similarly for~$M$ and~$M'$ exchanged. More precisely, we prove the following.

\begin{lm}\label{lm:matchingtying}
	Let $0<\frac{1}{m}\ll \varepsilon, \zeta \ll d\leq 1$ and suppose~$k\in \mathbb{N}^*$ is even. Let~$G$ be a graph on vertex set~$V$ and~$V_1, \dots, V_k$ be a partition of~$V$ into~$k$ clusters of size~$m$. Suppose that for any~$i\in [k]$, the pair~$G[V_i, V_{i+1}]$ is~$[\varepsilon, d]$-superregular (where~$V_{k+1}\coloneqq V_1$). 
	Suppose that~$\cP_1, \dots, \cP_\ell$ are sets of paths on~$V$ satisfying the following.
	\begin{enumerate}
		\item $\ell \leq \zeta m$.
		%\item The paths in~$\cP_1\cup \dots \cup \cP_\ell$ are edge-disjoint from each other and from~$G$. 
		\item For each~$i\in [\ell]$, there exists~$I_i\subseteq [k]_{\rm odd}$ or~$I_i\subseteq [k]_{\rm even}$ such that we can write~$\cP_i=\{P_{i,j}\mid j\in I_i\}$, where, for each~$j\in I_i$,~$P_{i,j}$ is a path of length at most~$\frac{dm}{10}$ with an endpoint in~$V_j$, an endpoint in~$V_{j+1}$ and $V(P_{i,j})\subseteq V_j\cup V_{j+1}$.
		\item Any vertex~$x\in V$ is an endpoint of at most~$4$ paths in~$\cP_1\cup \dots\cup \cP_\ell$. \label{lm:tyingmatchingendpoint}
	\end{enumerate}
	Then there exist disjoint $E_1, \dots, E_\ell\subseteq E(G)$ such that the following hold.
	\begin{enumerate}[label=\upshape(\alph*)]
		\item Each~$x\in V$ is an endpoint of at most~$6$ edges in~$E_1\cup \dots \cup E_\ell$.\label{lm:matchingtyingendpoints}
		\item For each~$i\in [\ell]$,~$\cP_i\cup E_i$ forms a cycle.
	\end{enumerate}
\end{lm}

\begin{proof}
	Assume inductively that for some~$0\leq i\leq \ell$, we have constructed disjoint sets~$E_1, \dots, E_i\subseteq E(G)$ such that the following hold.
	\begin{enumerate}[label=(\arabic*)]
		\item For each~$j\in [i]$,~$\cP_j\cup E_j$ forms a cycle~$C_j$.
		\item Each~$x\in V$ is an internal vertex of at most one link path in~$C_1\cup \dots\cup C_i$.\label{lm:matchingtyinglinkvertex}
		\item Let~$j\in [i]$ and~$i_1<\dots < i_s$ be an enumeration of~$I_j$. Let~$Q$ be a link path in~$C_j$. Then there exists~$t\in [s]$ such that the following hold. The path~$Q$ links~$P_{j,i_t}$ and~$P_{j,i_{t+1}}$, where~$i_{s+1}\coloneqq i_1$\COMMENT{i.e. the endpoints of~$Q$ are an endpoint of~$P_{j,i_t}$ and an endpoint of~$P_{j,i_{t+1}}$}. Moreover, $V(Q)\subseteq V_{i_t+1}\cup \dots \cup V_{i_{t+1}}$. Finally,~$Q$ contains at most~$3$ edges of~$G[V_{i_t+1}, V_{i_t+2}]$ and at most one edge of~$G[V_{j'}, V_{j'+1}]$ for each $j'= i_t+2, i_t+3, \dots, i_{t+1}-1$. \label{lm:matchingtyinglinkpath}
	\end{enumerate} 
	Observe that by \cref{lm:matchingtyinglinkvertex,lm:tyingmatchingendpoint}, a vertex~$x\in V$ is an endpoint of at most~$6$ edges in~$E_1\cup \dots \cup E_i$. 
	Thus, if~$i=\ell$, we are done. We may therefore assume that~$i<\ell$.  
	
	Let~$i_1<i_2<\dots < i_s$ be an enumeration of~$I_{i+1}$. For each~$t\in [s]$, denote by~$x_{i_t}$ and~$x_{i_t+1}$ the endpoints of~$P_{i_t}$ in~$V_{i_t}$ and~$V_{i_t+1}$, respectively.
	Define \NEW{$G'\coloneqq (G\setminus \bigcup_{j\in [i]}E_j)-(V(\cP_{i+1})\cup V(\bigcup_{j\in [i]}E_j)\setminus \{x_{i_t}, x_{i_t+1}\mid t\in [s]\})$}\OLD{$G'\coloneqq G-(V(\cP_{i+1})\cup V(E_1\cup \dots \cup E_i)\setminus \{x_{i_t}, x_{i_t+1}\mid t\in [s]\})$}.
	For any~$j\in [k]$, let~$V_j'$ be obtained from~$V_j$ by removing the vertices in $V(\cP_{i+1})\cup V(E_1\cup \dots \cup E_i)\setminus \{x_{i_t}, x_{i_t+1}\mid t\in [s]\}$ and note that, by \cref{lm:matchingtyinglinkpath,lm:tyingmatchingendpoint}, $|V_j\setminus V_j'|\leq \frac{dm}{10}+2\ell\leq \frac{dm}{5}$. \NEW{Thus, by \cref{lm:subgraph,lm:verticesedgesremoval}, $G'[V_j', V_{j+1}']$ is $\varepsilon^{\frac{1}{3}}$-regular%
		\COMMENT{First apply \cref{lm:verticesedgesremoval} to account for deleted edges then apply \cref{lm:subgraph} to account for deleted vertices.}.
	Moreover, each $x\in V_j'$ satisfies $|N_{G'}(x)\cap V_{j+1}'|\geq (d-\varepsilon)m-6-\frac{dm}{5}\geq \frac{3d}{5}|V_{j+1}'|$ and, similarly, each $x'\in V_{j+1}'$ satisfies $|N_{G'}(x')\cap V_j'|\geq \frac{3d}{5}|V_j'|$.}\OLD{Thus, since for each~$j\in [k]$, $G_{i+1}[V_j, V_{j+1}]$ is~$[\varepsilon, d]$-superregular, $G'[V_j', V_{j+1}']$ contains a path of length at most~$4$ between any two of its vertices.}
	
	Then, for each~$t\in [s]$, we find an~$(x_{i_t+1}, x_{i_{t+1}})$-path~$Q_{i_t}$ in~$G'$ as follows. First, we find a path $Q_{i_t}'=x_{i_t+2}\dots x_{i_{t+1}}$ in~$G'$ where~$x_j\in V_j'$ for each~$j=i_t+2, \dots, i_{t+1}$. \NEW{Then, we apply \cref{prop:supreg4} (with $G'[V_{i_t+1}', V_{i_t+2}'], \varepsilon^{\frac{1}{3}}$, and $\frac{3d}{5}$ playing the roles of $G, \varepsilon$, and $d$) to}\OLD{Then, by the above, we can} find an~$(x_{i_t+1},x_{i_t+2})$-path~$Q_{i_t}''$ of length at most~$3$ in~\NEW{$G'[V_{i_t+1}', V_{i_t+2}']$}\OLD{$G[V_{i_t+1}', V_{i_t+2}']$}. Let $Q_{i_t}\coloneqq x_{i_t+1}Q_{i_t}''x_{i_t+2}Q_{i_t}'x_{i_{t+1}}$. 
	Setting $E_{i+1}\coloneqq \bigcup_{t\in [s]}E(Q_{i_t})$ completes the proof.
\end{proof}

\onlyinsubfile{\bibliographystyle{abbrv}
	\bibliography{Bibliography/papers}}

%% file: Regularising.tex
\onlyinsubfile{
\setcounter{section}{4}
\setcounter{subsection}{3}
\addtocounter{subsection}{-1}
\setcounter{definition}{18}
}

\subsection{Making superregular pairs Eulerian and regular}\label{sec:regularising}

As discussed in the proof overview, in \cref{step:sketchGamma} of the proof of \cref{thm:n/2}\cref{thm:n/2cycledecomp}, we will need to decompose superregular pairs into Hamilton cycles. Thus we will need to ensure that our superregular pairs are Eulerian and regular. In this section, we introduce efficient tools for achieving this.

\begin{lm}\label{lm:Eulerianpairs}
	Let $0< \frac{1}{m}\ll \frac{1}{k} \ll \varepsilon \ll d\leq 1$. Then there exists a constant~$c=c(d,k)$ such that the following holds. \NEW{Let~$G$ be an Eulerian graph and~$V_0,V_1, \dots, V_k$ be an $(\varepsilon, \geq d, k,m, m', R)$-superregular equalised partition of~$G$. Suppose that~$V_0$ is a set of isolated vertices in~$G$.
	Then, there exists}\OLD{Let~$G$ be an Eulerian graph and~$V_0,V_1, \dots, V_k$ be an $(\varepsilon, \geq d, k,m, m', R)$-superregular equalised partition of~$G$. Suppose~$V_0$ is a set of isolated vertices in~$G$.
	Then, there exist a constant~$c=c(d,k)$ and} a spanning subgraph~$G'\subseteq G$ such that the following hold. For any~$ij\in E(R)$, $G'[V_i, V_j]$ is Eulerian. Moreover,~$G'$ can be obtained from~$G$ by removing at most~$c$ edge-disjoint cycles. \NEW{In particular, by \cref{lm:verticesedgesremoval},  $V_0,V_1, \dots, V_k$ is an $(2\sqrt{\varepsilon}, \geq d, k,m, m', R)$-superregular equalised partition of~$G'$.}\OLD{and, thus, $V_0,V_1, \dots, V_k$ is an $(2\sqrt{\varepsilon}, \geq d, k,m, m', R)$-superregular equalised partition of~$G'$.} 
\end{lm}

\begin{proof}
	For simplicity, we assume that~$R$ is connected. If~$R$ is not connected, we can proceed similarly, but apply our arguments to each component of~$R$ separately.
	For any~$AB\in E(R)$, we write~$A^B$ for the support cluster of~$A$ with respect to~$B$.
	Let~$H\subseteq G$,~$i\in [k]$, and~$x\in V_i$. We define the \textit{oddity} of~$x$ in~$H$, denoted~$\cO_H(x)$, as the number of indices~$j\in [k]$ such that~$|N_H(x)\cap V_j|$ is odd. The \textit{oddity} of~$H$ is defined as $\cO(H)\coloneqq \sum_{x\in V(H)} \cO_H(x)$. Let~$S(H)\coloneqq \{x\in V(H)\mid \cO_H(x)>0\}$ and $\cN(H)\coloneqq |S(H)|$. Thus,~$G[V_i,V_j]$ is Eulerian for all~$ij\in E(R)$ if and only if~$\cN(G)=0$, or, equivalently, if and only if~$\cO(G)=0$.	
	Our argument relies on the two following observations: 
	\begin{enumerate}[label=(\roman*)]
		\item any graph contains an even number of odd degree vertices, and, \label{evennumberodddegree}
		\item in an Eulerian graph, the oddity of each vertex is even.\label{evenoddity}
	\end{enumerate}
	Our proof splits into three steps.
	In \cref{step:oddity1}, we significantly reduce the number of vertices of positive oddity by removing cycles of linear length. Then, in \cref{step:oddity2}, we will proceed similarly but optimise the number of vertices whose oddity is reduced in order to decrease~$\cN(G)$ to a bounded number. Then, in \cref{step:oddity3}, we will be able to use a greedy approach.
	
	\begin{steps}
		\item \textbf{Decreasing the number of vertices with positive oddity to fewer than~$\frac{dm'}{2}$.}\label{step:oddity1} If~$\cN(G)<\frac{dm'}{2}$, let~$G_1\coloneqq G$ and go to the next step. Otherwise, we claim that there exists~$G_1\subseteq G$ such that $\cN(G_1)<\frac{dm'}{2}$ and~$G_1$ can be obtained from~$G$ by removing at most $c_1\coloneqq\frac{20 k^2}{d}$ cycles.
		
		Consider the following algorithm.
		Pick~$x_0\in S(G)$ and let~$P_0$ be the path~$x_0$ of length~$0$.
		Suppose that after~$i\geq 0$ steps, we have extended~$P_0$ to an~$(x_0, x_i)$-path~$P_i$. Let~$A_i$ be the cluster such that~$x_i\in A_i$. Let $G^i\coloneqq G-(V(P_i)\setminus\{x_i\})$ and $G^{i,0}\coloneqq G-(V(P_i)\setminus\{x_i, x_0\})$.
		\begin{case}
			\item\textbf{$P_i$ has length less than~$\frac{dm'}{4}$.}\label{algo1:1}
			\begin{enumerate}[label=(\alph*)]
				\item\label{algo1:1a} If there exist a cluster~$B_i$ and a vertex $x_{i+1}\in (A_i\cup B_i)\cap (S(G)\setminus V(P_i))$ such that both~$x_i$ and~$x_{i+1}$ have odd degree in $(G\setminus P_i)[A_i, B_i]$, pick such~$B_i$ and~$x_{i+1}$. 
				
				\NEW{Note that $|V(P_i)|\leq \frac{dm'}{4}$, so $\delta((G^i-V(P_i))[A_i^{B_i}, B_i^{A_i}])\geq \frac{2dm'}{3}$. Moreover, by \cref{lm:subgraph}, $(G^i-V(P_i))[A_i^{B_i}, B_i^{A_i}]$ is $\sqrt{\varepsilon}$-regular.
				Apply \cref{prop:supreg4} (with $(G^i-V(P_i))[A_i^{B_i}, B_i^{A_i}]$,  $\sqrt{\varepsilon}$, and $\frac{2d}{3}$ playing the roles of $G, \varepsilon$, and $d$) to find an $(x_i, x_{i+1})$-path $Q$ of length at most $4$ in $G^i[A_i, B_i]$.
				Let $P_{i+1}\coloneqq x_0P_ix_iQx_{i+1}$.}%
				\OLD{Let~$P_{i+1}$ be obtained by concatenating~$P_i$ with an~$(x_i, x_{i+1})$-path of~$G^i$ of length at most~$4$. Note that such a path exists. Indeed, $|V(P_i)|\leq \frac{dm'}{4}$, so $\delta(G^i[A_i^{B_i}, B_i^{A_i}])\geq \frac{2dm'}{3}$ and thus, since~$G[A_i^{B_i},B_i^{A_i}]$ is~$\varepsilon$-regular,~$G^i[A_i^{B_i},B_i^{A_i}]$ contains a path of length at most~$4$ between any two of its vertices.}
				Finally, observe that for any~$x\in V(G)$,
				\begin{subequations}\label{eqs:1a}
				\begin{align}
					\cO_{G\setminus P_{i+1}}(x)=\cO_{G\setminus P_i}(x)-1, \quad &\text{if~$x\in\{x_i, x_{i+1}\}$};\label{eq:1a}\\
					\cO_{G\setminus P_{i+1}}(x)=\cO_{G\setminus P_i}(x), \quad &\text{otherwise}.
				\end{align}
				\end{subequations}
				\item\label{algo1:1b} Otherwise, pick any $x_{i+1}\in S(G)\setminus V(P_i)$.
				\NEW{Let $A_{i+1}$ denote the cluster which contains $x_{i+1}$. We claim that there exists an $(x_i,x_{i+1})$-path $Q$ of length at most $2k$ in $G^i-(V(P_i)\setminus \{x_i, x_{i+1}\})$. Indeed, observe that for any $UW\in E(R)$, $\delta((G^i-(V(P_i)\setminus \{x_i, x_{i+1}\}))[U,W])\geq (d-\varepsilon)m'-\frac{dm'}{4}\geq \frac{2d m'}{3}$ and, by \cref{lm:subgraph}, $(G^i-(V(P_i)\setminus \{x_i, x_{i+1}\}))[U,W]$ is $\sqrt{\varepsilon}$-regular.
				Thus, if $A_i=A_{i+1}$, we can let $U\in N_R(A_i)$ and apply \cref{prop:supreg4} (with $(G^i-(V(P_i)\setminus \{x_i, x_{i+1}\}))[A_i,U], \sqrt{\varepsilon}$, and $\frac{2d}{3}$ playing the roles of $G, \varepsilon$, and $d$) to obtain an $(x_i, x_{i+1})$-path $Q$ of length at most $4$ in $G^i-(V(P_i)\setminus \{x_i, x_{i+1}\})$.
				Similarly, if $A_iA_{i+1}\in E(R)$, then we can apply \cref{prop:supreg4} (with $(G^i-(V(P_i)\setminus \{x_i, x_{i+1}\}))[A_i,A_{i+1}], \sqrt{\varepsilon}$, and $\frac{2d}{3}$ playing the roles of $G, \varepsilon$, and $d$) to obtain an $(x_i, x_{i+1})$-path $Q$ of length at most $4$ in $G^i-(V(P_i)\setminus \{x_i, x_{i+1}\})$.
				Suppose that $A_i\neq A_{i+1}$ and $A_iA_{i+1}\notin E(R)$. Let $Q'$ be an $(A_i, A_{i+1})$-path in $R$. (This is possible since, by assumption, $R$ is connected.) Note that $Q'$ is a path of length at most $k-1$.
				Denote $Q'=A_iU_1\dots U_{\ell}A_{i+1}$.
				(Note that, by assumption, $\ell\geq 1$.) By the above, there exists a path $x_iu_1\dots u_\ell$ in $G^i-(V(P_i)\setminus \{x_i,x_{i+1}\})$ such that, for each $j\in [\ell]$, $u_j\in U_j$. Apply \cref{prop:supreg4} (with $(G^i-(V(P_i)\setminus \{x_i,x_{i+1}\}))[U_\ell, A_{i+1}], \sqrt{\varepsilon}$, and $\frac{2d}{3}$ playing the roles of $G, \varepsilon$, and $d$) to obtain a $(u_\ell, x_{i+1})$-path $Q''$ of length at most $4$ in $(G^i-(V(P_i)\setminus \{x_i,x_{i+1}\}))[U_\ell, A_{i+1}]$. Then, $Q\coloneqq x_iu_1\dots u_\ell Q''x_{i+1}$ is a path of length at most $(k-2)+4\leq 2k$ in $G^i- (V(P_i)\setminus \{x_i, x_{i+1}\})$, as desired.} 
				
				\NEW{Let $P_{i+1}\coloneqq x_0P_ix_iQx_{i+1}$.}%				
				\OLD{Let~$P_{i+1}$ be obtained by concatenating~$P_i$ with an~$(x_i,x_{i+1})$-path of length at most~$4k$ in~$G^i$. Such a path exists by similar arguments as above. Moreover, $\cN(G)\geq \frac{dm'}{2}>\frac{dm'}{4}$ implies that~$S(G)\setminus V(P_i)$ is non-empty, as desired.}
				Note that \cref{evenoddity} implies that there exists a cluster~$B_i$ such that~$x_i$ has odd degree in $(G\setminus P_i)[A_i, B_i]$. Thus, \cref{algo1:1}\cref{algo1:1b} can only occur at most~$|E(R)|< k^2$ times in total\COMMENT{Since if \cref{algo1:1}\cref{algo1:1b} occurs, there are no vertices in $(A_i\cup B_i)\setminus V(P_i)$ with odd degree in $(G\setminus P_i)[A_i,B_i]$.}.
				Finally, observe that for any~$x\in V(G)$, 
				\begin{subequations}\label{eqs:1b}
				\begin{align}
				\cO_{G\setminus P_i}(x)-1\leq \cO_{G\setminus P_{i+1}}(x) \leq \cO_{G\setminus P_i}(x)+1, \quad &\text{if~$x\in \{x_i, x_{i+1}\}$};\label{eq:1b}\\
				\cO_{G\setminus P_i}(x)-2\leq \cO_{G\setminus P_{i+1}}(x)\leq \cO_{G\setminus P_i}(x)+2, \quad  &\text{if $x\in V(P_{i+1})\backslash (V(P_i)\cup \{x_{i+1}\})$};\\
				\cO_{G\setminus P_{i+1}}(x)=\cO_{G\setminus P_i}(x), \quad &\text{otherwise}.
				\end{align}
				\end{subequations}	 
			\end{enumerate}
			
			\item\label{algo1:2} \textbf{$P_i$ has length at least~$\frac{dm'}{4}$.}
			\NEW{Note that, by \cref{algo1:1}, $P_i$ has length at most $\frac{dm'}{4}+2k$. Thus, by similar arguments as above, there exists} an~$(x_i,x_0)$-path~\NEW{$Q$} in~$G^{i,0}$ of length at most~$2k$.\OLD{Note that such a path exists by similar arguments as above.} Output the cycle $C\coloneqq x_0P_ix_iQx_0$ and observe that~$C$ has length at least~$\frac{dm'}{4}$. Moreover,
			for any~$x\in V(G)$, 
			\begin{subequations}\label{eqs:2}
				\begin{align}
			\cO_{G\setminus P_i}(x)-1\leq \cO_{G\setminus C}(x) \leq \cO_{G\setminus P_i}(x)+1, \quad &\text{if~$x\in \{x_i, x_0\}$};\label{eq:2}\\
			\cO_{G\setminus P_i}(x)-2\leq \cO_{G\setminus C}(x)\leq \cO_{G\setminus P_i}(x)+2, \quad  &\text{if $x\in V(C)\backslash V(P_i)$};\\
			\cO_{G\setminus C}(x)=\cO_{G\setminus P_i}(x), \quad &\text{otherwise}.
			\end{align}\end{subequations}
		\end{case}
			
		We claim that $\cO(G\setminus C)\leq \cO(G)- \frac{dm'}{20}$. Indeed, as observed above, \cref{algo1:1}\cref{algo1:1b} can only occur fewer than~$k^2$ times, and, clearly, \cref{algo1:2} can only occur at most once. Thus, \cref{algo1:1}\cref{algo1:1a} occurs at least~$\frac{dm'}{32}$ times\COMMENT{at least $\left(\frac{dm'}{4} - k^2 \cdot (4k+1)\right) \cdot \frac{1}{4}=\frac{dm'}{16}-k^3-\frac{k^2}{4} \geq \frac{dm'}{32}$ times} and, therefore, \cref{eqs:1a,eqs:1b,eqs:2} imply \NEW{$\cO(G\setminus C)\leq \cO(G)-2\cdot \frac{dm'}{32}+2(2k+1) k^2\leq \cO(G)-\frac{dm'}{20}$}\OLD{$\cO(G\setminus C)\leq \cO(G)-2\cdot \frac{dm'}{32}+2(4k+1) k^2\leq \cO(G)-\frac{dm'}{20}$}.
		
		If $\cN(G\setminus C)<\frac{dm'}{2}$, let $G_1\coloneqq G\setminus C$. 
		Otherwise, repeatedly run the algorithm (where, in each iteration, the current graph plays the role of~$G$) and delete the resulting cycle until a graph~$G_1$ with $\cN(G_1)<\frac{dm'}{2}$ is obtained. 
		Note that we need to run the algorithm and delete the cycle obtained at most~$c_1= \frac{20 k^2}{d}$ times. Indeed, assume we repeatedly ran the algorithm and deleted the resulting cycle~$c_1$ times and let~$G_1$ be the graph obtained. First, observe that we have delete at most~$2c_1$ edges incident to each vertex, so~$\varepsilon$-regular pairs still have minimum degree at least~$(d-2\varepsilon)m'$ and, thus, in each iteration, the algorithm is always well defined. Since~$\cO(G)\leq m'k^2$, we have $\cO(G_1)\leq m'k^2-\frac{dm'}{20}\cdot \frac{20 k^2}{d}<\frac{dm'}{2}$ and in particular $\cN(G_1)<\frac{dm'}{2}$. Thus~$G_1$ can be obtained from~$G$ by removing at most~$c_1$ cycles, as desired. 	
		
		\item \textbf{Decreasing the number of vertices of positive oddity to fewer than~$100k^4$.}\label{step:oddity2} If~$\cN(G_1)<100k^4$, let~$G_2\coloneqq G_1$ and go to the next step. Otherwise, we claim that there exists~$G_2\subseteq G$ such that~$G_2$ can be obtained from~$G_1$ by removing at most $c_2\coloneqq \frac{21k}{2}$ cycles and such that~$\cN(G_2)<100k^4$.
		
		We proceed similarly as above, but since the number of vertices of positive oddity has now been significantly reduced, we can proceed more carefully. Indeed, we observe that, in the above algorithm, oddity may be created whenever \cref{algo1:1}\cref{algo1:1b} occurs (as well as in \cref{algo1:2}). Note that \cref{algo1:1}\cref{algo1:1b} occurs at stage~$i$ if, for all~$B_i$ as in \cref{algo1:1}\cref{algo1:1a}, all vertices of odd degree in~$G_1[A_i, B_i]$ already belong to~$V(P_i)$. Thus, in order to make our algorithm more efficient, we shall add the extra condition that the internal vertices of the short path used to extend the paths~$P_i$ have oddity~$0$. 	
		Namely, we now let $G^i\coloneqq G_1-((S(G_1)\cup V(P_i))\setminus\{x_i\})$ and $G^{i,0}\coloneqq G_1-((S(G_1)\cup V(P_i))\setminus\{x_i, x_0\})$.
		We note that this improvement could not have been implemented in \cref{step:oddity1} since~$\cN(G)$ was large\COMMENT{Otherwise, there would have been too few vertices available to find our short~$(x_i, x_{i+1})$-paths.}. Moreover, we observe that \cref{algo1:1}\cref{algo1:1b} may still occur. 
		We proceed as in \cref{algo1:1} of \cref{step:oddity1} if~$P_i$ has length less than~$\frac{dm'}{4}$ and $S(G_1)\not\subseteq V(P_i)$ and as in \cref{algo1:2} of \cref{step:oddity1} if~$P_i$ has length at least~$\frac{dm'}{4}$ (Case \hypertarget{algo2:2a}{2(a)}) or~$S(G_1)\subseteq V(P_i)$ (Case \hypertarget{algo2:2b}{2(b)}), with~$G_1$ playing the role of~$G$.
		
		By similar arguments as above, $|S(G_1)|\leq \frac{5dm'}{8}$ implies the desired short paths always exist and so the algorithm is well defined.
		Using similar arguments as in \cref{step:oddity1}, one can show that \NEW{$\cN(G_1\setminus C)<\cN(G_1)+(2k+1)k^2$}\OLD{$\cN(G_1\setminus C)<\cN(G_1)+(4k+1)k^2$}.\COMMENT{As mentioned in \cref{step:oddity1}, \namecrefs{algo1:1} \labelcref{algo1:1}\cref{algo1:1b} and \labelcref{algo1:2} occur at most~$k^2$ times in total, and, each of these times, we create at most~$2k+1$ vertices of positive oddity.}
		Moreover, if the algorithm terminates in Case~\hyperlink{algo2:2a}{2(a)} (i.e.\ if~$P_i$ has length at least $\frac{dm'}{4}$), then, as before, $\cO(G_1\setminus C)\leq \cO(G_1)-\frac{dm'}{20}$. If the algorithm terminates in Case~\hyperlink{algo2:2b}{2(b)} (i.e.\ if~$S(G_1)\subseteq V(C)$), then we note that, since \namecrefs{algo1:1} \labelcref{algo1:1}\cref{algo1:1b} and \labelcref{algo1:2} occur at most~$k^2$ times in total, by \cref{eq:1a,eq:1b,eq:2}, we have $\cO_{G_1\setminus C}(x)=\cO_{G_1}(x)-2$ for all but at most~$2k^2$ vertices~$x\in S(G_1)$.\COMMENT{We may only increase the oddity of vertices in~$S(G_1)$ when \cref{algo1:1}\cref{algo1:1b} or \cref{algo1:2} occurs. In such cases, we may only increase the oddity of two vertices in~$S(G_1)$ ($x_i$ and~$x_{i+1}$, or,~$x_i$ and~$x_0$, respectively).}
				
		If $\cN(G_1\setminus C)<100k^4$, then let $G_2\coloneqq G_1\setminus C$. Otherwise, repeatedly run the algorithm (where, in each iteration, the current graph plays the role of~$G_1$) and delete the resulting cycle until a graph~$G_2\subseteq G_1$ with $\cN(G_2)<100k^4$ is obtained. 
		We claim that we need to run the algorithm and delete the cycle obtained at most~$c_2= \frac{21k}{2}$ times. Indeed, assume we ran the algorithm and deleted the resulting cycle~$c_2$ times and let~$G_2$ be the graph obtained. Note that, in each iteration of the algorithm, the current graph has fewer than $\frac{dm'}{2}+5c_2k^3\leq\frac{5dm'}{8}$ vertices of positive oddity, so the algorithm is well defined in each of the iterations. 
		If the algorithm terminates in Case~\hyperlink{algo2:2b}{2(b)} in at least~$\frac{k}{2}$ of the iterations, then we note that all but at most~$2c_2k^2$ of the vertices in~$S(G_1)$ now have oddity~$0$. Therefore, $\cN(G_2)\leq 2c_2k^2+5c_2k^3<100k^4$, as desired.\COMMENT{All but at most~$c_2\cdot 2k^2$ vertices of positive oddity in~$G_1$ have their oddity decreased by~$2$ each time we obtain a cycle covering all vertices of positive oddity. Thus, all but at most $\frac{k}{2}\cdot 2k^2$ vertices of positive oddity in~$G_1$ have oddity~$0$ in~$G_2$. Moreover, we have created at most~$c_2\cdot 4k^3$ vertices of positive oddity in total.}
		Otherwise, the algorithm terminates in Case~\hyperlink{algo2:2a}{2(a)} in at least~$10k$ of the iterations. Therefore, $\cO(G_2)\leq \cO(G_1)-10k\cdot \frac{dm'}{20}\leq 0$ and so~$\cN(G_2)=0$.
		Thus~$G_2$ can be obtained from~$G_1$ by removing at most~$c_2$ cycles.
		
		\item \textbf{Removing all oddity.}\label{step:oddity3} If~$\cN(G_2)=0$, we set~$G'\coloneqq G_2$. Otherwise, we claim that there exists~$G'\subseteq G_2$ such that~$G'$ can be obtained from~$G_2$ by removing at most $c_3\coloneqq25k^4(k-1)$ cycles and such that~$\cN(G')=0$.
		
		Consider the following algorithm. Pick a vertex~$x_0\in S(G_2)$ and let~$P_0$ be the path~$x_0$ of length~$0$. Suppose that after~\NEW{$|S(G_2)|\geq i\geq 0$}\OLD{$i\geq 0$} steps we have extended~$P_0$ to an~$(x_0, x_i)$-path~$P_i$ \NEW{of length at most $4i$} such that $x_i\in S(G_2)$, $\cO_{G_2\setminus P_i}(x_i)=\cO_{G_2}(x_i)-1$ and $\cO_{G_2\setminus P_i}(x)\leq \cO_{G_2}(x)$ for all $x\in V(G)\setminus V_0$. Denote by~$A_i$ the cluster such that~$x_i\in A_i$. Let $x_{i+1}\in S(G_2)\setminus \{x_i\}$ be such that there exists a cluster~$B_i\neq A_i$ such that both~$x_i$ and~$x_{i+1}$ have odd degree in $(G_2\setminus P_i)[A_i,B_i]$. Observe that such cluster and vertex exist by \cref{evennumberodddegree,evenoddity}.
		
		\NEW{Since $|V(P_i)\cup S(G_2)|\leq 6|S(G_2)|\leq 600k^4$, \cref{lm:verticesedgesremoval} implies that $(G_2-((V(P_i)\cup S(G_2))\setminus\{x_i, x_{i+1}\}))[A_i^{B_i},B_i^{A_i}]$ is $[2\sqrt{\varepsilon}, \geq d]$-superregular. Apply \cref{prop:supreg4} (with $(G_2-((V(P_i)\cup S(G_2))\setminus\{x_i, x_{i+1}\})[A_i^{B_i},B_i^{A_i}], 2\sqrt{\varepsilon}$, and $d-2\sqrt{\varepsilon}$ playing the roles of $G, \varepsilon$, and $d$) to obtain an~$(x_i,x_{i+1})$-path~$Q_{i+1}$}\OLD{Let~$Q_{i+1}$ be an~$(x_i,x_{i+1})$-path} of length at most~$4$ in $(G_2-((V(P_i)\cup S(G_2))\setminus\{x_i, x_{i+1}\}))[A_i^{B_i},B_i^{A_i}]$.\OLD{Note that this path exists by similar arguments as in \cref{step:oddity1}.} Then,
		\begin{enumerate}[label=(\arabic*)]
			\item if~$x_{i+1}\in V(P_i)$, output the cycle $C\coloneqq x_{i+1}P_ix_iQ_{i+1}x_{i+1}$; \label{algo2:A}
			\item if~$x_{i+1}\notin V(P_i)$, let $P_{i+1}\coloneqq x_0P_ix_iQ_{i+1}x_{i+1}$. \label{algo2:B}
		\end{enumerate}
		\NEW{Note that if $i=S(G_2)$, then $S(G_2)\subseteq V(P_i)$. Thus, if we are in case \cref{algo2:B}, then $i+1\leq |S(G_2)|$, as desired.}
			
		Clearly, this algorithm eventually terminates, and, for each~$x\in V(G)$, we have 
		\[\cO_{G_2\setminus C}(x)=
		\begin{cases}
		\cO_{G_2}(x)-2, \quad &\text{if~$x\in V(C)\cap S(G_2)$,};\\
		\cO_{G_2}(x), \quad &\text{otherwise}.
		\end{cases}\]
		 Moreover,~$|V(C)\cap S(G_2)| \geq 2$ and, thus, $\cO(G_2\setminus C)\leq \cO(G_2)-4$.
		
		If $\cN(G_2\setminus C)=0$, then let $G'\coloneqq G_2\setminus C$. Otherwise, repeatedly run the algorithm (where, in each iteration, the current graph plays the role of~$G_2$) and delete the resulting cycle until a graph~$G'$ with $\cN(G')=0$ is obtained. By the above, we clearly need to run the algorithm and delete the cycle obtained at most $c_3= 25k^4(k-1)$ times\COMMENT{Since $\cO(G_2)\leq 100k^4(k-1)$.}. 
		Let $c\coloneqq c_1+c_2+c_3$. This completes the proof.\qed		
	\end{steps}
	\renewcommand{\qed}{}
\end{proof}

To regularise an Eulerian~$\varepsilon$-regular pair, we adapt an argument of \cite{glock2016optimal}. The idea is to repeatedly remove cycles covering all vertices of maximum degree. By ensuring that each vertex of minimum degree is covered by at most half of the cycles, we are able to regularise the pair by deleting only a few cycles.

\begin{lm} \label{lm:regularising}
	Suppose $0< \frac{1}{m} \ll \eta, \varepsilon \ll d \leq 1$ and let $\varepsilon'\coloneqq \max\{2\sqrt \varepsilon,4\sqrt{\eta}\}$. Let~$G$ be an Eulerian~$(\varepsilon,d)$-regular bipartite graph on vertex classes~$A, B$ of size~$m$. Let $\Theta \coloneqq \Delta(G)-\delta(G)$ and suppose~$\Theta\leq \eta m$.  Then there exists a spanning subgraph~$H\subseteq G$ such that~$H$ is regular,~$\varepsilon'$-regular, and can be obtained from~$G$ by removing at most~$2\Theta$ edge-disjoint cycles of length at least~$\frac{2m}{3}$. In particular,~$H$ is~$r$-regular for some $r\geq \Delta(G)-4\Theta$.\COMMENT{Let~$x$ be a vertex of maximum degree in~$G$. Then at most~$2\cdot 2\Theta$ edges incident to~$x$ are removed to obtain~$H$.}
\end{lm}

\begin{proof}
	First note that $\delta(G)\geq \frac{dm}{2}$.%
	\COMMENT{Indeed, since~$G$ is~$(\varepsilon,d)$-regular, we have $e_{G}(A, B)\geq dm^2$. Also, $e_{G}(A, B)\leq m\Delta(G)$. Thus $\Delta(G) \geq dm$ and $\delta(G) \geq dm- \eta m\geq \frac{dm}{2}$. }	
	Let~$G_0\coloneqq G$. We proceed inductively to build
	\begin{itemize}
		\item spanning subgraphs $G_1\supseteq G_2 \supseteq \dots \supseteq G_\ell~$ of~$G$;
		\item sets of vertices $A_0^\Delta\subseteq A_1^\Delta \subseteq \dots \subseteq  A_{\ell -1}^\Delta \subseteq A$ and $B_0^\Delta\subseteq B_1^\Delta \subseteq \dots B_{\ell -1}^\Delta \subseteq B$;
		\item sets of vertices $A_j^{\delta, 1}, A_j^{\delta, 2}\subseteq A$ and $B_j^{\delta, 1}, B_j^{\delta, 2}\subseteq B$ for each even $j\in \{0,1,\dots,\ell -1\}$;
		\item sets of vertices~$S_j^A\subseteq A$ and~$S_j^B\subseteq B$ for each $j\in \{0, 1,\dots, \ell-1\}$; and
		\item edge-disjoint cycles $C_0, C_1, \dots, C_{\ell -1}$;
	\end{itemize}   
	such that~$G_\ell~$ is regular,~$\ell \leq 2\Theta$, and, for each $i\in\{0,1,\dots, \ell -1\}$, the following hold.
	\begin{enumerate}[label=(\roman*)]
		\item $A_i^\Delta = \{a\in A \mid d_{G_i}(a)=\Delta(G_i)\}$ and $B_i^\Delta = \{b\in B \mid d_{G_i}(b)=\Delta(G_i)\}$. \label{Delta}
		\item If~$i$ is even,~$A_i^{\delta, 1}$ and~$A_i^{\delta, 2}$ are disjoint and such that $A_i^{\delta, 1} \cup A_i^{\delta, 2}=\{a\in A \mid d_{G_i}(a)=\delta(G_i)\}$, and, similarly,~$B_i^{\delta, 1}$ and~$B_i^{\delta, 2}$ are disjoint and such that $B_i^{\delta, 1} \cup B_i^{\delta, 2}=\{b\in B \mid d_{G_i}(b)=\delta(G_i)\}$. \label{delta}
		\item If~$i$ is even, then $A_i^\Delta \subseteq S_i^A \subseteq A\setminus A_i^{\delta, 1}$, $B_i^\Delta \subseteq S_i^B \subseteq B\setminus B_i^{\delta, 1}$ and $|S_i^A|=|S_i^B|\geq \frac{m}{3}$. \label{Seven}
		\item If~$i$ is odd, then $A_i^\Delta \subseteq S_i^A \subseteq A\setminus (A_{i-1}^{\delta, 2}\cap S_{i-1}^A)$, $B_i^\Delta \subseteq S_i^B \subseteq B\setminus (B_{i-1}^{\delta, 2}\cap S_{i-1}^B)$ and $|S_i^A|=|S_i^B|\geq \frac{m}{3}$.\COMMENT{We need $S_i^A \subseteq A\setminus (A_{i-1}^{\delta, 2}\cap S_{i-1}^A)$ instead of $S_i^A \subseteq A\setminus A_{i-1}^{\delta, 2}$ when $\Delta(G_i)-\delta(G_i)=2$.} \label{Sodd}
		\item $C_i$ is a Hamilton cycle of $G_i[S_i^A \cup S_i^B]$. \label{Ci}
		\item $G_{i+1}=G_i \setminus C_i$. \label{Gi+1}
	\end{enumerate}
	Assume that for some even~$0\leq i \leq 2\Theta$, we have already constructed subgraphs~$G_j$ for each~$j\in [i]$, sets~$A_j^\Delta$ and~$B_j^\Delta$ for each $j\in \{0, 1, \dots, i-1\}$, sets $A_j^{\delta, 1}, A_j^{\delta, 2}$ and $B_j^{\delta, 1}, B_j^{\delta, 2}$ for each even $j\in \{0, 2, \dots, i-2\}$, sets~$S_j^A$ and~$S_j^B$ for each $j\in \{0, 1, \dots, i-1\}$, and cycles~$C_j$ for each $j\in \{0, \dots, i-1\}$ such that \cref{Delta,delta,Seven,Sodd,Ci,Gi+1} are satisfied with~$j$ playing the role of~$i$ for all~$0\leq j\leq i$. If~$G_i$ is regular, let~$\ell \coloneqq i$. Otherwise, proceed as follows.
	
	Let $A_i^\Delta\coloneqq  \{a\in A\mid d_{G_i}(a)=\Delta(G_i)\}$ and $B_i^\Delta\coloneqq  \{b\in B\mid d_{G_i}(b)=\Delta(G_i)\}$, so that \cref{Delta} is satisfied for~$i$. 
	Also define $A_i^\delta\coloneqq  \{a\in A\mid d_{G_i}(a)=\delta(G_i)\}$ and $B_i^\delta\coloneqq  \{b\in B\mid d_{G_i}(b)=\delta(G_i)\}$. We will now construct $A_i^{\delta,1}, A_i^{\delta,2}$ and $B_i^{\delta,1}, B_i^{\delta,2}$, but, first note that, by \cref{delta,Gi+1,Seven,Sodd,Ci}, $\delta(G_i)\geq \delta(G_0)-i\geq \frac{dm}{2}-2\eta m\geq \frac{dm}{3}$. 
	A simple application of \cref{lm:Chernoff} shows that there exists a partition $A_i^{\delta,1}\cup A_i^{\delta,2}$ of~$A_i^\delta$ such that all of the following hold.
		\begin{enumerate}[label=\upshape(\alph*)]
			\item If $|A_i^\delta|\geq \frac{dm}{6}$, then $\frac{|A_i^\delta|}{3}\leq |A_i^{\delta,2}|\leq |A_i^{\delta,1}|\leq \frac{2|A_i^\delta|}{3}$. \label{Adeltalarge}
			\item If $|A_i^\delta|< \frac{dm}{6}$, then $ |A_i^{\delta,1}|=\left\lceil\frac{|A_i^\delta|}{2}\right\rceil$ and $|A_i^{\delta,2}|=\left\lfloor\frac{|A_i^\delta|}{2}\right\rfloor$. \label{Adeltasmall}
			\item For any~$b\in B$, $|N_{G_i}(b)\setminus A_i^{\delta,1}|,|N_{G_i}(b)\setminus A_i^{\delta,2}|\geq \frac{dm}{12}$. \label{Adeltadegree}
		\end{enumerate}
		Similarly, there exists a partition $B_i^{\delta,1}\cup B_i^{\delta,2}$ of~$B_i^\delta$ satisfying analogous properties.%
	Note that, in particular, \ref{delta} holds for~$i$.
	We will now construct sets~$S_i^A, S_i^B$ such that \ref{Seven} is satisfied and~$G_i[S_i^A\cup S_i^B]$ is Hamiltonian (in order to satisfy \ref{Ci}).
	
	We may assume without loss of generality that $|A\setminus A_i^{\delta,1}|\geq |B\setminus B_i^{\delta,1}|$ (the other case is similar).
	Clearly, $|A_i^\Delta|\leq m-|B_i^\delta|$.\COMMENT{Indeed, assume for a contradiction that $|A_i^\Delta |> m-| B_i^\delta |$. Then $e(G_i)=\sum_{a\in A}d_{G_i}(a)\geq |A_i^\Delta|\Delta(G_i)+(m-|A_i^\Delta|)\delta(G_i)>(m-|B_i^\delta|)\Delta(G_i)+|B_i^\delta|\delta(G_i)\geq \sum_{b\in B}d_{G_i}(b)=e(G_i)$, a contradiction.		
	%Let $C\subseteq B_i^\delta$ be a set of size $m-| A_i^\Delta|$. Note that $B_i^\delta \setminus C \neq \emptyset$. Then, we have $\sum_{a\in A\setminus A_i^\Delta} d_{G_i}(a) \geq \sum_{b\in C} d_{G_i}(b)$ and $\sum_{a\in A_i^\Delta} d_{G_i}(a) >\sum_{b\in B\setminus C} d_{G_i}(b)$. Therefore, $e(G_i)=\sum_{a\in A} d_{G_i}(a) > \sum_{b\in B} d_{G_i}(b)=e(G_i)$, a contradiction.
} 
	Thus,
	\begin{equation}\label{claim:stepi}
	m -| B_i^{\delta, 1}| - | A_i^\Delta |\geq|B_i^{\delta,2}|.
	\end{equation}
		
	Let $S_i^B\coloneqq B\setminus B_i^{\delta,1}$ and $T_i^A\subseteq A\setminus ( A_i^{\delta,1} \cup A_i^\Delta)$ be a set of size $|B\setminus B_i^{\delta,1}|-|A_i^\Delta|$ chosen uniformly at random.%
	\COMMENT{Note that this is possible since $|A\setminus ( A_i^{\delta,1} \cup A_i^\Delta)| = |A\setminus  A_i^{\delta,1} |-| A_i^\Delta| \geq |B\setminus B_i^{\delta,1}|-|A_i^\Delta|\overset{\text{\cref{claim:stepi}}}{\geq}0$.}
	Let $S_i^A\coloneqq  T_i^A \cup A_i^\Delta$. 
	By construction, \ref{Seven} holds for~$i$.%
	\COMMENT{We have $|S_i^A|=|S_i^B|=| B\setminus B_i^{\delta,1}| \geq \frac{m}{3}$, by construction.}
	Thus, by \cref{Gi+1}, \cref{lm:verticesedgesremoval,lm:subgraph} imply that~$G_i[S_i^A\cup S_i^B]$ is~$3\varepsilon'$-regular.%
	\COMMENT{by \ref{Gi+1} and \cref{lm:verticesedgesremoval},~$G_i$ is~$\varepsilon'$-regular so \cref{lm:subgraph} implies that~$G_i[S_i^A\cup S_i^B]$ is~$3\varepsilon'$-regular.}
	
	\begin{claim}
		With positive probability,~$G_i[S_i^A\cup S_i^B]$ is Hamiltonian.
	\end{claim}
	
	\begin{proofclaim}
		By \cref{lm:regHamcycle}, it suffices to show that $\delta (G_i[S_i^A\cup S_i^B])\geq \frac{d^2 m}{10^3}$ with positive probability.
		By construction, for any~$a\in S_i^A$, we have $d_{S_i^B}(a)\geq \frac{dm}{12}$, as needed. It remains to show that, with high probability, $d_{S_i^A}(b)\geq \frac{d^2 m}{10^3}$ for all~$b\in S_i^B$. 
		If $|B_i^{\delta,1}|-|A_i^{\delta, 1}|\leq \frac{dm}{24}$, then $|A\setminus (A_i^{\delta, 1}\cup S_i^A)|\leq \frac{dm}{24}$ and thus, by \cref{Adeltadegree}, $d_{S_i^A}(b)\geq \frac{dm}{24}$ for all~$b\in S_i^B$. 
		
		We may therefore assume that $|B_i^{\delta,1}|-|A_i^{\delta, 1}|\geq \frac{dm}{24}$. Then, $|B_i^{\delta,1}|\geq \frac{dm}{24}$, and, by \cref{Adeltalarge,Adeltasmall}, $| B_i^{\delta, 2}| \geq \frac{| B_i^{\delta, 1}|}{2}\geq \frac{dm}{48}$.\COMMENT{If $|B_i^\delta|\geq \frac{dm}{6}$, we have $|B_i^{\delta,2}|\geq \frac{|B_i^\delta|}{3}\geq \frac{|B_i^{\delta,1}|}{2}$. Otherwise, $|B_i^{\delta,2}|=\left\lfloor\frac{|B_i^\delta|}{2}\right\rfloor\geq |B_i^{\delta,1}|-1\geq\frac{|B_i^{\delta,1}|}{2}$ since $|B_i^{\delta,1}|\geq \frac{dm}{24}$ by assumption.}		
		Let~$b\in S_i^B$. Then,%
		\COMMENT{If $d_{A_i^\Delta}(b)\geq \frac{dm}{800}$, then $\mathbb{E}\left[d_{S_i^A}(b)\right]\geq
			d_{A_i^\Delta}(b)\geq \frac{dm}{800}$. Otherwise, by \cref{Adeltadegree}, $d_{A\setminus (A_i^{\delta ,1}\cup A_i^\Delta)}(b)\geq \frac{dm}{15}$ and thus, since			
			\[d_{T_i^A}(b)\sim \HGeom\left( | A\setminus (A_i^{\delta, 1}\cup A_i^\Delta)|, |B\setminus B_i^{\delta, 1}|-|A_i^\Delta|, d_{A\setminus ( A_i^{\delta, 1} \cup A_i^ \Delta)}(b)\right),\]
			we have
			\[\mathbb{E}\left[d_{S_i^A}(b)\right]\geq \mathbb{E}\left[d_{T_i^A}(b) \right]
			\geq \frac{\left(m -| B_i^{\delta, 1}| - \left| A_i^\Delta \right|\right)}{\left(m-|A_i^{\delta, 1}|-| A_i^\Delta |\right) }\cdot \frac{dm}{15}
			\geq \frac{\left(m -| B_i^{\delta, 1}| - | A_i^\Delta |\right)d}{15}.\]}
		\[\mathbb{E}\left[d_{S_i^A}(b)\right]\geq
		\mathbb{E}\left[d_{T_i^A}(b) \right]+d_{A_i^\Delta}(b)\overset{\text{\cref{Adeltadegree}, \cref{claim:stepi}}}{\geq}\frac{d^2 m}{800}.\]%
		Thus, by \cref{lm:Chernoff},%
		\COMMENT{\begin{equation*}
		\mathbb{P}\left[ d_{S_i^A}(b) < \frac{800}{10^3}\cdot \frac{d^2 m}{800}\right] \leq \exp \left(-\frac{1}{3}\cdot \left(\frac{800}{10^3}\right)^2\cdot \frac{d^2m}{800}\right)=\exp \left(-\frac{800 d^2m}{3\cdot 10^6}\right)\leq \frac{1}{m^2}.
		\end{equation*}}
		$\mathbb{P}\left[ d_{S_i^A}(b) < \frac{d^2 m}{10^3}\right] \leq \frac{1}{m^2}$, and,
		a union bound over all~$b\in S_i^B$ gives that, with positive probability, $d_{S_i^A}(b)\geq \frac{d^2 m}{10^3}$ for all~$b\in S_i^B$.
	\end{proofclaim}

	Thus, we can let~$C_i$ be a Hamilton cycle of~$G_i[S_i^A\cup S_i^B]$ and $G_{i+1}\coloneqq G_i\setminus C_i$. Then, \cref{Ci,Gi+1} are satisfied for~$i$.
	 
	If~$G_{i+1}$ is regular, let~$\ell \coloneqq i+1$. 
	Otherwise, proceed as follows. Let $A_{i+1}^\Delta\coloneqq  \{a\in A\mid d_{G_{i+1}}(a)=\Delta(G_{i+1})\}$, $B_{i+1}^\Delta\coloneqq  \{b\in B\mid d_{G_{i+1}}(b)=\Delta(G_{i+1})\}$, so that \ref{Delta} is satisfied for~$i+1$. We will now proceed similarly as above to construct $S_{i+1}^A, S_{i+1}^B, C_{i+1}$ and~$G_{i+2}$.
	
	If $\Delta(G_i)-\delta(G_i)=2$, we have $|A_i^\Delta|=|B_i^\Delta|$, $|A_i^\delta|=|B_i^\delta|$, $A=A_i^\delta \cup A_i^\Delta$, and $B=B_i^\delta \cup B_i^\Delta$. Then, by construction of~$G_{i+1}$, we have $A_{i+1}^\Delta = A_i^\Delta \cup A_i^{\delta,1}\cup  (A_i^{\delta,2}\setminus S_i^A)$ and $B_{i+1}^\Delta=B_i^\Delta \cup B_i^{\delta,1}$. Moreover, all vertices of $A\setminus A_{i+1}^\Delta$ and $B\setminus B_{i+1}^\Delta$ are vertices of minimum degree in~$G_{i+1}$. Thus we can let $S_{i+1}^A\coloneqq A_{i+1}^\Delta$ and $S_{i+1}^B\coloneqq B_{i+1}^\Delta$ in order to satisfy \ref{Sodd}\COMMENT{Clearly, $|S_{i+1}^A|=|S_{i+1}^B|$. Moreover, $|S_{i+1}^B|=|B_i^\Delta|+|B_i^{\delta, 1}|\geq |B_i^\Delta|+|B_i^{\delta, 2}|=|S_i^B|\geq \frac{m}{3}$.}. We can then proceed as above to define~$C_{i+1}$ and~$G_{i+2}$ satisfying \ref{Ci} and \ref{Gi+1} (and, in particular,~$G_{i+2}$ is regular). We may therefore assume that $\Delta(G_i)-\delta(G_i)>2$. Note that
	
	\begin{equation}\label{claim:stepi+1}
	    m\geq 
        \begin{cases}
	        | B_i^\delta|+|A_{i+1}^\Delta|, \quad &\text{if }|A\setminus A_i^{\delta, 2}|\geq |B\setminus B_i^{\delta, 2}|,\\
		    | A_i^\delta|+|B_{i+1}^\Delta|, \quad &\text{otherwise}.
        \end{cases}
	\end{equation}

	We construct $S_{i+1}^A, S_{i+1}^B, C_{i+1},$ and~$G_{i+2}$ similarly as above, but, we now let $A_{i+1}^\Delta \subseteq S_{i+1}^A \subseteq A\setminus A_i^{\delta, 2}$,  $B_{i+1}^\Delta \subseteq S_{i+1}^B \subseteq B\setminus B_i^{\delta, 2}$, and use \cref{claim:stepi+1} instead of \cref{claim:stepi}.
	
	We now show that we eventually obtain a regular graph~$G_\ell$, with~$\ell \leq 2\Theta$. Assume~$i$ is even and~$G_i$ is not regular. By \cref{Delta,Seven,Sodd,Ci,Gi+1}, $\Delta(G_i)=\Delta(G_0)-2i$. Moreover, by \cref{delta,Seven,Sodd,Ci,Gi+1}, $\delta(G_i)\geq \delta(G_0)-i$\COMMENT{Recall~$i$ is even.}. Thus,
	\begin{equation*}
	0< \Delta(G_i)-\delta(G_i)\leq (\Delta(G_0)-2i)- (\delta(G_0)-i)\leq \Theta - i,
	\end{equation*}
	and, therefore,~$i< \Theta$. Thus, $\ell \leq \Theta+1\leq 2\Theta $, as desired.
	
	Let~$H\coloneqq  G_\ell~$. Clearly,~$H$ is regular. Moreover, $G-H= \bigcup_{i=0}^{\ell -1} C_i$, with~$\ell \leq 2\Theta$, and so~$H$ can be obtained from~$G$ be removing at most~$2\Theta$ cycles. Moreover, by \cref{Seven,Sodd,Ci}, the cycles $C_0, \dots, C_{\ell -1}$ have length at least~$\frac{2m}{3}$, as desired. Finally, by \cref{lm:verticesedgesremoval},~$H$ is~$\varepsilon'$-regular. 
	This completes the proof.
\end{proof}

\onlyinsubfile{\bibliographystyle{abbrv}
	\bibliography{Bibliography/papers}}

%% file: Robust_Decomposition_Lemma.tex
\onlyinsubfile{
\setcounter{section}{4}
\setcounter{subsection}{4}
\addtocounter{subsection}{-1}
\setcounter{definition}{20}
}

\subsection{Robust decomposition lemma}\label{sec:robust}
\NEW{Note that the contents of this section will only be used in \cref{sec:Gamma} and so the reader may skip it and return to it later on.}

A key tool in our proofs will be the robust decomposition lemma of \cite{kuhn2013hamilton}, which implies the existence of a ``robust" Hamilton decomposition of superregular pairs.
More precisely, given a graph~$G$ consisting of suitable superregular pairs, it guarantees the existence of a spanning superregular graph~$G^{\rm rob}$ such that~$G^{\rm rob}\cup H$ has a Hamilton decomposition for any very sparse~$H$, i.e.\ $G^{\rm rob}$ is ``robustly" Hamilton decomposable. (The graph~$H$ will be the ``leftover" of an approximate decomposition of~$G\setminus G^{\rm rob}$.) Moreover, we can prescribe that a given Hamilton cycle in this decomposition contains a given set of edges. These edges will be the ``fictive edges" discussed in the proof overview (see \cref{fig:sketch} for example). To formalise the latter property, we need the notion of special path systems and special factors defined below. The fictive edges are extended into such special path systems prior to applying the robust decomposition lemma. It turns out to be more convenient to consider digraphs rather than graphs.

Given a digraph~$\overrightarrow{G}$ and a partition~$\cP$ of~$V(\overrightarrow{G})$ into~$k$ clusters~$V_1,\dots,V_k$ of equal size, a partition~$\cP'$ of~$V(\overrightarrow{G})$ is an \emph{$\ell$-refinement of~$\cP$}%
\COMMENT{page 130 of 1-factorisation paper}
if~$\cP'$ is obtained by splitting each~$V_i$ into~$\ell$ subclusters of equal size.\COMMENT{(So~$\cP'$ consists of~$\ell k$ clusters.)}
Let~$\overrightarrow{R}$ be the reduced digraph of~$\overrightarrow{G}$ with respect to~$\cP$ and assume that for any~$VW\in E(\overrightarrow{R})$, the pair~$\overrightarrow{G}[V, W]$ is~$[\varepsilon, d]$-superregular. We say~$\cP'$ is an \emph{$\varepsilon$-superregular~$\ell$-refinement of~$\cP$} if the following holds. For any~$V, W\in \cP$ and~$V', W'\in \cP'$ with~$V'\subseteq V$ and~$W'\subseteq W$, if~$VW\in E(\overrightarrow{R})$ then~$\overrightarrow{G}[V', W']$ is~$[\varepsilon, d]$-superregular.

We say $(\overrightarrow{G},\cP,\cP',\overrightarrow{R}, \overrightarrow{R}', C)$ is an \emph{$(\ell, k,m,\eps,d)$-bi-setup}%
\COMMENT{page 131 of 1-factorisation paper}
if the following properties are satisfied.
\begin{enumerate}[longlabel,label=\upshape(\text{BST}\arabic*)]
	\item $\overrightarrow{G}$ is a directed graph.
	\item $\cP$ is a partition of~$V(G)$ into~$k$ clusters of size~$m$, where~$k$ is even, and,~$\overrightarrow{R}$ is the reduced digraph of~$\overrightarrow{G}$ with respect to~$\cP$.
	\item $C$ is a Hamilton cycle of~$\overrightarrow{R}$.
	\item $\cP'$ is an~$\ell$-refinement of~$\cP$, and,~$\overrightarrow{R}'$ is the reduced digraph of~$\overrightarrow{G}$ with respect to~$\cP'$. 
	\item $R$ and~$R'$ are complete balanced bipartite digraphs.\label{def:BSTcompletebipartite}
	\item For each $VW\in E(\overrightarrow{R})\cup E(\overrightarrow{R}')$, the corresponding pair~$\overrightarrow{G}[V,W]$ is~$[\varepsilon, d]$-superregular.\label{def:bisetupsupreg}
	\COMMENT{In the original definition, a bi-universal walk~$U$ is also needed. But in our case such a walk is guaranteed by \cite[Lemma 4.5.1]{csaba2016proof} so no need to specify one. Moreover, the original definition requires~$\cP'$ to be~$\varepsilon$-uniform (see definition on page 130 of \cite{csaba2016proof}). However, \cref{def:bisetupsupreg} implies that~$\cP'$ is~$\sqrt{\varepsilon}$-uniform. (Indeed, if~$V\in \cP$ and~$x\in V(\overrightarrow{G})$ with~$|N^+(x)\cap V|>0$ then $|N^+(x)\cap V|=(d\pm \varepsilon)m$. Moreover, for any~$V'\in \cP'$ such that~$V'\subseteq  V$, we have $|N^+(x)\cap V'|=(d\pm \varepsilon)\frac{m}{\ell}$. Thus, $|N^+(x)\cap V'|\leq (d+\varepsilon)\frac{m}{\ell}\leq (1+ \sqrt{\varepsilon})(d-\varepsilon)\frac{m}{\ell}\leq (1+ \sqrt{\varepsilon})\frac{|N^+(x)\cap V'|}{\ell}$ and similarly, $|N^+(x)\cap V'|\geq (1- \sqrt{\varepsilon})\frac{|N^+(x)\cap V'|}{\ell}$.) Thus, since the exact value of~$\varepsilon$ is not important in the robust deceomposition lemma, we do not need to explicitly require~$\cP'$ to be~$\varepsilon$-uniform in our definition.}
\end{enumerate}
This is a special case of the setting in \cite{kuhn2013hamilton}, which also requires the existence of a ``universal walk"~$U$ in the reduced digraph~$\overrightarrow{R}$. This is trivially implied by assumption \cref{def:BSTcompletebipartite}.

Let~$\overrightarrow{G}$ be a digraph,~$\cP$ be a partition of~$V(\overrightarrow{G})$ into~$2k$ clusters $A_1, \dots, A_k, B_1, \dots, B_k$ of size~$m$, and~$\overrightarrow{R}$ be the corresponding reduced digraph of~$\overrightarrow{G}$. 
Let~$C\coloneqq A_1B_1\dots A_k B_k$ be a Hamilton cycle of~$\overrightarrow{R}$.
	Suppose~$f\in \mathbb{N}$ divides~$k$. 
	
The \emph{canonical interval partition~$\cI$ of~$C$ into~$f$ intervals}\COMMENT{page 124 of 1-factorisation paper} consists of the intervals \[I_i\coloneqq A_{(i-1)\frac{k}{f}+1} B_{(i-1)\frac{k}{f}+1} A_{(i-1)\frac{k}{f}+2} \dots B_{i\frac{k}{f}} A_{i\frac{k}{f}+1}\] for all~$i\in [f]$, where~$A_{k+1}\coloneqq A_1$.
		
Let~$\cP'$ be an~$\ell$-refinement of~$\cP$ and for each~$V\in \cP$ denote by~$V^1, \dots, V^\ell$ the partition of~$V$ induced by~$\cP'$.
Suppose~$\frac{2k}{f}\geq 3$.%
\COMMENT{Need~$\frac{2k}{f}\geq 3$ because fictive edges cannot have an endpoint in the endcluster of an interval. (See \cref{def:SPSfictive}.)}
Let~$i\in [f],h\in [\ell]$. Denote~$I_i$ by $V_{i,1}\dots V_{i,\frac{2k}{f}+1}$. 
A \emph{special path system~$SPS$ of style~$h$ in~$\overrightarrow{G}$ spanning the interval~$I_i$} consists of~$\frac{2k}{f}$ matchings~$M_1, \dots, M_{\frac{2k}{f}}$ such that the following hold.
\begin{enumerate}[label=\upshape(SPS\arabic*),longlabel]
	\item For all~$j\in [\frac{2k}{f}]$,~$M_j$ is a perfect matching between~$V_{i,j}^h$ and~$V_{i,j+1}^h$, with all edges oriented from~$V_{i,j}^h$ to~$V_{i,j+1}^h$.
	\item $SPS$ contains a \emph{fictive edge}~$f_{SPS}\in M_{\frac{2k}{f}-1}$ such that $E(SPS)\setminus \{f_{SPS}\}\subseteq E(\overrightarrow{G})$.\label{def:SPSfictive}
\end{enumerate}
A \emph{special factor~$SF$ with parameters~$(\ell, f)$ with respect to~$C$ and~$\cP'$ in~$\overrightarrow{G}$}
is a~$1$-regular digraph on~$V(\overrightarrow{G})$ consisting of~$\ell f$ special path systems~$SPS_{h,i}$, one for each~$(h,i)\in [\ell]\times [f]$, where~$SPS_{h,i}$ is a special path system of style~$h$ in~$\overrightarrow{G}$ spanning the interval~$I_i$. We denote the set of fictive edges of~$SF$ by $\Fict(SF)\coloneqq \{f_{SPS_{h,i}}\mid h\in[\ell], i\in [f]\}$.

\begin{figure}[htb]
	\centering
	\includegraphics[width=0.8\textwidth]{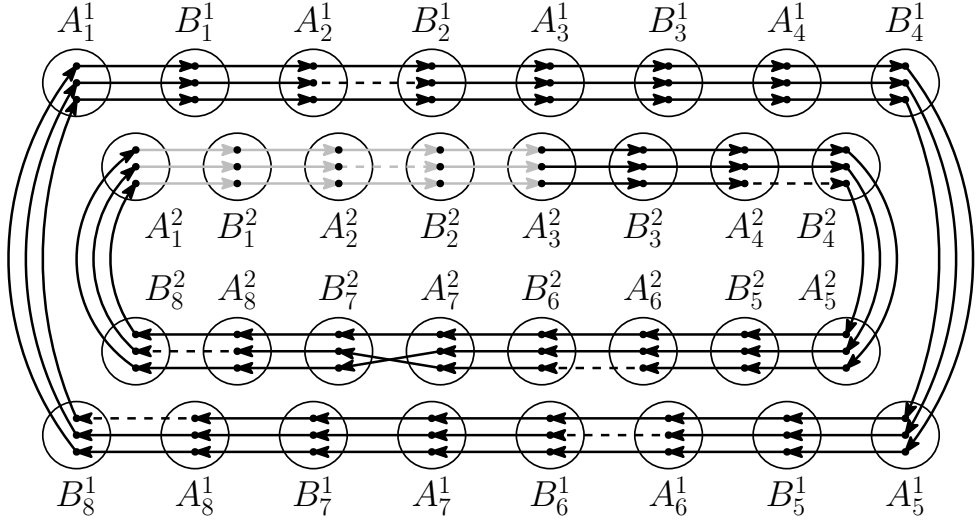}
	\caption{A special factor with parameters $(2,4)$ with respect to $C=A_1B_1A_2\dots B_8$ and $\cP'=\{A_1^1, A_1^2, A_2^1\dots, A_8^2, B_1^1, \dots, B_8^2\}$. The fictive edges are represented by dashed edges. The gray edges form a special path system of style $2$ spanning the first interval in the canonical interval partition of $C$ into $4$ intervals.}
	\label{fig:SF-SPS}
\end{figure}

More generally, we will use the term \emph{fictive edges} to refer to auxiliary edges that are artificially added to graphs. Whenever we add a set~$\cF$ of fictive edges to a (di)graph~$G$, we view them as being distinct from those in~$G$, even if they create multiple edges. Similarly, we also allow multiple edges within~$\cF$ and view these edges as being distinct from each other.
We are now ready to state the (bipartite version of the) robust decomposition lemma. As indicated above, it guarantees the existence of a robustly decomposable digraph~$\overrightarrow{G}^{\rm rob}$ which consists of a ``chord absorber"~$\overrightarrow{CA}(r)$ and a ``parity extended cycle absorber"~$\overrightarrow{PCA}(r)$, as well as a prescribed set of special factors (which contain the fictive edges).

\begin{lm}[Robust Decomposition Lemma \cite{kuhn2013hamilton}] \label{lm:robustdecomp}
\NEW{Let} $0<\frac{1}{m}\ll \frac{1}{k}\ll \eps \ll \frac{1}{q} \ll \frac{1}{f} \ll \frac{r_1}{m}\ll d\ll \frac{1}{\ell},\frac{1}{g}\ll 1$ and
that~$rk^2\le m$. Let
\[r_2\coloneqq 96\ell g^2kr, \quad r_3\coloneqq \frac{rfk}{q}, \quad r^\diamond\coloneqq r_1+r_2+r-(q-1)r_3, \quad s\coloneqq rfk+7r^\diamond\]
and suppose that $\frac{k}{14}, \frac{k}{f}, \frac{k}{g}, \frac{q}{f}, \frac{m}{4\ell}, \frac{fm}{q}, \frac{2fk}{3g(g-1)}, \frac{\ell}{2}\in \mathbb{N}$.
\NEW{Let} $(\overrightarrow{G},\cP,\cP',\overrightarrow{R},\overrightarrow{R}',C)$ \NEW{be} an $(\ell,k,m,\eps,d)$-bi-setup and~$C=V_1\dots V_k$.
Suppose that~$\cP^*$ is a~$\frac{q}{f}$-refinement
of~$\cP$ and that~$SF_1,\dots, SF_{r_3}$ are edge-disjoint special factors
with parameters~$(\frac{q}{f},f)$ 
with respect to~$C$ and~$\cP^*$ in~$\overrightarrow{G}$. Let $\mathcal{SF}\coloneqq SF_1\cup\dots \cup SF_{r_3}$.
Then there exists a digraph~$\overrightarrow{CA}(r)$ for which the following hold.
\begin{enumerate}
\item\label{lm:robustdecompCA} $\overrightarrow{CA}(r)$ is an~$(r_1+r_2)$-regular spanning subdigraph of~$\overrightarrow{G}$ which is edge-disjoint from~$\mathcal{SF}$.
\item\label{lm:robustdecompPCA} Suppose that~$SF'_1,\dots, SF'_{r^\diamond}$ are special factors
with parameters~$(1,7)$
with respect to~$C$ and~$\cP$ in~$\overrightarrow{G}$ which are edge-disjoint from each other and from $\overrightarrow{CA}(r)\cup \mathcal{SF}$.
Let $\mathcal{SF}'\coloneqq SF'_1\cup\dots \cup SF'_{r^\diamond}$.
Then there exists a digraph~$\overrightarrow{PCA}(r)$ for which the following hold.
\begin{enumerate}[label=\upshape(\alph*)]
\item $\overrightarrow{PCA}(r)$ is a~$5r^\diamond$-regular spanning subdigraph of~$\overrightarrow{G}$ which
is edge-disjoint from $\overrightarrow{CA}(r)\cup \mathcal{SF}\cup \mathcal{SF}'$.
\item Let~$\mathcal{SPS}$ be the set consisting of all the~$s$ special path systems
contained in~$\mathcal{SF}\cup \mathcal{SF}'$. Let~$V_{\rm even}$ denote the union of all~$V_i$ over all~$i\in [k]_{\rm even}$
and define~$V_{\rm odd}$ similarly.
Suppose that~$\overrightarrow{H}$ is an~$r$-regular bipartite digraph on~$V(\overrightarrow{G})$ with vertex classes~$V_{\rm even}$ and~$V_{\rm odd}$
which is edge-disjoint from $\overrightarrow{G}^{\rm rob}\coloneqq \overrightarrow{CA}(r)\cup \overrightarrow{PCA}(r)\cup\mathcal{SF}\cup \mathcal{SF}'$.
Then $\overrightarrow{H}\cup \overrightarrow{G}^{\rm rob}$ has a decomposition into~$s$
edge-disjoint Hamilton cycles~$C_1,\dots,C_s$.
Moreover,~$C_i$ contains one of the special path systems from~$\mathcal{SPS}$, for each~$i\in [s]$.
\end{enumerate}
\end{enumerate}
\end{lm}

\onlyinsubfile{\bibliographystyle{abbrv}
	\bibliography{Bibliography/papers}}

%% file: Szemeredi_and_Cleaning.tex
\onlyinsubfile{
\setcounter{section}{5}
\setcounter{subsection}{1}
\addtocounter{subsection}{-1}
}

\subsection{Applying Szemer\'{e}di's regularity lemma and setting aside random subgraphs \texorpdfstring{$\Gamma$}{Gamma} and \texorpdfstring{$\Gamma'$}{Gamma'}}\label{sec:Szemeredi}
	
This section corresponds to \cref{step:sketchcleaning} of the proof overview. The proof of \cref{lm:goodgraph} relies of a straightforward application of Szemer\'{e}di's regularity lemma and a cleaning procedure similar to the one used to \NEW{prove} the degree form of the regularity lemma.
\NOAPPENDIX{For details, see the appendix of the arXiv version.}%
\APPENDIX{For details, see the appendix.}

\begin{lm}\label{lm:goodgraph}
	Let $0< \frac{1}{M}\ll \varepsilon \ll \zeta \ll d \ll\beta \ll \alpha \leq 1$ and $\frac{1}{M}\leq \frac{1}{L}\ll d$. Let~$r\in \mathbb{N}^*$. Then there exist~$M', n_0\in \mathbb{N}^*$ such that the following holds. Let~$G$ be a graph on vertex set~$V$ with~$|V|=n\geq n_0$, and~$\delta(G)\geq \alpha n$. Then~$G$ can be decomposed into edge-disjoint graphs~$G', \Gamma, \Gamma', H$, and,~$V$ can be partitioned into~$k$ clusters~$V_1, \dots, V_k$ of size~$m$ and an exceptional set~$V_0$ such that the following properties are satisfied.
	\begin{enumerate}
		\item $M\leq k\leq M'$. \label{lm:goodk}
		\item $\frac{m'}{r}\in \mathbb{N}^*$.\label{lm:goodm'}
		\item $V_0, V_1, \dots, V_k$ is an \label{lm:goodsupregeqpart}
		\begin{itemize}
			\item $(\varepsilon, \geq d, k, m, m', R)$-superregular equalised partition of~$G'$;
			\item $(\varepsilon, \beta, k, m, m', R')$-superregular equalised partition of~$\Gamma$;
			\item $(\varepsilon, \zeta, k, m, m', R'')$-superregular equalised partition of~$\Gamma'$.
		\end{itemize}
		\item $R'$ and~$R''$ are edge-disjoint and~$R=R'\cup R''$. \label{lm:goodreduced}
		\item $G', \Gamma$ and,~$\Gamma'$ have the same support clusters. \label{lm:goodsamesupport}
		\item Each~$x\in V\setminus V_0$ belongs to at least~$\beta k$ superregular pairs of~$\Gamma$.\label{lm:goodbadinfewpairs}
		\item There exists a decomposition~$\cD_{R'}$ of~$R'$ into at most~$\frac{k}{2}$ cycles whose lengths are even and at least~$L$ and such that for any distinct~$i,j,j'\in [k]$, if~$V_jV_iV_{j'}$ is a subpath of a cycle in~$\cD_{R'}$ then the support clusters of~$V_i$ with respect to~$V_j$ and~$V_{j'}$ are the same.\label{lm:goodreducedpartition}
		\item $V_0$ is a set of isolated vertices in~$\Gamma, \Gamma'$, and~$H$.\label{lm:goodV0GammaH}
		\item $\Delta(H)\leq 4dn$. \label{lm:goodDeltabad}
	\end{enumerate}
\end{lm}

\onlyinsubfile{\bibliographystyle{plain}
	\bibliography{Bibliography/papers}}

%% file: Inside_V0.tex
\onlyinsubfile{
		\setcounter{section}{5}
		\setcounter{subsection}{2}
		\addtocounter{section}{-1}
		\setcounter{definition}{1}
	}
	
	\subsection{Covering the edges inside the exceptional set}\label{sec:insideV0} This section corresponds to \cref{step:sketchinsideV0} of the proof overview.

\begin{lm}\label{lm:insideV0cycles}
	Suppose $0<\frac{1}{n}\ll \frac{1}{k}\ll \varepsilon \ll d\ll \beta \ll \alpha\leq 1$.
	Let~$G$ be a graph on vertex set~$V$ with~$|V|=n$ and let~$\Gamma$ be edge-disjoint from~$G$. Assume~$G$ and~$\Gamma$ satisfy the following.
	\begin{enumerate}
		\item $V_0, V_1, \dots, V_k$ is an $(\varepsilon, \beta, k, m, m' ,R')$-superregular equalised partition of~$\Gamma$.\label{lm:insideV0Gammasupreg}
		\item Any~$x\in V\setminus V_0$ belongs to at least~$\beta k$ superregular pairs of~$\Gamma$.\label{lm:insideV0badinfewpairs}
		\item $V_0, V_1, \dots, V_k$ is an $(\varepsilon, \geq d, k, m, m', R)$-superregular equalised partition of~$G$.\label{lm:insideV0supreg}
		\item For any~$x\in V_0$,~$d_G(x)\geq \alpha n$. \label{lm:insideV0degree}
	\end{enumerate}
	Let $\varepsilon'\coloneqq\varepsilon ^{\frac{1}{73}}$. Then there exists~$H\subseteq G\cup \Gamma$ such that the following hold.
	\begin{enumerate}[label=\upshape(\alph*)]
		\item $G[V_0]=H[V_0]$.\label{lm:insideV0coverV0}
		\item $V_0, V_1, \dots, V_k$ is an $(\varepsilon', \geq d, k, m, m', R)$-superregular equalised partition of~$G\setminus H$.\label{lm:insideV0Gstillsupreg}
		\item $V_0, V_1, \dots, V_k$ is an $(\varepsilon', \beta, k, m, m', R')$-superregular equalised partition of~$\Gamma\setminus H$.\label{lm:insideV0Gammastillsupreg}
		\newcounter{insideV0}
		\setcounter{insideV0}{\value{enumi}}
		\item There exists a decomposition~$\cD\cup \cD'$ of~$H$ where~$\cD$ is a set of at most~$\beta n$ cycles and~$\cD'$ is a set of at most~$\beta^{-2}$ edges.\label{lm:insideV0cycledecomp}
	\end{enumerate}
\end{lm}

\begin{proof} 
	This can be proved in a similar way as \cref{lm:closingcycles}, so we only provide a sketch of \NEW{the} proof. Let~$\cD\coloneqq \emptyset$. Let~$C_1,\dots, C_c$ be the connected components of~$R'$, where, by \cref{lm:insideV0badinfewpairs},~$c\leq \beta^{-1}$. For each~$x\in V_0$, by \cref{lm:insideV0supreg,lm:insideV0degree}, $|N_{G}(x)\setminus V_0|\geq (\alpha-\varepsilon)n$ and thus there exists~$i\in [c]$ such that $|N_{G}(x)\cap V_\Gamma(C_i)|\geq \beta^2 n$. Colour each~$x\in V_0$ with such a colour~$i\in[c]$.
	
	Apply \cref{thm:Lovasz} to decompose~$G[V_0]$ into~$\ell$ paths~$P_1, \dots, P_\ell$ and~$\ell^*$ cycles. Add the~$\ell^*$ cycles to~$\cD$ and note that $|\cD|,\ell \leq \varepsilon n$. We apply the arguments of the proof of \cref{lm:closingcycles} with~$\cP_i=\{P_i\}$ for each~$i\in [\ell]$. The only difference is that we now have paths with endpoints in~$V_0$. We adapt to this setting as follows.
	
	Partition the paths~$P_1, \dots, P_\ell$ into monochromatic subpaths and bichromatic edges as in \cref{lm:closingcycles}.%
	\COMMENT{Each~$P_j$ admits a decomposition $\cD_{\rm mono}\cup \cD_{\rm bi}$,  where~$\cD_{\rm mono}$ is a set of at most~$c\leq\beta^{-1}$ monochromatic paths of distinct colours and~$\cD_{\rm bi}$ is a set of bichromatic edges such that, if~$e,e'\in \cD_{\rm bi}$ are distinct, then they are coloured with distinct pairs of colours. Thus we can decompose $P_1\cup \dots\cup P_\ell$ into
		\begin{itemize}
			\item $\ell'\leq \beta^{-1}\ell$ monochromatic subpaths~$P_1', \dots, P_{\ell'}'$; and
			\item for each~$1\leq i<i'\leq c$, a set~$\cQ_{ii'}$ of at most~$\ell$ bichromatic edges coloured with~$\{i,i'\}$.
	\end{itemize}
	By removing at most one edge from each~$\cQ_{ii'}$, we may assume~$|\cQ_{ii'}|$ is even for any~$1\leq i<i'\leq c$. Let~$\cD'$ be the set of deleted edges. Then, $|\cD'|\leq\beta^{-2}$, as desired.}
	Let~$1\leq i < i' \leq c$. We may assume that the set~$\cQ_{ii'}$ of bichromatic edges coloured with~$\{i,i'\}$ is a matching. Indeed, if~$\cQ_{ii'}$ contains distinct edges~$e$ and~$e'$ with a common endpoint, then we can delete~$e$ and~$e'$ from~$\cQ_{ii'}$ and consider~$e\cup e'$ as a monochromatic path instead. Since for each~$1\leq i < i' \leq c$, we have~$|\cQ_{ii'}|\leq \ell$, we obtain, in total, at most $\binom{c}{2}\frac{\ell}{2}\leq \sqrt{\varepsilon} n$ additional monochromatic paths.\COMMENT{I.e.\ we have $\ell'\leq \beta^{-1}\ell+\sqrt{\varepsilon} n$.}	
	
	For each~$i\in [c]$, extend each monochromatic path coloured~$i$ to a path with internal vertices in~$V_0$ and endpoints in~$V_\Gamma(C_i)$. Similarly, for each~$1\leq i<i'\leq c$, extend~$\cQ_{ii'}$ to a set of vertex-disjoint paths of length~$3$ with internal vertices in~$V_0$, an endpoint in~$V_\Gamma(C_i)$ and an endpoint in~$V_\Gamma(C_{i'})$\COMMENT{This is possible by construction of the colouring and since $\ell+\sqrt{\varepsilon}n\leq 2\sqrt{\varepsilon}n<\beta^2 n$.}. Then, one can easily show that we can proceed as in the proof of \cref{lm:closingcycles}.%

	One can easily verify that, in the end, we have covered all but at most~$\beta^{-2}$ edges of~$G[V_0]$ with at most~$\beta n$ cycles, as desired.
\end{proof}

\onlyinsubfile{\bibliographystyle{plain}
	\bibliography{Bibliography/papers}}

%% file: Exceptional_Edges.tex
\onlyinsubfile{
	\setcounter{section}{5}
	\setcounter{subsection}{3}
	\addtocounter{subsection}{-1}
	\setcounter{definition}{2}
}

\subsection{Main step of the decomposition}\label{sec:exceptionaledges}
This section corresponds to \cref{step:sketchexceptional} of the proof overview. \Cref{lm:exceptionaledges} will be used to obtain a cycle decomposition in the proof of \cref{thm:n/2}\cref{thm:n/2cycledecomp} and \cref{lm:fictive} will be used to obtain a path decomposition in the proof of \cref{thm:n/2}\cref{thm:n/2pathdecomp}. \cref{lm:connectedR} will play a similar role in the proof of \cref{thm:Delta/2}.

\begin{lm}\label{lm:exceptionaledges}
	Suppose $0<\frac{1}{n}\ll\frac{1}{k}\ll\varepsilon\ll\zeta\ll d\ll \beta \leq 1$. Let~$G, \Gamma, \Gamma'$ be edge-disjoint graphs on the same vertex set~$V$ of order~$n$. Assume~$V_0, V_1, \dots, V_k$ is a partition of~$V$ such that the following hold.
	\begin{enumerate}
		\item $V_0, V_1, \dots, V_k$ is \label{lm:exceptionaledgessupregpart}
		\begin{itemize}
			\item an $(\varepsilon, \geq d, k,m, m', R)$-superregular equalised partition of~$G$,\label{lm:exceptionaledgesG}
			\item an $(\varepsilon, \beta, k, m, m', R')$-superregular equalised partition of~$\Gamma$,
			\item an $(\varepsilon, \zeta, k, m, m', R'')$-superregular equalised partition of~$\Gamma'$.
		\end{itemize}
		\item $R'$ and~$R''$ are edge-disjoint and~$R'\cup R'' =R$. \label{lm:exceptionaledgesreduced}
		\item $G, \Gamma$, and~$\Gamma'$ have the same support clusters.\label{lm:exceptionaledgessamesupport}
		\item $V_0$ is a set of isolated vertices in~$\Gamma$ and~$\Gamma'$. Moreover,~$G[V_0]$ is empty.\label{lm:exceptionaledgesV0}
		\item Any~$x\in V\setminus V_0$ belongs to at least~$\beta k$ superregular pairs of~$\Gamma$. \label{lm:exceptionaledgesGammabadinfewpairs}
		
		\item For any~$x\in V_0$,~$d_G(x)$ is even.\label{lm:exceptionaledgeseven}
	\end{enumerate}
	Then,~$G\cup \Gamma\cup \Gamma'$ can be decomposed into edge-disjoint graphs~$G', \tGamma$, and~$H$ such that~$G, \Gamma'\subseteq G'\cup H$,~$\tGamma\subseteq \Gamma$, and the following hold.
	\begin{enumerate}[label=\upshape{(\alph*)}]
		\item $\Delta (H)\leq 13\zeta n$ and~$V_0$ is a set of isolated vertices in~$H$.\label{lm:exceptionaledgesDeltaH}
		\item $V_0, V_1, \dots, V_k$ is a~$(\zeta, \beta, k, m, m', R')$-superregular equalised partition of~$\tGamma$.\label{lm:exceptionaledgesGamma}
		\newcounter{exceptionaledges}
		\setcounter{exceptionaledges}{\value{enumi}}
		\item There exists a decomposition~$\cD\cup \cD_{\rm exc}$ of~$G'$ such that~$\cD$ is a set of at most~$\frac{n}{2}+2\beta n$ cycles and~$\cD_{\rm exc}$ is a set of at most~$\beta^{-2}$ exceptional edges.\label{lm:exceptionaledgesdecomp}
	\end{enumerate}
\end{lm}
$H$ should be thought of as a sparse ``leftover" and~$\Gamma$ is a graph which we want to use as little as possible.
An analogous result can be obtained for path decompositions.

\begin{lm}\label{lm:fictive}
	Suppose $0<\frac{1}{n}\ll\frac{1}{k}\ll\varepsilon\ll\zeta\ll d\ll \beta \leq 1$. Let~$G, \Gamma$, and~$\Gamma'$ be edge-disjoint graphs on the same vertex set~$V$ of order~$n$. Assume~$V_0, V_1, \dots, V_k$ is a partition of~$V$ such that properties \cref{lm:exceptionaledgessupregpart,lm:exceptionaledgesV0,lm:exceptionaledgessamesupport,lm:exceptionaledgeseven,lm:exceptionaledgesGammabadinfewpairs,lm:exceptionaledgesreduced} of \cref{lm:exceptionaledges} hold.
	Let~$U\subseteq V\setminus V_0$ have even size\COMMENT{$U$ will be the set of non-exceptional vertices of odd degree.}.
	Then,~$G\cup \Gamma\cup \Gamma'$ can be decomposed into edge-disjoint graphs~$G', \tGamma$, and~$H$ such that~$G, \Gamma'\subseteq G'\cup H$,~$\tGamma\subseteq \Gamma$, and properties \cref{lm:exceptionaledgesDeltaH,lm:exceptionaledgesGamma} of \cref{lm:exceptionaledges} are satisfied. Moreover,
	\begin{enumerate}[label=\upshape{(\alph*$'$)}]
		\setcounter{enumi}{\value{exceptionaledges}}
		\item there exists a path decomposition~$\cD$ of~$G'$ such that 
		$|\cD|\leq \frac{n}{2}+4\beta n$ and, 
		for any~$x\in V\setminus V_0$,~$\cD$ contains an odd number of paths with~$x$ as an endpoint if and only if~$x\in U$.\label{lm:fictivedecomp}
	\end{enumerate}
\end{lm}

The next lemma shows that stronger results can be obtained if the reduced graph~$R$ of~$G$ is assumed to be connected.

\begin{lm}\label{lm:connectedR}
	Suppose $0<\frac{1}{n}\ll\frac{1}{k}\ll\varepsilon\ll\zeta\ll d\ll \beta\leq 1$. Let~$G, \Gamma$, and~$\Gamma'$ be edge-disjoint graphs on the same vertex set~$V$ of order~$n$. Assume~$V_0, V_1, \dots, V_k$ is a partition of~$V$ such that properties \cref{lm:exceptionaledgessupregpart,lm:exceptionaledgesV0,lm:exceptionaledgessamesupport,lm:exceptionaledgeseven,lm:exceptionaledgesGammabadinfewpairs,lm:exceptionaledgesreduced} of \cref{lm:exceptionaledges} hold. Suppose furthermore that~$R$ is connected. Then, the following hold.
	\begin{enumerate}[label=\upshape{(\alph*)}]
		\item $G\cup \Gamma\cup \Gamma'$ can be decomposed into edge-disjoint graphs~$G', \tGamma$, and~$H$ such that~$G, \Gamma'\subseteq G'\cup H$,~$\tGamma\subseteq \Gamma$, and properties \cref{lm:exceptionaledgesDeltaH,lm:exceptionaledgesGamma} of \cref{lm:exceptionaledges} hold. Furthermore,~$G'$ can be decomposed into at most~$\frac{\Delta(G) }{2}+7\zeta n$ cycles.\label{lm:connectedRcycles}
		\item Let~$U\subseteq V\setminus V_0$ have even size. Then~$G\cup \Gamma\cup \Gamma'$ can be decomposed into edge-disjoint graphs~$G', \tGamma$, and~$H$ such that~$G, \Gamma'\subseteq G'\cup H$,~$\tGamma\subseteq \Gamma$, and properties \cref{lm:exceptionaledgesDeltaH,lm:exceptionaledgesGamma} of \cref{lm:fictive} hold. Furthermore, there exists a path decomposition~$\cD$ of~$G'$ such that $|\cD|\leq \max\left\{\frac{\Delta(G)}{2},\frac{|U|}{2}\right\}+8\zeta n$ and each vertex~$x\in V\setminus V_0$ is an endpoint of an odd number of paths in~$\cD$ if and only if~$x\in U$.\label{lm:connectedRpaths}
	\end{enumerate}
\end{lm}

\cref{lm:exceptionaledges,lm:fictive,lm:connectedR} will be proved simultaneously. 
To obtain a path decomposition, the idea is to insert suitable fictive edges and then construct a cycle decomposition such that each cycle in the decomposition contains exactly one fictive edge.

We will need the following result of \cite{frieze2005packing}\COMMENT{Adapted from {\cite[Theorem 3]{frieze2005packing}}.}.

\begin{thm}[{\cite{frieze2005packing}}]\label{thm:FK}
	Let~$0<\frac{1}{m}\ll \varepsilon \ll d<1$ and assume~$G$ is a bipartite graph on vertex classes~$A,B$ of size~$m$. If~$G$ is~$[\varepsilon,d]$-superregular then~$G$ contains ~$\frac{(d-19\varepsilon)m}{2}$ edge-disjoint Hamilton cycles.
	\COMMENT{\cite[Theorem 3]{frieze2005packing} with~$d-\varepsilon$ and~$3\varepsilon$ playing the roles of~$\alpha$ and~$\varepsilon$, respectively, to account for the different definitions of superregularity.}	
\end{thm}

\begin{proof}[Proof of \cref{lm:exceptionaledges,lm:fictive,lm:connectedR}]
	Define
	\[\varepsilon_1\coloneqq \varepsilon^{\frac{1}{12}}, \quad
	\varepsilon_2\coloneqq \varepsilon_1^{\frac{1}{7}}, \quad 
	\varepsilon_3\coloneqq \varepsilon_2^{\frac{1}{5}}, \quad
	\varepsilon_4\coloneqq \varepsilon_2^{\frac{1}{15}}, \quad
	\varepsilon_5\coloneqq \varepsilon_2^{\frac{1}{75}}, \quad
	\varepsilon_6\coloneqq \varepsilon_2^{\frac{1}{225}}, 
	\]
	\[\varepsilon_7\coloneqq \varepsilon_2^{\frac{1}{1125}},\quad 
	\varepsilon_8\coloneqq \varepsilon_2^{\frac{1}{3375}}, \quad
	\varepsilon_9\coloneqq \varepsilon_8^{\frac{1}{3375}}, \quad
	\zeta_1\coloneqq \sqrt{2\zeta}, \quad
	\zeta_2\coloneqq \sqrt{\zeta _1}.\]
	\NEW{Let $i,j\in [k]$ be distinct. Denote by~$V_{ij}$ and~$V_{ji}$ the support clusters of~$G[V_i, V_j]$. If~$G[V_i, V_j]$ is empty let~$d_{ij}\coloneqq 0$. Otherwise, by \cref{def:supregpartitionreduced}, \cref{def:supregpartitionsupreg}, and definition of the reduced graph, $G[V_{ij}, V_{ji}]$ is $[\varepsilon, \geq d]$-superregular and so we can let~$d_{ij}$ be a constant such that~$G[V_{ij}, V_{ji}]$ is~$[\varepsilon, d_{ij}]$-superregular.}\OLD{For any distinct~$i,j\in[k]$, denote by~$V_{ij}$ and~$V_{ji}$ the support clusters of~$G[V_i, V_j]$. If~$G[V_i, V_j]$ is empty let~$d_{ij}\coloneqq 0$, otherwise, let~$d_{ij}$ be a constant such that~$G[V_{ij}, V_{ji}]$ is~$[\varepsilon, d_{ij}]$-superregular.}
	We let~$G'$ and~$H$ be empty graphs on~$V$. Throughout this proof, we will repeatedly add edges to~$G'$ and~$H$, and, whenever we do so, these edges are deleted from~$G\cup \Gamma \cup \Gamma'$. 
	
	Let~$G_{\rm exc}\coloneqq \emptyset$. For each connected component~$C$ of~$R$ and~$x\in V_0$, if~$|N_G(x)\cap V_G(C)|$ is odd, add exactly one edge of~$G[\{x\}, V_G(C)]$ to~$G_{\rm exc}$. Delete the edges in~$G_{\rm exc}$ from~$G$. Observe that we may now assume that for each connected component~$C$ of~$R$, any $x\in V_0$ has even degree in~$G[V_0\cup V_G(C)]$. The graph~$G_{\rm exc}$ will be covered in \cref{step:Rdisconnected}. Observe that, by \cref{lm:exceptionaledgeseven},~$G_{\rm exc}$ is empty in the proof of \cref{lm:connectedR}.
	
	We now assume~$R$ is connected. If it is disconnected, we will apply \cref{step:fictive,step:exceptionaledges,step:partitioning,step:equalisingsupport,step:paths,step:setsofpaths,step:tyingpaths} to each connected component of~$R$ separately and then cover the potentially remaining edges (i.e.\ the edges of~$G_{\rm exc}$) in \cref{step:Rdisconnected}. Fix~$\Delta\coloneqq \Delta(G)$, $\Delta_0\coloneqq \max\{d_G(x)\mid x\in V_0\}$, and $\Delta'\coloneqq \max\{d_G(x)\mid x\notin V_0\}$. In particular, in what follows,~$\Delta, \Delta'$, and~$\Delta_0$ are left unchanged when we delete some edges from~$G$.
	 
	\begin{steps}
		\item \textbf{Partitioning the edges of~$G$ and constructing reservoirs of vertices.}\label{step:partitioning} We will partition each superregular pair of~$G$ into subgraphs of small comparable density. Each subgraph will be assigned a reservoir, that is a small number of vertices that will be set aside to tie paths together later on. To do so, we will partition each cluster into small subclusters of equal size and, in each subgraph, one of these subclusters will play the role of the reservoir.

		Let~$r\coloneqq \left\lfloor\zeta ^{-1}\right\rfloor$ ($r$ will be the number of reservoirs). For each~$ij\in E(R)$, let $\ell_{ij}\coloneqq \left\lfloor\zeta^{-1}d_{ij}\right\rfloor$ and apply \cref{lm:edgepartition} to partition~$G[V_{ij}, V_{ji}]$ into~$r\ell_{ij}$ spanning edge-disjoint~$[\varepsilon_1, \zeta^2]$-superregular graphs $G_{1,1}^{ij}, \dots,G_{1,\ell_{ij}}^{ij}, G_{2,1}^{ij}, \dots, G_{r, \ell_{ij}}^{ij}$ and a leftover graph~$G_0^{ij}$ which we add to~$H$. Note that, for each~$ij\in E(R)$, $\Delta(G_0^{ij})\leq (d_{ij}+\varepsilon)m-(\zeta^2-\varepsilon)mr\ell_{ij}\leq 3\zeta m$ and thus 
		\begin{equation}\label{eq:Delta1}
			\Delta(H)\leq 3\zeta n.
		\end{equation}

		For each~$i\in[k]$, randomly partition~$V_i$ into~$r$ subclusters~$V_i^1, \dots, V_i^r$ of size~$\zeta m$. (If $\zeta m\notin \mathbb{N}^*$, then the subclusters will only have sizes roughly~$\zeta m$, but this does not affect the argument below.) Let~$V^\ell\coloneqq \bigcup_{i\in[k]} V_i^\ell$. Also define~$U^\ell\coloneqq V^\ell \cap U$ for the proof of \cref{lm:fictive,lm:connectedR}\cref{lm:connectedRpaths}. We claim that the following hold with positive probability.
		\begin{enumerate}[label=(\arabic*)]
			\item For any~$ij\in E(R)$,~$\ell\in [r]$, and~$\ell'\in [\ell_{ij}]$,\label{subgraphssupreg}
			\begin{enumerate}[label=(\arabic{enumi}\alph*)]
				\item $G_{\ell, \ell'}^{ij}[X, Y]$ is~$[\varepsilon_2, \zeta^2]$-superregular for each \NEW{$X\in \{V_{ij}\setminus V_i^\ell, V_{ij}\cap V_i^\ell\}$ and $Y\in \{V_{ji}\setminus V_j^\ell, V_{ji}\cap V_j^\ell\}$}\OLD{$X\in \{V_{ij}\setminus V_i^\ell, V_{ij}\cap V^\ell\}$ and $Y\in \{V_{ji}\setminus V_i^\ell, V_{ji}\cap V^\ell\}$};\label{G}
				\item if~$ij\in E(R')$, then $\Gamma_{ij \ell}:=\Gamma[V_{ij}\cap V_i^\ell, V_{ji}\cap V_j^\ell]$ is~$[\varepsilon_2, \beta]$-superregular;\label{Gammaij}
				\item if~$ij\in E(R'')$, then $\Gamma_{ij\ell}':=\Gamma'[V_{ij}\cap V_i^\ell, V_{ji}\cap V_j^\ell]$ is~$[\varepsilon_2, \zeta]$-superregular.\label{Gamma'ij}
			\end{enumerate}
			\item For any~$ij\in E(R)$ and~$\ell \in [r]$, $|V_{ij}\cap V_i^\ell|=(1\pm \varepsilon_1)\zeta m'$.%
			\COMMENT{Recall~$m'=|V_{ij}|$.}\label{reservoirsupport}
			\item For each~$\ell\in [r]$ and~$x\in V_0$, $\left|N_G(x)\cap V^\ell\right|=\frac{|N_G(x)|}{r}\pm \varepsilon n$.\label{V0neighbours}
		\end{enumerate}
		Additionally, for the proof of \cref{lm:fictive,lm:connectedR}\cref{lm:connectedRpaths}, the following holds with positive probability.
		\begin{enumerate}[resume,label=(\arabic*)]
			\item For each~$\ell\in [r]$, $|U^\ell|=\frac{|U|}{r}\pm \varepsilon n$.\label{odd}
		\end{enumerate}
		Indeed, by \cref{lm:vertexpartition}, \cref{subgraphssupreg,reservoirsupport} hold with high probability, and, a simple application of \cref{lm:Chernoff}%
		shows that \cref{V0neighbours,odd} hold with high probability. Therefore, by a union bound, we may assume that \cref{subgraphssupreg,odd,reservoirsupport,V0neighbours} are satisfied.

		For each~$ij\in E(R), \ell\in [r]$, and~$\ell'\in [\ell_{ij}]$, define $\tG_{\ell, \ell'}^{ij}\coloneqq G_{\ell, \ell'}^{ij}[V_{ij}\setminus V_i^\ell, V_{ji}\setminus V_i^\ell]$ and add all edges of $G_{\ell, \ell'}^{ij}[V_{ij}\cap V_i^\ell,V_{ji}\cap V_j^\ell]$ to~$H$. By \cref{G}, we add, in total, at most $(\zeta^2+\varepsilon_2)\zeta m \cdot \zeta^{-1}\cdot k\leq 2\zeta^2 n$ edges incident to each vertex, so, by \cref{eq:Delta1},
		\begin{equation}\label{eq:Delta1'}
		\Delta(H)\leq 4\zeta n.
		\end{equation}

		\item \textbf{Equalising the support cluster sizes.}\label{step:equalisingsupport} 
		For any distinct~$i,j\in [k]$ and~$\ell\in [r]$ in turn, we now construct a subset $V_{ij\ell}\subseteq  V_{ij}\setminus V_i^{\ell}$ of size $m''\coloneqq (1-\zeta-\varepsilon)m'$ by removing exactly $|V_{ij}\setminus V_i^{\ell}|-m''$ vertices.
		We build these sets one by one and, in each step, we only remove vertices which have already been removed fewer than~$\sqrt{\varepsilon_1}k$ times in the construction so far.
		This is possible since in each step, by \cref{reservoirsupport}, we need to remove at most~$2\varepsilon_1 m$ vertices%
		\COMMENT{Since by \cref{reservoirsupport}, we have $|V_{ij}\setminus V_i^{\ell}|\leq (1-\zeta+\varepsilon_1\zeta)m'$. Thus, we need to remove at most $(\varepsilon_1\zeta+\varepsilon)m'\leq 2\varepsilon_1 m$ vertices.},
		and so, in each step, there are at most\COMMENT{$\leq \frac{2\varepsilon_1 m \cdot k \cdot r}{\sqrt{\varepsilon_1}k}=\frac{2\sqrt{\varepsilon_1}m}{\zeta}\leq \zeta m$}~$\zeta m$ vertices which \NEW{we} are not allowed to remove anymore. On the other hand, by \cref{reservoirsupport}, $|V_{ij}\setminus V_i^\ell|\geq (1-2\zeta)m$ for any distinct~$i,j\in [k]$, and~$\ell\in [r]$%
		\COMMENT{We have $|V_{ij}\setminus V_i^\ell|\geq m'-\zeta m' (1+\varepsilon_1)\geq (1-\zeta -\zeta\varepsilon_1)(1-\varepsilon)m=(1-\varepsilon-\zeta+\zeta \varepsilon-\zeta\varepsilon_1+\zeta\varepsilon_1\varepsilon)m\geq (1-2\zeta)m$\label{comment:sizesupportreservoir}.}

		For any $ij\in E(R), \ell\in [r]$, and $\ell'\in [\ell_{ij}]$, define $\hG_{\ell, \ell'}^{ij}\coloneqq \tG_{\ell, \ell'}^{ij}[V_{ij\ell}, V_{ji\ell}]$. By \cref{G} and \cref{lm:verticesedgesremoval},~$\hG_{\ell, \ell'}^{ij}$ is $[2\sqrt{\varepsilon_2}, \zeta^2]$-superregular%
		\COMMENT{We remove at most $2\varepsilon_1 m\leq \varepsilon_2 |V_{ij}\setminus V_i^\ell|$ vertices from each cluster (as $|V_{ij}\setminus V_{ij}^\ell|\geq (1-2\zeta)m$ (see \cref{comment:sizesupportreservoir})).}.
		Add to~$H$ all edges of $\tG_{\ell, \ell'}^{ij}\setminus \hG_{\ell, \ell'}^{ij}$. Since~$\hG_{\ell, \ell'}^{ij}$ is obtained from $\tG_{\ell, \ell'}^{ij}$ by deleting at most~$2\varepsilon_1 m$ vertices from each cluster, for each $x\in V_{ij\ell}$,~$x$ has degree at most~$2\varepsilon_1 m$ in $\tG_{\ell, \ell'}^{ij}\setminus \hG_{\ell, \ell'}^{ij}$. Moreover, for each~$i\in [k]$ and~$x\in V_i$, there are at most~$\sqrt{\varepsilon_1}k$ pairs $(j, \ell)\in [k]\times [r]$ such that $x\in V(\tG_{\ell, \ell'}^{ij})\setminus V(\hG_{\ell, \ell'}^{ij})$. Thus, we have added to~$H$ at most~$2\sqrt{\varepsilon_1}n$ edges incident to each vertex and, thus, by \cref{eq:Delta1'}, 
		\begin{equation}\label{eq:Delta2}
		\Delta(H)\leq 5\zeta n.
		\end{equation}

		\item \textbf{Decomposing non-exceptional edges of~$G$ into long paths with endpoints in reservoirs.}\label{step:paths} For each~$ij\in E(R), \ell\in [r]$, and~$\ell'\in [\ell_{ij}]$, apply \cref{thm:FK} with $\hG_{\ell, \ell'}^{ij}, \zeta^2$, and~$2\sqrt{\varepsilon_2}$ playing the roles of~$G, d$, and~$\varepsilon$ to obtain a set~$\cH_{\ell,\ell'}^{ij}$ of $h\coloneqq \frac{(\zeta^2-38\sqrt{\varepsilon_2})m''}{2}$ edge-disjoint Hamilton cycles of~$\hG_{\ell, \ell'}^{ij}$. We turn each cycle in~$\cH_{\ell,\ell'}^{ij}$ into a path one by one by deleting an edge~$xy$ such that no edge incident to~$x$ or~$y$ has already been deleted from~$\cH_{\ell,\ell'}^{ij}$. This is possible since $|\cH_{\ell,\ell'}^{ij}|\leq \frac{\zeta^2 m''}{2}$ and each cycle in~$\cH_{\ell,\ell'}^{ij}$ is of length~$2m''$. We add all these edges as well as all the edges in $E(\hG_{\ell, \ell'}^{ij})\setminus E(\cH_{\ell, \ell'}^{ij})$ to~$H$. Thus, we add at most $(1+40\sqrt{\varepsilon_2}m'') r\zeta^{-1}k\leq \zeta^2 n$ edges incident to each vertex, and so, by \cref{eq:Delta2},\COMMENT{Recall~$G[V_i, V_j]$ is decomposed into $r\ell_{ij}=r\zeta^{-1}d_{ij}\leq r\zeta^{-1}$ subgraphs. Since~$\hG_{\ell, \ell'}^{ij}$ is $[2\sqrt{\varepsilon_2}, \zeta^2]$-superregular, there are at most $(\zeta^2+2\sqrt{\varepsilon_2})m''-(\zeta^2-38\sqrt{\varepsilon_2})m''$ edges incident to each vertex in $E(\hG_{\ell, \ell'}^{ij})\setminus E(\cH_{\ell, \ell'}^{ij})$.}
		\begin{equation}\label{eq:Delta3}
		\Delta(H)\leq 6\zeta n.
		\end{equation}	
		
		We now extend the paths in~$\cH_{\ell,\ell'}^{ij}$ to paths with internal vertices in~$V_{ij\ell}\cup V_{ji\ell}$ and endpoints in~$V_i^\ell\cup V_j^\ell$ one by one as follows. Given an~$(x,y)$-path~$P$ in~$\cH_{\ell,\ell'}^{ij}$ with~$x\in V_{ij\ell}$ and~$y\in V_{ji\ell}$, pick~$x'\in V_j^\ell$ and~$y'\in V_i^\ell$ such that $xx', yy'\in E(G_{\ell, \ell'}^{ij})$ and~$\cH_{\ell,\ell'}^{ij}$ contains fewer than $\varepsilon^2 m$ paths with~$x'$ as endpoint and similarly for~$y'$.
		Replace~$P$ in~$\cH_{\ell,\ell'}^{ij}$ by the path~$x'xPyy'$ and delete~$xx', yy'$ from~$G_{\ell, \ell'}^{ij}$. Note that the existence of~$x'$ and~$y'$ is guaranteed by \cref{G,reservoirsupport}, and the fact that $|\cH_{\ell,\ell'}^{ij}|\leq \frac{\zeta^2 m}{2}$%
		\COMMENT{By \cref{G,reservoirsupport},~$G_{\ell, \ell'}^{ij}$ originally contained at least $(\zeta^2-\varepsilon_2)(1-\varepsilon_1)\zeta m'\geq\frac{\zeta^3 m}{2}$ edges incident to~$x$ with an endpoint in~$V_j^\ell$. Since $|\cH_{\ell,\ell'}^{ij}|\leq \frac{\zeta^2 m}{2}$, at most \[\frac{2\frac{\zeta^2 m}{2}}{\varepsilon^2 m}<\frac{\zeta^3 m}{2}\] of these edges cannot be used anymore. The same holds for~$y$ by \cref{G}.}. 

		Once all paths in~$\cH_{\ell,\ell'}^{ij}$ have been extended as above, add all remaining edges of $G_{\ell, \ell'}^{ij}[V_{ij}\cap V_i^\ell, V_{ji}\setminus V_j^\ell]$ and $G_{\ell, \ell'}^{ij}[V_{ij}\setminus V_i^\ell,  V_{ji}\cap V_j^\ell]$ to~$H$. Note that by \cref{G}, if~$x\in V_i^\ell$, then $G_{\ell, \ell'}^{ij}[V_{ij}\cap V_i^\ell, V_{ji}\setminus V_j^\ell]$ contains at most\COMMENT{$(\zeta^2+\varepsilon_2)m$}~$2\zeta^2 m$ edges incident to~$x$, and if~$x\in V_i\setminus V_i^\ell$, $G_{\ell, \ell'}^{ij}[V_{ij}\setminus V_i^\ell,  V_{ji}\cap V_j^\ell]$ contains at most\COMMENT{$(\zeta^2+\varepsilon_2)\zeta m$}~$2\zeta^3 m$ edges incident to~$x$. This holds for any~$ij\in E(R), \ell\in[r]$, and~$\ell'\in [\ell_{ij}]$ so, in total, we have added to~$H$ at most~$4\zeta n$ edges incident to each vertex%
		\COMMENT{at most $2\zeta^2m\cdot k\cdot \zeta^{-1}+2\zeta^3m\cdot k \cdot (r-1)\cdot \zeta^{-1}\leq 4\zeta n$, since~$\ell_{ij}\leq \zeta^{-1}$ for each~$ij\in E(R)$}
		and thus, by \cref{eq:Delta3} 
		\begin{equation}\label{eq:Delta3'}
		\Delta(H)\leq 10\zeta n.
		\end{equation}
		 Moreover, for each~$ij\in E(R)$,~$\ell\in [r]$, and~$\ell'\in [\ell_{ij}]$, add all edges of~$\cH_{\ell, \ell'}^{ij}$ to~$G'$. Note that all (remaining) edges of~$G$ now have exactly one endpoint in~$V_0$. 
		
		\item \textbf{Combining the paths into sets of vertex-disjoint paths.}\label{step:setsofpaths} Let~$\hR$ be the multigraph obtained from~$R$ by replacing each edge~$ij\in E(R)$ by~$\ell_{ij}$ parallel edges denoted $e_{ij}^1, \dots, e_{ij}^{\ell_{ij}}$.	
		\begin{claim}\label{claim:DeltatR}
			$\Delta(\hR)\leq \frac{\Delta'}{\zeta m}+\sqrt{\varepsilon}k$.\COMMENT{Similar to \cite[Lemma 13]{christofides2012edge}.}
		\end{claim}
		
		\begin{proofclaim}
			Let~$i\in [k]$ and \NEW{recall that the graphs $G$, $G'$, and $H$ at the end of \cref{step:partitioning} form a decomposition of the original graph G.}\OLD{recall that~$G, G', H$ forms a decomposition of the original graph~$G$.} Clearly, we have
			\[e_{G\cup G'\cup H}(V_i, V\setminus V_0)\leq \Delta' m.\]
			Moreover, by \cref{lm:exceptionaledgessupregpart}, 
			\[e_{G\cup G'\cup H}(V_i, V\setminus V_0)\geq \sum_{j\neq i}(d_{ij}-\varepsilon)(1-\varepsilon)^2m^2\geq  \sum_{j\neq i}(\zeta \ell_{ij}-4\varepsilon)m^2\geq\zeta d_{\hR}(V_i) m^2-4\varepsilon m^2 k \]
			\COMMENT{so $d_{\hR}(V_i)\leq \frac{\Delta '}{\zeta m}+\frac{4\varepsilon k}{\zeta}$}and the claim holds.
		\end{proofclaim}
		Note that there are at most~$\zeta^{-1}$ parallel edges between any two vertices of~$\hR$. Thus, we can apply \cref{claim:DeltatR} and \cref{thm:Vizing} to fix a decomposition~$\cD_{\hR}$ of~$\hR$ into at most $\frac{\Delta'}{\zeta m}+ 2\sqrt{\varepsilon}k$ matchings. 
		For each~$M\in \cD_{\hR}$ and each~$\ell\in [r]$, we decompose $\bigcup_{e_{ij}^{\ell'}\in M} \cH_{\ell, \ell'}^{ij}$ into~$h$\COMMENT{Recall that~$h$ is the number of paths in each~$\cH_{\ell, \ell'}^{ij}$.} disjoint sets of paths containing exactly one path of~$\cH_{\ell, \ell'}^{ij}$ for each~$e_{ij}^{\ell'}\in M$. 

		Let~$\sP$ be the collection of all the~$|\cD_{\hR}|rh$ linear forests obtained. Note that for each $x\in V\setminus V_0$, all non-exceptional edges incident to~$x$ are covered by paths in $\bigcup \sP$, apart from those lying in~$H$. Thus,\COMMENT{The first inequality follows from definition of~$\Delta'$,~$|V_0|\leq \varepsilon n$, and~$\Delta(H)\leq 10\zeta n$.} $\frac{\Delta'}{2}-6\zeta n\leq |\sP|\leq rh\left(\frac{\Delta'}{\zeta m}+ 2\sqrt{\varepsilon}k\right)\leq \frac{\Delta'}{2}+\sqrt{\varepsilon}n$. We also note that for each~$\cP\in \sP$ there exists~$r(\cP)\in [r]$ such that~$\cP$ is a set of vertex-disjoint paths with endpoints in~$V^{r(\cP)}$ and internal vertices in~$V\setminus V^{r(\cP)}$. For each~$\ell\in [r]$, let $\sP_\ell\coloneqq \{\cP\in \sP \mid r(\cP)=\ell\}$. \NEW{By construction, $|\sP_1|=\dots=|\sP_r|=\frac{|\sP|}{r}$.}
		
		\item \textbf{Including exceptional edges.}\label{step:exceptionaledges} For each~$\ell\in [r]$, we add exceptional edges to the linear forests in~$\sP_\ell$ as follows. If possible, pick~$\cP\in  \sP_\ell$ such that we have not yet added exceptional edges to~$\cP$ and such that~$G$ contains a set of paths~$\cP_{\rm exc}$ satisfying the following.
		\begin{enumerate}[label=(\Roman*)]
			\item $\cP_{\rm exc}$ is a set of vertex-disjoint paths of~$G$ of length of~$2$.\label{exclength}
			\newcounter{exc}
			\setcounter{exc}{\value{enumi}}
			\item The paths in~$\cP_{\rm exc}$ have their endpoints in~$V^\ell\setminus V(\cP)$ and internal vertex in~$V_0$. \label{excpaths}
			\item $V(\cP_{\rm exc})\cap V_0$ is the set of vertices~$x\in V_0$ such that~$|N_G(x)\cap V^\ell|$ is maximum.\label{excmaxdegree}
			\newcounter{exc2}
			\setcounter{exc2}{\value{enumi}}
			\item $|V(\cP_{\rm exc})\cap V_i^\ell|\leq \zeta^2 m$ for each~$i\in [k]$.\label{excroom}
		\end{enumerate} 
		Fix such a set~$\cP_{\rm exc}$ and add the paths in~$\cP_{\rm exc}$ to~$\cP$. Add the edges of~$\cP_{\rm exc}$ to~$G'$. We repeat this procedure until there is no such~$\cP$. Then, we claim that the following holds.
		
		\begin{claim}\label{claim:excdegree}
			For each~$x\in V_0$, $d_{G}(x)\leq \max\{\Delta_0-\Delta'+13\zeta n, \zeta n\}$.
		\end{claim}
		
		\begin{proofclaim}
			\NEW{Note that it is enough to show that for each $\ell\in [r]$ and $x\in V_0$, we have $|N_G(x)\cap V^\ell|\leq \max\{\frac{1}{r}(\Delta_0-\Delta'+13\zeta n), \sqrt{\varepsilon}n\}$.}
			
			\NEW{Let $\ell\in [r]$. Suppose first that we have added a set $\cP_{\rm exc}$ satisfying \cref{exclength,excmaxdegree,excpaths,excroom} to each~$\cP\in \sP_\ell$.
			By \cref{V0neighbours}, each vertex in $V_0$ initially had at most $\frac{\Delta_0}{r}+\varepsilon n$ neighbours in $V^\ell$. Thus, by \cref{excpaths,excmaxdegree}, we now have $|N_G(x)\cap V^\ell|\leq \frac{\Delta_0}{r}+\varepsilon n-2|\sP_\ell|\leq \frac{\Delta_0}{r}+\varepsilon n-\frac{2}{r}(\frac{\Delta'}{2}-6\zeta n)\leq \frac{1}{r}(\Delta_0-\Delta'+13\zeta n)$ for each $x\in V_0$.}
			
			\NEW{Suppose now that there exists $\cP\in \sP_\ell$ which does not contain any exceptional edges.
			We claim that for any~$x\in V_0$, we have $|N_G(x)\cap V^\ell|\leq \sqrt{\varepsilon} n$. Suppose not. Let $x_1, \dots, x_s$ be an enumeration of the vertices $x\in V_0$ such that $|N_G(x)\cap V^\ell|$ is maximum. By assumption, we have $|N_G(x_i)\cap V^\ell|>\sqrt{\varepsilon}n$. Suppose inductively that, for some $0\leq i\leq s$, we have constructed a set of paths $\cP_{\rm exc}^i$ satisfying \cref{exclength,excpaths,excroom} and such that $|V(\cP_{\rm exc}^i)\cap V_0|=\{x_j\mid j\in [i]\}$.
			If $i<s$, we construct $\cP_{\rm exc}^{i+1}$ as follows. Let $X$ be the set of indices $j\in [\ell]$ such that $|V(\cP_{\rm exc}^i)\cap V_j^\ell|=\zeta^2m$.
			Let $Y\coloneqq (N_G(x_{i+1})\cap V_\ell)\setminus (V(\cP)\cup V(\cP_{\rm exc}^i)\cup \bigcup_{j\in X}V_j^\ell)$.
			Then, $|Y|\geq \sqrt{\varepsilon}n-k-2i-\frac{2i}{\zeta^2m}\cdot \zeta m\geq 2$. Let $y,y'\in Y$ be distinct. Set $\cP_{\rm exc}^{i+1}\coloneqq \cP_{\rm exc}^i\cup \{yx_{i+1}y'\}$. By construction, $\cP_{\rm exc}^{i+1}$ satisfies \cref{exclength,excpaths,excroom} and $|V(\cP_{\rm exc}^{i+1})\cap V_0|=\{x_j\mid j\in [i+1]\}$.
			Therefore, we can construct a set $\cP_{\rm exc}^s$ satisfying \cref{exclength,excmaxdegree,excpaths,excroom}, a contraction. Thus, $|N_G(x)\cap V^\ell|\leq \sqrt{\varepsilon} n$ for all~$x\in V_0$.}%
		\end{proofclaim} 
			 
		We now split most of the remaining exceptional edges into sets of vertex-disjoint paths, in a similar way as above. Let~$\sP'\coloneqq \emptyset$. Assume there exists a set of paths~$\cP_{\rm exc}$ satisfying \cref{exclength} and the following.
		\begin{enumerate}[label=(\Roman*$'$)]
			\setcounter{enumi}{\value{exc}}
			\item The paths in~$\cP_{\rm exc}$ have their endpoints in~$V\setminus V_0$ and internal vertex in~$V_0$. \label{excpaths'}
			\item $V(\cP_{\rm exc})\cap V_0$ is the set of vertices~$x\in V_0$ such that~$|N_G(x) \cap V\setminus V_0|$ is maximum.\label{excmaxdegree'}
			\item $|V(\cP_{\rm exc})\cap V_i|\leq \zeta m$ for each~$i\in [k]$.\label{excroom'}
		\end{enumerate}
		Then add such a set~$\cP_{\rm exc}$ to~$\sP'$ and add the edges of~$\cP_{\rm exc}$ to~$G'$.
		We repeat this procedure until we cannot find any~$\cP_{\rm exc}$ as above. Then, one can show using similar arguments as in the proof of \cref{claim:excdegree} that~$d_G(x)\leq \sqrt{\varepsilon}n$ for each~$x\in V_0$.%
		\COMMENT{Assume~$G$ still contains some exceptional edges. Assume for a contradiction that there exists $x\in V_0$ with $d_G(x)\geq \sqrt{\varepsilon}n$. Then, we can greedily construct a set~$\cP_{\rm exc}$ satisfying \cref{exclength}, \cref{excpaths',excmaxdegree',excroom'}. (At most $2\varepsilon n+\frac{2\varepsilon n}{\zeta m}\cdot m\leq \frac{\sqrt{\varepsilon }n}{2}$ forbidden vertices at each stage.)}

		We will now form at most~$4\sqrt{\varepsilon} n$ sets of linear forests which cover all remaining exceptional edges of~$G$. Assume that~$\cP_1, \dots, \cP_{4\sqrt{\varepsilon} n}$ are (possibly empty) edge-disjoint sets of paths such that for each~$i\in[4\sqrt{\varepsilon}n]$ the following are satisfied and such that $\sum_{i\in [4\sqrt{\varepsilon} n]}|\cP_i|$ is maximal.
		\begin{enumerate}[label=(\Roman*$''$)]
			\item $\cP_i$ is a set of vertex-disjoint paths of~$G$ of length~$2$. \label{exclength''}
			\item The paths in~$\cP_i$ have their endpoints in~$V\setminus V_0$ and internal vertex in~$V_0$. \label{excpaths''}
			\setcounter{enumi}{\value{exc2}}
			\item $|V(\cP_i)\cap V_j|\leq \zeta m$ for each~$j\in [k]$.\label{excroom''}
		\end{enumerate}
		Add all edges of $\cP_1\cup \dots\cup \cP_{4\sqrt{\varepsilon} n}$ to~$G'$. If~$d_G(x)=0$ for all~$x\in V_0$, add~$\cP_1, \dots, \cP_{4\sqrt{\varepsilon} n}$ to~$\sP'$ and we are done. We may therefore assume that there is~$x\in V_0$ with~$d_G(x)\geq 2$. Pick distinct~$y,z\in N_G(x)$ and let~$i,i'\in [k]$ be such that~$y\in V_i$ and~$z\in V_{i'}$. By maximality, we only need to find~$j\in[4\sqrt{\varepsilon} n]$ such that \cref{exclength'',excpaths'',excroom''} are still satisfied if we add~$yxz$ to~$\cP_j$ and thus obtain a contradiction. By construction,~$x$ belongs to fewer than~$\sqrt{\varepsilon} n$ of the~$\cP_j$\COMMENT{as~$d_G(x)\leq \sqrt{\varepsilon} n$}, and, since~$|V_0|\leq \varepsilon n$, each of~$y$ and~$z$ belong to fewer than~$\varepsilon n$ of the~$\cP_j$. Moreover, there are at most $\frac{\varepsilon n  m}{\zeta m -2}\leq \sqrt{\varepsilon} n$ indices~$j\in [4\sqrt{\varepsilon}n]$ such that~$|V(\cP_j)\cap V_i|> \zeta m -2$ and similarly, there are at most~$\sqrt{\varepsilon} n$ indices~$j\in [4\sqrt{\varepsilon} n]$ such that $|V(\cP_j)\cap V_{i'}|> \zeta m -2$. Thus, there are at least $4\sqrt{\varepsilon} n- \sqrt{\varepsilon} n- 2\varepsilon n-2\sqrt{\varepsilon}n>0$ indices~$j$ such that we can add the path~$yxz$ to~$\cP_j$ and \cref{exclength'',excpaths'',excroom''} are still satisfied.

		To summarise, we have constructed sets~$\sP_\ell$, for each~$\ell\in [r]$, and a set~$\sP'$ satisfying the following. 
		\begin{enumerate}[label=(\Alph*)]
			\item For each~$\ell\in [r]$ and~$\cP\in \sP_\ell$,~$\cP$ is a set of vertex-disjoint paths with endpoints in~$V^\ell$ and internal vertices in~$V\setminus V_\ell$. Moreover, $|\cP|\leq \frac{k}{2}+\varepsilon n$ and $|V(\cP)\cap V_i^\ell|\leq k+\zeta^2 m$ for each~$i\in [k]$.\label{sP}
			\item For each~$\cP\in \sP'$,~$\cP$ is a set of vertex-disjoint paths with endpoints in~$V\setminus V_0$ and internal vertices in~$V_0$. Moreover,~$|\cP|\leq \varepsilon n$  and~$|V(\cP)\cap V_i|\leq \zeta m$ for each~$i\in [k]$. \label{sP'}
			\item For each~$x\in V\setminus V_0$, there are at most~$2\varepsilon n$ paths in $\bigcup (\sP\cup \sP')$ which have~$x$ as an endpoint (recall $\sP=\sP_1 \cup \dots \cup \sP_r $)\COMMENT{Recall that in \cref{step:paths},~$\cH_{\ell, \ell'}^{ij}$ contains, for each~$x\in V_i^\ell$, at most~$\varepsilon^2 m$ paths with~$x$ as an endpoint. This is for any~$ij\in E(R)$ and $1\leq \ell'\leq \ell_{ij}=\zeta^{-1}d_{ij}\leq \zeta^{-1}$ so we have at most~$\varepsilon n$ non-exceptional paths with~$x$ as an endpoint. Additionally, since~$|V_0|\leq \varepsilon n$, there are at most~$\varepsilon n$ exceptional paths with~$x$ as an endpoint.}.\label{endpoints}
			\item By \cref{step:setsofpaths}, \cref{claim:excdegree}, and the above construction, $|\sP\cup \sP'|\leq \frac{\Delta}{2}+7\zeta n$%
			\COMMENT{By \cref{claim:excdegree} and the above procedure, we have $|\sP\cup \sP'|\leq \frac{\Delta'}{2}+\sqrt{\varepsilon}n+\frac{\Delta_0-\Delta'+13\zeta n}{2}+4\sqrt{\varepsilon} n\leq \frac{\Delta_0}{2}+7\zeta n\leq \frac{\Delta}{2}+7\zeta n$, or $|\sP\cup \sP'|\leq \frac{\Delta'}{2}+\sqrt{\varepsilon}n+\frac{\zeta n}{2}+4\sqrt{\varepsilon} n\leq \frac{\Delta'}{2}+\zeta n\leq \frac{\Delta}{2}+7\zeta n$.}.
			Moreover,~$G$ is now empty.\label{numbersets}
		\end{enumerate}
		Indeed, in order to check \cref{endpoints}, recall from \cref{step:paths} that, for each $x\in V_i^\ell$, each $\cH_{\ell, \ell'}^{ij}$ contains at most $\varepsilon^2 m$ paths with~$x$ as an endpoint.
		
		\item \textbf{Including fictive edges.}\label{step:fictive} We ignore this step for the proof of \cref{lm:exceptionaledges,lm:connectedR}\cref{lm:connectedRcycles}. For the proof of \cref{lm:fictive,lm:connectedR}\cref{lm:connectedRpaths}, we construct a multiset~$E_{\rm fict}$ of fictive edges on~$V\setminus V_0$. As discussed in \cref{sec:robust}, we view edges in~$E_{\rm fict}$ as distinct from each other and from edges in~$G, G', \Gamma, \Gamma'$, and~$H$.
		We will add a fictive edge to each linear forest in~$\sP\cup \sP'$. Moreover, in order to satisfy \cref{lm:fictivedecomp}, we make sure that for any~$x\in V\setminus V_0$,~$E_{\rm fict}$ contains an odd number of edges incident to~$x$ (counting multiplicity) if and only if~$x\in U$.

		Start with~$E_{\rm fict}=\emptyset$. In what follows, we denote by~$U_{\rm even}$ the set of vertices in~$U$ which are incident to an even number of edges in~$E_{\rm fict}$ and, for each~$\ell\in [r]$, we denote by~$U_{\rm even}^\ell$ the set~$U_{\rm even}\cap U^\ell$.\COMMENT{Our goal is to have~$U_{\rm even}=\emptyset$ at the end of the procedure.}
		In what follows, we will update~$U_{\rm even}$ and~$U_{\rm even}^\ell$ at each step of our algorithm.
		For each~$\ell\in [r]$, we add a fictive edge to each linear forest in~$\sP_\ell$ as follows. Assume~$\cP\in \sP_\ell$ does not contain a fictive edge yet. If there exist distinct $x,y\in U_{\rm even}^\ell\setminus V(\cP)$, add the edge~$xy$ to~$E_{\rm fict}$ and to~$\cP$. If there are no such~$x$ and~$y$, then note that, by \cref{sP}, $|U_{\rm even}^\ell|\leq |V(\cP)\cap V^\ell|+1\leq 2|\cP|+1\leq 3\varepsilon n$ and proceed as follows. If~$\cP$ is the only linear forest in~$\sP_\ell$ which does not contain a fictive edge, we remove~$\cP$ from~$\sP_\ell$, add all its edges to~$H$ and we are done. We note that this increases the maximum degree of~$H$ by at most~$2$. Otherwise, pick $\cP'\in \sP_\ell\setminus\{\cP\}$ such that~$\cP'$ does not contain a fictive edge. Note that, by \cref{sP}, $|V^\ell\cap V(\cP\cup \cP')|\leq \sqrt{\varepsilon} n$. Moreover, by \cref{numbersets}, there are at most $\frac{2|\sP_\ell|}{\varepsilon n}\leq 2\varepsilon^{-1}$ vertices in~$V_\ell$ which are incident to at least~$\varepsilon n$ edges in~$E_{\rm fict}$. Thus,\COMMENT{since~$|V^\ell|=\zeta mk$} we can choose distinct $x,y\in V^\ell\setminus V(\cP\cup \cP')$ such that~$E_{\rm fict}$ contains fewer than~$\varepsilon n$ edges incident to~$x$ and fewer than~$\varepsilon n$ edges incident to~$y$. Add the edge~$xy$ to both~$\cP$ and~$\cP'$, and, add two edges between~$x$ and~$y$ to~$E_{\rm fict}$. We repeat this procedure until each linear forest in~$\sP_\ell$ contains a fictive edge.

		We then proceed similarly to add a fictive edge to each linear forest in~$\sP'$%
		\COMMENT{Assume~$\cP\in \sP'$ does contain a fictive edge. If there exist distinct $x,y\in U_{\rm even}\setminus  V(\cP)$, add the edge~$xy$ to~$E_{\rm fict}$ and~$\cP$. If there are no such~$x$ and~$y$ then note that, by \cref{sP'}, $|U_{\rm even}|\leq 2|\cP|+1\leq 3\varepsilon n$ and proceed as follows. If~$\cP$ is the only set in~$\sP'$ which does not contain a fictive edge, we remove~$\cP$ from~$\sP'$ and add its edges to~$H$ and we are done. Otherwise, pick~$\cP'\in \sP'\setminus\{\cP\}$ such that~$\cP'$ does not contain a fictive edge. Note that, by \cref{sP'}, $|(V\setminus V_0)\cap V(\cP\cup \cP')|\leq 4\varepsilon n$. Moreover, there are at most $\frac{2|\sP\cup \sP'|}{\varepsilon n}\leq 2\varepsilon^{-1}$ vertices in~$V\setminus V_0$ which are incident to at least~$\varepsilon n$ edges in~$E_{\rm fict}$. Thus, (since~$|V\setminus V_0|= mk$) we can choose distinct $x,y\in (V\setminus V_0)\setminus V(\cP\cup \cP')$ such that~$E_{\rm fict}$ contains fewer than~$\varepsilon n$ edges incident to~$x$ and fewer than~$\varepsilon n$ edges incident to~$y$. Add the edge~$xy$ to both~$\cP$ and~$\cP'$, and, add two edges between~$x$ and~$y$ to~$E_{\rm fict}$. We repeat this procedure until we have added a fictive edge to each linear forest in~$\sP'$.},
		now using \cref{sP'} instead of \cref{sP} and allowing fictive edges to have endpoints in~$V\setminus V_0$ instead of~$V^\ell$ for some~$\ell\in [r]$%
		\COMMENT{Recall that no reservoir is assigned to the linear forests in~$\sP'$.}.	
		Once each linear forest in~$\sP'$ contains a fictive edge, observe that we have added at most~$r+1$ linear forests to~$H$, so, by \cref{eq:Delta3'},
		\begin{equation}\label{eq:Delta6}
		\Delta(H)\leq 11\zeta n.
		\end{equation}
		Moreover, the following holds.
		
		\begin{claim}\label{claim:U}
		 The set~$U_{\rm even}$ has even size. Moreover, $|U_{\rm even}| \leq \max \{|U|-2(|\sP\cup \sP'|-(r+1)), \sqrt{\varepsilon}n\}$.
		\end{claim}
	
		\begin{proofclaim}
			By construction and since~$|U|$ is even,~$|U_{\rm even}|$ is even. For the second part of the claim, we distinguish three cases.
			
			Firstly, assume that for any distinct~$\cP, \cP'\in \sP\cup \sP'$, the fictive edge of~$\cP$ is vertex-disjoint from the fictive edge of~$\cP'$.
			Then we clearly have $|U_{\rm even}|\leq |U|-2(|\sP\cup \sP'|-(r+1))$.
					
			Secondly, assume there exists~$\ell\in [r]$ such that there exist distinct~$\cP, \cP'\in \sP_\ell$ and~$x,y\in V^\ell$ such that both~$\cP$ and~$\cP'$ contain a fictive edge between~$x$ and~$y$.
			Then, by construction, $|U_{\rm even}^\ell|\leq 3\varepsilon n$, and so $|U^\ell|\leq 2|\sP_\ell|+3\varepsilon n$. 
			By \cref{odd}, for any~$\ell'\in [r]$, we have $|U^{\ell'}|\leq  2|\sP_{\ell'}|+5\varepsilon n$ and, thus, $|U_{\rm even}^{\ell'}|\leq 5\varepsilon n$.\COMMENT{Either we have decreased~$|U_{\rm even}^{\ell'}|$ by~$2$ at each stage and $|U_{\rm even}^{\ell'}|\leq 5\varepsilon n$, or we have  $|U_{\rm even}^{\ell'}|\leq 3\varepsilon n$.}
			Therefore, $|U_{\rm even}|\leq \sqrt{\varepsilon} n$, as desired.
		
			Thirdly, assume that there exist distinct~$\cP, \cP'\in \sP'$ and~$x,y\in V\setminus V_0$ such that~$\cP$ and~$\cP'$ both contain a fictive edge between~$x$ and~$y$.
			Then by construction $|U_{\rm even}|\leq 3\varepsilon n$. 
		\end{proofclaim}
		
		Pair all vertices in~$U_{\rm even}$ and for each pair~$(x,y)$, add~$xy$ to~$E_{\rm fict}$ and~$\{xy\}$ to~$\sP'$. By construction,~\cref{odd},~\cref{sP,sP',endpoints,numbersets}, and \cref{claim:U}, the following hold.
		\begin{enumerate}[label=(\Alph*$'$)]
			\item For each~$\ell\in [r]$ and~$\cP\in \sP_\ell$,~$\cP$ is a set of vertex-disjoint paths with endpoints in~$V^\ell$ and internal vertices in~$V\setminus V_\ell$. Moreover, $|\cP|\leq \frac{k}{2}+\varepsilon n+1$ and $|V(\cP)\cap V_i^\ell|\leq 2\zeta^2 m$ for each~$i\in [k]$\COMMENT{$|V(\cP)\cap V_i^\ell|\leq \zeta^2 m+2\leq 2\zeta^2 m$}.\label{sPfictive}
			\item For each~$\cP\in \sP'$,~$\cP$ is a set of vertex-disjoint paths with endpoints in~$V\setminus V_0$ and internal vertices in~$V_0$. Moreover,~$|\cP|\leq \varepsilon n+1$ and~$|V(\cP)\cap V_i|\leq 2\zeta m$ for each~$i\in [k]$.
			\item For each~$x\in V\setminus V_0$, there are at most~$4\varepsilon n$ paths in~$\bigcup(\sP\cup \sP')$ which have~$x$ as an endpoint%
			\COMMENT{At most $2\varepsilon n+\varepsilon n+1\leq 4\varepsilon n$.}.\label{endpointsfictive}
			\item By \cref{numbersets}, \cref{claim:U}, and construction, $|\sP\cup \sP'|\leq \max\left\{\frac{\Delta}{2},\frac{|U|}{2}\right\}+8\zeta n$%
			\COMMENT{If in \cref{claim:U},~$|U_{\rm even}|\leq |U|-2(|\sP\cup\sP'|-(r+1))$, then (after having paired up the vertices in~$U_{\rm even}$ and having added those to~$\cP'$), we now have~$|\sP\cup\sP'|\leq \frac{|U|}{2}+(r+1)\leq \frac{|U|}{2} +8\zeta n$. Otherwise, we now have $|\sP\cup\sP'|\leq \frac{\Delta}{2}+7\zeta n+\frac{\sqrt{\varepsilon}n}{2}\leq \frac{\Delta}{2}+8\zeta n$.}.
			Moreover,~$G$ is now empty. \label{numbersetsfictive}
			\item Each set in~$\cP\in \sP\cup \sP'$ contains exactly one edge of~$E_{\rm fict}$. Moreover, for each~$x\in V\setminus V_0$,~$E_{\rm fict}$ contains an odd number of edges incident to~$x$ if and only if~$x\in U$.\label{setsfictive}
		\end{enumerate}
			
		\item \textbf{Tying each set of paths into a cycle.}\label{step:tyingpaths} We now tie each linear forest~$\cP\in \sP\cup\sP'$ into a cycle using edges of~$\Gamma\cup \Gamma'$. This is achieved by successively applying \cref{lm:tyingmanypaths,lm:tyingfewpaths,lm:closingcycle} several times as follows.

		For each~$\ell\in[r]$, we tie the paths of each linear forest in~$\sP_\ell$ as follows. Let~$\Gamma_\ell$ be the graph on vertex set~$V_0\cup V^\ell$ and edge set $\bigcup_{ij\in E(R')} E(\Gamma_{ij\ell})$. Note that by \cref{Gammaij},~$V_0, V_1^\ell, \dots, V_k^\ell$ is an $(\varepsilon_2, \beta, k, \zeta m, R_\ell')$-superregular partition of~$\Gamma_\ell$%
		\COMMENT{Note that $|V_0|\leq \varepsilon n\leq \frac{\varepsilon_2\zeta n}{2}\leq \frac{\varepsilon_2 \zeta mk}{(1-\varepsilon_2)}$ so $|V_0|\leq \varepsilon_2 (\zeta mk+|V_0|)=\varepsilon_2 |V(\Gamma^\ell)|$, as needed. Moreover, if $\zeta m\notin \mathbb{N}$ and so the $V_i^\ell$ only have size $\zeta m\pm 1$, then we can still choose the $V_i^\ell$ such that for all $\ell\in [r]$ the sizes of $V_1^\ell, \dots, V_k^\ell$ are equal.},
		where $V_i^\ell V_j^\ell \in E(R_\ell')$ if and only if~$V_i V_j \in E(R')$. 
		Define~$\Gamma_\ell'$ similarly. Moreover, by \cref{Gamma'ij},~$V_0, V_1^\ell, \dots, V_k^\ell$ is an $(\varepsilon_2, \zeta, k, \zeta m, R_\ell'')$-superregular partition of~$\Gamma_\ell'$, where $V_i^\ell V_j^\ell \in E(R_\ell'')$ if and only if~$V_i V_j \in E(R'')$.
		
		For each~$\ell\in [r]$, we can then successively apply \cref{lm:tyingmanypaths,lm:tyingfewpaths,lm:closingcycle} as follows. Write $\sP_\ell =\{\cP_1, \dots, \cP_{\ell'}\}$. 
		First, we apply \cref{lm:tyingmanypaths} with $\Gamma_\ell, \zeta m, V_1^\ell, \dots, V_k^\ell, \varepsilon_2, R_\ell',\ell',\cP_1, \dots, \cP_{\ell'}$, and~$2\zeta$ playing the roles of $\Gamma, m, V_1, \dots, V_k, \varepsilon, R, \ell,\cP_1, \dots, \cP_\ell$, and~$\zeta$, respectively%
		\COMMENT{Properties \cref{lm:tyingmanypathssupreg,lm:tyingmanypathsbadinfewpairs} of \cref{lm:tyingmanypaths} hold by construction and properties \cref{lm:exceptionaledgessupregpart,lm:exceptionaledgesGammabadinfewpairs} of \cref{lm:exceptionaledges}. Properties \cref{lm:tyingmanypathsvertexdisjoint,lm:tyingmanypathsroom} of \cref{lm:tyingmanypaths} follow from \cref{sP} and \cref{sPfictive}. Property \cref{lm:tyingmanypathsendpoint} of \cref{lm:tyingmanypaths} follows from \cref{endpoints} and \cref{endpointsfictive}.}.
		We thus obtain disjoint $E_1, \dots, E_{\ell'}\subseteq E(\Gamma_\ell)$ such that the following hold. 
		For any distinct~$i,j\in [k]$, and~$x\in V_i^\ell$, the set $E\coloneqq E_1\cup \dots\cup E_{\ell'}$ contains at most~$\varepsilon_3\zeta m$ edges of $\Gamma_\ell [V_i^\ell, V_j^\ell]$ which are incident to~$x$, and, thus, by \cref{lm:verticesedgesremoval}, it follows that~$V_0, V_1^\ell, \dots, V_k^\ell$ is an $(\varepsilon_4, \beta, k, \zeta m, R_\ell')$-superregular partition of~$\Gamma_\ell \setminus E$. Moreover, for any~$i\in [\ell']$ and~$j\in [k]$, $\lvert V(\cP_i\cup E_i)\cap V_j^\ell\rvert \leq \zeta_1 \zeta m$. Finally, for any~$i\in [\ell']$, by using each edge in~$E_i$ exactly once, we can tie some of the paths in~$\cP_i$ to form a set of vertex-disjoint paths~$\cQ_i$ such that, for any~$j\in [k]$, at most~$2\beta^{-2}$ paths in~$\cQ_i$ have an endpoint in~$V_j^\ell$.
		
		We then apply \cref{lm:tyingfewpaths} with $\Gamma_\ell\setminus E, \zeta m, V_1^\ell, \dots, V_k^\ell, \varepsilon_4, R_\ell', \ell',\cQ_1, \dots, \cQ_{\ell'}$, and~$\zeta_1$ playing the roles of $\Gamma, m, V_1, \dots, V_k, \varepsilon, R, \ell,\cP_1, \dots, \cP_\ell$, and~$\zeta$, respectively. We thus obtain disjoint $E_1', \dots, E_{\ell'}'\subseteq E(\Gamma_\ell)\setminus E$ satisfying the following.  
		For any distinct~$i,j\in [k]$, and~$x\in V_i^\ell$, the set $E'\coloneqq E_1'\cup \dots\cup E_{\ell'}'$ contains at most~$\varepsilon_5\zeta m$ edges of $\Gamma_\ell [V_i^\ell, V_j^\ell]$ which are incident to~$x$, and, thus, by \cref{lm:verticesedgesremoval}, it follows that $V_0, V_1^\ell, \dots, V_k^\ell$ is an $(\varepsilon_6, \beta, k, \zeta m, R_\ell')$-superregular partition of $\Gamma_\ell \setminus (E\cup E')$. Moreover, for any~$i\in [\ell']$ and~$j\in [k]$, $| V(\cQ_i\cup E_i')\cap V_j^\ell| \leq \zeta_2 \zeta m$. Finally, for any~$i\in [\ell']$, by using each edge in~$E_i'$ exactly once, we can tie the paths in~$\cQ_i$ to form a set of vertex-disjoint paths~$\cQ_i'$ such that, for any component~$C$ of~$R_\ell'$,~$\cQ_i'$ contains at most one path with an endpoint in~$V_{\Gamma_\ell}(C)$.
		
		We then apply \cref{lm:closingcycle} with $\Gamma_\ell\setminus (E\cup E'), \zeta m, V_1^\ell, \dots, V_k^\ell, \varepsilon_6, R_\ell', \ell',\cQ_1', \dots, \cQ_{\ell'}', \zeta_2, \Gamma_\ell'$, and~$R_\ell''$ playing the roles of $\Gamma, m, V_1, \dots, V_k, \varepsilon, R, \ell,\cP_1, \dots, \cP_\ell, \zeta, \Gamma'$, and~$R'$, respectively. We thus obtain disjoint $E_1'', \dots, E_{\ell'}''\subseteq E(\Gamma_\ell\cup \Gamma_\ell')\setminus (E\cup E')$ satisfying the following.  
		For any distinct~$i,j\in [k]$, and~$x\in V_i^\ell \cup V_j^\ell$, the set $E''\coloneqq E_1''\cup \dots\cup E_{\ell'}''$ contains at most~$\varepsilon_7\zeta m$ edges of $\Gamma_\ell [V_i^\ell, V_j^\ell]$ which are incident to~$x$ and at most~$\varepsilon_7\zeta m$ edges of $\Gamma_\ell' [V_i^\ell, V_j^\ell]$ which are incident to~$x$. Moreover,~$\cQ_i'\cup E_i''$ forms a cycle, i.e.\ $\sP_\ell \cup (E\cup E'\cup E'')$ admits a cycle decomposition~$\cD_\ell$ of size~$|\sP_\ell|$. Add all edges in~$E\cup E'\cup E''$ to~$G'$. Proceed in this way for each~$\ell\in [r]$. Observe that by construction and since the reservoirs are pairwise disjoint, \cref{lm:verticesedgesremoval} implies that~$V_0, V_1, \dots, V_k$ is now an $(\varepsilon_8, \beta, k, m, m', R')$-superregular equalised partition of~$\Gamma$ and an $(\varepsilon_8, \zeta, k, m, m', R'')$-superregular equalised partition of~$\Gamma'$. 		
		\COMMENT{For any~$i,j\in[k]$, and~$x\in V_i^\ell$,~$E\cup E'\cup E''$ contains at most~$(\varepsilon_3+\varepsilon_5+\varepsilon_7)\zeta m\leq \varepsilon_7 m$ edges of $\Gamma_\ell[V_i^\ell, V_j^\ell]$ incident to~$x$ and at most~$\varepsilon_7 \zeta m\leq \varepsilon_7 m$ edges of $\Gamma_\ell'[V_i^\ell, V_j^\ell]$ incident to~$x$.}

		We proceed similarly to obtain \NEW{a set} $E^* \subseteq E(\Gamma\cup \Gamma')$ such that~$\sP'\cup E^*$ admits a cycle decomposition~$\cD'$ of size~$|\sP'|$.
		Add the edges in~$E^*$ to~$G'$. Then, \cref{lm:verticesedgesremoval} implies that~$V_0, V_1, \dots, V_k$ is an $(\varepsilon_9,\beta, k, m, m', R')$-superregular equalised partition of~$\Gamma$. Add all remaining edges of~$\Gamma'$ to~$H$. By \cref{lm:exceptionaledgessupregpart}, we add at most~$(\zeta +\varepsilon)n$ edges incident to each vertex. Thus, by \cref{eq:Delta6},~$\Delta(H)\leq 13\zeta n$, as desired for~\cref{lm:exceptionaledgesDeltaH}.
		
		Let $\cD\coloneqq\bigcup_{\ell\in [r]} \cD_\ell \cup \cD'$.
		Observe that for the proof of \cref{lm:exceptionaledges,lm:connectedR}\cref{lm:connectedRcycles}, by \cref{numbersets} and construction,~$\cD$ is a cycle decomposition of~$G'$ of size at most~$|\sP\cup \sP'|\leq \frac{\Delta}{2}+7\zeta n$.
		For the proof of \cref{lm:fictive,lm:connectedR}\cref{lm:connectedRpaths}, remove all fictive edges from~$\cD$. Then, by \cref{numbersetsfictive,setsfictive},~$\cD$ is now a path decomposition of~$G'$ of size at most
		$\max\left\{\frac{\Delta}{2},\frac{|U|}{2}\right\}+8\zeta n$.
		Moreover, each vertex~$x\in V\setminus V_0$ is an endpoint of an odd number of paths in~$\cD$ if and only if~$x\in U$.
		This completes the proof of \cref{lm:connectedR} (where we set $\tGamma\coloneqq \Gamma$).
		
		\item \textbf{Covering the remaining exceptional edges.}\label{step:Rdisconnected} 
		If~$R$ is disconnected, we apply the above argument to each component of~$R$.
		More precisely, for each connected component~$C$ of~$R$, we apply \cref{step:equalisingsupport,step:exceptionaledges,step:fictive,step:partitioning,step:partitioning,step:paths,step:setsofpaths,step:tyingpaths} with~$R[C]$ and~$G[V_0\cup V_G(C)]$ playing the roles of~$R$ and~$G$, and, for the proof of \cref{lm:fictive},~$U\cap V_G(C)$ playing the role of~$U$.
		In particular, observe that for each component~$C$ of~$R$, $\Delta(G[V_0\cup V_G(C)])\leq |V_G(C)|+\varepsilon n$. Moreover, by \cref{lm:exceptionaledgesGammabadinfewpairs,lm:exceptionaledgesreduced},~$R$ has at most~$\beta^{-1}$ components. Therefore, we obtain, in the proof of \cref{lm:exceptionaledges} (\cref{lm:fictive}), a cycle (path) decomposition~$\cD$ of~$G$ of size at most $\frac{n}{2}+ \frac{8\zeta n}{\beta}+\frac{\varepsilon n}{\beta}\leq \frac{n}{2}+\beta n$.		
		Then, there only remains to decompose~$G_{\rm exc}$ into at most~$\beta n$ cycles and~$\beta^{-2}$ exceptional edges for \cref{lm:exceptionaledges}, or,~$3\beta n$ paths for \cref{lm:fictive}.
		
		Recall that \NEW{$G_{\rm exc}$ was introduced at the beginning of the proof. By construction,} all edges of~$G_{\rm exc}$ are exceptional and for any~$x\in V_0$, if~$xy,xy'\in E(G_{\rm exc})$ are distinct then there exist distinct components~$C$ and~$C'$ of~$R$ such that~$y\in V_G(C)$ and~$y'\in V_G(C')$. Decompose~$G_{\rm exc}$ into~$s \leq \zeta n$ paths of length~$2$ with endpoints in~$V\setminus V_0$ and an internal vertex in~$V_0$.\COMMENT{$s \leq \frac{\varepsilon n }{2\beta}\leq \zeta n$ since~$|V_0|\leq \varepsilon n$ and~$R$ has at most~$\beta^{-1}$ connected components.} 
		Note that, by construction, each path has endpoints in clusters which lie in different connected components of~$R$. Apply \cref{lm:closingcycles}, with~$\varepsilon_9$ and~$s$ playing the roles of~$\varepsilon$ and~$\ell$, and, each~$\cP_i$ consisting of exactly one of the paths constructed above. We thus obtain~$E^\diamond\subseteq E(\Gamma)$ such that~$G_{\rm exc}\cup E^\diamond$ admits a decomposition~$\cD''\cup \cD_{\rm exc}$ where~$\cD''$ is a set of at most~$\beta n$ cycles and~$\cD_{\rm exc}$ is a set of at most~$\beta^{-2}$ exceptional edges. 
		Add all edges in~$E^\diamond$ and~$G_{\rm exc}$ to~$G'$. By part~\cref{lm:closingcyclesneighbours} of \cref{lm:closingcycles},~$V_0, V_1, \dots, V_k$ is a~$(\zeta, \beta, k, m, m', R')$-superregular equalised partition of~$\Gamma$.
		
		Set~$\tGamma\coloneqq\Gamma$.
		For the proof of \cref{lm:exceptionaledges}, add all cycles in~$\cD''$ to~$\cD$.
		By construction,~$\cD\cup \cD_{\rm exc}$ satisfies the desired properties.
		For the proof of \cref{lm:fictive}, split each cycle in~$\cD''$ into two paths 
		and add them to~$\cD$. Add the edges in~$\cD_{\rm exc}$ to~$\cD$. This completes the proof of \cref{lm:exceptionaledges,lm:fictive}.\qed
	\end{steps}
	\renewcommand{\qed}{}
\end{proof}

\onlyinsubfile{\bibliographystyle{plain}
	\bibliography{Bibliography/papers}}

%% file: Bad_Graph.tex
\onlyinsubfile{
	\setcounter{section}{5}
	\setcounter{subsection}{4}
	\addtocounter{section}{-1}
	\setcounter{definition}{6}
}

\subsection{Covering the leftovers}\label{sec:bad} This section corresponds to \cref{step:sketchbad} of the proof overview.
We will need the following fact.

\begin{fact}\label{fact:euleriancycle}
	Assume~$G$ is Eulerian and~$e\in E(G)$. Then~$G$ contains a cycle~$C$ such that~$e\in C$. 
\end{fact}

\begin{lm}\label{lm:bad}
	Suppose $0<\frac{1}{n}\ll\frac{1}{k}\ll d\ll \beta \leq 1$ and~$d'\coloneqq d^{\frac{1}{10^5}}$. Let~$G$ and~$\Gamma$ be edge-disjoint graphs on vertex set~$V$ of size~$n$ such that~$G\cup\Gamma$ is Eulerian and~$\Delta(G)\leq dn$. Assume~$V_0, V_1, \dots, V_k$ is an $(d, \beta, k, m, m', R)$-superregular equalised partition of~$\Gamma$ such that any~$x\in V\setminus V_0$ belongs to at least~$\beta k$ superregular pairs of~$\Gamma$. Moreover, suppose that~$V_0$ is a set of isolated vertices in~$G$.
	Then there exists~$E\subseteq E(\Gamma)$ such that~$G\cup E$ can be decomposed into at most~$2\beta n$ cycles. Moreover,~$V_0, V_1, \dots, V_k$ is a~$(d', \beta, k, m, m', R)$-superregular equalised partition of~$\Gamma\setminus E$.
\end{lm}

\begin{proof}
	Fix an additional constant~$\zeta$ such that~$d\ll\zeta \ll \beta$.
	The idea is to decompose~$G$ into matchings and then apply \cref{lm:tyingmanypaths,lm:tyingfewpaths,lm:closingcycles} to tie the edges in each matching together to form cycles using~$\Gamma$.
	
	By Vizing's theorem (\cref{thm:Vizing}), we can decompose~$G$ into~$\ell\leq d n+1\leq 2d n$ matchings~$M_1$, $\dots, M_\ell$. Randomly split each matching~$M_i$ into~$2\zeta^{-1}$ submatchings $M_{i,1}, \dots, M_{i,2\zeta^{-1}}$, by including each edge of~$M_i$ to~$M_{i,j}$ independently with probability~$\frac{\zeta}{2}$ for each~$j\in \left[2\zeta^{-1}\right]$. By \cref{lm:Chernoff}, we may assume that for each~$i\in[\ell]$ and~$i'\in \left[2\zeta^{-1}\right]$, we have~$|M_{i,i'}|\leq \zeta n$ and for each~$j\in [k]$, we have $|V(M_{i,i'})\cap V_j|\leq \zeta m$.%
	For simplicity, set $\ell'\coloneqq \frac{2\ell}{\zeta}\leq \zeta n$ and relabel $M_{1,1}, \dots, M_{1,2\zeta^{-1}}, \dots,M_{\ell, 1}$, $\dots, M_{\ell,2\zeta^{-1}}$ to~$M_1, \dots, M_{\ell'}$. We successively apply \cref{lm:tyingmanypaths,lm:tyingfewpaths,lm:closingcycles}, starting with~$d, M_1, \dots, M_{\ell'}$ playing the roles of $\varepsilon, \cP_1, \dots, \cP_\ell$ in \cref{lm:tyingmanypaths}.
	We thus obtain~$E_1\subseteq E(\Gamma)$ such that~$G\cup E_1$ admits a decomposition~$\cD\cup \cD'$ where~$\cD$ is a set of at most~$\beta n$ cycles and~$\cD'$ is a set of at most~$\beta^{-2}$ edges. Moreover, by \cref{lm:verticesedgesremoval} and part \cref{lm:tyingmanypathsneighbours,lm:tyingfewpathsneighbours,lm:closingcyclesneighbours} of \cref{lm:tyingmanypaths,lm:tyingfewpaths,lm:closingcycles}, $V_0, V_1, \dots, V_k$ is a~$(d^{\frac{1}{2\cdot 10^4}}, \beta, k, m, m',R)$-superregular equalised partition of~$\Gamma\setminus E_1$. 
	By \cref{fact:euleriancycle,lm:verticesedgesremoval}, there exists $E_2\subseteq \Gamma \setminus E_1$ such that~$E(\cD')\cup E_2$ can be decomposed into at most~$\beta^{-2}$ cycles and~$V_0, V_1, \dots, V_k$ is a~$(d', \beta, k, m, m', R)$-superregular equalised partition of~$\Gamma\setminus (E_1\cup E_2)$. Let~$E\coloneqq E_1\cup E_2$. This completes the proof.
\end{proof}

\onlyinsubfile{\bibliographystyle{plain}
	\bibliography{Bibliography/papers}}

%% file: Decomposing_Gamma.tex
\onlyinsubfile{
\setcounter{section}{5}
\setcounter{subsection}{5}
\setcounter{definition}{8}
}

\subsection{Fully decomposing \texorpdfstring{$\Gamma$}{Gamma}}\label{sec:Gamma}

This section corresponds to \cref{step:sketchGamma} of the proof overview. 

\begin{lm}\label{lm:cycleGamma}
	Let $0< \frac{1}{m}\ll \frac{1}{k}\ll \frac{1}{K}\ll\varepsilon \ll \frac{1}{q}\ll \frac{1}{f}\ll d\ll \frac{1}{\ell}, \frac{1}{g}\ll 1$ and suppose that $\frac{K}{7}, \frac{2K}{f}, \frac{2K}{g}, \frac{q}{f}, \frac{m}{4\ell K}, \frac{fm}{qK}, \frac{4fK}{3g(g-1)},\frac{\ell}{2}\in \mathbb{N}^*$. 
	Let~$G$ be an Eulerian graph on~$n$ vertices. Assume \NEW{that} $V_0,V_1, \dots, V_k$ is a partition of~$V(G)$ into an exceptional set~$V_0$ consisting of at most~$\varepsilon n$ isolated vertices and~$k$ clusters $V_1, \dots, V_k$ of size~$m$ such that the corresponding reduced graph~$R$ of~$G$ is a cycle of even length, and for each~$ij\in E(R)$, the pair~$G[V_i,V_j]$ is~$[\varepsilon, d]$-superregular.
	Then~$G$ admits a cycle decomposition~$\cD$ of size at most~$dm+\varepsilon^{\frac{1}{16}}m$.
\end{lm}

\NEW{To prove \cref{lm:cycleGamma}, we will use the robust decomposition lemma (\cref{lm:robustdecomp}). In order to apply this result, $m$ needs to satisfy certain divisibility conditions and we need to find several refinements of the partition $V_0, V_1, \dots, V_k$. This would not be possible if, for example, $m$ was prime. This explains why it is necessary to introduce the parameters $K, q, f, \ell$, and $g$ in the statement of \cref{lm:cycleGamma}.}

\begin{cor}\label{cor:Gammadecomp}
	Let $0< \frac{1}{m}\ll \frac{1}{k}\leq \frac{1}{L}\ll \frac{1}{K}\ll\varepsilon \ll \frac{1}{q}\ll \frac{1}{f}\ll d\ll \frac{1}{\ell}, \frac{1}{g}\ll 1$ and suppose that $\frac{K}{7}, \frac{2K}{f}, \frac{2K}{g}, \frac{q}{f}, \frac{m'}{4\ell K}, \frac{fm'}{qK}, \frac{4fK}{3g(g-1)},\frac{\ell}{2}\in \mathbb{N}^*$.
	Let~$G$ be an~$n$-vertex Eulerian graph and assume~$V_0, V_1, \dots, V_k$ is an~$(\varepsilon, d, k,m,m',R)$-superregular equalised partition of~$G$ such that~$V_0$ is a set of isolated vertices in~$G$. Assume~$R$ admits a decomposition~$\cD_R$ satisfying the following properties.~$\cD_R$ consists of at most~$\frac{k}{2}$ cycles whose lengths are even and at least~$L$. Moreover, for any distinct~$i,j,j'\in [k]$, if~$V_j V_i V_j'$ is a subpath of a cycle in~$\cD_R$, then the support clusters of~$V_i$ with respect to~$V_j$ and~$V_{j'}$ are the same. Then~$G$ admits a cycle decomposition~$\cD$ of size at most $\frac{dn}{2}+ \varepsilon^{\frac{1}{33}} n$.
\end{cor}

\begin{proof}
	Let $\cD\coloneqq \emptyset$. First apply \cref{lm:Eulerianpairs} and add the cycles obtained to~$\cD$ and delete their edges from~$G$. Then, for each cycle $C=V_{i_1}\dots V_{i_{k'}}$ in~$\cD_R$, apply \cref{lm:cycleGamma} with $2\sqrt{\varepsilon}, V_{i_1},\dots, V_{i_{k'}}, i_{k'}, m', |V_0|+i_{k'}m'$ and $G[V_0\cup V_G(C)]$ playing the roles of $\varepsilon, V_1, \dots, V_k, k, m, n$ and~$G$, respectively\COMMENT{Note that this is possible since~$i_k\geq L$ and is even, and, by construction, the superregular pairs of~$G'$ are all Eulerian.}, and add the cycles obtained to~$\cD$.
\end{proof}

To prove \cref{lm:cycleGamma}, we will find an approximate decomposition of~$G$ using \cref{lm:approxdecomp} and cover the leftover using the robust decomposition lemma (introduced in \cref{sec:robust}). 
The approximate decomposition will be obtained be repeatedly applying the following lemma, which is a special case of \cite[Lemma 6.4]{kuhn2013hamilton}.

	\begin{lm}\label{lm:Hamcycle}
		Let $0<\frac{1}{m}\ll d'\ll \varepsilon \ll d \ll \zeta, \frac{1}{t}\leq \frac{1}{2}$ and $k\geq 3$. Let~$G$ be a graph and~$V_1, \dots , V_k~$ be a partition of~$V(G)$ into~$k$ clusters of size~$m$. Suppose that the following hold.
		\begin{itemize}
			\item For each~$i\in[k -1]$,~$G[V_i, V_{i+1}]$ is a perfect matching~$M_i$.
			\item $G[V_1, V_k]$ is $(\varepsilon, d', \zeta d', \frac{t d'}{d})$-superregular.
		\end{itemize} 
		Then,~$G[V_1, V_k ]$ contains a perfect matching~$M$ such that $M\cup \left(\bigcup_{i=1}^{k -1}M_i\right)$ is a Hamilton cycle of~$G$.
	\end{lm}

\begin{lm}\label{lm:approxdecomp}
	Suppose $0< \frac{1}{m}\ll \frac{1}{k}\ll d' \ll \varepsilon \ll d \leq 1$. Let
	\[r\coloneqq \frac{18d' m}{d} \quad {\rm and} \quad h\coloneqq dm-r,\]
	and assume that $r, \frac{h}{k},dm\in \mathbb{N}^*$.
	Let~$G$ be an~$n$-vertex graph with vertex set~$V$. Assume \NEW{that} $V_0,V_1, \dots, V_k$ is a partition of~$V$ into an exceptional set~$V_0$ consisting of at most~$\varepsilon n$ isolated vertices and~$k$ clusters $V_1, \dots, V_k$ of size~$m$ such that the corresponding reduced graph~$R$ of~$G$ is a cycle, and, for each~$ij\in E(R)$, the pair~$G[V_i,V_j]$ is~$\varepsilon$-regular and~$dm$-regular.
	Then there exists~$H\subseteq G$ such that, for each $ij\in E(R)$,~$H[V_i, V_j]$ is~$r$-regular and $G'\coloneqq G\setminus H$ admits a decomposition~$\cD$ into~$h$ Hamilton cycles of~$G'-V_0$.
\end{lm}

\begin{proof}
	Let~$H$ be the empty graph on~$V$. Let \[
	\varepsilon_1\coloneqq \varepsilon^{\frac{1}{12}}, \quad
	\varepsilon_2\coloneqq \varepsilon^{\frac{1}{25}}.\]
	We may assume without loss of generality that $E(R)\coloneqq \{V_i V_{i+1} \mid i\in [k]\}$, where $V_{k+1}\coloneqq V_1$. For each~$i\in [k]$, denote $G_i\coloneqq G[V_i, V_{i+1}]$.
	
	Let~$i\in [k]$. Apply \cref{lm:sparsesubgraph} to obtain an $(\varepsilon_1, d', \frac{d'}{2}, \frac{3d'}{2d})$-superregular spanning subgraph~$\Gamma_i\subseteq G_i$. Let $G'_i\coloneqq G_i\setminus \Gamma_i$. One can easily verify that~$G'_i$ is~$\varepsilon_1$-regular and that, for each~$x\in V_i\cup V_{i+1}$, we have $d_{G_i'}(x)=(d\pm \frac{3d'}{2d})m$.%
	\COMMENT{
	Let~$S\subseteq V_i$ and~$T\subseteq V_{i+1}$ be such that~$|S|, |T|\geq \varepsilon_1 m$. Then we have
	\begin{equation*}
	\begin{aligned}
	d_{G'_i}(S, T)-d_{G'_i}(V_i,V_{i+1})&=(d_{G_i}(S, T)-d_{G_i}(V_i, V_{i+1}))-(d_{\Gamma_i}(S,T)-d_{\Gamma_i}(V_i, V_{i+1}))\\
	&\geq -\varepsilon - (1+\varepsilon_1)d'\\
	& \geq- \varepsilon_1;
	\end{aligned}	
	\end{equation*}
	and similarly
	\begin{equation*}
	d_{G'_i}(S,T)-d_{G'_i}(V_i, V_{i+1})\leq \varepsilon - (1-\varepsilon_1)d'\leq \varepsilon_1,	
	\end{equation*}
	so~$G'_i$ is~$\varepsilon_1$-regular.}
	
	In order to apply \cref{lm:Hamcycle}, we need to decompose each~$G_i'$ into perfect matchings. Thus, we will first ensure that the pairs~$G_i'$ are Eulerian and then apply \cref{lm:regularising} to regularise them.
	
	Let~$i\in [k]$. Apply \cref{lm:regHamcycle} to obtain a Hamilton cycle~$x_1 \dots x_{2m}$ of~$G'_i$. Let~$i_1<\dots < i_\ell$ be the indices of the odd-degree vertices of~$G_i'$. For each~$s\in \{1, 3, \dots, \ell-1\}$, add the edges of the path $x_{i_s}x_{i_{s}+1}\dots x_{i_{s+1}}$ to~$H$ and delete them from~$G_i'$. By construction and \cref{lm:verticesedgesremoval},~$G'_i$ is now Eulerian and~$\varepsilon_2$-regular. Moreover, we have $d_{G_i'}(x)=(d\pm \frac{2d'}{d})m$ for each~$x\in V_i\cup V_{i+1}$.  
	Apply \cref{lm:regularising} to obtain a regular spanning subgraph~$G_i''$ of~$G_i'$. By removing perfect matchings if necessary, we may assume~$G''_i$ is~$(\frac{k-1}{k}h)$-regular.%
	\COMMENT{By \cref{lm:regularising}, the degree of~$G''_i$ is at least $\Delta(G_i')-4(\Delta(G_i')-\delta(G_i'))\geq (d-\frac{2d'}{d})m-\frac{16d'}{d}m=(d-\frac{18d'}{d})m=h\geq \frac{k-1}{k}h$.}
	Apply Hall's theorem to obtain a decomposition of~$G_i''$ into edge-disjoint sets~$\cD_i^s$, with $s\in [k]\setminus \{i\}$, each containing~$\frac{h}{k}$ edge-disjoint perfect matchings.
	Add all edges of~$G_i'\setminus G''_i$ to~$H$.
	
	Let~$\ell\in [k]$. Let~$G^\ell$ be the graph on vertex set $V\setminus V_0$ with $E(G^\ell)\coloneqq (\bigcup_{i\in [k]\setminus \{\ell\}}E(\cD_i^\ell)) \cup \Gamma_\ell$. We construct~$\frac{h}{k}$ edge-disjoint Hamilton cycles~$C_1, \dots, C_\frac{h}{k}$ of~$G^\ell$ such that, for each~$s\in [\frac{h}{k}]$ and $i\in [k]\setminus \{\ell\}$,~$C_s$ contains a perfect matching in~$\cD_i^\ell$. In particular, observe that this implies that $C_s[V_{\ell} , V_{\ell+1}]$ is a perfect matching of~$G^\ell[V_\ell, V_{\ell+1}]=\Gamma_\ell$.
	
	Assume inductively that we have already constructed~$C_1, \dots, C_s$ for some~$0\leq s<\frac{h}{k}$. Delete from~$G^\ell$ all edges of~$C_1, \dots, C_s$. Note that since~$\frac{h}{km}\ll \varepsilon_1, d'$, by \cref{lm:edgesremovalsparse}, the pair $G^\ell [V_\ell, V_{\ell+1}]$ is still $(2\varepsilon_1, d', \frac{d'}{4}, \frac{3d'}{2d})$-superregular. Let~$F$ be a spanning subgraph of~$G^\ell$ such that, for each $i\in [k]\setminus \{\ell\}$, $F[V_i, V_{i+1}]$ is a \NEW{perfect} matching in~$\cD_i^\ell$ which has not been used for $C_1, \dots, C_s$ and $F[V_\ell, V_{\ell+1}]=G^\ell[V_\ell, V_{\ell+1}]$. Then, \cref{lm:Hamcycle} gives a Hamilton cycle~$C_{s+1}$ of~$F\subseteq G^\ell$ satisfying the desired properties.
	
	Proceed as above for each~$\ell\in [k]$ and add all cycles obtained to~$\cD$. Add to~$H$ all remaining edges of $\bigcup_{i\in [k]} \Gamma_i$. This completes the proof.
\end{proof}

We are now ready to prove \cref{lm:cycleGamma}, using the robust decomposition lemma. Recall the terminology introduced in \cref{sec:robust}.

\begin{proof}[Proof of \cref{lm:cycleGamma}]
	Let~$\cD\coloneqq \emptyset$. We will repeatedly add cycles to~$\cD$. Whenever a cycle is added to~$\cD$, it is removed from~$G$ so that all cycles in~$\cD$ are always pairwise edge-disjoint, as desired.
	In \cref{step:leftover,step:approx}, we will construct at most~$dm$ edge-disjoint Hamilton cycles of~$G$. The additional cycles will be created during the regularising step (see \cref{step:regularising}).
	
	Note that~$k$ is even. We may assume without loss of generality that $E(R)=\{V_iV_{i+1} \mid i\in [k]\}$, where~$V_{k+1}\coloneqq V_1$.
	For each~$i\in [k]$, denote~$G_i\coloneqq G[V_i, V_{i+1}]$.
	Decompose~$R$ into two perfect matchings $M\coloneqq \{V_i V_{i+1} \mid i\in[k]_{\rm odd}\}$ and $M'\coloneqq \{V_i V_{i+1} \mid i\in[k]_{\rm even}\}$. 
	As explained in \cref{step:sketchGamma} of the proof overview, the idea is decompose each superregular pair of~$G$ into Hamilton paths using the robust decomposition lemma and suitable fictive edges. We will then form Hamilton cycles of~$G-V_0$ by tying together a Hamilton path of each pair in~$M$ using an edge from each pair in~$M'$, and similarly for~$M$ and~$M'$ exchanged.
	
	\begin{steps}
		\item \textbf{Choosing the constants.}
		Fix additional constants such that \[0<\frac{1}{m}\ll \frac{1}{k}\ll d''\ll \frac{1}{K}\ll \varepsilon\ll \frac{1}{q}\ll \frac{1}{f}\ll \frac{r_1 K}{m}\ll d \ll \frac{1}{\ell}, \frac{1}{g}\ll 1.\]
		Let
		\[\varepsilon_1\coloneqq \varepsilon^{\frac{1}{12}}, \quad
		\varepsilon_2\coloneqq \varepsilon^{\frac{1}{84}}, \quad
		\varepsilon_3\coloneqq \varepsilon^{\frac{1}{588}},\quad
		\varepsilon_4\coloneqq \varepsilon^{\frac{1}{1177}},
		\]
		and, 
		\[\varepsilon_1^*\coloneqq \varepsilon^{\frac{1}{3}},\quad
		\varepsilon_2^*\coloneqq \varepsilon^{\frac{1}{7}}, \quad
		\varepsilon_3^*\coloneqq \varepsilon^{\frac{1}{15}},\quad
		\varepsilon_4^*\coloneqq \varepsilon^{\frac{1}{31}}, \quad
		\varepsilon_5^* \coloneqq \varepsilon^{\frac{1}{125}}, \quad
		\varepsilon_6^* \coloneqq \varepsilon^{\frac{1}{875}}, \quad
		\varepsilon_7^* \coloneqq \varepsilon^{\frac{1}{1751}}.\]
		Let~$c\coloneqq c(d,k)$ be the constant in \cref{lm:Eulerianpairs} and define
		\[d'\coloneqq d-11\varepsilon_2^*, \quad r\coloneqq \frac{9d''m}{d'}.\]
		Observe that~$r(2K)^2\leq \frac{m}{K}$.
		Let \[r_2\coloneqq 192\ell g^2Kr, \quad r_3\coloneqq \frac{2rfK}{q}, \quad r^\diamond\coloneqq r_1+r_2+r-(q-1)r_3, \quad s\coloneqq 2rfK+7r^\diamond.
		\]
		Note that~$r,r_2, r_3\leq r_1$ and~$r^\diamond \leq 2r_1$.\COMMENT{$r\leq \frac{m}{4K^3}\leq \frac{m}{K^2}\ll r_1$; $\frac{1}{K}\ll \frac{r_1 K}{m}\ll \frac{1}{\ell}, \frac{1}{g}$ implies $\frac{192 \ell g^2}{4K}\ll\frac{r_1 K}{m}$ so $r_2\leq 192\ell g^2 \frac{m}{4K^2}\ll r_1$;~$r_3\leq \frac{2fm}{4K^2}\ll r_1$.}
		By adjusting $\varepsilon, d$, and $d''$ slightly, we may assume that $\frac{(d-10\varepsilon_2^*)m}{2}, \frac{d'm}{2}, \frac{d'm-2r}{k}, r, r_3\in \mathbb{N}^*$.
		
		\item \textbf{Constructing the bi-setups.}\label{step:bisetup}
		Let~$i\in [k]$. 		
		Apply \cref{cor:supregdigraph} to obtain an orientation~$\overrightarrow{G}_i$ of~$G_i$ such that both $\overrightarrow{G}_i[V_i, V_{i+1}]$ and $\overrightarrow{G}_i[V_{i+1}, V_i]$ are~$[\varepsilon_1,\frac{d}{2}]$-superregular (here and below, the index is taken modulo~$k$).
		
		For each~$i\in [k]$, randomly partition~$V_i$ into~$K$ subclusters~$V_{i,1}, \dots, V_{i,K}$ of equal size. 
		This induces, for each~$i\in [k]$, a partition~$\cP_i$ of~$V(\overrightarrow{G}_i)$ into~$2K$ clusters of size~$\frac{m}{K}$. By \cref{lm:vertexpartition}, we may assume that for each~$i\in [k]$, the partition~$\cP_i$ is an~$\varepsilon_2$-superregular~$K$-refinement of~$\{V_i, V_{i+1}\}$. Let~$\overrightarrow{R}_i$ be the reduced digraph of~$\overrightarrow{G}_i$ with respect to~$\cP_i$. (Thus,~$\overrightarrow{R}_i$ is the complete bipartite digraph with vertex classes of size~$K$.)
		Proceed similarly to obtain, for each~$i\in [k]$, an~$\varepsilon_3$-superregular~$\ell$-refinement~$\cP_i'$ of~$\cP_i$.%
		\COMMENT{For each~$i\in [k]$ and~$j\in [K]$, randomly partition~$V_{i,j}$ into~$\ell$ subclusters of equal size. This induces, for each~$i\in [k]$, a partition~$\cP_i'$ of~$\overrightarrow{G}_i$ into~$K\ell$ clusters of size~$\frac{m}{K\ell}$. By \cref{lm:vertexpartition}, we may assume that~$\cP_i'$ is an~$\varepsilon_3$-superregular~$\ell$-refinement of~$\cP_i$.}
		Let~$\overrightarrow{R}_i'$ be the reduced digraph of~$\overrightarrow{G}_i$ with respect to~$\cP_i'$.
		
		For each~$i\in [k]$, let $C_i\coloneqq V_{i,1} V_{i+1,1}V_{i,2} \dots V_{i+1,K}$ and observe that~$C_i$ is a Hamilton cycle of~$\overrightarrow{R}_i$.			
		Thus, by construction, for each~$i\in [k]$, $(\overrightarrow{G}_i, \cP_i, \cP_i', \overrightarrow{R}_i, \overrightarrow{R}_i', C_i)$ is an $(\ell, 2K, \frac{m}{K}, \varepsilon_3, \frac{d}{2})$-bi-setup.
		
		\item \textbf{Selecting the fictive edges.}\label{step:fict}
		\NEW{In this step, we will set aside a set $\cE$ of edges which will enable us to tie together the Hamilton path obtained with the robust decomposition lemma. Then, we will construct a corresponding set $\cF$ of fictive edges which will prescribe the endpoints of the Hamilton paths (recall \cref{fig:sketch}).
		These fictive edges will then be incorporated in the special path systems. Thus, in order to satisfy \cref{def:SPSfictive}, we will need to ensure that the endpoints of the edges in $\cE$ lie in the appropriate subclusters (see \cref{fig:SF-SPS}).}
		\OLD{In order to apply \cref{lm:robustdecomp}, we need to select appropriate fictive edges.}% 
		We start by choosing the fictive edges which will be included in the special factors required for finding the chord absorbers (see part \cref{lm:robustdecompCA} of \cref{lm:robustdecomp}).
		
		First, proceed as in \cref{step:bisetup} to construct, for each~$i\in [k]$, an~$\varepsilon_3$-superregular~$\frac{q}{f}$-refinement~$\cP_i^*$ of~$\cP_i$.%
		\COMMENT{For each~$i\in [k]$ and~$j\in [K]$, randomly partition~$V_{i,j}$ into~$\frac{q}{f}$ parts $V_{i,j,1}, \dots, V_{i,j,\frac{q}{f}}$. This induces, for each~$i\in [k]$, a partition~$\cP_i^*$ of~$\overrightarrow{G}_i$ into~$\frac{qK}{f}$ clusters of size~$\frac{fm}{qK}$. By \cref{lm:vertexpartition}, we may assume that for each~$i\in [k]$, the partition~$\cP_i^*$ is an~$\varepsilon_3$-superregular~$\frac{q}{f}$-refinement of~$\cP_i$.}
		For each~$i\in [k]$ and~$V_{i,j}\in \cP_i$, denote by $V_{i,j,1}, \dots, V_{i,j,\frac{q}{f}}$ the partition of~$V_{i,j}$ induced by~$\cP_i^*$.
		
		For each~$i\in [k]$, denote by $\cI_i\coloneqq \{I_{i,1}, \dots, I_{i,f}\}$ the canonical interval partition of~$C_i$ into~$f$ intervals of length~$\frac{2K}{f}$. 
		For each~$i\in [k]$,~$j\in [f]$ and~$h\in [\frac{q}{f}]$, 
		apply \cref{lm:regperfectmatching} to obtain a set~$\cE_{i,j,h}^{CA}$ of~$r_3$ vertex-disjoint edges of $\overrightarrow{G}_i[V_{i,j',h}, V_{i+1,j',h}]$, where~$j'\coloneqq \frac{jK}{f}$.%
		\COMMENT{Since $\frac{1}{K}\ll 1$, we may assume $\frac{2K}{f}\geq 3$ and so the fictive edges lie in $\overrightarrow{G}_i[V_{i,j',h}, V_{i+1,j',h}]$, i.e.\ in the penultimate pair of the interval~$I_{i,j}$ (see \cref{def:SPSfictive} in \cref{sec:robust}).}%
		\COMMENT{Find a perfect matching and delete an appropriate number of edges. This is possible since, by assumption, $r_3=\frac{rfK}{q}\leq\frac{fm}{qK^2}\leq \varepsilon|V_{i,j\frac{K}{f},h}|$.\label{com:r3}}
		Let $e_{i,j,h,1}, \dots, e_{i,j,h,r_3}$ be an enumeration of the edges in~$\cE_{i,j,h}^{CA}$.
		
		We construct fictive edges as follows. Let~$i\in [k], j\in [f]$, and~$h\in \left[\frac{q}{f}\right]$. For each~$t\in [r_3]$,
		let~$f_{i,j,h,t}$ be a fictive edge from~$x$ to~$y$, where~$x$ is the endpoint of~$e_{i-1,j,h,t}$ which belongs to~$V_i$ and~$y$ is the endpoint of~$e_{i+1, j, h,t}$ which belongs to~$V_{i+1}$. Let~$\cF_{i,j,h}^{CA}$ be the set of all these fictive edges. 
		\NEW{Observe that each edge in $\cF_{i,j,h}^{CA}$ will be a suitable fictive edge for a special path system of style $h$ spanning the interval $I_{i,j}$. Indeed, by construction of $\cE_{i-1,j,h}^{CA}$ and $\cE_{i+1,j,h}^{CA}$, $\cF_{i,j,h}^{CA}\subseteq E(\overrightarrow{G}_i[V_{i,j',h}, V_{i+1,j',h}])$, where~$j'\coloneqq \frac{jK}{f}$. I.e.\ each edge in $\cF_{i,j,h}^{CA}$ lies in the ``$h^{\rm th}$ subpair" of the penultimate pair along the interval $I_{i,j}$, as desired for \cref{def:SPSfictive}.}		
		Let~$\cE^{CA}$ be the union of the sets~$\cE_{i,j,h}^{CA}$ for each~$i\in [k], j\in [f],$ and~$h\in [\frac{q}{f}]$. Let~$\cE_i^{CA}$ be the union of the sets~$\cE_{i,j,h}^{CA}$ for each~$j\in [f]$ and~$h\in [\frac{q}{f}]$. Define~$\cF^{CA}$ and~$\cF_i^{CA}$ similarly.
		
		Note that for all~$j\in [f]$,~$h\in [\frac{q}{f}]$, and~$t\in [r_3]$, the graph $\left(\bigcup_{i\in [k]_{\rm odd}}f_{i,j,h,t}\right)\cup\left(\bigcup_{i\in [k]_{\rm even}}e_{i,j,h,t}\right)$ is a (directed) cycle of length~$k$ which intersects each of the clusters~$V_1, \dots, V_k$. The same holds with~$[k]_{\rm odd}$ and~$[k]_{\rm even}$ exchanged.	
		Therefore, in particular, the following property is satisfied.
		\begin{property}{\dagger\dagger}\label{property:fictiveedges}
			$\cE^{CA}\cup \cF^{CA}$ can be decomposed into edge-disjoint (directed) cycles of length~$k$, \NEW{each containing either}\OLD{either containing}
			\begin{itemize}
				\item an edge of~$\cE^{CA}$ between~$V_i$ and~$V_{i+1}$ for each~$i\in [k]_{\rm odd}$ and an edge of~$\cF^{CA}$ between~$V_i$ and~$V_{i+1}$ for each~$i\in [k]_{\rm even}$, or
				\item an edge of~$\cF^{CA}$ between~$V_i$ and~$V_{i+1}$ for each~$i\in [k]_{\rm odd}$ and an edge of~$\cE^{CA}$ between~$V_i$ and~$V_{i+1}$ for each~$i\in [k]_{\rm even}$.
			\end{itemize}
		\end{property}
		Property \cref{property:fictiveedges} will eventually enable us to construct Hamilton cycles of~$G-V_0$ by tying together a Hamilton path of each pair in~$M$ using an edge from each pair in~$M'$, or vice versa (recall \cref{fig:sketch}).
		
		We now select the fictive edges which will be included in the special factors required for finding the parity chord absorbers (see part \cref{lm:robustdecompPCA} of \cref{lm:robustdecomp}). We proceed as above to construct, for each~$i\in [k]$ and $j\in [7]$, a set~$\cE_{i,j}^{PCA}$ of~$5r^\diamond$ vertex-disjoint edges of $\overrightarrow{G}_i[V_{i,j'}, V_{i+1, j'}]$, where $j'\coloneqq \frac{j K}{7}$. For each~$i\in [k]$, since~$\cE_i^{CA}$ is a matching, by \cref{lm:verticesedgesremoval}, we can ensure that~$\cE_i^{CA}$ and $\cE_{i,j}^{PCA}$ are edge-disjoint. Let~$\cE^{PCA}$ be the union of the sets~$\cE_{i,j}^{PCA}$, for $i\in [k]$ and $j\in [7]$. For each~$i\in [k]$ and $j\in [7]$, we construct a set~$\cF_{i,j}^{PCA}$ of fictive edges as above and let~$\cF^{PCA}$ be the union of these sets. Importantly, the edges in~$\cE^{PCA}$ and~$\cF^{PCA}$ satisfy (the analogue of) property \cref{property:fictiveedges}. 
		
		Define $\cE\coloneqq \cE^{CA}\cup \cE^{PCA}$ and $\cF\coloneqq \cF^{CA}\cup \cF^{PCA}$. Delete from~$\overrightarrow{G}$ (and~$G$) all the edges in~$\cE$. For each~$i\in [k]$, note that we have deleted form~$\overrightarrow{G}_i$ at most~$2$ edges incident to each vertex so, by \cref{lm:verticesedgesremoval}, $(\overrightarrow{G}_i, \cP_i, \cP_i', \overrightarrow{R}_i, \overrightarrow{R}_i', C_i)$ is still an $(\ell, 2K, \frac{m}{K}, \varepsilon_4, \frac{d}{2})$-bi-setup \NEW{and $\cP_i^*$ is still a $\varepsilon_4$-superregular $\frac{q}{f}$-refinement of $\cP_i$.}
		
		\item\textbf{Constructing the special factors.}\label{step:SF}
		\NEW{In this step, we will construct, for $i\in [k]$, edge-disjoint special factors~$SF_{i,1}, \dots, SF_{i,r_3}$ with parameters~$(\frac{q}{f},f)$ with respect to~$C_i$,~$\cP^*_i$ in~$\overrightarrow{G}_i$ such that, for all~$t\in [r_3]$, $\Fict(SF_{i,t})=\{f_{i,j,h,t}\mid j\in [f], h\in [\frac{q}{f}]\}$. (Recall \cref{fig:SF-SPS}.)}
		
		\NEW{Let $i\in [k]$, $j\in [f]$, and $h\in [\frac{q}{f}]$. Suppose inductively that, for some $0\leq t\leq r_3$, we have constructed edge-disjoint special path systems $SPS_{i,j,h,1}, \dots, SPS_{i,j, h, t}$ of style $h$ in $\overrightarrow{G}_i$ spanning the interval $I_j$ and such that, for each $i'\in [t]$, $f_{i,j,h,i'}$ is the fictive edge contained in $SPS_{i,j,h,i'}$.
		If $t<r_3$, we construct $SPS_{i,j,h,t+1}$ as follows. For simplicity, denote $I_{i,j}=U_1\dots U_{\frac{2K}{f}+1}$ and, for each $i'\in [\frac{2K}{f}+1]$, let $U_{i',h}$ denote the $h^{\rm th}$ subcluster of $U_{i'}$ in $\cP_i^*$.
		Let $\overrightarrow{G}_i'\coloneqq \overrightarrow{G}_i\setminus \bigcup_{i'\in [t]}SPS_{i,j,h,i'}$. By \cref{lm:verticesedgesremoval} and since $r_3\leq \varepsilon_4\frac{fm}{qK}=\varepsilon_4|U_{i',h}|=\varepsilon_4|U_{i'+1,h}|$, $\overrightarrow{G}_i'[U_{i',h}, U_{i'+1,h}]$ is $[2\sqrt{\varepsilon_4}, \geq \frac{d}{2}]$-superregular for each $i'\in [\frac{2K}{f}]$ and $\overrightarrow{G}_i'[U_{\frac{2K}{f}-1,h}\setminus V(f_{i,j,h,t+1}), U_{\frac{2K}{f},h}\setminus V(f_{i,j,h,t+1})]$ is $[2\sqrt{\varepsilon_4}, \geq \frac{d}{2}]$-superregular.
		For each $i'\in [\frac{2K}{f}]\setminus \{\frac{2K}{f}-1\}$, apply \cref{lm:regperfectmatching} with $\overrightarrow{G}_i'[U_{i',h}, U_{i'+1,h}], 2\sqrt{\varepsilon_4}$, and $\frac{d}{3}$ playing the roles of $G, \varepsilon$, and $\alpha$ to obtain a perfect matching $M_{i'}$ in $\overrightarrow{G}_i'[U_{i',h}, U_{i'+1,h}]$.
		Apply \cref{lm:regperfectmatching} with $\overrightarrow{G}_i'[U_{\frac{2K}{f}-1,h}\setminus V(f_{i,j,h,t+1}), U_{\frac{2K}{f},h}\setminus V(f_{i,j,h,t+1})], 2\sqrt{\varepsilon_4}$, and $\frac{d}{3}$ playing the roles of $G, \varepsilon$, and $\alpha$ to obtain a perfect matching $M_{\frac{2K}{f}-1}'$ in $\overrightarrow{G}_i'[U_{\frac{2K}{f}-1,h}\setminus V(f_{i,j,h,t+1}), U_{\frac{2K}{f},h}\setminus V(f_{i,j,h,t+1})]$. Let $M_{\frac{2K}{f}-1}\coloneqq M_{\frac{2K}{f}-1}'\cup \{f_{i,j,h,t+1}\}$. Then, $SPS_{i,t+1}\coloneqq \bigcup_{i'\in [\frac{2K}{f}]} M_{i'}$ is a special path system of style $h$ in $\overrightarrow{G}_i$ spanning the interval $I_j$ which is edge-disjoint from $SPS_{i,1}, \dots, SPS_{i,t}$ and which contains the fictive edge $f_{i,j,h,t+1}$.}
		
		\NEW{Thus, we can construct, for each $i\in [k]$, $j\in [f]$, and $h\in [\frac{q}{f}]$, edge-disjoint special path systems $SPS_{i,j,h,1}, \dots, SPS_{i,j, h, r_3}$ of style $h$ in $\overrightarrow{G}_i$ spanning the interval $I_j$ such that, for each $t\in [r_3]$, $f_{i,j,h,t}$ is the fictive edge contained in $SPS_{i,j,h,t}$.
		For each $i\in [k]$ and $t\in [r_3]$, let $SF_{i,t}\coloneqq \bigcup_{j\in [f]}\bigcup_{h\in [\frac{q}{f}]}SPS_{i,j,h,t}$. Then, for each $i\in [k]$, $SF_{i,1}, \dots, SF_{i,r_3}$ are edge-disjoint special factors with parameters $(\frac{q}{f},f)$ with respect to~$C_i$,~$\cP^*_i$ in~$\overrightarrow{G}_i$ such that, for all~$t\in [r_3]$, $\Fict(SF_{i,t})=\{f_{i,j,h,t}\mid j\in [f], h\in [\frac{q}{f}]\}$. For each $i\in [k]$, let $\cS\cF_i\coloneqq SF_{i,1} \cup \dots \cup SF_{i,r_3}$.} 
		\OLD{Let~$i\in [k]$. Recall that $\cF_i^{CA}=\{f_{i,j,h,t} \mid j\in [f], h\in [\frac{q}{f}], t\in [r_3]\}$ and, for each $j\in [f], h\in [\frac{q}{f}]$, and $t\in [r_3]$, the fictive edge~$f_{i,j,h,t}$ has an endpoint in~$V_{i,j',h}$ and an endpoint in~$V_{i+1,j', h}$, where $j'\coloneqq \frac{jK}{f}$. 
		Thus, by \cref{lm:verticesedgesremoval} and since $r_3\leq \varepsilon\frac{fm}{qK}$%
		%\COMMENT{ $r_3\leq\varepsilon|V_{i,j,h}|$ (see \cref{com:r3})},
		one can repeatedly apply \cref{lm:regperfectmatching} to construct edge-disjoint special factors~$SF_{i,1}, \dots, SF_{i,r_3}$ with parameters~$(\frac{q}{f},f)$ with respect to~$C_i$,~$\cP^*_i$ in~$\overrightarrow{G}_i$ such that, for all~$t\in [r_3]$, $\Fict(SF_{i,t})=\{f_{i,j,h,t}\mid j\in [f], h\in [\frac{q}{f}]\}$.  Let $\cS\cF_i\coloneqq SF_{i,1} \cup \dots \cup SF_{i,r_3}$.}
		
		\item\textbf{Finding the robustly decomposable digraphs.}
		For any~$i\in [k]$, we apply \cref{lm:robustdecomp} with $\frac{m}{K}$, $2K$, $\varepsilon_4$, $\frac{d}{2}$, $\overrightarrow{G}_i$, $\cP_i$, $\cP_i'$, $\overrightarrow{R}_i$, $\overrightarrow{R}_i'$, $C_i$, and~$\cP_i^*$ playing the roles of $m$, $k$, $\varepsilon$, $d$, $\overrightarrow{G}$, $\cP$, $\cP'$, $\overrightarrow{R}$, $\overrightarrow{R}'$, $C$, and~$\cP^*$ to obtain a digraph~$\overrightarrow{CA}_i(r)$ satisfying the properties described in \cref{lm:robustdecomp}. 
		
		Since $\overrightarrow{CA}_i(r)\cup\cS\cF_i$ is~$(r_1+r_2+r_3)$-regular and~$r_1+r_2+r_3\leq 3r_1$, \cref{lm:verticesedgesremoval} implies that one can proceed similarly as in \cref{step:SF} to construct special factors $SF_{i,1}', \dots, SF_{i,r^\diamond}'$ with parameters~$(1, 7)$ with respect to~$C_i, \cP_i$ in~$\overrightarrow{G}_i$ which are edge-disjoint from each other and from $\overrightarrow{CA}_i(r)\cup \cS\cF_i$. 
		
		Let~$\overrightarrow{PCA}_i(r)$ and~$\overrightarrow{G}_i^{\rm rob}$ be as in \cref{lm:robustdecomp}.
		For each~$i\in [k]$, delete the edges of~$\overrightarrow{G}_i^{\rm rob}$ from~$G_i$ \NEW{(and $G$)}. Since~$\overrightarrow{G}_i^{\rm rob}$ is $(r_1+r_2+r_3+5r^\diamond+r^\diamond)$-regular, we have deleted at most\COMMENT{We have ``at most'' here since~$\overrightarrow{G}_i^{\rm rob}$ is~$(r_1+r_2+r_3+6r^\diamond)$-regular when ``including" the fictive edges contained in the special factors, but of course these fictive edges are not deleted from~$G$.} $2(r_1+r_2+r_3+5r^\diamond+r^\diamond)\leq 30 r_1\leq \varepsilon m$ edges incident to each vertex in~$V_i\cup V_{i+1}$. Moreover, recall that, at the end of \cref{step:fict}, we have already deleted from~$G$ the edges in~$\cE$, which contains at most two edges incident to each vertex in~$G_i$, for each~$i\in [k]$. Thus, \cref{lm:verticesedgesremoval} implies that~$G_i$ is still~$[\varepsilon_1^*,d]$-superregular. Furthermore, \cref{property:fictiveedges} and its analogue for~$\cE^{PCA}$ and~$\cF^{PCA}$ ensure that~$G$ is still Eulerian.

		\item \textbf{Regularising the superregular pairs.}\label{step:regularising}
		In order to apply \cref{lm:approxdecomp}, we need to regularise each superregular pair of~$G$. We will first apply the tools of \cref{sec:regularising} to each~$G_i$ separately, but, we will see that this yields too many cycles. We will therefore use a few further edges of~$G$ to tie together some of the cycles obtained to form longer cycles. We make sure that, when tying some of the cycles together, we use only a bounded number of edges incident to each vertex. Thus, applying the tools of \cref{sec:regularising} once again will only yield a few additional cycles.
		
		First, apply \cref{lm:Eulerianpairs} \NEW{(to the current graph $G$)} with~$\varepsilon_1^*$ and $m$ playing the roles of~$\varepsilon$ and~$m'$. Add the resulting cycles to~$\cD$ and delete their edges from~$G$. Note that we have added at most~$c$ cycles to~$\cD$. Moreover, for each~$i\in [k]$, the pair~$G_i$ is now Eulerian and~$[\varepsilon_2^*, d]$-superregular.
		
		For each~$i\in [k]$, apply \cref{lm:regularising} with~$G_i, \varepsilon_2^*$, and $2\varepsilon_2^*m$ playing the roles of~$G, \varepsilon$, and $\Theta$ in order to obtain a set~$\sC_i$ of at most~$4\varepsilon_2^*m$ edge-disjoint cycles of length at least $\frac{2m}{3}$ such that the following holds. Delete the edges of~$\sC_i$ from~$G_i$. Then,~$G_i$ is regular and~$[\varepsilon_3^*,d]$-superregular. By adding additional edge-disjoint Hamilton cycles to~$\sC_i$ if necessary, we may assume that~$G_i$ is~$(d-10\varepsilon_2^*)m$-regular and $|\sC_i|\leq 10\varepsilon_2^*m$.\COMMENT{The graph we obtain is Eulerian, has degree at most~$(d+\varepsilon_2^*)m$, and degree at least $(d-2\varepsilon_2^*)m-8\varepsilon_2^*m=(d-10\varepsilon_2^*)m$. Thus, we need to remove at most~$\frac{11\varepsilon_2^*m}{2}$ Hamilton cycles. This is possible since, by assumption,~$(d-10\varepsilon_2^*)m$ is even.}
		We observe that~$\bigcup_{i\in [k]} \sC_i$ may contain up to~$10\varepsilon_2^*mk$ cycles, so we need to split each of these cycles into paths and tie them together to form fewer cycles.
		
		Let~$i\in [k]$. Split one by one each cycle in~$\sC_i$ into at most~$\frac{30}{d}$ paths of length at most~$\frac{d m}{10}$, each with an endpoint in~$V_i$ and an endpoint in~$V_{i+1}$, and such that each vertex in~$V_i\cup V_{i+1}$ is an endpoint of at most~$2$ paths. This is possible since the cycles in~$\sC_i$ have length at least~$\frac{2m}{3}$ while, on the other hand, in each step there are at most $\frac{30}{d}|\sC_i|\leq \varepsilon_3^* m$ vertices in each cluster which are already endpoints of~$2$ paths.
		%Thus, given a cycle~$C\in \sC_i$, we can select two distinct vertices~$x\in V(C)\cap V_i$ and~$y\in V(C)\cap V_{i+1}$ such that~$x$ and~$y$ are not yet endpoints of any paths and split~$C$ into two~$(x,y)$-paths. 
		Let~$\sP_i$ be the set of paths obtained at the end of this procedure and observe that $|\sP_i|\leq \varepsilon_3^* m$.
		
		Decompose $\bigcup_{i\in [k]_{\rm odd}}\sP_i$ into at most~$\varepsilon_3^* m$ sets of paths, each containing at most one path in~$\sP_i$ for each~$i\in [k]_{\rm odd}$. Decompose $\bigcup_{i\in [k]_{\rm even}} \sP_i$  similarly. Let~$\sP_1', \dots, \sP_{\ell'}'$ be the sets of paths obtained. Thus, $\ell'\leq 2\varepsilon_3^* m$. 
		Apply \cref{lm:matchingtying} with $\varepsilon_3^*, 2\varepsilon_3^*,\ell'$, and~$\sP_1', \dots, \sP_{\ell'}'$ playing the roles of~$\varepsilon, \zeta, \ell$, and~$\sP_1, \dots, \sP_\ell$ to obtain~$E\subseteq E(G)$ such that $(\sP_1'\cup \dots \cup \sP_{\ell'}')\cup E$ can be decomposed into~$\ell'$ cycles. Add these cycles to~$\cD$ and delete from~$G$ all the edges in~$E$.
		Note that, for each~$i\in [k]$, by \cref{lm:verticesedgesremoval} and part \cref{lm:matchingtyingendpoints} of \cref{lm:matchingtying},~$G_i$ is now~$[\varepsilon_4^*, d]$-superregular with maximum degree at most~$(d-10\varepsilon_2^*)m$ and minimum degree at least~$(d-10\varepsilon_2^*)m-6$.
		
		We now need to regularise the superregular pairs of~$G$ once again. First, we apply \cref{lm:Eulerianpairs}\COMMENT{with $\varepsilon_4^*$ and $m$ playing the roles of $\varepsilon$ and $m'$} and add the resulting cycles to~$\cD$. Then, we apply \cref{lm:regularising} to~$G_i$, for each $i\in [k]$, and add all cycles obtained to~$\cD$. Using similar arguments as above, we may assume that, for each~$i\in [k]$, the pair~$G_i$ is now~$[\varepsilon_5^*, d]$-superregular and~$d' m$-regular.%
		\COMMENT{We apply \cref{lm:regularising} with~$\Theta\leq 6+2c$ and maximum degree at least~$(d-10\varepsilon_2^*)m-6-2c$ so by \cref{lm:regularising} the resulting graph has degree at least $(d-10\varepsilon_2^*)m-5(6+2c)>d'm$. Also recall that~$d'm$ is assumed to be even.}
		We note that $|\cD|\leq 3\varepsilon_3^* m\leq \varepsilon^{\frac{1}{16}}m$, as desired.\COMMENT{$|\cD|\leq c+2\varepsilon_3^*m+ \frac{(d-10\varepsilon_2^*)m-d'm}{2}$.}			
		
		\item \textbf{Approximately decomposing (the remainder of)~$G$.}\label{step:approx}
		Apply \cref{lm:approxdecomp} with $\varepsilon_5^*, d',d''$, and $2r$ playing the roles of $\varepsilon, d, d'$, and $r$ to obtain $H\subseteq G$ such that, for each $i\in [k]$, $H_i\coloneqq H[V_i, V_{i+1}]$ is $2r$-regular and $G\setminus H$ can be decomposed into $d'm-2r$ cycles. Add these cycles to~$\cD$.
		
		\item \textbf{Decomposing the leftover and robustly decomposable graphs.}\label{step:leftover}
		Let \NEW{$i\in [k]$}\OLD{$i\in [r]$}. Since $H_i$ is $2r$-regular, there exists an orientation $\overrightarrow{H}_i$ of $H_i$ such that $\overrightarrow{H}_i[V_i\cup V_{i+1}]$ is an $r$-regular bipartite oriented graph with vertex classes $V_i$ and $V_{i+1}$. Let~$\cD_i$ be the Hamilton decomposition of $\overrightarrow{H}_i\cup \overrightarrow{G}_i^{\rm rob}$ guaranteed by \cref{lm:robustdecomp}. Note that, in particular, each cycle in~$\cD_i$ contains exactly one fictive edge and thus corresponds to a Hamilton path of the original graph $G[V_i, V_{i+1}]$. Moreover, $|\cD_i|=s$.
		
		We form~$2s$ cycles by removing the fictive edges in~$\cF$ and inserting back the edges in~$\cE$ as follows. Fix a decomposition of~$\cE\cup \cF$ into edge-disjoint cycles of length~$k$ satisfying the property described in \cref{property:fictiveedges}. Let~$C$ be a cycle in this decomposition and assume without loss of generality that the fictive edges in~$C$ lie between~$V_i$ and~$V_{i+1}$ for~$i\in [k]_{\rm odd}$. Let~$f_1,e_2,f_3,\dots, e_k$ be an enumeration of the edges of~$C$ where, for each~$i\in [k]$, the edge~$f_i$ (respectively~$e_i$) lies between~$V_i$ and~$V_{i+1}$. For each~$i\in [k]_{\rm odd}$, let~$C_i\in \cD_i$ be the cycle which contains~$f_i$. Then, by construction, $(\bigcup_{i\in [k]_{\rm odd}}C_i\setminus \{f_i\} )\cup (\bigcup_{i\in [k]_{\rm even}} e_i)$ is a cycle and we add this cycle to~$\cD$.
		We proceed in this way for every cycle~$C$ in the cycle decomposition of $\cE\cup \cF$. This gives a cycle decomposition~$\cD$ of our original graph~$G$ of size at most $(d'm-2r)+2s+\varepsilon^{\frac{1}{16}}m\leq dm+\varepsilon^{\frac{1}{16}}m$.
		\COMMENT{$(d'm-2r)+2s+\varepsilon^{\frac{1}{16}}m=dm-11\varepsilon_2^*m-2r+4rfK+14r^\diamond+\varepsilon^{\frac{1}{16}}m\leq dm +\varepsilon^{\frac{1}{16}}m + 4rfK+28r_1-11\varepsilon^{\frac{1}{7}}m\leq dm +\varepsilon^{\frac{1}{16}}m +\frac{fm}{K^2}+\frac{28 m}{K}-11\varepsilon^{\frac{1}{7}}m\leq dm+\varepsilon^{\frac{1}{16}}m$}
		\qed
	\end{steps}
	\renewcommand{\qed}{}
\end{proof}

\onlyinsubfile{\bibliographystyle{plain}
	\bibliography{Bibliography/papers}}

%% file: Proof_of_Main_Theorems.tex
\onlyinsubfile{
\setcounter{section}{5}
\setcounter{subsection}{6}
\addtocounter{subsection}{-1}
\setcounter{definition}{13}
}

\subsection{Proof of the main theorems}\label{sec:thms}

We are now ready to prove Theorems \labelcref{thm:n/2}\labelcref{thm:n/2pathdecomp}, \labelcref{thm:n/2}\labelcref{thm:n/2cycledecomp}, \labelcref{thm:Delta/2,thm:epspathdecomp}. 

\begin{proof}[Proof of \cref{thm:n/2}\cref{thm:n/2cycledecomp}]
	Let~$\cD\coloneqq \emptyset$. We will repeatedly add cycles to~$\cD$. The set~$\cD$ will eventually be the set of cycles for our final decomposition of~$G$.
	The proof is structured as described in \cref{sec:sketch}.
	
	Fix additional constants such that $0<\frac{1}{n_0}\ll \frac{1}{M}\ll \varepsilon\ll \zeta \ll d \ll \beta \ll \alpha, \delta \leq 1$ and $0<\frac{1}{M}\leq \frac{1}{L}\ll \frac{1}{K}\ll d\ll \frac{1}{q}\ll \frac{1}{f}\ll \beta \ll \frac{1}{\ell},\frac{1}{g}\ll1$, with $\frac{K}{7}. \frac{2K}{f}, \frac{2K}{g}, \frac{q}{f}, \frac{4fK}{3g(g-1)}, \frac{\ell}{2}\in \mathbb{N}^*$.
	Let~$G$ be a graph on~$n\geq n_0$ vertices with $\delta(G)\geq \alpha n$. Let~$V\coloneqq V(G)$ and 
	\[\varepsilon'\coloneqq \varepsilon^{\frac{1}{75}},
	\quad \zeta'\coloneqq \zeta^{\frac{1}{3}},
	\quad d_1\coloneqq d^{\frac{1}{10^7}},
	\quad d_2\coloneqq d_1^{\frac{1}{33}}.
	\]
	\begin{steps}
		\item \textbf{Applying Szemer\'{e}di's regularity lemma and setting aside \NEW{some} random subgraphs~$\Gamma$ and~$\Gamma'$.}\label{step:cleaning}
		Apply \cref{lm:goodgraph} with parameters $M, L, \varepsilon, \zeta, d, \beta, \alpha$ and with $4q\ell K$ playing the role of~$r$ to obtain parameters~$M', m'\in \mathbb{N}^*$, a decomposition of~$G$ into four edge-disjoint graphs~$G^*, \Gamma, \Gamma'$, and~$H$, and a partition of~$V$ into~$k$ clusters~$V_1, \dots, V_k$ and an exceptional set~$V_0$ satisfying the properties described in \cref{lm:goodgraph}. In particular, the following property is satisfied. 
		\begin{property}{\ddagger}\label{property:reduceddecomp}
			The reduced graph~$R'$ of~$\Gamma$ admits a decomposition~$\cD_{R'}$ such that the following hold.~$\cD_{R'}$ consists of at most~$\frac{k}{2}$ cycles whose lengths are even and at least~$L$. Moreover, for any distinct~$i,j,j'\in [k]$, if~$V_jV_iV_{j'}$ is a subpath of a cycle in~$\cD_{R'}$ then the support clusters of~$V_i$ with respect to~$V_j$ and~$V_{j'}$ are the same.
		\end{property}

		\item \textbf{Covering the edges of~$G[V_0]$.}\label{step:insideV0}
		Apply \cref{lm:insideV0cycles} with~$G^*$ playing the role of~$G$ to obtain a graph~$H_0\subseteq G^*\cup \Gamma$ satisfying properties \cref{lm:insideV0coverV0,lm:insideV0Gstillsupreg,lm:insideV0Gammastillsupreg,lm:insideV0cycledecomp} of \cref{lm:insideV0cycles}. In particular, there exists a decomposition~$\cD_0\cup \cD_0'$ of~$H_0$ such that~$\cD_0$ is a set of at most~$\beta n$ cycles and~$\cD_0'$ is a set of at most~$\beta^{-2}$ edges. Add the cycles in~$\cD_0$ to~$\cD$. Since~$G$ is Eulerian, by \cref{fact:euleriancycle}, we can cover the edges in~$\cD_0'$ with at most~$\beta^{-2}$ edge-disjoint cycles. Add these cycles to~$\cD$ and delete the edges in all these cycles from $G, G^*, \Gamma$, and $\Gamma'$. Observe that by \cref{lm:verticesedgesremoval,lm:insideV0cycles},~$V_0, V_1, \dots, V_k$ is  
		\begin{itemize}
			\item an $(\varepsilon', \geq d, k, m, m', R)$-superregular equalised partition of~$G^*$;
			\item an $(\varepsilon', \beta, k, m, m', R')$-superregular equalised partition of~$\Gamma$; and
			\item an $(\varepsilon', \zeta, k, m, m', R'')$-superregular equalised partition of~$\Gamma'$,
		\end{itemize}
		where $R'$ and $R''$ are edge-disjoint and satisfy $R'\cup R''=R$. Moreover,~$G[V_0]$ is now empty, as desired.
	
		\item \textbf{Covering most of~$G^*$ with at most roughly~$\frac{n}{2}$ cycles.}\label{step:exceptional}
		We now apply \cref{lm:exceptionaledges} with~$G^*$ and~$\varepsilon'$ playing the roles of~$G$ and~$\varepsilon$ to obtain a decomposition of~$G^*\cup \Gamma \cup \Gamma'$ into edge-disjoint graphs~$G', \tGamma$, and~$H'$ such that $G^*, \Gamma' \subseteq G' \cup H'$,~$\tGamma \subseteq \Gamma$, and properties \cref{lm:exceptionaledgesDeltaH,lm:exceptionaledgesGamma,lm:exceptionaledgesdecomp} of \cref{lm:exceptionaledges} are satisfied. In particular, there exists a decomposition~$\cD'\cup \cD'_{\rm exc}$ of~$G'$ such that~$\cD'$ is a set of at most~$\frac{n}{2}+2\beta n$ cycles and~$\cD'_{\rm exc}$ is a set of at most~$\beta^{-2}$ edges.
		Add all cycles in~$\cD'$ to~$\cD$. Apply \cref{fact:euleriancycle} with $\tGamma\cup H\cup H'\cup \cD_{\rm exc}'$ playing the role of~$G$ to cover the edges in~$\cD_{\rm exc}$ with at most~$\beta^{-2}$ edge-disjoint cycles. Add these cycles to~$\cD$ and delete the edges in all these cycles from $\tGamma, H$, and $H'$. By \cref{lm:verticesedgesremoval,lm:exceptionaledges},~$V_0, V_1, \dots, V_k$ is a~$(\zeta', \beta, k, m, m', R)$-superregular equalised partition of~$\tGamma$.
		Also note that $\Delta(H\cup H')\leq 4dn + 13\zeta n \leq 5dn$.

		\item \textbf{Covering the leftovers.}\label{step:bad}
		We now cover the edges of~$H\cup H'$ by applying \cref{lm:bad} with~$H\cup H', \tGamma$, and~$5d$ playing the roles of~$G, \Gamma$, and~$d$ to obtain a subgraph~$\tH\subseteq \tGamma$ and a decomposition~$\widetilde{\cD}$ of~$H\cup H'\cup \tH$ into at most~$2\beta n$ cycles. Add all cycles in~$\widetilde{\cD}$ to~$\cD$ and let $\tGamma'\coloneqq \tGamma \setminus \tH$. By \cref{lm:bad},~$V_0, V_1, \dots, V_k$ is a~$(d_1, \beta, k, m, m', R')$-superregular equalised partition of~$\tGamma'$.
		
		\item \textbf{Fully decomposing~$\Gamma$.}\label{step:Gamma}
		Finally, observe that~$\tGamma'$ is an Eulerian subgraph of~$\Gamma$ with the same reduced graph and the same support clusters, so property \cref{property:reduceddecomp} holds for~$\tGamma'$. Thus, we can apply \cref{cor:Gammadecomp} with~$\tGamma', R', d_1$, and~$\beta$ playing the roles of~$G, R, \varepsilon$, and~$d$ to obtain a decomposition of~$\tGamma'$ into at most~$\frac{\beta n}{2}+d_2 n$ cycles.
		Add these cycles to~$\cD$. Then,~$\cD$ forms a cycle decomposition of~$G$ and~$|\cD|\leq \frac{n}{2}+\delta n$, as desired\COMMENT{$|\cD|\leq \beta n+\beta^{-2}+\frac{n}{2}+2\beta n+\beta^{-2}+2\beta n+\frac{\beta n}{2}+d_2n$}. \qed
	\end{steps}
	\renewcommand{\qed}{}
	\end{proof}

	\begin{proof}[Proof of \cref{thm:n/2}\cref{thm:n/2pathdecomp}]
		We modify the proof of \cref{thm:n/2}\cref{thm:n/2cycledecomp} to \NEW{get} a path decomposition as follows. \Cref{step:cleaning} is identical. For \cref{step:insideV0}, we simply apply \cref{thm:Lovasz} to obtain a path decomposition of~$G[V_0]$ into at most~$\varepsilon n$ paths. 
		For \cref{step:exceptional}, first remove an edge incident to each exceptional vertex of odd degree in~$G^*$ so that property \cref{lm:exceptionaledgeseven} of \cref{lm:fictive} holds. We view these edges as individual paths in our decomposition. We can thus apply \cref{lm:fictive} instead of \cref{lm:exceptionaledges}, with the set of odd degree vertices of $G^*\cup H\cup \Gamma \cup \Gamma'$ playing the role of~$U$.
		Then, by \cref{lm:fictive}\cref{lm:fictivedecomp},~$H\cup H'\cup \tGamma$ is Eulerian at the end of \cref{step:exceptional}. 
		Thus, we can apply the arguments of \cref{step:Gamma,step:bad} and split each cycle obtained in these steps into two paths in order to obtain a path decomposition of~$H\cup H'\cup \tGamma$. One can easily verify that we obtain at most~$\frac{n}{2}+\delta n$ paths in total.
	\end{proof}

	The weak quasirandomness assumed in the next two proofs allows for a more efficient decomposition. The critical property implied by weak quasirandomness is that the reduced graph~$R$ is connected.

	\begin{proof}[Proof of \cref{thm:Delta/2}]
		We use the same arguments as in the proof of Theorems~\labelcref{thm:n/2}\labelcref{thm:n/2pathdecomp} and \labelcref{thm:n/2}\cref{thm:n/2cycledecomp} with $\beta\ll \alpha, \delta, p$ and applying \cref{lm:connectedR} instead of \cref{lm:exceptionaledges,lm:fictive}.
		This is possible since the reduced graph~$R$ of~$G^*$ is connected.
		
		Indeed, assume for a contradiction that~$R$ is disconnected and let~$C$ be a component of~$R$. Let~$A\coloneqq V_{G^*}(C)$ and~$B\coloneqq V(G)\setminus A$. Since~$\delta(G)\geq \alpha n$, it is easy to see that $|A|, |B|\geq \frac{\alpha n}{2}$.
		But, by \cref{lm:goodgraph}, 
		\begin{align*}
		e_G(A, B)&=e_H(A, B)+e_\Gamma(A, B)+e_{\Gamma'}(A,B)+e_{G^*}(A, V_0)\\
		&\leq |A|(4dn+(\beta+\varepsilon)n+(\zeta+\varepsilon)n+\varepsilon n)\\
		&< p|A||B|,
		\end{align*} 
		contradicting the fact that~$G$ is weakly-$(\frac{\alpha}{2},p)$-quasirandom.
		
		Also observe that in the path decomposition case, we have $\odd(G^*\cup H\cup \Gamma \cup \Gamma')\leq \odd(G)+|V_0|$\COMMENT{We could have a vertex~$x\in V_0$ with even degree in the original graph~$G$ but with $|N_{G^*}(x)\setminus V_0|$ odd, and so after covering the edges inside~$V_0$ with paths and after deleting our edge at~$x$ going to $V\setminus V_0$ from~$G^*$, we might have created an additional odd degree vertex for~$x$.} at the point where we apply \cref{lm:connectedR} (recall the proof of \cref{thm:n/2}\cref{thm:n/2pathdecomp}), so we obtain a decomposition of the desired size.
	\end{proof}

	\begin{proof}[Proof of \cref{thm:epspathdecomp}]
		First observe that since~$G$ is weakly-$(\varepsilon,p)$-quasirandom,~$G$ has fewer than~$\varepsilon n$ vertices of degree less that~$\frac{p n}{2}$. Let~$X$ be the set of these vertices. We modify the proof \cref{thm:Delta/2}\cref{thm:Delta/2pathdecomp} as follows.	
		Apply the arguments of \cref{step:cleaning} with~$G-X$ and~$\frac{p}{3}$ playing the roles of~$G$ and~$\alpha$. Add the vertices in~$X$ to the exceptional set and all edges incident to these vertices to~$G^*$. The remainder of the proof is identical.
	\end{proof}

\onlyinsubfile{\bibliographystyle{abbrv}
	\bibliography{Bibliography/papers}}

%% file: Concluding_Remarks.tex
	\onlyinsubfile{
		\setcounter{section}{6}
		\addtocounter{section}{-1}
\section{Concluding remarks}}

We conclude by deriving \cref{thm:n/2}\cref{thm:ErdosGallai} and providing some remarks on our results.

\subsection{Proof of Theorem~{\ref{thm:n/2}\ref{thm:ErdosGallai}}}\label{sec:ErdosGallai}

We now show how \cref{thm:n/2}\cref{thm:ErdosGallai} can be derived from \cref{thm:n/2}\cref{thm:n/2cycledecomp}. Let~$G$ be a graph. We saw in the introduction that one can remove at most~$n-1$ edges of~$G$ to obtain an Eulerian graph. However, in order to apply \cref{thm:n/2}\cref{thm:n/2cycledecomp}, we also need to make sure that the resulting Eulerian graph still has linear minimum degree.

\begin{proof}[Proof of \cref{thm:n/2}\cref{thm:ErdosGallai}]
	Fix~$\ell\coloneqq \odd(G)$. Let~$V_{\rm odd}$ be the set of odd-degree vertices of~$G$. We repeatedly remove short paths with endpoints in~$V_{\rm odd}$ (but~$\ell$ is left unchanged). Fix a maximum matching~$M$ of~$G[V_{\rm odd}]$. Delete the edges of~$M$ from~$G$ and remove the vertices in~$V(M)$ from~$V_{\rm odd}$. We observe that~$V_{\rm odd}$ is now an independent set of~$G$.
	
	If there exists a path~$xyz$ in~$G$ such that~$x,z\in V_{\rm odd}$ are distinct and fewer than~$\frac{\alpha n}{4}$ edges incident to~$y$ have been deleted so far, remove the edges~$xy$ and~$yz$ from~$G$ and the vertices~$x, z$ from~$V_{\rm odd}$. 
	We repeat this procedure until there exists no such path of length~$2$. 
	
	Then, we claim that $|V_{\rm odd}|\leq \frac{2}{\alpha}$. Indeed, at each stage, there are at most\COMMENT{$\leq \frac{2\cdot 2\cdot \frac{n}{2}}{\frac{\alpha n}{4}}=\frac{8}{\alpha}$}
	$\frac{8}{\alpha}$ vertices~$y\in V\setminus V_{\rm odd}$ such that we have deleted at least~$\frac{\alpha n}{4}$ edges incident to~$y$.
	By construction, for each~$x\in V_{\rm odd}$, no edge incident to~$x$ has been deleted from~$G$ and, thus,~$x$ has more than~$\frac{\alpha n}{2}$ neighbours~$y$ such that fewer than~$\frac{\alpha n}{4}$ edges incident to~$y$ have been deleted so far. Thus, we must have $|V_{\rm odd}|\leq \frac{2}{\alpha}$\COMMENT{Since, by construction, vertices in~$V_{\rm odd}$ have disjoint neighbourhoods}.
	
	Pair all remaining vertices of~$V_{\rm odd}$ such that, in each pair, the vertices belong to a same connected component of~$G$. By construction, $\delta(G)\geq \frac{3\alpha n}{4}-1$.
	\COMMENT{$\geq \alpha n-(\frac{\alpha n}{4}+1)$}%
	Thus, for each pair~$(x,y)$ in turn, we can find a path of length at most~$\frac{5}{\alpha}$ between~$x$ and~$y$, which we delete from~$G$\COMMENT{Since at each stage $\delta(G)\geq \frac{3\alpha n}{4}-1-\frac{2}{\alpha}\geq \frac{3\alpha n}{5}$.}.%
	\COMMENT{Indeed, assume for a contradiction that ${\rm diam}(G)\geq \frac{5}{\alpha}$. Then, there exist~$x,y\in V(G)$ with~$d(x,y)=\frac{5}{\alpha}$. Fix an~$(x,y)$-path $P\coloneqq x_0, x_1, \dots, x_{\ell'}$ of length~$\ell'=\frac{5}{\alpha}$. By minimality of~$\ell'$, for any~$0< i+2<j\leq \ell'$, we have~$N(x_i)\cap N(x_j)=\emptyset$. Moreover, for each $i\in \{0,3,6, \dots, \frac{5}{\alpha}\}$, we have $|N(x_i)|\geq \delta(G)\geq \frac{3\alpha n}{5}$. Then, $|V(G)|\geq \frac{5}{3\alpha}(\frac{3\alpha n}{5}+1)>n$, a contradiction.} 
	Let $\cP$ be the set of edge-disjoint paths deleted. Note that $|E(\cP)|\leq \frac{5}{\alpha^2}$.
	Moreover, we have deleted at most $\ell+\frac{5}{\alpha^2}\leq n+\frac{\delta n}{2}$ edges in total. Finally,~$G$ is Eulerian and $\delta(G)\geq \frac{\alpha n}{2}$. 
	Applying \cref{thm:n/2}\cref{thm:n/2cycledecomp} with~$\frac{\alpha}{2}$ and~$\frac{\delta}{2}$ playing the roles of~$\alpha$ and~$\delta$ completes the proof.
\end{proof}

\subsection{Some remarks on Theorem~\ref{thm:Delta/2}}

As discussed in the introduction, we now show that neither the linear minimum degree condition (or even the stronger assumption of linear connectivity), nor the weakly-$\left(\frac{\alpha}{2}, p\right)$-quasirandom property is sufficient on its own to obtain the bounds in \cref{thm:Delta/2}.

\begin{prop}
	For any odd integer~$n\geq 20$, there exists an~$\left\lfloor\frac{n}{10}\right\rfloor$-connected Eulerian graph~$G$ on~$2n$ vertices such that the following hold.
	\begin{enumerate}
		\item $G$ cannot be decomposed into fewer than $\max\left\{ \frac{\odd(G)}{2},\frac{\Delta(G)}{2}\right\} + \frac{n}{10}$ paths.\label{prop:mindegpathdecomp}
		\item $G$ cannot be decomposed into fewer than $\frac{\Delta(G)}{2}+ \frac{n}{10}$ cycles.\label{prop:mindegcycledecomp}
	\end{enumerate}
\end{prop}

\begin{proof}
	Assume~$G_1, G_2$ are two vertex-disjoint cliques of size~$n$ and let~$V_1\subseteq V(G_1)$ and~$V_2\subseteq V(G_2)$ with~$|V_1|, |V_2|= \left\lfloor\frac{n}{10}\right\rfloor$. 
	Let~$G$ be obtained from~$G_1\cup G_2$ by adding two edge-disjoint perfect matchings between~$V_1$ and~$V_2$. 
	Note that~$G$ is an~$\left\lfloor\frac{n}{10}\right\rfloor$-connected Eulerian graph on~$2n$ vertices with~$\Delta(G)=n+1$. 
	
	Since there are at most~$\frac{n}{5}$ edges between~$G_1$ and~$G_2$, any cycle decomposition of~$G$ will contain at most~$\frac{n}{10}$ cycles with edges of both~$G_1$ and~$G_2$ and these will cover at most~$\frac{n}{5}$ edges incident to each vertex of~$G$. Thus any cycle decomposition of~$G$ will contain at least~$\frac{4n+5}{10}$ cycles of~$G_1$\COMMENT{$\geq \frac{n+1-\frac{n}{5}}{2}$} and at least~$\frac{4n+5}{10}$ cycles of~$G_2$. Therefore, any cycle decomposition of~$G$ will contain at least $\frac{4n}{5}+1>\frac{\Delta(G)}{2}+\frac{n}{10}$ cycles\COMMENT{$\frac{4n}{5}+1=\frac{n+1}{2}+\frac{3n}{10}+\frac{1}{2}$}.
	Similar arguments show that~$G$ cannot be decomposed into fewer than $\frac{3n}{5}+1>\max\{\frac{\odd(G)}{2},\frac{\Delta(G)}{2}\}+\frac{n}{10}$ paths.\COMMENT{At least $2\cdot \frac{n+1-\frac{2n}{5}}{2}=\frac{3n}{5}+1=\frac{n+1}{2}+\frac{n}{10}+\frac{1}{2}$.}
\end{proof}	

\begin{prop}\label{prop:random}
	For all~$0<\alpha \leq 1$, and all~$n_0\in \mathbb{N}^*$, the following hold.
	\begin{enumerate}
		\item There exists a weakly-$\left(\frac{\alpha}{2}, \frac{\alpha^2}{100}\right)$-quasirandom graph~$G$ on~$n\geq n_0$ vertices such that~$G$ cannot be decomposed into fewer than $\max\left\{ \frac{\odd(G)}{2},\frac{\Delta(G)}{2}\right\} +\frac{\alpha n}{10}$ paths.\label{prop:randompathdecomp}
		\item There exists an Eulerian weakly-$\left(\frac{\alpha}{2}, \frac{\alpha^2}{100}\right)$-quasirandom graph~$G$ on~$n\geq n_0$ vertices such that~$G$ cannot be decomposed into fewer than $\frac{\Delta(G)}{2}+\frac{\alpha n}{10}$ cycles.\label{prop:randomcycledecomp}
	\end{enumerate}
\end{prop}	

\begin{proof}
	Let~$m$ be a sufficiently large odd integer, $\delta \coloneqq \frac{\alpha}{10}$, and $\ell\coloneqq \frac{2\delta m+4}{1-2\delta}$. 
	
	For part \cref{prop:randompathdecomp}, let~$S_\ell$ be a star with~$\ell$ leaves and~$K_m$ be a complete graph on~$m$ vertices such that~$V(S_\ell)\cap V(K_m)=\{x\}$ for some leaf~$x$ of~$S_\ell$. Let~$G\coloneqq K_m\cup S_\ell$. Then~$G$ is graph of order~$n\coloneqq m+\ell$, with~$\Delta(G)=m$ and at least~$\ell$ vertices of odd degree. Let~$A\cup B$ be a partition of~$V(G)$ with $|A|, |B|\geq \frac{\alpha n}{2}$. Then, both~$A$ and~$B$ contain at least~$\frac{\alpha n}{10}$ vertices of~$K_m$. Thus, $e_G(A, B)\geq \frac{\alpha^2 n^2}{100}\geq \frac{\alpha^2}{100}|A||B|$ and~$G$ is weakly-$(\frac{\alpha}{2},\frac{\alpha^2}{100})$-quasirandom. But, one can easily show that~$G$ cannot be decomposed into fewer than $\frac{m+1}{2}+\frac{\ell -2}{2}>\max\left\{ \frac{\odd(G)}{2},\left\lceil\frac{\Delta(G)+1}{2}\right\rceil\right\} +\frac{\alpha n}{10} $ paths%
	\COMMENT{The inequality follows since~$\ell=2\delta n+4$.}.	
	
	For part \cref{prop:randomcycledecomp}, let~$G$ be obtained from~$K_m$ by appending~$\frac{\ell}{2}$ vertex-disjoint triangles with exactly one endpoint in~$V(K_m)$. Clearly,~$G$ is an Eulerian graph on~$n\coloneqq m+\ell$ vertices with~$\Delta(G)=m+1$. Now let~$A\cup B$ be a partition of~$V(G)$ with~$|A|,|B|\geq \frac{\alpha n}{2}$. Then, similarly as before, it follows that~$G$ is weakly-$(\frac{\alpha}{2},\frac{\alpha^2}{100})$-quasirandom. But~$G$ cannot be decomposed into fewer than $\frac{m-1}{2}+\frac{\ell}{2}>\frac{\Delta (G)}{2}+\frac{\alpha n}{10}$ cycles.
\end{proof}

\subsection{Some remarks on Conjecture~\ref{conj:epscycledecomp}}

As discussed in the introduction, we show that the Erd\H{o}s-Gallai conjecture is equivalent to \cref{conj:epscycledecomp}.

\begin{prop}\label{prop:equivErdosGallai}
	\cref{conj:epscycledecomp} is equivalent to the Erd\H{o}s-Gallai conjecture (\cref{conj:Erdos-Gallai}). 
\end{prop}

\begin{proof}
	~$(\Leftarrow)$ Assume \cref{conj:Erdos-Gallai} holds and let~$c$ be a constant such that any~$N$-vertex graph can be decomposed into at most~$cN$ cycles and edges, for each~$N\in \mathbb{N}^*$. Let~$\varepsilon\ll c^{-1},p$. Let~$G$ be as in \cref{conj:epscycledecomp} and $\cD\coloneqq \emptyset$. We repeatedly add cycles to~$\cD$ until it forms a cycle decomposition of~$G$. 	
	Weak-$(\varepsilon,p)$-quasirandomness implies that fewer than~$\varepsilon n$ vertices of~$G$ have degree less that~$\frac{pn}{2}$. Let~$S$ be the set of these vertices and apply the arguments of \cref{step:sketchcleaning} of the proof \cref{thm:n/2}\cref{thm:n/2cycledecomp} with~$G-S$ and~$\frac{p}{2}$ playing the roles of~$G$ and~$\alpha$ to obtain a decomposition of~$G$ into $G^*, \Gamma, \Gamma'$, and $H$.    
	Add the vertices in~$S$ to the exceptional set~$V_0$ and all edges incident to these vertices to~$G^*$. Note that we now have~$|V_0|\leq 2\varepsilon n$.
	Moreover, by similar arguments as in the proof of \cref{thm:Delta/2}, the reduced graph~$R$ of~$G^*$ is connected. 
	
	Decompose~$G[V_0]$ into at most~$2c\varepsilon n$ cycles and edges. Add the cycles obtained to~$\cD$ and delete their edges from~$G$. By choosing a decomposition where the number of edges is minimal, we may assume~$G[V_0]$ is now a forest and, thus, contains at most~$2\varepsilon n$ edges.
	
	Then, we can decompose~$G[V_0]$ into $\ell \leq 2\varepsilon n$ edge-disjoint paths $P_1, \dots, P_\ell$ such that the following hold. For each $i\in [\ell]$, the endpoints~$x_i$ and~$y_i$ of~$P_i$ have odd degree in~$G[V_0]$ and, moreover, each~$x\in V_0$ is an endpoint of at most one of the~$P_i$.
	Let~$i\in [\ell]$. Since~$G$ is Eulerian, there exist $x_i'\in N_G(x_i)\setminus V_0$ and $y_i'\in N_G(y_i)\setminus V_0$.  If $x_i'=y_i'$, add the cycle $x_i'x_iP_iy_i$ to~$\cD$. Otherwise, let $P_i'\coloneqq x_i'x_iP_iy_iy_i'$. 
	
	Apply \cref{lm:closingcycle} with~$2\varepsilon$ playing the role of~$\varepsilon$ and each~$\cP_i$ consisting of exactly one of the paths~$P_j'$ constructed above.
	Add all the cycles obtained to~$\cD$ and delete their edges from $G^*, \Gamma,\Gamma'$, and $H$. Thus,~$G^*[V_0]$ is now empty.
 	By \cref{lm:verticesedgesremoval} and \cref{lm:closingcycle}\cref{lm:closingcycleneighbours}, $V_0, V_1, \dots, V_k$ is now an $(\varepsilon^{\frac{1}{9}}, \beta, k, m, m', R')$-superregular equalised partition of~$\Gamma$ and an $(\varepsilon^{\frac{1}{9}}, \zeta, k, m, m', R'')$-superregular equalised partition of~$\Gamma'$.
 	Decompose the remainder of~$G$ as in \cref{thm:Delta/2}\cref{thm:Delta/2cycledecomp} (see \cref{step:sketchGamma,step:sketchbad,step:sketchexceptional} of the proof of \cref{thm:n/2}\cref{thm:n/2cycledecomp}).
	
	~$(\Rightarrow)$ Assume \cref{conj:Erdos-Gallai} does not hold and assume for a contradiction that \cref{conj:epscycledecomp} is true. Fix~$\delta>0$ and~$\frac{1}{4}\geq p>0$. Let~$1>\varepsilon >0$ and~$n_0$ be as in \cref{conj:epscycledecomp} and fix a constant~$c$ such that $c\geq \delta (1+\frac{2}{\varepsilon})$. 
	Let~$H$ be an Eulerian graph on~$m\geq \varepsilon n_0$ vertices such that any cycle decomposition of~$H$ contains more than~$cm$ cycles. Note that such graph exists since, as mentioned in the introduction, the Erd\H{o}s-Gallai conjecture is equivalent to the problem of decomposing Eulerian graphs of order~$n$ in~$O(n)$ cycles, and, by assumption, the Erd\H{o}s-Gallai conjecture is false.
	
	Assume without loss of generality that~$\frac{2m}{\varepsilon}$ is an odd integer. Let~$G$ be the disjoint union of~$H$ and~$K_{\frac{2m}{\varepsilon}}$. Note that~$G$ is a graph on $n=(1+\frac{2}{\varepsilon})m\geq n_0$ vertices. Moreover, $\Delta(G)=\Delta(K_{\frac{2m}{\varepsilon}})$. Thus, any cycle decomposition of~$G$ will contain more than $\frac{\Delta(G)}{2}+cm \geq  \frac{\Delta(G)}{2}+\delta n$ cycles.	
	But, for any partition~$A, B$ of~$V(G)$ with~$|A|, |B|\geq \varepsilon n$, we have $|A|, |B|\geq \varepsilon(1+\frac{2}{\varepsilon})m\geq 2m$ and thus $|A\cap V(K_{\frac{2m}{\varepsilon}})|\geq \frac{|A|}{2}$ and $|B\cap V(K_{\frac{2m}{\varepsilon}})|\geq \frac{|B|}{2}$. Therefore, $e_G(A, B)\geq \frac{1}{4}|A||B|\geq p|A||B|$ and~$G$ is weakly-$(\varepsilon, p)$-quasirandom, a contradiction.
\end{proof}

Using \cref{thm:loglognErdosGallai} and the arguments of the proof of \cref{prop:equivErdosGallai}, one can show the following.

\begin{prop}\label{prop:loglogn}
    For any~$\delta, p>0$, there exist~$\varepsilon,n_0>0$ such that the following hold. If~$G$ is an Eulerian weakly-$(\frac{\varepsilon}{\log\log n},p)$-quasirandom graph on~$n\geq n_0$ vertices, then~$G$ can be decomposed into at most~$\frac{\Delta(G)}{2} +\delta n$ cycles.
\end{prop}

\onlyinsubfile{\bibliographystyle{abbrv}
	\bibliography{Bibliography/papers}}

%% file: Appendix.tex
\onlyinsubfile{\appendix
\section*{Appendix: Proof of Lemma~\ref{lm:goodgraph}}}

\renewcommand{\thedefinition}{\Alph{definition}}
\setcounter{definition}{0}

\begin{lm}[{Szemer\'edi's regularity lemma \cite{szemeredi1978regular}}]\label{lm:Szemeredi}
	For every~$\varepsilon>0$ and~$M\in \mathbb{N}^*$, there exist $n_0,M'\in \mathbb{N}^*$ such that the following holds. For any graph~$G$ on~$n\geq n_0$ vertices, there exists a partition of~$V(G)$ into~$k$ clusters~$V_1, \dots, V_k$ and an exceptional set~$V_0$ satisfying the following properties.
	\begin{enumerate}
		\item $M\leq k\leq M'$.
		\item $|V_1|=\dots =|V_k|=:m$.
		\item $|V_0|\leq \varepsilon n$.
		\item All but at most~$\varepsilon k^2$ pairs~$G[V_i, V_j]$ with~$1\leq i<j\leq k$ are~$\varepsilon$-regular.
	\end{enumerate}
\end{lm}

We will also need the following easy observation about $\varepsilon$-regular pairs.

\begin{lm}[{Similar to \cite[Proposition 4.2]{kuhn2013hamilton}}]\label{lm:badvertices}
	Let~$0< \frac{1}{m_A}, \frac{1}{m_B} < \varepsilon \leq d \leq 1$ and~$G$ be an~$(\varepsilon, d)$-regular bipartite graph on vertex classes~$A$ and~$B$ of size~$m_A$ and~$m_B$, respectively. Then, fewer than~$\varepsilon m_A$ vertices in~$A$ have degree at least~$(d+\varepsilon)m_B$ and fewer than~$\varepsilon m_A$ vertices in~$A$ have degree at most~$(d-\varepsilon)m_B$.
\end{lm}

\begin{proof}[Proof of \cref{lm:goodgraph}]
	We start by applying Szemer\'edi's regularity lemma (\cref{step:Szemeredi}) and a cleaning procedure (\cref{step:cleaningexceptionaledges,step:smalldensitypairs,step:splitting,step:nonregular,step:nonexceptionaldeg,step:degrees}) similar to the one used to prove the degree form of the regularity lemma (see for instance \cite{taylor2014regularity}) to obtain a graph which admits a superregular partition. In \cref{step:randomslices,step:lovasz}, we obtain~$\Gamma$ and~$\Gamma'$ satisfying the desired properties by taking appropriate random subgraphs of~$G$. Finally, in \cref{step:support}, we equalise all support cluster sizes.
	\begin{steps}
		\item \textbf{Applying the regularity lemma.}\label{step:Szemeredi} Fix additional constants $\varepsilon_1, \varepsilon_2 >0$ and~$N\in \mathbb{N}^*$ such that $\frac{1}{N}\ll\varepsilon_1\ll \frac{1}{M}\ll \varepsilon_2 \ll \varepsilon$. Let 
		\[\varepsilon_1'\coloneqq \varepsilon_1^{\frac{1}{3}}, 
		\quad \varepsilon_1''\coloneqq \varepsilon_1^{\frac{1}{7}},
		\quad \varepsilon_1^*\coloneqq \varepsilon_1^{\frac{1}{14}},
		\quad \varepsilon_1^{**}\coloneqq \varepsilon_1^{\frac{1}{29}},
		\quad\varepsilon_2'\coloneqq \varepsilon_2^{\frac{1}{12}},
		\quad \varepsilon_2''\coloneqq \varepsilon_2^{\frac{1}{84}}.\]
		 Let~$N'$ be the constant obtained by applying \cref{lm:Szemeredi} with~$\varepsilon_1$ and~$N$ playing the role of~$\varepsilon$ and~$M$, respectively. Let~$n_0\in \mathbb{N}^*$ be such that $\frac{1}{n_0}\ll \frac{1}{N}$, and let~$G$ be a graph on~$n\geq n_0$ vertices. By \cref{lm:Szemeredi}, we can partition~$V$ into~$k_1$ clusters~$V_1^1, \dots, V_{k_1}^1$ and an exceptional set~$V_0^1$ such that the following hold.
		\begin{enumerate}[label=(\alph*)]
			\item $N\leq k_1\leq N'$; \label{k_1}
			\item $|V_1^1|=\dots =|V_k^1|\eqqcolon m_1$; \label{m_1}
			\item $|V_0^1|\leq \varepsilon_1 n$; \label{V0'}
			\item all but at most~$\varepsilon_1 k_1^2$ pairs~$G[V_i^1, V_j^1]$ with~$1\leq i<j\leq k_1$ are~$\varepsilon_1$-regular. \label{irregularpairs}
		\end{enumerate}
		 Let~$H$ be the empty graph on~$V$. We will repeatedly add some edges of~$G$ to~$H$. Whenever we do so, we also delete the corresponding edges from~$G$ so that at all stages,~$G$ and~$H$ are edge-disjoint.
	
		\item \textbf{Deleting edges between non-regular pairs and within clusters.}\label{step:nonregular} For any distinct~$i,j\in [k_1]$ such that~$G[V_i^1, V_j^1]$ is not~$\varepsilon_1$-regular, add all the edges of~$G[V_i^1, V_j^1]$ to~$H$. Additionally, for any~$i\in [k_1]$, add the edges of~$G[V_i^1]$ to~$H$. Note that by \crefrange{m_1}{irregularpairs}, we have added at most $\varepsilon_1 n^2+k_1\binom{m_1}{2}\leq 2\varepsilon_1 n^2$ edges to~$H$. 
		
		\item \textbf{Removing edges incident to vertices of small or high degree in~$\varepsilon_1$-regular pairs.}\label{step:degrees} 
		For any distinct~$i,j\in [k_1]$, denote by~$d_{ij}$ the density of~$G[V_i^1, V_j^1]$ and let $V_{ij}^1\coloneqq\{x\in V_i^1\mid d_{G[V_i^1,V_j^1]}(x)=(d_{ij}\pm \varepsilon_1)m_1\}$. Add to~$H$ all edges of~$G[V_i^1,V_j^1]$ which do not have both endpoints in~$ V_{ij}^1\cup V_{ji}^1$. Note that by 
		\cref{lm:badvertices}, $|V_{ij}^1|,|V_{ji}^1| \geq (1-2\varepsilon_1)m_1$ and so we have added at most~$4\varepsilon_1 n^2$ edges to~$H$. Moreover,~$V_{ij}^1$ and~$V_{ji}^1$ are now the support clusters of~$G[V_i^1, V_j^1]$. Finally, if~$d_{ij}\geq 2d$, \cref{lm:verticesedgesremoval} implies that~$G[V_{ij}^1, V_{ji}^1]$ is~$[\varepsilon_1', d_{ij}]$-superregular\COMMENT{\(\varepsilon_1'\)-regularity follows by \cref{lm:verticesedgesremoval} \cref{lm:verticesedgesremovalreg} and superregularity follows by construction of~$V_{ij}^1, V_{ji}^1$}, otherwise, note that~$G[V_{ij}', V_{ji}']$ contains at most~$(2d+\varepsilon_1)m_1$ edges incident to each vertex. 
			
		\item \textbf{Removing non-exceptional vertices of high degree in~$H$.}\label{step:nonexceptionaldeg} Move any~$x\in V\setminus V_0^1$ such that~$d_H(x) > \varepsilon n$ to the exceptional set. If we have removed at least~$\sqrt{\varepsilon_1} m_1$ vertices from a cluster, we add the entire cluster to the exceptional set. Denote by~$V_0^2, V_1^2, \dots, V_{k_2}^2$ the partition obtained. We may assume without loss of generality that~$V_i^2\subseteq V_i^1$ for each~$i\in [k_2]$. For any distinct~$i,j\in[k_2]$, let~$V_{ij}^2, V_{ji}^2$ be the support clusters of~$G[V_i^2, V_j^2]$. Note that if~$d_{ij}\geq 2d$, then, by \cref{step:degrees}, \(|V_{ij}^2|, |V_{ji}^2|\geq (1-2\varepsilon_1-\sqrt{\varepsilon_1})m_1\geq (1-\varepsilon_1'')m_1\) and, by \cref{lm:verticesedgesremoval},~$G[V_{ij}^2, V_{ji}^2]$ is~$[\varepsilon_1'', d_{ij}]$-superregular. 
		Note that the total number of edges removed in \cref{step:nonregular,step:degrees} is at most~$6\varepsilon_1 n^2$. Thus the number of clusters we fully add to~$V_0^1$ is at most~$\varepsilon_1' k_1$\COMMENT{at most $\frac{2\cdot 6\varepsilon_1 n^2}{\varepsilon n \cdot \sqrt{\varepsilon_1}m_1}=\frac{12 \sqrt{\varepsilon_1}n}{\varepsilon m_1}\leq \varepsilon_1' k_1$} and $|V_0^2|\leq |V_0^1|+\sqrt{\varepsilon_1} n+ \varepsilon_1'n \leq \varepsilon_1'' n$. Moreover, by \cref{k_1}, $k_1\geq k_2\geq \left(1-\varepsilon_1'\right)k_1\geq \left(1-\varepsilon_1'\right)N\geq M$.
		
		\item \textbf{Removing small density pairs.}\label{step:smalldensitypairs} For any~$i,j\in [k_2]$ such that~$d_{ij}<2d$\COMMENT{We remove pairs of density less than~$2d$ instead of less than~$d$ because we will take out a random subgraph later on.}, add all edges of~$G[V_i^2, V_j^2]$ to~$H$. By \cref{step:nonexceptionaldeg,step:degrees}, for any~$x\in V\setminus V_0^2$, we have $d_H(x)\leq(2d+\varepsilon_1)n +\varepsilon n\leq 3dn$\COMMENT{Since~$G[V_{ij}^2, V_{ji}^2]$ is a subgraph of~$G[V_{ij}^1, V_{ji}^1]$ and every vertex had degree at most~$\varepsilon n$ in~$H$ at the end of \cref{step:nonexceptionaldeg}.}.
		
		\item \textbf{Removing vertices of degree zero in many non-empty pairs.}\label{step:degreezerovertices}
		Let~$R_2$ be the reduced graph of~$G$ with respect to the partition~$V_0^2, V_1^2, \dots, V_{k_2}^2$. For any~$ij\in E(R_2)$, recall that $|V_i^2\setminus V_{ij}^2|\leq \varepsilon_1''m_1$. Thus, for each~$i\in [k_2]$, there are at most~$\varepsilon_1^*m_1$ vertices~$x\in V_i^2$ such that~$x\notin V_{ij}^2$ for at least~$\varepsilon_1^*k_2$ indices~$j$ such that~$ij\in E(R_2)$.%
		\COMMENT{$\frac{\varepsilon_1''m_1\cdot k_2}{\varepsilon_1^* k_2}$}
		Add all these vertices to the exceptional set~$V_0^2$. We now have $|V_0^2|\leq \varepsilon_1'' n+\varepsilon_1^*n\leq 2\varepsilon_1^*n$. Moreover, by \cref{lm:verticesedgesremoval}, for any distinct~$i,j\in [k_2]$, the pair~$G[V_{ij}^2, V_{ji}^2]$ is~$[\varepsilon_1^{**}, d_{ij}]$-superregular. Finally, the following holds.
		
		\begin{property}{\dagger}\label{property:notinsupport}
			For any~$i\in [k]$ and~$x\in V_i^2$, there are at most~$\varepsilon_1^*k_2$ indices~$j$ such that~$ij\in E(R_2)$ but~$x\notin V_{ij}^2$.
		\end{property}
		
		\item \textbf{Splitting clusters into subclusters of equal size.}\label{step:splitting} 
		\cref{step:nonexceptionaldeg,step:degreezerovertices} have the effect that the cluster sizes are no longer equal. To equalise them again, we proceed as follows.
		Randomly split~$V_1^2, \dots, V_{k_2}^2$ into subclusters of size exactly $m_3\coloneqq \left\lfloor\frac{\varepsilon_2 m_1}{4}\right\rfloor$ and put leftover vertices in the exceptional set. Let~$V_0^3, V_1^3, \dots, V_{k_3}^3$ be the resulting partition and~$R_3$ be the corresponding reduced graph of~$G$. Note that by \cref{step:nonexceptionaldeg,k_1},~$k_3\geq k_2\geq M$ and $k_3\leq \frac{4k_2}{\varepsilon_2}\leq \frac{4N'}{\varepsilon_2}$. 
		Moreover, $|V_0^3|\leq \left|V_0^2\right|+\frac{\varepsilon_2 m_1}{4}\cdot k_2\leq \varepsilon_2 n$. We claim that 
		\begin{property}{\ddagger}\label{property:supregpart}
			$V_0^3, V_1^3, \dots, V_{k_{3}}^3$ is an $(\varepsilon_2, \geq 2d, k_3, m_3, R_3)$-superregular partition of~$G$, with support clusters induced by the partition $V_0^2, V_1^2, \dots, V_{k_2}^2$.
		\end{property}
		Indeed, \cref{def:supregpartitionkm,def:supregpartitionV0,def:supregpartitionreduced} are clearly satisfied, while \cref{def:supregpartitioninside} holds by \cref{step:nonregular}. One can show that \cref{def:supregpartitionsupreg} holds with positive probability using \cref{lm:vertexpartition,lm:Chernoff}%
		\COMMENT{Recall from \cref{step:nonregular,step:nonexceptionaldeg,step:smalldensitypairs} that for any distinct \(i,j\in[k_2]\), if~$d_{ij}<2d$, then~$G[V_i^2, V_j^2]$ is empty, and otherwise,~$G[V_{ij}^2, V_{ji}^2]$ is \([\varepsilon_1'', \geq d]\)-superregular, with $|V_{ij}^2|, |V_{ji}^2|\geq (1-\varepsilon_1'')m_1$.
		Now assume \(V_{i'}^3\subseteq V_i^2\) and \(V_{j'}^3\subseteq V_j^2\) ($i',j'\in [k_3]$), and let \(V_{i'j'}^3=V_{i'}^3\cap V_{ij}^2\) and \(V_{j'i'}^3=V_{j'}^3\cap V_{ji}^2\). By \cref{lm:vertexpartition},~$G[V_{i'j'}^3, V_{j'i'}^3]$ is~$[\varepsilon_2, d_{ij}]$-superregular with high probability. Note that in particular, since~$d_{ij}\geq 2d$,~$V_{i'j'}^3, V_{j'i'}^3$ are the support clusters of the pair~$G[V_{i'}^3, V_{j'}^3]$. Moreover, we have \(\mathbb{E}[|V_{i'j'}^3|]\geq (1-\varepsilon_1'')m_3\) so by \cref{lm:Chernoff}
			\begin{equation*}
			\begin{aligned}
			\mathbb{P}\left[|V_{i'j'}|< (1-\varepsilon_2)m_3\right]&\leq\mathbb{P} \left[|V_{i'j'}<(1-\varepsilon_2^2)\mathbb{E}[|V_{i'j'}|]\right]
			\leq \exp \left(-\frac{\varepsilon_2^4(1-\varepsilon_1'')m_3}{3}\right).
			\end{aligned}
			\end{equation*}
			Similarly, \(\mathbb{P}[|V_{j'i'}|< (1-\varepsilon_2)m_3]\leq \exp \left(-\frac{\varepsilon_2^4(1-\varepsilon_1'')m_3}{3}\right)\).
			A union bound over all~$1\leq i'<j'\leq k$ gives that \cref{def:supregpartitionsupreg} holds with positive probability.}, %
		so we may assume that \cref{property:supregpart} is satisfied. For any distinct~$i,j\in [k_3]$, denote by~$V_{ij}^3$ and~$V_{ji}^3$ the support clusters of the pair~$G[V_i^3, V_j^3]$. Then, the following holds.
		\begin{property}{\dagger\dagger}\label{property:badinfewpairs}
			For each~$i\in[k_3], x\in V_i^3$, there are at least~$2\beta k_3$ indices~$j$ such that~$x\in V_{ij}^3$ and~$G[V_{ij}^3, V_{ji}^3]$ is $[\varepsilon_2, \geq \frac{\alpha}{2}]$-superregular. Moreover, there are at most~$\varepsilon k_3$ indices~$j$ such that~$ij\in E(R_3)$ but~$x\notin V_{ij}^3$.
		\end{property}
		Indeed, if the first part of \cref{property:badinfewpairs} is not satisfied then, by \cref{step:smalldensitypairs}, we have \(\delta (G\cup H)<2\beta n+(\frac{\alpha}{2}+\varepsilon_2)n+3dn+\varepsilon_2 n<\alpha n\), a contradiction. The second part of the statement follows from \cref{property:notinsupport,property:supregpart}.%
		\COMMENT{at most $\frac{4\varepsilon_1^*k_2}{\varepsilon_2}\leq \varepsilon k_3$ indices~$j$ such that~$ij\in E(R_3)$ but~$x\notin V_{ij}^3$, since support clusters are induced by the partition $V_0^2, V_1^2, \dots, V_{k_2}^2$.}
		
		\item \textbf{Finding~$\Gamma$ and~$\Gamma'$.}\label{step:randomslices}
		For simplicity, for any distinct~$i,j\in [k]$, let $G_{ij}\coloneqq  G[V_{ij}^3, V_{ji}^3]$. Also define~$d_{ij}'$ as follows. If~$G_{ij}$ is empty, let~$d_{ij}'\coloneqq 0$, otherwise, define~$d_{ij}'$ as the largest constant such that~$G_{ij}$ is~$[\varepsilon_2, d_{ij}']$-superregular.\COMMENT{Needed since some pairs could be both $[\varepsilon_2, <\frac{\alpha}{2}]$-superregular and $[\varepsilon_2, \geq\frac{\alpha}{2}]$-superregular but we want to proceed differently for each of these cases (a~$\zeta$ density slice of the former pairs will form~$\Gamma'$ and a~$\beta$ density slice of the latter pairs will form~$\Gamma$). We take~$d_{ij}'$ as large as possible so that by \cref{property:badinfewpairs} we will get a nice lower bound on the minimum degree of~$\Gamma$.} 
	
	We construct~$\Gamma$ as follows. Let~$S$ be the set of pairs~$(i,j)$ such that~$1\leq i<j\leq k$ and~$d_{ij}'\geq \frac{\alpha}{2}$. For any~$(i,j)\in S$, apply \cref{lm:edgepartition} to obtain a spanning subgraph $\Gamma_{ij}\subseteq G_{ij}$ such that~$\Gamma_{ij}$ is~$[\varepsilon_2', \beta]$ superregular and~$G_{ij}\setminus\Gamma_{ij}$ is $[\varepsilon_2', d_{ij}'-\beta]$-superregular. Let~$\Gamma$ be the graph with vertex set~$V$ and edge set $\bigcup_{(i,j)\in S} E(\Gamma_{ij})$.
	
	We construct~$\Gamma'$ using similar arguments as above, so that for any~$1\leq i<j\leq k$ such that~$d_{ij}'<\frac{\alpha}{2}$,~$\Gamma'[V_{ij}^3, V_{ji}^3]$ is~$[\varepsilon_2', \zeta]$-superregular and~$G_{ij}\setminus\Gamma'$ is $[\varepsilon_2', d_{ij}'-\zeta]$-superregular. In particular,~$\Gamma$ and~$\Gamma'$ are edge-disjoint. Set $G'\coloneqq G\setminus (\Gamma \cup \Gamma')$.
		
		\item \textbf{Decomposing~$R'$ into long cycles of even length.}\label{step:lovasz}
		Apply \cref{thm:Lovasz} to obtain a decomposition~$\cD$ of the reduced graph of~$\Gamma$ into at most~$\frac{k_3}{2}$ paths and cycles. Let~$V_0, V_1, \dots, V_k$ be the partition obtained from~$V_0^3, V_1^3, \dots, V_{k_3}^3$ by randomly splitting~$V_i^3$ into~$2^L$ subclusters of size $m\coloneqq \left\lfloor \frac{m_3}{2^L}\right\rfloor$ and adding the leftover to the exceptional set, for each~$i\in [k_3]$.\COMMENT{We split into~$2^L$ parts rather than~$2L$ for simplicity.} Thus, $|V_0|\leq \varepsilon_2'' n$. Denote by~$R, R'$, and~$R''$ the corresponding reduced graphs of~$G', \Gamma$, and~$\Gamma'$. Then, by \cref{lm:vertexpartition}, \cref{step:randomslices}, and \cref{property:badinfewpairs}, we may assume that the following hold.
		\begin{itemize}
			\item $V_0, V_1, \dots, V_k$ is
			\begin{itemize}
				\item an $(\varepsilon_2'', \geq d, k, m, R)$-superregular partition of~$G'$,
				\item an $(\varepsilon_2'', \beta, k, m, R')$-superregular partition of~$\Gamma$,
				\item an $(\varepsilon_2'', \zeta, k, m, R'')$-superregular partition of~$\Gamma'$.
			\end{itemize}
			\item Each~$x\in V\setminus V_0$ belongs to at least~$2\beta k$ superregular pairs of~$\Gamma$.
			\item For each~$i\in [k]$ and~$x\in V_i$, there are at most~$\varepsilon k$ indices~$j\in [k]$ such that~$ij\in E(R')$ but~$d_{\Gamma[V_i,V_j]}(x)=0$.\COMMENT{Actually such that~$ij\in E(R)$ but~$d_{G[V_i,V_j]}(x)=0$ by \cref{property:badinfewpairs}. But we are only interested in pairs of~$\Gamma$ for \cref{step:support}.}
			\item Properties \cref{lm:goodsamesupport,lm:goodreduced} are satisfied. Moreover,~$V_0$ is a set of isolated vertices in~$\Gamma$ and~$\Gamma'$, as desired for \cref{lm:goodV0GammaH}.
			\item Let~$\cD_{R'}$ be the decomposition of~$R'$ induced by~$\cD$. Then~$\cD_{R'}$ is a set of at most~$2^L\frac{k_3}{2}=\frac{k}{2}$ cycles whose lengths are even and at least~$L$.
			\COMMENT{Splitting in~$2^L$ parts is the same as repeatedly splitting in two~$L$ times. \\
			Assume we have some clusters $A_1,\dots, A_\ell$ and say each~$A_i$ is partitioned into~$A_i^1, A_i^2$.
			A path~$A_1$, $\dots$, $A_{2\ell+1}$ gives two cycles $A_1^1$, $A_2^2$, $A_3^1$, $A_4^2$, $\dots$, $A_{2\ell+1}^1$, $A_{2\ell}^1$, $A_{2\ell-1}^2$, $\dots$, $A_2^1$ and $A_1^2$, $A_2^2$, $A_3^2$, $\dots$, $A_{2\ell+1}^2$, $A_{2\ell}^1$, $\dots$, $A_2^1$ of length~$4\ell$.
			A path~$A_1$, $\dots$, $A_{2\ell}$ gives two cycles $A_1^1$, $A_2^2$, $A_3^1$, $A_4^2$, $\dots$, $A_{2\ell}^2$, $A_{2\ell-1}^2$, $A_{2\ell-2}^1$, $\dots$, $A_2^1$ and $A_1^2$, $A_2^2$, $A_3^2$, $\dots$, $A_{2\ell-1}^2$, $A_{2\ell}^1$, $A_{2\ell -1}^1$, $\dots$, $A_2^1$ of length~$4\ell-2$. 
			A cycle~$A_1$, $\dots$, $A_{2\ell+1}$  gives two cycles $A_1^1$, $A_2^2$, $A_3^1$, $A_4^2$, $\dots$, $A_{2\ell+1}^{1}$, $A_1^2$, $A_{2\ell+1}^2$, $A_{2\ell}^1$, $A_{2\ell-1}^2$, $\dots$, $A_2^1$ and $A_1^2$, $A_2^2$, $\dots$, $A_{2\ell+1}^2$, $A_1^1$, $A_{2\ell+1}^1$, $A_{2\ell}^1$, $\dots$, $A_2^1$ of length~$4\ell+2$.
			A cycle~$A_1$, $\dots$, $A_{2\ell}$  gives two cycles $A_1^1$, $A_2^2$, $A_3^1$, $A_4^2$, $\dots$, $A_{2\ell}^{2}$, $A_1^2$, $A_{2\ell}^1$, $A_{2\ell-1}^2$, $\dots$, $A_2^1$ and $A_1^2$, $A_2^2$, $\dots$, $A_\ell^2$, $A_1^1$, $A_{2\ell}^1$, $A_{2\ell-1}^1$, $\dots$, $A_2^1$ of length~$4\ell$.}
		\end{itemize}
	 	 Also observe that~$k=2^Lk_3$ and so by \cref{step:splitting}, we may set~$M'\coloneqq \frac{2^{L+2}N'}{\varepsilon_2}$ in order to satisfy \cref{lm:goodk}.
		
		\item \textbf{Removing from~$H$ all edges with exactly one endpoint in~$V_0$.}	\label{step:cleaningexceptionaledges} Delete from~$H$ all edges with an endpoint in~$V_0$ in order to satisfy \cref{lm:goodV0GammaH}. Add these edges to~$G'$.
		
		\item \textbf{Equalising the support cluster sizes.}\label{step:support} Let~$ij\in E(R)$. Denote by~$V_{ij}$ the support cluster of~$V_i$ with respect to~$V_j$ (for the graph~$G$). We define~$V_{ij}'$ as follows. Let~$j'$ be such that~$V_jV_iV_{j'}$ is a subpath of a cycle in~$\cD_{R'}$\COMMENT{(so~$j'$ is uniquely defined if it exists)}, and, let $V_{ij}'\coloneqq V_{ij}\cap V_{ij'}$ if~$j'$ exists,~$V_{ij}'\coloneqq V_{ij}$ otherwise. Note that, by \cref{step:lovasz}, for each~$i\in [k]$ and~$x\in V_i$, there are at most~$\varepsilon k$ indices~$j$ such that~$x\in V_{ij}\setminus V_{ij}'$. Moreover, $|V_{ij}'|\geq (1-2\varepsilon_2'')m$. Pick $(1-\varepsilon^4)m\leq m'\leq (1-2\varepsilon_2'')m$ such that~$\frac{m'}{r}\in \mathbb{N}^*$, as desired for \cref{lm:goodm'}. Note that $|V_{ij}'|-m'\geq 0$. 
		
		For each~$ij\in E(R)$ in turn, we now construct subset~$V_{ij}''$ of~$V_{ij}'$ by removing exactly~$|V_{ij}'|-m'$ vertices. We build these sets one by one and, in each step, we only remove vertices which have already been removed fewer than~$\varepsilon^3k$ times in the construction so far. This is possible since we need to remove at most~$\varepsilon^4 m$ vertices from each~$V_{ij}'$. Thus, in each step, there at most~$\varepsilon m$ vertices\COMMENT{at most $\frac{\varepsilon^4 m k}{\varepsilon^3 k}$} that cannot be removed anymore. For each~$ij\in E(R)$, delete from~$G'[V_i, V_j], \Gamma[V_i, V_j]$, and~$\Gamma'[V_i, V_j]$ all edges with an endpoint in~$(V_{ij}\setminus V_{ij}'')\cup (V_{ji}\setminus V_{ji}'')$, and add all these edges to~$H$. For any~$i\in [k]$ and~$x\in V_i$, if~$x\in V_{ij}''$ then we have removed at most~$\varepsilon^4 m$ edges of~$G[V_i, V_j]$ which are incident to~$x$, and, by construction, there are at most~$\varepsilon^3k$ indices~$j$ such that~$x\in V_{ij}\setminus V_{ij}''$. Thus, we have added to~$H$ at most~$2\varepsilon^3 n$ edges incident to each vertex and, by \cref{step:smalldensitypairs}, \cref{lm:goodDeltabad} holds. Moreover, by \cref{step:lovasz,lm:verticesedgesremoval}, \cref{lm:goodsupregeqpart,lm:goodbadinfewpairs,lm:goodsamesupport,lm:goodreducedpartition} are satisfied, and this finishes the proof.
		\qed
	\end{steps}
	\renewcommand{\qed}{}
\end{proof}

\onlyinsubfile{\bibliographystyle{plain}
	\bibliography{Bibliography/papers}}

%% file: main.bbl
\begin{thebibliography}{10}

\bibitem{alon1994algorithmic}
N.~Alon, R.~Duke, H.~Lefmann, V.~R{\"o}dl, and R.~Yuster.
\newblock The algorithmic aspects of the regularity lemma.
\newblock {\em Journal of Algorithms}, 16(1):80--109, 1994.

\bibitem{alspach2008wonderful}
B.~Alspach.
\newblock The wonderful {W}alecki construction.
\newblock {\em Bulletin of the Institute of Combinatorics and its
  Applications}, 52:7--20, 2008.

\bibitem{bienia1986partitions}
W.~Bienia and H.~Meyniel.
\newblock Partitions of digraphs into paths or circuits.
\newblock {\em Discrete Mathematics}, 61(2-3):329--331, 1986.

\bibitem{bollobas1996proof}
B.~Bollob{\'a}s and A.~Scott.
\newblock A proof of a conjecture of {B}ondy concerning paths in weighted
  digraphs.
\newblock {\em Journal of Combinatorial Theory, Series B}, 66(2):283--292,
  1996.

\bibitem{bonamy2019gallai}
M.~Bonamy and T.~J. Perrett.
\newblock Gallai's path decomposition conjecture for graphs of small maximum
  degree.
\newblock {\em Discrete Mathematics}, 342(5):1293--1299, 2019.

\bibitem{bondy2008graph}
J.~A. Bondy and U.~S.~R. Murty.
\newblock Graph theory.
\newblock {\em Graduate Texts in Mathematics}, 2008.

\bibitem{botler2017path}
F.~Botler and A.~Jim{\'e}nez.
\newblock On path decompositions of $2k$-regular graphs.
\newblock {\em Discrete Mathematics}, 340(6):1405--1411, 2017.

\bibitem{conlon2014cycle}
D.~Conlon, J.~Fox, and B.~Sudakov.
\newblock Cycle packing.
\newblock {\em Random Structures \& Algorithms}, 45(4):608--626, 2014.

\bibitem{csaba2016proof}
B.~Csaba, D.~K{\"u}hn, A.~Lo, D.~Osthus, and A.~Treglown.
\newblock Proof of the 1-factorization and {H}amilton decomposition
  conjectures.
\newblock {\em Memoirs of the American Mathematical Society}, 244(1154), 2016.

\bibitem{dean1986smallest}
N.~Dean.
\newblock What is the smallest number of dicycles in a dicycle decomposition of
  an {E}ulerian digraph?
\newblock {\em Journal of Graph Theory}, 10(3):299--308, 1986.

\bibitem{dean2000gallai}
N.~Dean and M.~Kouider.
\newblock Gallai's conjecture for disconnected graphs.
\newblock {\em Discrete Mathematics}, 213(1-3):43--54, 2000.

\bibitem{donald1980upper}
A.~Donald.
\newblock An upper bound for the path number of a graph.
\newblock {\em Journal of Graph Theory}, 4(2):189--201, 1980.

\bibitem{erdos1983some}
P.~Erd{\H{o}}s.
\newblock On some of my conjectures in number theory and combinatorics.
\newblock In {\em Proceedings of the Fourteenth Southeastern Conference on
  Combinatorics, Graph Theory and Computing}, volume~39, pages 3--19, 1983.

\bibitem{erdos1966representation}
P.~Erd{\H{o}}s, A.~Goodman, and L.~P{\'o}sa.
\newblock The representation of a graph by set intersections.
\newblock {\em Canadian Journal of Mathematics}, 18:106--112, 1966.

\bibitem{fan2002subgraph}
G.~Fan.
\newblock Subgraph coverings and edge switchings.
\newblock {\em Journal of Combinatorial Theory, Series B}, 84(1):54--83, 2002.

\bibitem{fan2003covers}
G.~Fan.
\newblock Covers of {E}ulerian graphs.
\newblock {\em Journal of Combinatorial Theory, Series B}, 89(2):173--187,
  2003.

\bibitem{fan2005path}
G.~Fan.
\newblock Path decompositions and {G}allai's conjecture.
\newblock {\em Journal of Combinatorial Theory, Series B}, 93(2):117--125,
  2005.

\bibitem{fan2002hajos}
G.~Fan and B.~Xu.
\newblock Haj{\'o}s' conjecture and projective graphs.
\newblock {\em Discrete Mathematics}, 252(1-3):91--101, 2002.

\bibitem{frieze2005packing}
A.~Frieze and M.~Krivelevich.
\newblock On packing {H}amilton cycles in {$\epsilon$}-regular graphs.
\newblock {\em Journal of Combinatorial Theory, Series B}, 94(1):159--172,
  2005.

\bibitem{glock2016optimal}
S.~Glock, D.~K{\"u}hn, and D.~Osthus.
\newblock Optimal path and cycle decompositions of dense quasirandom graphs.
\newblock {\em Journal of Combinatorial Theory, Series B}, 118:88--108, 2016.

\bibitem{harding2014gallai}
P.~Harding and S.~McGuinness.
\newblock Gallai's conjecture for graphs of girth at least four.
\newblock {\em Journal of Graph Theory}, 75(3):256--274, 2014.

\bibitem{jackson1980decompositions}
B.~Jackson.
\newblock Decompositions of graphs into cycles.
\newblock {\em Regards sur la Th\'{e}orie des Graphes}, pages 259--261, 1980.

\bibitem{janson2011random}
S.~Janson, T.~\L{}uczak, and A.~Ruci\'{n}ski.
\newblock {\em Random Graphs}, volume~45.
\newblock John Wiley \& Sons, 2011.

\bibitem{jimenez2017path}
A.~Jim{\'e}nez and Y.~Wakabayashi.
\newblock On path-cycle decompositions of triangle-free graphs.
\newblock {\em Discrete Mathematics \& Theoretical Computer Science}, 19, 2017.

\bibitem{knierim2019long}
C.~Knierim, M.~Larcher, A.~Martinsson, and A.~Noever.
\newblock Long cycles, heavy cycles and cycle decompositions in digraphs.
\newblock {\em arXiv preprint arXiv:1911.07778}, 2019.

\bibitem{korandi2015decomposing}
D.~Kor{\'a}ndi, M.~Krivelevich, and B.~Sudakov.
\newblock Decomposing random graphs into few cycles and edges.
\newblock {\em Combinatorics, Probability and Computing}, 24(6):857--872, 2015.

\bibitem{kuhn2013hamilton}
D.~K{\"u}hn and D.~Osthus.
\newblock Hamilton decompositions of regular expanders: a proof of {K}elly's
  conjecture for large tournaments.
\newblock {\em Advances in Mathematics}, 237:62--146, 2013.

\bibitem{kuhn2014hamilton}
D.~K{\"u}hn and D.~Osthus.
\newblock Hamilton decompositions of regular expanders: applications.
\newblock {\em Journal of Combinatorial Theory, Series B}, 104:1--27, 2014.

\bibitem{lovasz1968covering}
L.~Lov{\'a}sz.
\newblock On covering of graphs.
\newblock In {\em Theory of Graphs (Proceedings of the Colloquium held at
  {T}ihany, {H}ungary, {S}eptember 1966)}, pages 231--236. Academic Press New
  York, 1968.

\bibitem{lucas1883recreationsII}
{\'E}.~Lucas.
\newblock {\em R{\'e}cr{\'e}ations Math{\'e}matiques}, volume~II.
\newblock Gauthier-Villars, 1883.

\bibitem{pyber1985erdHos}
L.~Pyber.
\newblock An {E}rd{\H{o}}s-{G}allai conjecture.
\newblock {\em Combinatorica}, 5(1):67--79, 1985.

\bibitem{pyber1996covering}
L.~Pyber.
\newblock Covering the edges of a connected graph by paths.
\newblock {\em Journal of Combinatorial Theory, Series B}, 66(1):152--159,
  1996.

\bibitem{szemeredi1978regular}
E.~Szemer\'{e}di.
\newblock Regular partitions of graphs.
\newblock In {\em Probl\`{e}mes Combinatoires et Th\'{e}orie des Graphes
  ({C}olloque {I}nternational du {CNRS}, {U}niversit\'{e} d'{O}rsay, {O}rsay,
  1976)}, volume 260 of {\em Colloques Internationaux du CNRS}, pages 399--401.
  CNRS, Paris, 1978.

\bibitem{taylor2014regularity}
A.~Taylor.
\newblock {\em The regularity method for graphs and digraphs}.
\newblock M{S}ci thesis, University of Birmingham, 2014.

\bibitem{yan1998path}
L.~Yan.
\newblock {\em On path decompositions of graphs}.
\newblock PhD thesis, Arizona State University, 1998.

\end{thebibliography}
